\def\today{\ifcase\month\or January\or February\or March\or
April\or May\or June\or July\or August\or September\or
October\or November\or December\fi \space\number\day,
\number\year}

\def\a{\alpha}
\def\b{\beta}
\def\d{\delta}

\def\g{\gamma}

\def\l{\lambda}
\def\om{\omega}
\def\r{\rho}
\def\s{\sigma}
\def\t{\theta}

\def\ve{\varepsilon}
\def\vp{\varphi}

\def\z{\zeta}

\def\D{\Delta}
\def\G{\Gamma}
\def\L{\Lambda}
\def\O{\Omega}
\def\Si{\Sigma}
\def\T{\Theta}

\let\nn=\noindent

\font\tenbb=msym10
\font\sevenbb=msym8
\font\fivebb=msym5
\newfam\bbfam
\textfont\bbfam=\tenbb \scriptfont\bbfam=\sevenbb
\scriptscriptfont\bbfam=\fivebb
\def\bb{\fam\bbfam}

\def\CC{{\bb C}}

\def\NN{{\bb N}}

\def\RR{{\bb R}}

\def\ZZ{{\bb Z}}

\font\titre=cmbx12
\font\tenfm=eufm10

\def\part{\partial}
\def\Inf{\infty}
\def\bl{\backslash}

\def\bu{\bullet}

\def\bgt{\nabla}
\def\mps{\mapsto}
\def\ts{\times}

\def\sbs{\subset}

\def\ra{\rightarrow}

\def\lra{\leftrightarrow}

\def\lbc{\lbrace}
\def\rbc{\rbrace}
\def\lbk{\lbrack}
\def\rbk{\rbrack}

\def\oo{\overline}
\def\ww{\widetilde}
\def\hh{\widehat}

\def\ds{\displaystyle }

\def\Da{{\cal D}}

\def\Fa{{\cal F}}

\def\Ja{{\cal J}}
\def\Ka{{\cal K}}
\def\La{{\cal L}}
\def\Ma{{\cal M}}

\def\Pa{{\cal P}}

\def\Sa{{\cal S}}
\def\Ta{{\cal T}}
\def\Ua{{\cal U}}

\def\and{\mathop{\rm and}\nolimits}

\def\by{\mathop{\rm by}\nolimits}

\def\cos{\mathop{\rm cos}\nolimits}

\def\det{\mathop{\rm det}\nolimits}

\def\diag{\mathop{\rm diag}\nolimits}

\def\Det{\mathop{\rm Det}\nolimits}

\def\Dev{\mathop{\rm Dev}\nolimits}

\def\Eucl{\mathop{\rm Eucl}\nolimits}
\def\Eucl.{\mathop{\rm Eucl.}\nolimits}

\def\exp{\mathop{\rm exp}\nolimits}

\def\or{\mathop{\rm or}\nolimits}

\def\sin{\mathop{\rm sin}\nolimits}

\def\supp{\mathop{\rm supp}\nolimits}
\def\supp.{\mathop{\rm supp.}\nolimits}

\def\space{\mathop{\rm space}\nolimits}

\def\Tr{\mathop{\rm Tr}\nolimits}
\def\tr{\mathop{\rm tr}\nolimits}

\def\Var{\mathop{\rm Var}\nolimits}

\catcode`\@=11
\def\displaylinesno #1{\displ@y\halign{
\hbox to\displaywidth{$\@lign\hfil\displaystyle##\hfil$}&
\llap{$##$}\crcr#1\crcr}}

\def\ldisplaylinesno #1{\displ@y\halign{
\hbox to\displaywidth{$\@lign\hfil\displaystyle##\hfil$}&
\kern-\displaywidth\rlap{$##$}
\tabskip\displaywidth\crcr#1\crcr}}
\catcode`\@=12

\def\buildrel#1\over#2{\mathrel{
\mathop{\kern 0pt#2}\limits^{#1}}}

\def\build#1_#2^#3{\mathrel{
\mathop{\kern 0pt#1}\limits_{#2}^{#3}}}

\def\hfl#1#2{\smash{\mathop{\hbox to 6mm{\rightarrowfill}}
\limits^{\scriptstyle#1}_{\scriptstyle#2}}}

\def\hfll#1#2{\smash{\mathop{\hbox to 6mm{\leftarrowfill}}
\limits^{\scriptstyle#1}_{\scriptstyle#2}}}

\def\up#1{\raise 1ex\hbox{\sevenrm#1}}

\def\per{|\!\raise -4pt\hbox{$-$}}
 \def\cqfd{\unskip\kern
6pt\penalty 500 \raise
-2pt\hbox{\vrule\vbox to10pt{\hrule
width 4pt \vfill\hrule}\vrule}\par}

\def\trait{\hbox to 12mm{\hrulefill}}
\def\2{{\mathop{\rm 
I }\nolimits}\!{\mathop{\rm  
I}\nolimits}}

\def\1{{\mathop{\rm I }\nolimits}}

\def\og{\leavevmode\raise.3ex\hbox{$
\scriptscriptstyle\langle\!\langle$}}
\def\fg{\leavevmode\raise.3ex\hbox{$
\scriptscriptstyle \,\rangle\!\rangle$}}

\def\picture #1 by #2 (#3)
{\vcenter{\vskip #2
\special{picture #3}
\hrule width #1 height 0pt depth 0pt
\vfil}}

\def\scaledpicture #1 by #2 (#3 scaled #4){{
\dimen0=#1 \dimen1=#2
\divide\dimen0 by 1000 \multiply\dimen0 by #4
\divide\dimen1 by 1000 \multiply\dimen1 by #4
\picture \dimen0 by \dimen1 (#3 scaled #4)}}

\def\un{{\rm 1\mkern-4mu l}}
\def\[{{[\mkern-3mu [}}
\def\]{{]\mkern-3mu ]}}

\def\tvi{\vrule height 12pt depth 5pt width 0pt}

\def\tv{\tvi\vrule}

\def\cc#1{\hfill\kern .7em#1\kern .7em\hfill}

\def\TeX{T\kern-.1667em\lower.5ex\hbox{E}\kern-.125em X}

\def\ins{{ \raise -2mm\hbox{$<$}  
\atop \raise 2mm\hbox{$\sim$}J}
}

\def\sus{{ \raise -2mm\hbox{$>$}  
\atop \raise 2mm\hbox{$\sim$}J}
}

\def\lta{\hbox{\raise.5ex\hbox{$<$}
\kern-1.1em\lower.5ex\hbox{$\sim$}}}
\def\gta{\hbox{\raise.5ex\hbox{$>$}
\kern-1.1em\lower.5ex\hbox{$\sim$}}}

\magnification=1200
\overfullrule=0mm
\baselineskip=12pt
 
\hsize=118mm 
\hoffset=7mm
\vsize=185mm
\voffset=10mm

\vglue 15mm

\centerline{\titre A new perspective on}
\smallskip
\centerline{\titre  Functional Integration}

\vglue 5mm
{\baselineskip=11pt
\centerline{\bf Pierre Cartier}

\centerline{ Ecole Normale Sup\'erieure }

\centerline{45 rue d'Ulm}

\centerline{F-75230 Paris C\'edex 05}

\centerline{France} 
\smallskip
\centerline{\it and}
\smallskip

\centerline{\bf C\'ecile DeWitt-Morette}

\centerline{Center for Relativity}

\centerline{and Department of Physics}

\centerline{The University of Texas at Austin}

\centerline{Austin, Texas 78712-1081 USA}
\bigskip

\centerline{Journal of Mathematical Physics, {\bf 36}, 2137-2340 (1995)}
\par }

\vglue 20mm

\nn
\centerline {\bf Abstract}

\bigskip

The core of this article is a general theorem with a large number of specializations. Given a
manifold $N$ and a finite number of one-parameter groups of point transformations on $N$ with
generators $Y, X_{(1)}, \cdots, X_{(d)} $, we obtain, via
functional integration over spaces of pointed paths on $N$ (paths with one fixed point), a
one-parameter group of functional operators acting on tensor or spinor fields on $N$. The
generator of this group is a quadratic form in the Lie derivatives $\La_{X_{(\a)}}$ in the 
$X_{(\a)}$-direction plus a term linear in $\La_Y$. 

\smallskip
The basic functional integral is over $L^{2,1}$ paths $x: {\bf T} \ra N$ (continuous paths
with square integrable first derivative). Although the integrator is invariant under time
translation, the integral is powerful enough to be used for systems which are not time
translation invariant. We give seven non trivial applications of the basic formula, and we
compute its semiclassical expansion.

\smallskip
The methods of proof are rigorous and combine Albeverio H\o egh-Krohn oscillatory integrals
with Elworthy's parametrization of paths in a curved space. Unlike other approaches we solve
 Schr\"odinger type equations directly, rather than solving first diffusion equations
and then using analytic continuation.

\vfill\eject

\nn
\centerline {\titre I - Introduction}

\bigskip
\bigskip

We have studied many applications of functional integration looking for its
{\it substantifique moelle}\footnote{$^{1}$ }{
F. Rabelais.
Literally ``bone marrow as a producer of substance''
in {\it Gargantua}, Prologue de l'auteur.}. Little by little, several ideas have taken shape 
concerning the domains of integration, the integrators, and the integrands.

\bigskip
\nn
{\bf 1. The domain of integration is a function space.} 

\smallskip
Working with an infinite-dimensional
space is easier than working with the limit for large $n$ of the product of $n$ copies of a
finite-dimensional space. For example, a space of pointed paths (paths with one fixed
point) is contractible even when the paths take their values in a non contractible space.
Over the years the advantages -- often the necessity -- of working with spaces of paths rather
than with the discretized version of the paths have become increasingly apparent:

\medskip
-- Semiclassical approximations, even in the presence of caustics, are obtained by expanding
a functional on the space of paths around a dominating contribution [1,2].

\medskip
-- If the paths take their values in a multiply-connected space the topology of the
space of paths plays a central role [3].

\medskip
-- A change of variable of integration, regarded as a map on the domain of integration
gives, in two lines,  the Cameron-Martin formula [4] obtained from discretized paths via a
lengthy derivation, and generalizes its applications [5].

\medskip
-- Computing functional determinants using properties of linear maps on Banach spaces is
simpler than computing limits of finite-dimensional determinants [6].

\medskip
--  A major progress in the definition and computation of functional integrals was
achieved by Elworthy [7] when he parametrized the space of paths on a Riemannian manifold
$M$, with fixed initial point $a$, by the space of paths on $T_a M$ starting at the
origin of $T_a M$.

\medskip
The importance of the domain of integration, noted in the above examples, is even more
striking in the formulation presented here. A key point is a generalization of 
 Elworthy's
idea.
 Consider a finite-dimensional manifold $N$, and denote by ${\bf T}$ some finite time interval.
 Here
we parametrize a space $\Pa_{x_0} N$ of pointed paths $x: {\bf T} \ra N$ by a space $\Pa_0
\RR^d$ of pointed paths $z: {\bf T} \ra \RR^d$. A general construction for maps
$$\Pa_0 \RR^d \ra \Pa_{x_0} N \qquad  {\mathop{\rm by }\nolimits}
\qquad z \mps x$$
is given in Section II by solving suitable  first-order {\it differential} equations, not by
constructing stochastic processes. It is stated in terms of vector fields $Y$, $X_{(1)}, \cdots , X_{(d)}$ on $N$. By
specializing $N $ and the vector fields, {\it one basic functional integral yields, with rigor and no Ansatz, a great
variety of functional integrals which are solutions of complex problems.}
A number of them are pre\-sented in Section IV.

\medskip

In general, to define the domain of integration, one needs to specify:

-- the analytic nature of the paths\footnote{\nn$^2$}{ \nn A path $z: {\bf
T} \ra \RR^{\bf d}$ is said to be $L^{2, 1}$ if 
$\int_{\bf T} dt \, |\dot z (t)|^2  < \Inf$ where 
$\dot z (t) = dz (t) / dt$.
 }  (continuous, $L^{2}, L^{2,1}, \cdots $);

-- the domains of the paths  and their ranges;

-- the behaviour of the paths  at the boundaries of their domains.

\bigskip
\nn
{\bf 2. Integrators cannot be expected to be universal.}
\smallskip

The naive approach to the definition of an integral in an infinite number of variables is to
take a limit $d = \Inf$ in a $d$-dimensional integral. Because of scaling problems this
procedure is well known to abort. For instance, if we wish to evaluate (for $a > 0$) the
integral $$I_{\Inf}: =  \int_{ \RR^\Inf }^{ }  d^\Inf x \,
\exp 
\left({- { \pi \over 	a} |x|^2 }\right) ,
\eqno  ({\mathop{\rm I.1 }\nolimits})$$
where 
$|x|^2: = \ds \sum_{ \a = 1 }^{ \Inf }
\left({ x^\a }\right)_{}^{2}$,
we may first evaluate the corresponding integral
$$I_d: = 
\int_{ \RR^d }^{ } d x \,
\exp 
\left({- { \pi \over 	a} |x|^2 }\right) 
\eqno  ({\mathop{\rm I.2 }\nolimits})$$
for finite $d$ and then set $d = \Inf$.
Since $I_d = a^{d/2}$ we get
$$I_\Inf =
\left\{\matrix{ 
0 \hfill& {\mathop{\rm if }\nolimits} & 0 < a < 1 \hfill\cr
1 \hfill& {\mathop{\rm if }\nolimits} & a = 1 \hfill\cr
\Inf \hfill& {\mathop{\rm if }\nolimits} & 1 < a \hfill\cr}\right.
\eqno  ({\mathop{\rm I.3 }\nolimits})$$
and this fails to be continuous in the parameter $a$, as should be reasonably desired.

\medskip
A way out of this difficulty is to introduce, for each value of the scaling parameter $a >
0$, an integrator $\Da_a x$ in the $d$-dimensional space $\RR^d$, namely
$$\Da_a x = a^{- d/2} d x^1 \cdots dx^d
\eqno  ({\mathop{\rm I.4 }\nolimits})$$
and to remark that it is characterized by the following integration formula
$$
\int_{ \RR^d }^{ }  \Da_a x \cdot
\exp 
\left({- { \pi \over 	a} |x|^2  - 2 \pi i 
\left\langle{ x', x }\right\rangle
}\right)=
\exp 
\left({  - \pi a |x'|^2 }\right) .
\eqno  ({\mathop{\rm I.5 }\nolimits})$$

\medskip
\nn
Here $x'$ runs over the  space $\RR_d$ dual to $\RR^d$ and the scalar product is
given\footnote{$^3$}{\nn Here and in the rest of
this paper we use the Einstein summation convention over repeated indices.  } by 
$\left\langle{ x', x }\right\rangle = x'_\a x^\a$.

\medskip

This formula is dimension-independent and hence suitable for the generalization from $\RR^d$
to a (real) Banach space ${\bf X}$.
Let ${\bf X'}$ be its dual and consider two continuous quadratic forms $Q$ on 
${\bf X}$ and $W$ on ${\bf X'}$.
Assume that $Q$ {\it and  $W$ are inverse to each other}
in the following sense. There exist continuous linear maps
$$D: {\bf X} \ra {\bf X'} \ ,
\qquad
G: {\bf X'} \ra {\bf X}$$
such that\footnote{$^4$}{\nn
\nn In standard applications, $D$ is a differential operator and $G$ the corresponding Green
operator, taking into account the boundary conditions of the domain of $D$. }
$$\left\{\matrix{
DG = GD = \un
  \hfill\cr
\noalign{\medskip}
\left\langle{ Dx, y }\right\rangle =
\left\langle{ Dy, x }\right\rangle 
   \hfill \cr
\noalign{\medskip}
Q(x) = 
\left\langle{ Dx, x }\right\rangle \  ,
\  
W(x') =
\left\langle{ x', Gx' }\right\rangle \, .
  \hfill \cr}
\right. 
\eqno  ({\mathop{\rm I.6 }\nolimits})$$
Here $\left\langle{ x', x }\right\rangle $ denotes the duality between 
 ${\bf X}$ (elements $x$) and its dual $
{\bf X'}$ (elements $x'$).

\medskip
Then we define the integrator 
$\Da_{s, Q} x$ (also denoted $\Da x$ for simplicity) by the following requirement 
$$\int_{ {\bf X} }^{ } \Da_{s, Q} x \cdot 
\exp \!
\left({- { \pi \over 	s} Q(x)  - 2 \pi i 
\left\langle{ x', x }\right\rangle
}\right)  =
\exp 
 \! \left({ - \pi s W (x') }\right) \, .
\eqno  ({\mathop{\rm I.7 }\nolimits})$$

\medskip
\nn
Here $x'$ runs over ${\bf X'}$ and there are two cases:

\smallskip
i) $s = 1$ and  $Q$ is positive definite, namely $Q(x) > 0$ for $x \ne 0 $;

ii) $s = i$ and there is  no restriction on $Q$ except it be real.

\medskip
\nn
{\it Introduction of the parameter $s$ enables us to treat in a unified way the diffusion and
the Schr\"odinger equations}\footnote{$^5$}{ \nn
It is not true, as is still often stated, that the case $s = i$ has no mathematical
foundation. See for example references [8, 9, 10, 11].  }.
We can introduce a suitable space $\Fa ({\bf X})$ of functionals on ${\bf X}$ integrable  by $\Da_{s,
Q}$ and  a norm on $\Fa ({\bf X}) $, and then compute integrals of the type
$$I = \int_{ {\bf X} }^{ } \Da_{s, Q} x \cdot F(x) 
\eqno  ({\mathop{\rm I.8 }\nolimits})$$
for $F $ in $\Fa ({\bf X})$. In both cases,
$s = 1$ and $s = i$,  we shall call $\Da_{s, Q} x$ a {\it Gaussian integrator}.

\medskip

 Gaussian integrators have the following properties:
$$\Da (x + x_0) = \Da x, \qquad
x_0
\ \ 
{\mathop{\rm  a \ fixed \ element \ of }\nolimits}
\ \ {\bf X} \eqno  ({\mathop{\rm I.9 }\nolimits})$$
$$ \Da (L x) = | \Det L | \cdot \Da x, \qquad L: {\bf X} \ra {\bf X}
\eqno  ({\mathop{\rm I.10 }\nolimits})$$
where $L$ is in a suitable class of linear changes of variables of integration,
including the obvious case  where $Q(Lx) = Q(x) $ and 
$| \Det L | = 1$ (cf. formula (A.140) in Appendix A).

\medskip
Gaussian integrators are not the only possible integrators. 
In reference [11] we  have developed an axiomatic for functional integrals on a Banach space
$\Phi$ expressed in terms of integrators 
$\Da_{\T, Z}$ defined by 
$$\int_\Phi \Da_{\T, Z} \vp \cdot \T (\vp, J)
= Z(J) 
\eqno  ({\mathop{\rm I.11 }\nolimits})$$
for $\vp$ in $\Phi$, $J$ in the dual $\Phi'$ of $\Phi$,
where $\T $ and $Z$ are two given continuous bounded functionals
$$\T: \Phi \ts \Phi' \ra \CC \quad ,
\quad
Z: \Phi' \ra \CC.$$

\smallskip
\nn
In quantum field theory, we interpret $\vp$ as a field and $J$ as a source,
$Z(J)$ is then the Schwinger generating functional for the $n$-point functions.

\bigskip
\nn
{\bf 3. Integrands and integrators. }

\smallskip
Splitting the quantity inside the integral sign into ``integrator'' and
``integrand'' belongs to the art of integration, but rules of thumb apply:

\smallskip
-- When the functional integral has its origin in physics try not to break up the action into,
say, kinetic and potential contributions. On the other hand, do not hesitate to work with a
potential which is a functional of a path rather than a function of its {\it value} (e.g. in
equation (III.1), $V$ is a functional of $z$, not a function of $z(t)$).

\smallskip

-- Look for a possible change of variable of integration;  this may suggest
 a practical choice for the integrand.

\smallskip
-- Gaussian integrators have a wealth of simple, powerful properties;
 look for exponentials of quadratic forms,  this suggests a practical choice for
the integrator. \medskip

\bigskip
\nn
{\bf 4. A basic functional integral. }
\smallskip

The core of this article is a theorem which provides the mathematical
underpinning for a great variety of functional integrals. 
It consists of two parts: the definition of a functional integral, and the partial
differential equation satisfied by 
the {\it value} of the functional integral, as a function of a set of parameters. Given a
manifold $N$  consider the  $L^{2,1}$ paths over a finite time interval ${\bf T}$ with values
in $N$,  $x: {\bf T} \ra  N$, and with a fixed point 
$x_0 \in N$; for instance, if 
${\bf T} = \left\lbk{ t_a ,  t_b  }\right\rbk$, the fixed point
can be either 
$x(t_a) $ or $x(t_b)$. Given $d + 1$ vector fields $Y$ and 
$ X_{(\a)}  $ on $N$,  define a map $P$ from a space $\Pa_0 \RR^d $ of $L^{2, 1}$
paths into a space $\Pa_{x_0} N$ of $L^{2, 1} $ paths, 
$$P: \Pa_0 \RR^d \ra \Pa_{x_0} N,$$
by $P(z) = x $, where
$$
\left\{\matrix{
 dx (t, z)= 
 X_{(\a)}
(x (t, z)) dz^\a +
Y
(x(t, z)) dt 
\hfill\cr
\noalign{\medskip}
   x(t_0, z)  = x_0 .\hfill \cr}
\right. \eqno ({\mathop{\rm I.12 }\nolimits} )
$$
Here $t_0$ and $x_0$ are fixed, with $t_0$ in ${\bf T}$ and $x_0$ in $N$, and the paths
$z : {\bf T} \ra \RR^d$  satisfy $z \!\left({ t_0 }\right) = 0$. 
In general, the vector fields do not commute, that is: 
$$\left\lbk{  X_{(\a)} ,  X_{(\b)}  }\right\rbk
\ne 0 \quad , \quad \left\lbk{  Y ,  X_{(\a)}  }\right\rbk
\ne 0\, ;$$
therefore the solution of (I.12) is of the form
$$x(t, z) = x_0 \cdot \Si(t, z)
\eqno ({\mathop{\rm I.13 }\nolimits} )$$
where $x$ is a {\it function of}  $x_0$ {\it and  $t$, and} {\it a functional of}  $z$.
Here $\Si (t, z)$ is a transformation in $N$, depending on $t$ and $z$ as stated.
Only when 
$\left\lbk{  X_{(\a)} ,  X_{(\b)}  }\right\rbk = \left\lbk{  Y ,  X_{(\a)}  }\right\rbk = 0$ can one express
$\Si(t, z) $ as a function of $t$ and $z(t)$.

\medskip
In Appendix C the familiar theorems of differential equations are extended to cover the properties of (I.12).
\medskip

Given (I.12) one can express an integral over $\Pa_{x_0} N $ as an integral over
$\Pa_0 \RR^d$, which is an integral that one can manipulate and compute. The partial
differential equation satisfied by the integral on 
$\Pa_{x_0} N$ is expressed in terms of Lie derivatives along the vector fields $Y$ and
$  X_{(\a)}  $. Conversely, given a parabolic 
partial differential equation, one can construct in many cases the path integral representation of its solutions
(see for instance in Section IV, the lift of a covariant Laplacian at a point of  a Riemannian manifold in terms
of  Lie derivatives at a point of its frame bundle [13]). 

\medskip

The basic equations are given in the first
paragraph of Section II, followed by a summary of notations used throughout the paper.

\bigskip
\nn
{\bf 5. Examples.}

\smallskip
In Section III we compute semiclassical approximations of the basic functional integral
(II.1), and in Section IV we specialize the basic integral by making particular choices of the
manifold $N$.
 We treat
in detail the following examples:
\smallskip
-- $N = \RR^d$, in cartesian or polar coordinates. 
\smallskip
-- $N$ is a frame bundle $O(M) $ over a Riemannian manifold $M$.  We give
the explicit functional integral representing  the solution $\Psi$ of the Schr\"{o}dinger
equation on $M$, with initial wave function $\phi$.

\smallskip

-- $N$ is  a multiply-connected space.

\smallskip
-- $N$ is a $U(1)$-bundle; the basic integral solves the Schr\"odinger equation for a
particle in an electromagnetic field.

\smallskip
-- $N$ is a symplectic manifold; coherent-state transitions can be
obtained from the basic integral.

\smallskip
\nn
We conclude this section by an analysis of the Bohm-Aharonov effect, where all the previous techniques are
brought to bear.

\bigskip
\nn
{\bf 6. Techniques.}
\smallskip

In the course of computing our basic functional integral  in various situations, we have
used the properties of linear changes\footnote{\nn$^6$}{  \nn
Linear changes of variable of integration in an infinite-dimensional space are sufficiently
powerful and varied for the purposes of this paper. In a later publication we shall present
nonlinear changes, simplifying and generalizing earlier works such as [12].  } of
variable of integration and properties of functional determinants; they are given in two
appendices.  The transformation $\Si$ in formula (I.13) is related to
the Cartan development map; the properties of this map are discussed in the
third appendix.

\vfill\eject
\vglue 1cm
\nn
\centerline {\titre II - A general theorem}

\bigskip
\bigskip

The primary goal of this section is to define functional integrals over $L^{2,1}$ pointed paths
(paths with a fixed point) taking their values in a manifold $N$, by reducing them to functional
integrals over paths taking their values in a flat space $\RR^d$. The main definition is given as
follows
$$\left({ U_T \phi }\right)
\left({ x_0 }\right) : = \int_{ {\bf Z }_0 }^{ }
\Da_{s, Q_0 } z\cdot \exp 
\left({ - {\pi  \over s } Q_0 (z)  }\right)
\phi \left({ x_0 \cdot \Si (T, z) }\right) .
\eqno ({\mathop{\rm II.1 }\nolimits} )$$
All the notations are given in paragraph 1. The previous integral, being also denoted $\Psi \! \left({ T, x_0
}\right)$, is a solution of the {\it generalized Schr\"odinger equation }
$$
{ \part \Psi  \over \part T} =
{ s \over 4 \pi } h^{\a \b } 
\La_{X_{(\a) } }
\La_{X_{(\b) } } \Psi + \La_{Y }\Psi
\eqno ({\mathop{\rm II.2 }\nolimits} )$$
with initial condition $\Psi \! \left({ 0, x_0 }\right)
= \phi \left({ x_0 }\right)$.

\medskip
We give also a general construction of {\it time-ordered products} in the form of a functional
integral generalizing equation (II.1), namely
$$
\left({ U_T^F \phi  }\right)
\left({ x_0  }\right)
= \int_{ {\bf Z }_0J}^{ } \Da_{s, Q_0 } z \cdot 
\exp \left({ - {\pi  \over s } Q_0 (z)  }\right)
F ({\bf T} , z) \phi \left({ x_0 \cdot \Si (T, z)  }\right).
\eqno ({\mathop{\rm II.3 }\nolimits} )$$
In paragraph 2, we construct the simplest functional integral of type (II.1) and prove that it
satisfies equation (II.2). In paragraph 3, we study the general case.

\medskip
\nn
{\bf 1. The setup, and a summary of notations.}
\smallskip
\nn 1.1. {\it A manifold and vector fields.}

\smallskip
-- A finite-dimensional manifold $N$.

\smallskip
-- One-parameter groups acting on $N$, denoted $\s_{(\a)} (r) $; 
here $\a$ takes the values  $0, 1, 2, \cdots, d$ and $r$ is a real parameter; the transform of a
point $x$ in $N$ under $\s_{(\a) }  (r) $ is denoted by 
$x \cdot \s_{(\a) }  (r)$, and
$$\s_{(\a) }  (r) \circ \s_{(\a) } (s) = \s_{(\a) }
(r + s) \,  . \eqno ({\mathop{\rm II.4 }\nolimits} )$$

\smallskip
-- The generator of 
$\left\lbc{ \s_{(\a) } (r) }\right\rbc$ is the vector field 
$X_{(\a) } $ in $N$ such that 
$${ d \over dr } \left({ x \cdot  \s_{(\a) } (r)  }\right)
= X_{(\a) } \left({ x \cdot \s_{(\a) }  (r) }\right) 
\eqno ({\mathop{\rm II.5 }\nolimits} )$$
and in particular
$$X_{(\a) } (x) = 
\left.{ { d \over dr } \left({ x \cdot \s_{(\a) }  (r) }\right)  }\right\vert_{r = 0}
\eqno ({\mathop{\rm II.6 }\nolimits} )$$
for any point $x$ in $N$.
We do not assume that the vector fields 
$X_{(\a)}$ commute, hence
$\left\lbk{ X_{(\a)}, X_{(\b)} }\right\rbk \ne 0$
in general. We often write $Y$ for $X_{(0)}$ emphasizing its special role.

\smallskip

-- $\La_X$ denotes the Lie derivative w.r.t. the vector field $X$.

\bigskip
\nn 1.2. {\it Pointed paths on the flat space $\RR^d$.}

\smallskip

-- ${\bf T} $ is a time interval of length $T$, hence
$$ {\bf T } 
= \left\lbk{ t_a, t_b }\right\rbk \ ,
\qquad T = t_b - t_a \, .
\eqno ({\mathop{\rm II.7 }\nolimits} )$$ 

\smallskip

-- $t_0$ is a chosen element of ${\bf T}$; the standard choices are $t_0 = t_a$ or $t_0 = t_b$.

\smallskip

-- ${\bf Z}_0$ (or ${\bf Z}_{0, {\bf T}}$ if we need to specify ${\bf T}$)
consists of the real vector-valued functions $z = 
\left({ z^1 , \cdots,  z^d }\right)$ whose components 
$t \mps z^\a (t) $  are continuous functions with square-integrable derivatives $\dot z^\a$. We
assume the normalization  $z^\a \!\left({ t_0 }\right) = 0$.

\smallskip

-- ${\bf Z}_0'$ is the space dual to ${\bf Z}_0$. Its elements are interpreted as vector-valued
distributions
$z' = \left({ z'_1 , \cdots, z'_d }\right)$, each component being the derivative of an
$L^2$-function. The duality is given by
$$\left\langle{ z', z }\right\rangle =
\int_{ {\bf T }J}^{ } dt \, z'_\a (t) z^\a (t) 
\eqno ({\mathop{\rm II.8 }\nolimits} )$$ 
(summation over $\a$ and integration over $t$).

\smallskip
--
$s$ is a parameter equal to $1$ or $i$; its square root 
$\sqrt{ s }$ is normalized as follows
$$\sqrt{ s } = 
\left\{\matrix{
 1 \hfill& {\mathop{\rm  if }\nolimits} \hfill& s = 1
   \hfill \cr
\noalign{\medskip}
e^{\pi i / 4} \hfill& {\mathop{\rm  if }\nolimits} \hfill& s = i \, .
  \hfill \cr}
\right. \eqno ({\mathop{\rm II.9 }\nolimits} )$$ 

\smallskip
--
$h = \left({  h_{\a \b } }\right)$ is a constant invertible symmetric real matrix of size $d \ts d$. 
We denote by $\left({  h^{\a \b } }\right)$ the inverse matrix and we assume that $h$ is positive
definite in case $s$ is equal to $1$. By a suitable linear change of coordinates, we may take $h$
into a diagonal form
$$h = \diag (1, \cdots, 1, - 1, \cdots, - 1)
\eqno ({\mathop{\rm II.10 }\nolimits} )$$ 
with $p$ elements $+ 1$, $q$ elements $- 1$ and $p + q = d$. 
Hence $p = d$, $q = 0$ in the case $s = 1$.

\smallskip

-- We introduce a quadratic form $Q_0 $ on ${\bf Z}_0$ as follows
$$Q_0 (z) = \int_{ {\bf T }J}^{ }
dt \, h_{\a \b } \, 
\dot z^\a (t) \dot z^\b (t) \, .
\eqno ({\mathop{\rm II.11 }\nolimits} )$$
The corresponding kernel is given by 
$$D_{\a \b } (u, r) =  
{ {}_{1} \over {}^{2} } 
{\d^2 Q_0 (z)  \over \d z^\a (u) \d z^\b (r) } \, , 
\eqno ({\mathop{\rm II.12 }\nolimits} )$$ 
hence the representation 
$$Q_0 (z) = 
\int_{ {\bf T }   }^{ }du
\int_{ {\bf T }J}^{ }
 dr\, D_{\a \b } (u, r) \, z^\a(u) \, z^\b (r) \, .
\eqno ({\mathop{\rm II.13 }\nolimits} )$$

\smallskip

-- On ${\bf Z}'_0$ we consider a quadratic form $W_0$, with kernel 
$G^{\a \b } (u, r)$. As above, we have the relations 
$$
G^{\a \b } (u, r) = 
{ {}_{1} \over {}^{2} } 
{\d^2 W_0 (z')  \over \d z'_\a (u) \d z'_\b (r) } \, ,
\eqno ({\mathop{\rm II.14 }\nolimits} )$$ 
$$
W_0 (z') = \int_{ {\bf T}  }^{ } du
\int_{ {\bf T }J}^{ }
 dr\, G^{\a \b } (u, r) \, z'_\a (u)  z'_\b (r) \, .
\eqno ({\mathop{\rm II.15 }\nolimits} )$$

\smallskip

-- We assume that the quadratic forms
$Q_0 $ on ${\bf Z}_0$ and $W_0$ on ${\bf Z}'_0$ are inverse to each other in the sense of relation
(I.6). In terms of kernels, this is expressed as follows
$$\int_{ {\bf T }J}^{ } dt \, D_{\a \b }
\left({  s_1, t }\right)
G^{\b \g }
\left({ t, s_2  }\right)
= \d_\a^\g \, 
\d \left({ s_1 - s_2  }\right).
\eqno ({\mathop{\rm II.16 }\nolimits} )$$
Here are explicit formulas for the kernels:
$$ D_{\a \b } (u, r) = 
\int_{ {\bf T }J}^{ } dt \, h_{\a \b } \, 
\d' (t - u) \d' (t - r) = - h_{\a \b } \, 
\d'' (u - r) \, , 
\eqno ({\mathop{\rm II.17 }\nolimits} )$$
that is $D : {\bf Z }_0 \ra {\bf Z }'_0 $ is the differential operator with matrix 
$\left({ - h_{\a \b } { d^2 \over dt^2} }\right) $.
Hence, 
$G^{\a \b }(u, r)$ is the corresponding Green's function  taking into account the boundary condition 
$z^\a \!\left({ t_0  }\right) = 0$, that is
$$G^{\a \b } (u, r) 
=
\left\{\matrix{ h^{\a \b }
\inf 
\left({ u - t_0, r - t_0 }\right)
\hfill & {\mathop{\rm for \  }\nolimits} u \geq t_0, r \geq t_0 ,
  \hfill\cr
\noalign{\medskip}
h^{\a \b } \inf \left({ t_0 - u , t_0 - r  }\right) \hfill &{\mathop{\rm for \  }\nolimits}
 u \leq t_0, r \leq t_0 ,
   \hfill \cr
\noalign{\medskip}
0 \hfill &
{\mathop{\rm otherwise. }\nolimits}
 \hfill \cr}
\right. \eqno ({\mathop{\rm II.18 }\nolimits} )$$

\medskip
\nn
{\it Remark.}
We shall refrain from integrating by parts in formulas like (II.11) in order not to have to make
explicit statements about boundary conditions.

\bigskip
\nn
1.3. {\it  Functional integrals on $\Pa_0 \RR^d $. }

\smallskip
 
-- The space of paths ${\bf Z}_0$ shall also be denoted by 
$\Pa_0 \RR^d$ to remind us of the flat space $\RR^d$ where the paths lie and of the fixed point $0$
(origin in $\RR^d$) of the paths $z$. We shall consider a variety of pairs $(Q, W)$ consisting of
a quadratic form $Q$ on a space ${\bf Z}$ and a quadratic form $W$ on its dual ${\bf Z}'$ satisfying the
analogues of relations (II.12), (II.14) and (II.16).

\smallskip

-- For such a pair $(Q, W)$ we have a translation invariant integrator 
$\Da_{s, Q} z$ on ${\bf Z}$ characterized by 
$$
\int_{ {\bf Z }J}^{ }
\Da_{s, Q} z \cdot 
\exp 
\left({ -{ \pi  \over s}  Q(z) - 2 \pi i \left\langle{ z', z }\right\rangle }\right)
=
\exp 
(-  \pi s W (z')) 
\eqno ({\mathop{\rm II.19 }\nolimits} )$$ 
for $z'$ in ${\bf Z'}$.
The normalizations are chosen so that Gaussian integrators and Fourier transforms do not include
powers of $\pi$ depending on the dimension of their domain of definition.

\medskip
We also write
$$\Da \om_{s, Q} (z) = \Da_{s, Q} z \cdot 
\exp 
\left({ - { \pi  \over s}  Q(z)  }\right)
 \eqno ({\mathop{\rm II.20 }\nolimits} )$$ 
and refer to both 
$\Da_{s, Q}$ and $\Da \om_{s, Q} $ as 
``Gaussian integrators''. When working with the basic pair
$ \left({ Q_0, W_0  }\right) $ we simply write 
$\Da_s$ and 
$\Da \om_s$. In the context of an application where $s$ has been chosen once for all, we omit it
in the notations.

\bigskip
\nn
1.4. {\it  Functional integrals on $\Pa_{x_0} N $. }

\smallskip

-- We fix a point $x_0$ in $N$ and consider continuous paths 
$x : {\bf T} \ra N$ with the fixed point $x \!\left({ t_0 }\right) = x_0$ and square-integrable
velocity\footnote{$^7$}{More precisely, for every smooth function $f$ in $C^\Inf (N)$, we assume
that the continuous function $t \mps f(x(t))$ on ${\bf T}$ is the primitive of a function  in 
$L^2 ({\bf T })$. }. The set of all such paths is denoted by
$\Pa_{x_0} N$.

\smallskip

-- The time interval ${\bf T}$ being given, consider an element $z$ of ${\bf Z}_0$. As we shall see
in Appendix C, the differential equation 
$$
\left\{\matrix{
  dx(t) = X_{(\a)}
( x(t) ) dz^\a (t) + Y (x(t)) dt 
   \hfill \cr
\noalign{\medskip}
x \!\left({ t_0 }\right) = x_0
  \hfill \cr}
\right. 
\eqno ({\mathop{\rm II.21 }\nolimits} )$$
admits a unique solution $x(\cdot)$ in 
$\Pa_{x_0} N$.

\smallskip

-- The previous construction defines a parametrization $P$ of the space
$\Pa_{x_0} N$ of pointed paths on $N$ by the space $\Pa_0 \RR^d$ of pointed paths in $\RR^d$, namely
$$P : \Pa_0 \RR^d \ra \Pa_{x_0} N$$
by taking $z$ into $x$. If necessary we shall denote by $x(t, z)$ the solution of the differential
equation (II.21) for given $z$ in $\Pa_0 \RR^d$. Hence $x(t, z)$ is {\it a function of $t$ and a
functional of $z$}.

\smallskip

--  Assume now that ${\bf T} = \left\lbk{ 0, T }\right\rbk$ and $t_0 = 0$.
 With the previous definitions, define 
$\Si (T, z)$ as the transformation taking a point $x_0$ in $N$ into the point
$x(T, z)$. 

\smallskip

-- Take a good\footnote{$^8$}{Any function in $C^\Inf (N)$ with compact support will do.}
 function $\phi$ on $N$. Define a functional $\ww \phi$ on ${\bf Z}_0$ by 
$$\ww \phi (z) = \phi
\left({  x_0 \cdot \Si (T, z)}\right) .
\eqno ({\mathop{\rm II.22 }\nolimits} )$$
It is integrable under the integrator $\Da \om_s$ on ${\bf Z}_0$.
By integrating we get
$$
\eqalign{ I \! \left({  \phi, T,  x_0}\right)
&= \int_{ {\bf Z}_0  }^{ }
\Da \om_s (z)  \ww \phi (z) \cr
& = \int_{ {\bf Z}_0  }^{ }
\Da_s z \cdot \exp 
\left({ - {\pi  \over s }  Q_0 (z) }\right)
\cdot \phi 
\left({ x_0 \cdot \Si (T, z)  }\right) .
\cr}
\eqno ({\mathop{\rm II.23 }\nolimits} )$$
We could replace $Q_0$ by another quadratic form on ${\bf Z}_0$.

\smallskip

-- 
The functional operator $U_T$ associates to the function $\phi$ on $N $ the function 
$x_0 \mps I \!\left({  \phi, T, x_0 }\right)$ on $N$.

\bigskip
\nn
1.5. {\it  Variational techniques.}
\smallskip

A frequently used technique consists in introducing one-parameter variations in the space of paths,
or in the space of functionals on the space of paths. We work in particular with the following
variations.

\smallskip
(i)  Fix a time interval 
${\bf T} = \left\lbk{ t_a, t_b  }\right\rbk$ and consider the space 
$\Pa \RR^d$ of paths $x : {\bf T} \ra \RR^d$, say of class\footnote{$^9$}{That is, continuous with
continuous first-order derivatives.} $C^1$. 
Assuming an action functional $S : \Pa \RR^d \ra \RR$ (for instance, the time integral of a
Lagrangian 
$L(x, \dot x, t)$), the critical points  of $S$ will form a $2d$-dimensional manifold 
$\Ua^{2d}$, the so-called {\it space of } (classical) {\it motions}, denoted by 
$x_{{\mathop{\rm cl }\nolimits} }$. We can parametrize them by a set of parameters $\mu =
\left({ \mu^1, \cdots, \mu^{2d}  }\right)$; for a given $k$ between $1$ and $2d$, the derivative 
$\part x_{{\mathop{\rm cl }\nolimits} }
 (\mu) / \part \mu^k$ defines a variation through classical paths, and we get $2d$ of
them.

\smallskip

(ii)
If $t_0$ is a given epoch in ${\bf T}$, the pointed paths are defined by the condition 
$x \!\left({ t_0 }\right) = 0$. They form the space $\Pa_0 \RR^d$. In paragraph III.1 we shall use a
one-parameter family of pointed paths $x (\l)$ in $\Pa_0 \RR^d$ (for $\l$ running over $[0, 1]$).

\smallskip

(iii) Introducing $d$ independent boundary conditions at $t_a$, and similarly at $t_b$, we define the
subspace
$\Pa_{a, b} \RR^d $ of $\Pa \RR^d$.
In paragraph B.2, we consider a one-parameter family of paths 
$x(\l) $ in $\Pa_{a, b} \RR^d$, with $\l$ in $[0, 1]$, such that $x(0)$ belongs to 
$\Ua^{2d} \cap \Pa_{a, b} \RR^d$. The presence and nature of the caustics is conveniently analyzed in
terms of this intersection.

\smallskip

(iv) In paragraph  III.2  we consider a one-parameter family of
action functionals $S(\nu)$ defined on a space $\Pa_0 \RR^d $ of pointed paths. 

\bigskip
\nn
{\bf 2. The one-dimensional case.}
\smallskip

We begin by the simple case $d = 1$. Here $X = X_{(1)}$ is a vector field on $N$ generating the one-parameter
group of transformation $\s(r)$ on $N$. Hence $\s(r)$ obeys the differential equation
$$d \! \left({ x_0 \cdot \s (r)  }\right) =
X \! \left({ x_0 \cdot \s (r)  }\right) dr 
\eqno ({\mathop{\rm II.24 }\nolimits})$$
for a fixed $x_0$ in $N$. We set $Y  = 0$ and consider the time interval
${\bf T}  = [0, T]$. The differential equation (II.21) reads now as 
$$
\left\{\matrix{
dx (t) = X (x(t)) \cdot dz (t)
  \hfill\cr
\noalign{\medskip}
x(0) = x_0     \, . 
  \hfill \cr}
\right. \eqno ({\mathop{\rm II.25 }\nolimits})$$
Substituting $r = z(t)$ in (II.24) we see that the solution to the previous equation is given by $x(t) =
x_0 \cdot \s (z(t))$. Hence the transformation $\Si (T, z)$ is simply $\s (z(T))$.

\bigskip
The path space ${\bf Z}_0$ consists of the $L^{2,1}$ functions 
$z : [0, T] \ra \RR$ such that $z(0) = 0$. It is endowed with the quadratic form
$$Q_0 (z) = \int_{0 }^{T } dt \, \dot z (t)^2 \, . 
\eqno ({\mathop{\rm II.26 }\nolimits})$$
The corresponding integrator is, for simplicity, denoted by $\Da_s z$, hence our basic path integral
specializes to
$$\Psi (T, x) =
\int_{ {\bf Z}_0 }^{ }
\Da_s z \cdot \exp 
\left({- {\pi  \over s } \int_{ 0}^{ T}  dt \, \dot z (t)^2}\right)
\phi (x \cdot \s (z(T))) \, . 
\eqno ({\mathop{\rm II.27 }\nolimits})$$
Our goal is to show that we have solved the differential equation
$$
{\part  \over \part T}
\Psi (T, x) = {s \over 4 \pi }
\La_X^2 \Psi (T, x) \, .
\eqno ({\mathop{\rm II.28 }\nolimits})$$

\bigskip
Fix a point $x$ in $N$ and define a function of a real variable 
$h(r) = $ \break  $\phi (x \cdot \s (r))$. From the group property 
$\s(r_1) \circ \s(r_2) = \s(r_1 + r_2)$, we get
$$\phi (x \cdot \s (r) \cdot \s (z(T))) =
h (r + z (T)) \, . 
\eqno ({\mathop{\rm II.29 }\nolimits})$$
Then define
$$H  (T, r) : = \int_{ {\bf Z}_0 }^{ } \Da_s z \cdot 
e^{- \pi Q_0 (z) / s }
h (r + z (T)) 
 \, . 
\eqno ({\mathop{\rm II.30 }\nolimits})$$

\medskip
\nn
Since the integrand $h(r + z(T))$ depends on the path $z$ through $z(T)$, a  linear change of variable $z \mps
z(T)$ transforms immediately this functional integral over ${\bf Z}_0$ into an ordinary integral over $\RR$
(see formula (A.38) in Appendix A). It is easier to use directly the properties of Fourier transforms of
Gaussian integrators underpinning (A.38). Denoting by $\hh h$ the Fourier transform of $h$, we obtain
$$h(r + z(T)) = \int_{ \RR }^{ }
d \r \, \hh h(\r) \, e^{2 \pi i \r (r + z (T))}
 \, . 
\eqno ({\mathop{\rm II.31 }\nolimits})$$
Therefore, after changing the order of integration, we get
$$H(T, r) = \int_{ \RR }^{ }
d \r \, \hh h(\r) \, e^{2 \pi i \r r}
\int_{ {\bf Z}_0 }^{ } \Da_s z \cdot \exp \!
\left({- {\pi  \over s } Q_0 (z) + 2 \pi i \r z (T)  }\right)
 \, . 
\eqno ({\mathop{\rm II.32 }\nolimits})$$
By the definition (II.19) and the duality 
$\r z (T) = 
\left\langle{ \r \d_T, z }\right\rangle$,
we get
$$\int_{ {\bf Z}_0 }^{ } \Da_s z \cdot \exp \! 
\left({- {\pi  \over s } Q_0 (z) + 2 \pi i \r z (T)  }\right)
=
\exp \! \left({ - \pi s W_0 (\r \d_T) }\right)
\, . 
\eqno ({\mathop{\rm II.33}\nolimits})$$

\medskip
\nn
By (II.15) and (II.18) we obtain
$$ W_0(\r \d_T) = \r^2 \, G (T, T) = \r^2 T
\, . 
\eqno ({\mathop{\rm II.34}\nolimits})$$
Collecting equations (II.32) to (II.34) we conclude
$$H (T, r) = \int_{ \RR }^{ }
d \r \, \hh h (\r) e^{2 \pi i \r r -  \pi s \r^2 T}
\eqno ({\mathop{\rm II.35}\nolimits})$$
and by derivation under this integral sign, we conclude
$$
{\part H \over \part T} 
= { s \over 4 \pi }  
{\part^2 H \over \part r^2} 
\, . 
\eqno ({\mathop{\rm II.36}\nolimits})$$

\bigskip
The vector field $X$ is the generator of the one-parameter group of transformations
$\s(r)$.
Hence, for every integer $m \geq 0$, we get
$$
\left({ {\part  \over \part r } }\right)_{}^{m} 
h (r + z(T)) =
\left({ \La_X^m \phi }\right)
(x \cdot \s (r + z (T))) .
\eqno ({\mathop{\rm II.37 }\nolimits})$$

\smallskip
\nn
By differentiating under the integral sign in formula (II.30), we deduce
$$\left({ {\part  \over \part r } }\right)_{}^{m} 
H(T, r) =
\int_{ {\bf Z}_0 }^{ }
\Da_{s} z \cdot \exp \! 
\left({ -  {\pi \over s}  Q_0 (z)  }\right) 
 \La_X^m \phi
(x \cdot \s (r + z (T))) \, .
\eqno ({\mathop{\rm II.38 }\nolimits})$$
Using definition (II.27), we get therefore
$$
\left.{ \left({ {\part  \over \part r } }\right)_{}^{m} 
H(T, r)  }\right\vert_{r = 0 } = 
 \La_X^m \Psi
(T, x)
\eqno ({\mathop{\rm II.39 }\nolimits})$$
by applying $m$ times the differential operator 
$\La_X$ acting on functions of $x$. In particular, for $m = 0$, we get
$$
\left.{ H(T, r)  }\right\vert_{r = 0}
=
\Psi (T, x)
\eqno ({\mathop{\rm II.40 }\nolimits})$$
and both sides can be derivated with respect to $T$, giving
$$
\left.{  {\part  \over \part T }  H (T, r)  }\right\vert_{r = 0}
= {\part  \over \part T }
\Psi (T, x) .
\eqno ({\mathop{\rm II.41 }\nolimits})$$
Setting $r = 0$ in the relation 
$\ds {\part H \over \part T } =
{s  \over 4 \pi }
{\part^2  \over \part r^2 } H$ and using relation 
(II.39) for $m = 2$, and relation (II.41), we conclude the proof of the differential
equation 
$${\part  \over \part T }
\Psi (T, x) =
{s \over 4 \pi }
\La_X^2 \Psi (T, x) .$$

\bigskip
\nn
{\it Remark. }
The previous proof can be recast in the following symbolical way [5].
We begin with a consequence of equations (II.33) and (II.34), namely
$$\int_{ {\bf Z}_0 }^{ }
\Da_{s} z \cdot \exp \!
\left({ -  {\pi \over s}  Q_0 (z) }\right)  \exp  ( 2 \pi i \r z (T) )
= \exp ( - \pi s \r^2 T )\, ,
\eqno ({\mathop{\rm II.42 }\nolimits})$$
one of many characterizations of our integrator 
$\Da_s z$. Make the formal substitution
$$\r = 
{1 \over 2 \pi i}
\La_X,$$
and get
$$\int_{ {\bf Z}_0 }^{ }
\Da_{s} z \cdot \exp \!
\left({ -  {\pi \over s}  Q_0 (z) }\right) 
\exp (z (T) \La_X) =
\exp \! \left({   {sT  \over 4 \pi }   \La_X^2 }\right) .
\eqno ({\mathop{\rm II.43 }\nolimits})$$

\smallskip
\nn
Apply this operator identity to the function $\phi (x)$ and notice that the operator 
$\exp (r \La_X)$ transforms $\phi (x) $ into 
$\phi (x \cdot \s (r))$.
Hence we get
$$\int_{ {\bf Z}_0 }^{ }
\Da_{s} z \cdot \exp \!
\left({ -  {\pi \over s}  Q_0 (z) }\right) 
\phi (x \cdot \s (z (T)) =
\left({  \exp  {sT  \over 4 \pi }   \La_X^2 }\right)
\phi (x) 
\eqno ({\mathop{\rm II.44 }\nolimits})$$
that is 
$$\Psi (T, x) = 
\exp \!
\left({   {sT  \over 4 \pi }   \La_X^2 }\right)
\phi (x) .
\eqno ({\mathop{\rm II.45 }\nolimits})$$
This integrated form is equivalent to the differential equation (II.28) together with the
initial condition 
$$\Psi (0, x) = \phi (x) \, . 
\eqno ({\mathop{\rm II.46 }\nolimits})$$

\bigskip
\nn
{\bf 3. The general case.}

\smallskip
\nn
3.1. {\it The group property.}
\smallskip
The functional operator $U_T$ defined by formula (II.1) satisfies the group property
$$U_T \circ U_{T'} = U_{T + T'}
\eqno ({\mathop{\rm II.47 }\nolimits})$$
for $T > 0$, $T' > 0$. The proof rests on three facts.

\medskip
a) {\it Group property for the point-transformations} $\Si (T, z)$:

The relevant function space ${\bf Z}_{0, T}$ consists of the paths 
$z : [0, T] \ra \RR^d$ with $L^2$ derivative.
For $z$ in ${\bf Z}_{0, T}$ and $z'$ in ${\bf Z}_{0, T'}$, we define a new path $z \ts z'$ by the rule
$$(z \ts z') (t) = 
\left\{\matrix{
z(t) \hfill & {\mathop{\rm for  }\nolimits} &
0 \leq t \leq T
  \hfill\cr
\noalign{\medskip}
   z(T) + z'(t- T) \hfill & {\mathop{\rm for  }\nolimits} &
T \leq t \leq T + T' \, .
  \hfill \cr}
\right. 
\eqno ({\mathop{\rm II.48 }\nolimits})$$
It is obvious that $z \ts z'$ is an element of ${\bf Z}_{0, T + T'}$.
Furthermore, the map $(z, z') \mps z \ts z'$ is an isomorphism of 
${\bf Z}_{0, T} \ts {\bf Z}_{0, T'}$
onto ${\bf Z}_{0, T+T'}$ and,  by the uniqueness of the solution of the differential equation (II.21), we obtain
$$x_0 \cdot \Si (T + T', z \ts z') = x_0 \cdot \Si (T, z) \cdot \Si (T', z') 
\eqno ({\mathop{\rm II.49 }\nolimits})$$
(see also formula (C.23) in Appendix C).

\medskip
b) {\it Quadratic forms}:

We denote by $Q_{0, T}$ the basic quadratic form on ${\bf Z}_{0, T}$ given by equation (II.11), namely:
$$Q_{0,T}(z) = \int_{0 }^{ T} dt \, 
h_{\a \b} \dot z^\a (t) \dot z^\b (t) \, .
\eqno ({\mathop{\rm II.50 }\nolimits})$$
Using similar definitions for $Q_{0, T'}$ and $Q_{0, T + T'}$, one obtains immediately
$$Q_{0, T + T'} (z \ts z') =
Q_{0,T} (z) + Q_{0,T'} (z')
\eqno ({\mathop{\rm II.51 }\nolimits})$$
and by exponentiating
$$\exp \! 
\left({ - { \pi \over s } Q_{0,T+T'} (z \ts z') }\right)
=
\exp \! 
\left({ - { \pi \over s } Q_{0,T} (z ) }\right)
\exp \! 
\left({ - { \pi \over s } Q_{0,T'} (z' ) }\right) \, . 
\eqno ({\mathop{\rm II.52 }\nolimits})$$

\medskip
We can identify 
${\bf Z}_{0, T+T'}$ to ${\bf Z}_{0, T} \ts {\bf Z}_{0, T'}$. This identification entails an identification of
the dual space ${\bf Z}'_{0, T+T'}$ to the product space 
${\bf Z}'_{0, T} \ts {\bf Z}'_{0, T'}$.
We use the notations $\z$ for elements of ${\bf Z}_{0, T}$, $\z'$ in ${\bf Z}_{0, T'}$
and denote by $\z \ts \z' $ the corresponding element in 
${\bf Z}'_{0, T+T'}$.
Using simple algebra, one derives the identity
$$W_{0, T + T'} (\z \ts \z') 
= W_{0, T} (\z) + W_{0, T'} (\z') 
\eqno ({\mathop{\rm II.53 }\nolimits})$$
from (II.51).
Here $W_{0,T}$ denotes the quadratic form on ${\bf Z}'_{0, T}$ inverse to 
$Q_{0,T}$, etc.

\medskip
c) {\it Integrators}:

The following form of Fubini's theorem holds
$$
\int_{ {\bf Z}_{0, T+T'} }^{ }
\Da_s u \cdot F(u) = 
\int_{ {\bf Z}_{0, T} }^{ } \Da_s z
\int_{ {\bf Z}_{0, T'} }^{ } \Da_s z' 
\cdot F(z \ts z') \, . 
\eqno ({\mathop{\rm II.54 }\nolimits})$$
According to the general method explained in Appendix A, 
it is enough to check this formula for a function of the form
$$F (z \ts z') =
\exp \!
\left({ - 2 \pi i \!
\left({ \left\langle{ \z, z }\right\rangle + \left\langle{ \z', z' }\right\rangle  }\right)
 }\right)
\, . \eqno ({\mathop{\rm II.55 }\nolimits})$$
But then our contention follows from the characterization (II.19) of the integrator, and the relations (II.51)
and (II.53).

\medskip
We combine this formula with equation (II.52) and obtain another form of Fubini's theorem
$$\int_{ {\bf Z}_{0, T+T'} }^{ }
\Da \om_s (u) \, G(u) =
\int_{ {\bf Z}_{0, T} }^{ }
\Da \om_s  (z)  
\int_{ {\bf Z}_{0, T'} }^{ }
\Da \om_s (z') \, G (z \ts z') \, . 
\eqno ({\mathop{\rm II.56 }\nolimits})$$

\medskip
We conclude the proof of equation (II.47). Indeed
$$
\eqalign{
\left({ U_T \! \left({ U_{T'} \phi }\right) }\right)
\left({x_0 }\right)
&= \int_{ {\bf Z}_{0, T} }^{ }
\Da \om_s (z) \left({ U_{T'} \phi }\right)
\left({x_0 \cdot \Si (T, z) }\right)\cr
&=
\int_{ {\bf Z}_{0, T} }^{ }
\Da \om_s (z)
\int_{ {\bf Z}_{0, T'} }^{ }
\Da \om_s (z')
\phi \! \left({ x_0 \cdot \Si (T, z) \cdot \Si (T', z') }\right)\cr
&=
\int_{ {\bf Z}_{0, T} }^{ }
\Da \om_s (z)
\int_{ {\bf Z}_{0, T'} }^{ }
\Da \om_s (z')
\phi \! \left({ x_0 \cdot \Si (T+ T' , z \ts z') }\right)\cr
&=
\int_{ {\bf Z}_{0, T + T'} }^{ }
\Da \om_s (u)
\phi \! \left({ x_0 \cdot \Si (T+ T' ,u) }\right)\cr
&
=
\left({ U_{T + T'} \phi }\right)   \left({x_0 }\right) \, .\cr}$$

\vfill\eject
\nn
3.2. {\it The differential equation.}

\smallskip
>From the group property (II.47) it follows that we need to establish the partial differential equation (II.2)
at the time $T = 0$, and that the general case will follow. That is, we want to prove
$$
\left({ U_T \phi }\right)
\left({ x_0 }\right) =
\phi \! \left({ x_0 }\right) + 
T \! \left({ {s \over 4 \pi} h^{\a \b} \La_{X_{(\a)}} \La_{X_{(\b)}} \phi 
\! \left({ x_0 }\right)
+ \La_Y \phi \! \left({ x_0 }\right) }\right)
 + o(T) 
\, . 
\eqno ({\mathop{\rm II.57 }\nolimits})$$
This relation implies also the initial condition
$$\build {\lim }_{ T = 0 }^{ }
\! \left({ U_T \phi }\right) \! \left({ x_0 }\right) =
\phi \! \left({ x_0 }\right)
\, . 
\eqno ({\mathop{\rm II.58 }\nolimits})$$

\medskip
For the proof, we rely on the scaling properties of paths as described in paragraph A.3.6. For $z$ in 
${\bf Z}_{0,1}$ we define the scaled path $z_T $ in ${\bf Z}_{0,T}$ by
$$z_T (t) = T^{ 1/2} z (t/T) 
\, . 
\eqno ({\mathop{\rm II.59 }\nolimits})$$
By scaling the differential equation (II.21), we get
$$dx (t/T) =
T^{1/2} X_{(\a)} (x (t/T)) 
dz_T^\a (t) + TY 
( x (t/T)) d (t/T) 
\, . 
\eqno ({\mathop{\rm II.60 }\nolimits})$$
Hence the transformation $\Si (T, z_T)$ in $N$ defined using the vector fields $X_{(\a)}$ and $Y$ is the same
as the transformation $\Si (1, z)$ using the vector fields $T^{1/2} X_{(\a)}$ and $TY$.
Moreover we can transfer the path integration from the variable space 
${\bf Z}_{0,T}$ to the fixed space ${\bf Z}_{0,1}$, that is 
$$\left({ U_T \phi }\right) \! \left({ x_0 }\right)  =
\int_{ {\bf Z}_{0,1} }^{ } \Da \om_s (z) \, \phi \! 
\left({ x_0 \cdot \Si \! \left({ T, z_T }\right) }\right)
\, . 
\eqno ({\mathop{\rm II.61 }\nolimits})$$
We use then a limited expansion of $\phi$ around $\phi (x_0)$, namely
$$
\eqalign{
\phi \! 
\left({ x_0 \cdot \Si \! \left({ T, z_T }\right) }\right)
&
= \phi (x_0) + T^{1/2} \La_{X_{(\a)}} \phi (x_0) z^\a (1) \cr
&
\quad +  {T \over 2}  \La_{X_{(\a)}}   \La_{X_{(\b)}} \phi (x_0) 
z^\a (1) z^\b (1) \cr
& 
\quad + T \La_Y \phi (x_0) + O \! \left({ T^{3/2} }\right)\, .\cr}
\eqno ({\mathop{\rm II.62 }\nolimits})$$

\smallskip
\nn
Recall the integration formulas
$$\int_{ {\bf Z}_{0,1} }^{ } \Da \om_s (z) \, z^\a (1) = 0
\eqno ({\mathop{\rm II.63 }\nolimits})$$
 and
$$
\int_{ {\bf Z}_{0,1} }^{ } \Da \om_s (z) \, z^\a (1) z^\b (1)  = sh^{\a \b }/2 \pi \, .
\eqno ({\mathop{\rm II.64 }\nolimits})$$
Collecting equations (II.61) to (II.64), we conclude the proof of (II.57).

\medskip
The previous calculation can be extended to give the complete Taylor expansion of 
$\left({ U_T \phi }\right) (x_0) $ around $T = 0$.
Since the vector fields $X_{(\a)}$ and $Y$ do not commute, we have to use a time-ordered exponential to express
the solution of the differential equation (II.21).

\bigskip
\nn
3.3. {\it About the construction of } $\Si (T, z)$.
\smallskip
For simplicity, assume $Y = 0 $ and $h_{\a \b} = \d_{\a \b }$.
Consider again a path $z : [0, T] \ra \RR^d $ of class $L^{2,1}$ with 
$z(0) = 0$.
For each $t$ in $[0, T]$, denote by $z^1 (t) , \cdots, z^d (t)$ the components of the vector $z(t)$ and define
the transformation 
${\bf s} (z(t)) $ on $N$ by
$$x_0 \cdot {\bf s} (z(t)) =
x_0 \cdot \s_{(1)} (z^1 (t)) \cdot \ldots \cdot \s_{(d)} (z^d (t)) \, . 
\eqno ({\mathop{\rm II.65 }\nolimits})$$
Here $\left\lbc{ \s_{(\a)} (r) }\right\rbc$ denotes the one-parameter group of transformations in $N$ with
generator $X_{(\a)}$.
In general, the vector fields $X_{(\a)}$ do not commute, hence the transformations 
${\bf s}(z(t))$ do not form a group.
In the general case, we can use a {\it multistep method} to solve the differential equation (II.21), hence
$$x_0 \cdot \Si (T, z) =
\build { \lim }_{n = \Inf }^{ }
x_0 \cdot {\bf s} (z(T/n)) \cdot {\bf s} (z (2T/n) - z(T/n))
\cdot \ldots \cdot {\bf s} (z(T) - z(T - T/n)) \, . 
\eqno ({\mathop{\rm II.66 }\nolimits})$$
Putting this into the integral (II.1), we obtain after some calculations the following variant of the Lie,
Trotter, Kato, Nelson formula
$$U_T \phi = \build { \lim }_{n = \Inf }^{ }
\left({ U_{T/n}^{(1)} \cdots U_{T/n}^{(d)}}\right)^n \phi \, . 
\eqno ({\mathop{\rm II.67 }\nolimits})$$

\nn
Here $U_T^{(\a)}$ for $\a$ in $\{ 1, \cdots, d \}$ corresponds to the one-dimensional case studied in paragraph
II.2, hence
$$ U_T^{(\a)} = \exp \! 
\left({ { sT\over 4 \pi } \La_{X_{(\a)}}^2 }\right)\, . 
\eqno ({\mathop{\rm II.68 }\nolimits})$$
These relations are in agreement with
$$ U_T = \exp \! 
\left({ { sT\over 4 \pi } \sum_{ \a }^{ }\La_{X_{(\a)}}^2 }\right)\, . 
\eqno ({\mathop{\rm II.69 }\nolimits})$$

\bigskip
\nn
{\bf 4. Some generalizations.}
\smallskip
\nn
4.1. {\it Including a potential.}
\smallskip

Consider the Schr\"odinger equation with potential
$${ \part \Psi\over \part t}
= {s \over 4 \pi } h^{\a \b }\, \La_{X_{(\a)}} \La_{X_{(\b)}}\Psi + V \Psi \, .
\eqno ({\mathop{\rm II.70 }\nolimits})$$

\nn
A path integral solution is obtained as follows
$$\Psi \! \left({ T, x_0 }\right)
=
\int_{ {\bf Z}_{0,T} }^{ } \Da_sz \cdot \exp \!  
\left({ -{\pi \over s} Q_0 (z) + \int_{{\bf T} }^{ }
dt \, V\! \left({ x_0 \cdot \Si (t, z) }\right)
 }\right) \phi \! \left({ x_0 \cdot \Si (T, z) }\right)
\,
.
\eqno ({\mathop{\rm II.71 }\nolimits})$$

\nn
The proof can be obtained by a slight generalization of the arguments presented in paragraph II.3.
We can also use the following trick:
we add one variable $\T$, considering the manifold $N \ts \RR$. The vector fields $X_{(\a)}$
in
 $N$ give vector fields, also denoted by $X_{\a}$, in 
$N \ts \RR$, which have a zero component on the factor $\RR$.
Moreover 
$Y = V(x) \part / \part \T$.
The function $\Psi (T, x) $ satisfies the equation (II.70) if, and only if, the function $\Psi (T, x)
\exp \T$ satisfies the equation (II.2) on $N \ts \RR$.
The differential system (II.21) is now written as
$$
\left\{\matrix{
dx = X_{(\a) } (x) d z^\a
  \hfill\cr
\noalign{\medskip}
    d \T = V(x) dt \, . 
\hfill \cr}
\right. 
\eqno ({\mathop{\rm II.72 }\nolimits})$$
Hence the transformation $\Si (T, z) $ takes $
\left({ x_0, \T_0 }\right)$ into 
$$\left({  x_0 \cdot \Si (T, z) , \ 
\T_0 + \int_{ {\bf T} }^{ } dt \  V 
(x_0 \cdot \Si (t, z))
}\right)\, ,$$
hence
$$
\phi \! \left({ (x_0, 0) \cdot \Si (T, z) }\right)
= \phi \! \left({ x_0 \cdot \Si (T, z) }\right)
\exp \! \left({ \int_{ {\bf T}}^{ } dt \, 
V \! \left({ x_0 \cdot \Si (t, z) }\right)
 }\right)
\, .\eqno ({\mathop{\rm II.73 }\nolimits})$$
Equation (II.71) follows easily from these remarks.
Notice that the exponent in this formula does not have a simple dynamical interpretation in general -- but
see our illustrations in section IV.

\bigskip
\nn
4.2. {\it Time-ordered product.}

\smallskip
Suppose that $F$ is a suitable functional on the space ${\bf Z}_{0,T}$.
We generalize equation (II.1) as follows
$$\left({ U_T^F \phi }\right)(x_0) =
\int_{ {\bf Z}_{0, T} }^{ } \Da \om_s (z) \, F(z) \, \phi  
\! \left({ x_0 \cdot \Si (T, z) }\right)
\, . \eqno ({\mathop{\rm II.74 }\nolimits})$$
By imitating the proof of (II.47), we get
$$U_{T+T'}^G = U_T^F U_{T'}^{F'}
\eqno ({\mathop{\rm II.75 }\nolimits})$$
for functionals $F$ on ${\bf Z}_{0, T}$ and $F'$ on ${\bf Z}_{0, T'}$, where the functional 
$G$ on 
${\bf Z}_{0, T+T'}$ is defined by
$$G( z \ts z') = F(z) F' (z') \,.
\eqno ({\mathop{\rm II.76 }\nolimits})$$

\vfill\eject
\vglue 1cm
\nn
\centerline {\titre III - Semiclassical Expansions}

\bigskip
\bigskip
We shall compute the semiclassical approximation of the expression in the basic equation
(II.1). The scalar potential $V$ is included via coupled equations rather than via an
additional variable, so that the equations will take readily their familiar forms; i.e. we
compute the wave function:
$$\eqalign{
\Psi \!
\left({ t_b, x_b  }\right)
= \int_{ {\bf Z }  }^{ }
&\Da_{s, Q_0 } z \cdot \exp \! 
\left({ - { \pi \over s } Q_0 (z)  }\right)  . \cr
&\exp \!
\left({ { 1 \over s\hbar  }
\int_{ t_a}^{ t_b } dt \, V \!
\left({ x_b \cdot \Si (t, z) }\right)
 }\right) \cdot 
\phi \left({ x_b \cdot \Si (t_a, z) }\right).\cr} 
\eqno ({\mathop{\rm III.1 }\nolimits} )
$$

\medskip
\nn
For ready use of our results in quantum physics, we compute $\Psi \!\left({ t_b, x_b  }\right)$.
Therefore in the general set up, we set $t_0 = t_b$, $x_0 = x_b$, hence ${\bf Z}$ is the {\it
space  ${\bf Z}_b$ of paths vanishing at} $t_b$ and $x(t_b, z) = x_b$. In paragraph 3, we give the
corresponding results for $\Psi \!\left({ t_b, x_a  }\right)$.

\medskip
We choose an initial wave function given by
$$\phi (x) = \exp \!
\left({-  { 1 \over s\hbar  } \Sa_0 (x)  }\right) \cdot 
\Ta (x), 
\eqno ({\mathop{\rm III.2 }\nolimits} )
$$
where $\Ta$ is a smooth function on $N$ with compact support, and 
$\Sa_0$ is an arbitrary, but reasonable function on $N$. 

\medskip
The initial wave function given by (III.2) generalizes plane waves on $\RR^d$, but is,
obviously, not a momentum eigenstate. We can nevertheless call the semiclassical
expansion of $\Psi \! \left({ t_b, x_b  }\right)$ with this initial wave function a
``{\it momentum-to-position  transition amplitude}'' for three reasons:

\smallskip
(i) In the limit $h = 0$, assuming $s = i$, the current density corresponding to the
initial wave function $\phi$ is
$$
\build { \lim }_{h = 0 }^{ }
{ \hbar \over 2 i  } 
\left({ \phi^* \bgt \phi - (\bgt \phi)^* \phi  }\right)
= 
\left\vert{ \Ta  }\right\vert^2 p\ , 
\quad 
{\mathop{\rm where  }\nolimits}
\ \ p(x) = \bgt \Sa_0 (x) \, .
\eqno ({\mathop{\rm III.3 }\nolimits} )
$$

\medskip
\nn
Consequently, $\Psi \! \left({ t_b, x_b  }\right)$ is the amplitude corresponding to the
transition from momentum $\bgt \Sa_0 (x)$ to position $x_b$.

\smallskip
(ii) We shall expand (III.1) around a ``classical path'' $z_{
{\mathop{\rm  cl }\nolimits}
 }$ characterized by initial momentum and final position.

\smallskip
(iii) The leading terms (henceforth labeled WKB) of semiclassical approximations combine
as if they were transitions between eigenstates, 
e.g.
$$
\eqalign{
\left\langle{  {\mathop{\rm position }\nolimits}
\bigm| {\mathop{\rm  position}\nolimits}  }\right\rangle_{
{\mathop{\rm  WKB}\nolimits}}
&= \int d \ {\mathop{\rm momentum  }\nolimits}
\left\langle{  {\mathop{\rm position }\nolimits}
\bigm| {\mathop{\rm  momentum }\nolimits}  }\right\rangle_{
{\mathop{\rm  WKB}\nolimits}}\cr
&\ \ \ \ \ \  \left\langle{  {\mathop{\rm momentum }\nolimits}
\bigm| {\mathop{\rm  position}\nolimits}  }\right\rangle_{
{\mathop{\rm  WKB}\nolimits}} \cr} $$
[see reference 15, \S 7 in the Appendix].

\medskip
Another useful initial wave function is
$$\phi (x) = \d_{x_a} (x) \, .
\eqno ({\mathop{\rm III.4 }\nolimits} )
$$
The corresponding expression 
 $\psi \!\left({ t_b, x_b  }\right)$
gives a {\it position-to-position transition} amplitude.
In general, this case is more complicated than the previous one, because now the initial
wave function restricts the domain of integration ${\bf Z}$.
Moreover the paths $z \in {\bf Z}$ such that $x \!
\left({t_a, z }\right) = x_a$ do not usually have a common origin.
Even at the very best, when $N = \RR^d$, and 
$z \left({t_a }\right) \simeq x_a - x_b$, one needs to be careful because the domain of
integration, say ${\bf Z}_{a, b} \sbs {\bf Z}$, is not a  vector space but an affine
space. Therefore, we shall compute position-to-position transitions  in section IV where
we specialize the basic equation (III.1). We shall treat in detail the case $N = \RR^d$
and give the necessary indications and references when $N = O (M)$, a frame bundle over a
Riemannian space $M^d$.

\medskip
Semiclassical expansions are best analyzed in the broader context of spaces of paths with
no requirement on their boundary values. For instance let $\Pa \RR^d$ be a space of paths
$z$ with no requirement on $z \left({t_a }\right) $, nor 
$z \left({t_b }\right) $, and let 
$\Pa_{a, b} \RR^d$ be the subspace of $\Pa \RR^d$ with $d$ requirements at 
$t_a$ and $d$ requirements at $t_b$. Let 
$\left\lbc{ S(\nu, z) }\right\rbc_{\nu}^{}$ be a one-parameter family of actions, and let
$U^{2d}(\nu)$ be the 2$d$-dimensional space of motions made of the critical points of
$S(\nu, z)$. A classical paths 
$z_{ {\mathop{\rm  cl }\nolimits} } \in \Pa_{a, b} \RR^d$ for the action $S(\nu, z)$ is at
the intersection of 
$\Pa_{a, b} \RR^d$ with $U^{2d} (\nu)$. If the intersection is transversal, no Jacobi field
of  $z_{ {\mathop{\rm  cl }\nolimits} }$ is in the tangent space to 
$\Pa_{a, b } \RR^d$. In this paper we assume that such is the case. Otherwise there are
caustics and we refer the reader to the literature [e.g. 2].

\bigskip
\nn
{\bf 1. General strategy.}
\smallskip
We introduce a Lagrangian
$$L(t, z) = {{}_1 \over {}^2} 
h_{\a \b } \dot z^\a (t) \dot z^\b (t) - V(x (t, z))
\eqno ({\mathop{\rm III.5 }\nolimits} )
$$
where as usual $x(t, z)$ is the solution of the differential equation (II.21) with
boundary condition $x \! \left({ t_b, z }\right) = x_b$. It is a function of $\dot z (t)$
and a functional of $z$. From this Lagrangian we deduce the {\it action functional}
$$S(z) = 
\int_{ {\bf T}J}^{ }
dt \, L(t, z) + \Sa_0 (x (t_a, z))
\eqno ({\mathop{\rm III.6 }\nolimits} )
$$
on the space ${\bf Z}_b = \Pa_0 \RR^d$ of paths $z$ obeying the boundary condition 
$z \!\left({ t_b }\right) = 0$.
Let 
$z_{ {\mathop{\rm  cl }\nolimits}} $
in ${\bf Z}_b$ be a critical point of the action functional $S$.

\medskip
As it is customary in the calculus of variations, we take a one-parameter variation
$$z(\l) = z_{ {\mathop{\rm  cl }\nolimits}}  + \l \z
\eqno ({\mathop{\rm III.7 }\nolimits} )
$$
with $\z$ in ${\bf Z}_b$
and the equation
$$
\left.{ {d\over d \l  } S (z (\l)) }\right\vert_{\l = 0} = 0
\eqno ({\mathop{\rm III.8 }\nolimits} )
$$
has to be satisfied for all $\z$. That yields, after integrating by parts, a functional
differential equation for 
$ z_{ {\mathop{\rm  cl }\nolimits} } (t_a)$.

\bigskip
Under the affine change of variable from $z $ to $\z$ given
by (III.7), we obtain
$$Q_0 (z) = Q_0 \!
\left({ z_{ {\mathop{\rm  cl }\nolimits}} + \l \z  }\right)
= Q_0 \!
\left({ z_{ {\mathop{\rm  cl }\nolimits}}  }\right)
+ 2 \l 
Q_0 \!\left({ z_{ {\mathop{\rm  cl }\nolimits}} , \z  }\right)
+ \l^2 Q_0 (\z) 
\eqno ({\mathop{\rm III.9 }\nolimits} )
$$
and
$$\Da_{s, Q_0} z = 
\Da_{s, \l^2 Q_0} \z .
\eqno ({\mathop{\rm III.10 }\nolimits} )
$$

\medskip
\nn
The expansion of $x(\cdot, z(\l))$ around 
$x(\cdot, z_{ {\mathop{\rm  cl }\nolimits}})$ reads 
$$P
\left({  z_{ {\mathop{\rm  cl }\nolimits}} + \l \z }\right)
=
P
\left({  z_{ {\mathop{\rm  cl }\nolimits}} }\right)
+ \l P' 
\left({  z_{ {\mathop{\rm  cl }\nolimits}} }\right)
\cdot \z +
{1\over 2}\l^2 P'' 
\left({  z_{ {\mathop{\rm  cl }\nolimits}}  }\right)
\cdot \z \z + O(\l^3).
\eqno ({\mathop{\rm III.11 }\nolimits} )$$

\medskip
\nn
Here 
$P' \left({  z_{ {\mathop{\rm  cl }\nolimits}} }\right)
$ and $P'' \left({  z_{ {\mathop{\rm  cl }\nolimits}} }\right)$
are the first and the second derivative mappings of $P$ at
$  z_{ {\mathop{\rm  cl }\nolimits}} $, where 
$P : {\bf Z}_b \ra \Pa_{x_b} N$ takes $z$ into 
$x (\cdot, z)$. We abbreviate $x(t,z_{{\mathop{\rm cl}\nolimits}})$ to $x_{ {\mathop{\rm cl }\nolimits}}(t)$. They are of the form
$$
\eqalign{
\left({ P'  \left({  z_{ {\mathop{\rm  cl }\nolimits}} }\right) \cdot \z  }\right)^\a (t) 
&= - \int_{t }^{t_b  } ds \, 
{ \d x_{ {\mathop{\rm  cl }\nolimits} }^\a  (t) \over
\d z_{ {\mathop{\rm  cl }\nolimits} }^\b  (s)  } \, 
\z^\b (s)\cr
&= - \int_{t }^{t_b  } ds \, 
k^\a {}_\b (t, s) \, \z^\b (s) = : \xi^\b (t)\cr}
\eqno ({\mathop{\rm III.12 }\nolimits} )$$
$$\eqalignno{
\left({ P''  \left({  z_{ {\mathop{\rm  cl }\nolimits}} }\right) \cdot \z \z  }\right)^\a
(t)  &= - \int_{t }^{t_b  } ds 
\left({  - \int_{t }^{t_b  } du }\right)
{ \d^2 x_{ {\mathop{\rm  cl }\nolimits} }^\a  (t) \over
\d z_{ {\mathop{\rm  cl }\nolimits} }^\b  (s)  
\d z_{ {\mathop{\rm  cl }\nolimits} }^\g  (u)
}
\z^\b (s) \z^\g (u) \cr
&=  \int_{t }^{t_b  } ds  \int_{t }^{t_b  } du \, k^\a 
 {}_{\b \g}  (t, s, u ) \z^\b (s) \z^\g (u) .
&({\mathop{\rm III.13 }\nolimits} )
 \cr} 
$$

\bigskip
We could obtain the short time propagator by expanding the quantity
$\phi \left({  x_b \cdot \Si (t, z) }\right)$ in equation (III.1) around 
$\phi \left({ x_b }\right)$ with $x_b \in N$.
Here, we shall expand 
$\phi \left({  x_b \cdot \Si (t, z) }\right)$ around
$\phi \left({  x_b \cdot \Si (t, z_{ {\mathop{\rm  cl }\nolimits} } ) }\right)
$
with $z_{ {\mathop{\rm  cl }\nolimits} } \in {\bf Z}_b$.
When we choose for $z_{ {\mathop{\rm  cl }\nolimits} }$  a critical point of the action, we
obtain the WKB approximation.
These are two different expansions, and in general, the short time propagator is
different from the WKB approximation [14]. But in some simple cases, they are equal [1].

\bigskip
\nn
{\it Remark.}
A simple and powerful argument of Stephen A. Fulling clarifies this issue: The Schr\"odinger operator 
$\exp (- it H / \hbar)$ can be written in terms of dimensionless terms 
${tH \over \hbar } = - {1 \over 2}  A \D + B V $, where 
$A = {\hbar t \over m} $ and $B = {\l t \over \hbar} $ with $\l$ a coupling constant. There are 
4 {\it different} possible expansions of physical interest:

-- expansion in $B$, equivalently expansion in $\l$;

-- expansion in $A$, equivalently expansion in $m^{- 1}$;

-- expansion in $A \cdot B$, equivalently expansion in $t$;

-- expansion in $A/B$, equivalently expansion in $\hbar$.

\bigskip
\nn
{\bf 2. Momentum-to-position transitions.}
\smallskip
In this paragraph we take $s = i$ and write 
$\Da_{s, Q_0}$ simply as $\Da_{Q_0}$. We expand the integrand in 
(III.1) by taking $z$ in the form
$z_{ {\mathop{\rm  cl }\nolimits} } + \z$ where 
$z_{ {\mathop{\rm  cl }\nolimits} }
$ is a critical point of the action functional $S$, 
 with boundary condition
$z_{ {\mathop{\rm  cl }\nolimits} }
\left({ t_b }\right) = 0$.
The initial wave function is given by (III.2).

\bigskip
The terms independent of $\z$ combine to make the action function 
$$\Sa \! \left({t_b, x_b   }\right) =
{{}_1 \over {}^2}
Q_0 \!
\left({ z_{ {\mathop{\rm  cl }\nolimits}}  }\right)
- \int_{t_a  }^{t_b  } dt \, 
V \!\left({ x_{ {\mathop{\rm  cl }\nolimits}} (t)  }\right)
+ \Sa_0 \!
\left({ x_{ {\mathop{\rm  cl }\nolimits}} (t_a)  }\right)
\eqno ({\mathop{\rm III.14 }\nolimits} )$$
where $x_{ {\mathop{\rm  cl }\nolimits}} (t) $ is by definition 
$x
\left({  t, z_{ {\mathop{\rm  cl }\nolimits}} .
}\right) $.
The terms linear in $\z$
$$
\eqalign{ \int_{\bf T }  dt & 
\left({  h_{\a \b } 
\dot z_{ {\mathop{\rm  cl }\nolimits} }^\a (t) 
\dot z_{ {\mathop{\rm  cl }\nolimits} }^\b (t) 
- \bgt_\a 
V \left({ x_{ {\mathop{\rm  cl }\nolimits}} (t) }\right)
\left({ P' \left({ z_{ {\mathop{\rm  cl }\nolimits}}   }\right) \cdot \z  }\right)^\a
(t)
}\right)
\cr 
& \ \ \ 
+ \bgt_\a \,  \Sa_0 
\left({ x_{ {\mathop{\rm  cl }\nolimits}} (t_a) }\right)
\left({ P' \left({ z_{ {\mathop{\rm  cl }\nolimits}} }\right) \cdot  \z  }\right)^\a 
(t_a)
\cr} \eqno ({\mathop{\rm III.15 }\nolimits} ) $$
vanish since $z_{ {\mathop{\rm  cl }\nolimits}}$ is a critical point of $S$ such that
$z_{ {\mathop{\rm  cl }\nolimits}} (t_b) = 0$. The  quadratic terms  can 
be written in the form
$\left({ Q_0 + Q }\right) (\z)$. Thanks to the equations  (A.127) and (A.132), the
corresponding integration is:
$$
\int_{{\bf Z} }^{ }
\Da_{Q_0 } \z \cdot
\exp \! \left({ -{\pi  \over s } \left({ Q_0 + Q }\right) (\z) }\right)
=
\Det \left({   Q_0 / (Q_0 + Q)  }\right)_{}^{1/2} \, . 
\eqno ({\mathop{\rm III.16 }\nolimits} )$$

\medskip
\nn
Thus the dominating term of the semiclassical expansion of 
$\Psi \!\left({ t_b, x_b  }\right)$ is:
$$\Psi_{ {\mathop{\rm  WKB }\nolimits}
 } \left({ t_b, x_b  }\right) =
\exp \!
\left({ {i \over \hbar} \Sa \! \left({ t_b, x_b  }\right) }\right) \cdot 
\Det \left({   Q_0 / (Q_0 + Q)  }\right)_{}^{1/2} \cdot
\Ta \!
\left({ x_{ {\mathop{\rm  cl }\nolimits}} (t_a) }\right) \, . 
\eqno ({\mathop{\rm III.17 }\nolimits} )$$

\medskip
We proceed to calculate 
$\Det (Q_0 / Q_\nu)$, with 
$Q_\nu = Q_0 + \nu Q$, in the case $\nu = 1$.
Physical arguments as well as previous calculations (see [6] and the appendix of [15])
suggest that $\Det (Q_0/Q_1)$ is equal to 
$\det K^\a {}_\b \left({ t_b, t_a  }\right)$, where the Jacobi matrix
$ K^\a {}_\b \left({ t_b, t_a  }\right)$ is defined by:
$$ K^\a {}_\b \left({ t_b, t_a  }\right)
= { \part z_{ {\mathop{\rm  cl }\nolimits}}^\a (t_b) \over
\part z_{ {\mathop{\rm  cl }\nolimits}}^\b (t_a) }\, .
\eqno ({\mathop{\rm III.18 }\nolimits} )$$

\medskip
\nn
Therefore we shall construct a one-parameter family of actions 
$S(\nu, z)$ and a one-parameter family of Jacobian matrices:
$$K^\a{}_\b 
\left({ \nu ; t_b, t_a }\right)
: = { \part z_{{\mathop{\rm  cl }\nolimits} }^\a  
\left({ \nu ; t_b  }\right)\over 
\part z_{{\mathop{\rm  cl }\nolimits} }^\b  
\left({ \nu ; t_a }\right) } \, .
\eqno ({\mathop{\rm III.19 }\nolimits} )$$

\medskip
\nn
We shall show that:
$$ \det \! \left({  K \!\left({ \nu ; t_b, t_a }\right)
\cdot  K\! \left({ 0 ; t_b, t_a }\right)^{- 1}  }\right) 
= : c (\nu) 
\eqno ({\mathop{\rm III.20 }\nolimits} )$$
satisfies the same boundary condition as 
$\Det \left({ Q_0 / Q_\nu }\right) $, namely
$$c(0) = \Det \left({ Q_0 / Q_0 }\right)  = 1, 
\eqno ({\mathop{\rm III.21 }\nolimits} )$$
and the same differential equation as 
$\Det \left({ Q_0 / Q_\nu  }\right) $, namely,
$${d  \over d \nu}  \ell n \, 
\Det \left({ Q_\nu / Q_0  }\right)   =
\Tr \! \left({ Q_\nu^{- 1}{ d \over d\nu } Q_\nu}\right)
= \Tr \! \left({ Q_\nu^{- 1} Q  }\right)
\eqno ({\mathop{\rm III.22 }\nolimits} )$$
(See equation (A.124) in Appendix A).

\bigskip
\nn
2.1. {\it  Differential equation satisfied by $c(\nu)$.}

\smallskip
We choose the family $S(\nu, z)$ such that its second variation evaluated at
 $z_{ {\mathop{\rm  cl }\nolimits} }$ be
$$S'' \left({ \nu, z_{ {\mathop{\rm  cl }\nolimits} } }\right)
 \cdot \z \z =
Q_0 (\z) + \nu Q (\z) = : Q_\nu (\z) \, .
\eqno ({\mathop{\rm III.23 }\nolimits} )$$

\medskip
\nn
It follows from its definition (III.18) that the $\b$-column 
$K_\b \left({ \nu ; t, t_a }\right)$ of the Jacobi matrix
$K \! \left({ \nu ; t, t_a }\right)$ is the following Jacobi field along the classical path
$z_{ {\mathop{\rm  cl }\nolimits} } (\nu, t)$ of the system governed by the action $S(\nu, z)$
$$Q_\nu K_{(\b)}  \! \left({ \nu ; t, t_a }\right) = 0
\eqno ({\mathop{\rm III.24 }\nolimits} )$$
$$
 K^{\a}_{\ (\b)} \! \left({ \nu ; t_a, t_a }\right) = \d_\b^\a \, .
\eqno ({\mathop{\rm III.25 }\nolimits} )$$

\medskip
\nn
Relation (III.24) encodes a functional differential equation for $K$ and an implicit
equation for 
$d K \! \left({ \nu ; t, t_a }\right) / dt$ at $t = t_a$.
It also gives, after differentiating w.r.t. $\nu$
$$Q_\nu {d \over d \nu } K(\nu) = - Q K (\nu)
\eqno ({\mathop{\rm III.26 }\nolimits} )$$
which can be solved with the {\it retarded Green's function}
$G_\nu^{ {\mathop{\rm  ret}\nolimits} }  $ of $Q_\nu$, namely:
$${d \over d \nu }
 K \!\left({ \nu ; t, t_a }\right)
= - \int_{ {\bf T } }^{ }
G_\nu^{ {\mathop{\rm  ret}\nolimits} } (t, s) Q(s) 
K \!\left({ \nu ; s, t_a }\right) ds \, .
\eqno ({\mathop{\rm III.27 }\nolimits} )$$

\medskip
The retarded Green's function can be expressed in terms of Jacobi fields, namely:
$$G_\nu^{ {\mathop{\rm  ret}\nolimits} } (t, s)
= \t (t - s) J (\nu ; t, s)
\eqno ({\mathop{\rm III.28 }\nolimits} )$$
where $\t$ is the step function, equal to $1$ for $t > s$ and $0$ for $t < s$; moreover
$J(\nu ; t, s)$ is a Jacobi matrix, {\it i.e.}  each column is a Jacobi field. The
$\a$-component of the $\b$-Jacobi field
$J^{(\b)} \left({ \nu ; t, t_a }\right)$ is
$$J^{\a (\b)}  \left({ \nu ; t, t_a }\right)
= {\part 
z_{ {\mathop{\rm  cl }\nolimits} }^\a (\nu, t) 
\over \part 
p_{ {\mathop{\rm  cl }\nolimits}, \b } \left({ \nu , t_a }\right)
}
\eqno ({\mathop{\rm III.29 }\nolimits} )$$
where
$
p_{ {\mathop{\rm  cl }\nolimits}  } (\nu, t) : =
\d \Sa (\nu, z) /
\d z_{ {\mathop{\rm  cl }\nolimits}} (\nu, t)$, here
$p_{ {\mathop{\rm  cl }\nolimits}, \a } (\nu, t) =
h_{\a \b } \dot z_{ {\mathop{\rm  cl }\nolimits}}^{\b}
(\nu, t).$

\medskip

>From (III.20), we get the equation for $c(\nu)$:
$${ d\over d \nu }J\ell n \, c (\nu) =
\tr \!
\left({ N  \! \left({ \nu ; t_a, t_b }\right) { d\over d \nu }J 
K  \! \left({ \nu ; t_b, t_a }\right) 
}\right)\eqno ({\mathop{\rm III.30 }\nolimits} )$$
where the matrix $N$ is the inverse of $K$
$$N \! \left({ \nu ; t_a, t_b }\right)  
K \! \left({ \nu ; t_b, t_a }\right) 
= \un . 
\eqno ({\mathop{\rm III.31 }\nolimits} )$$

\medskip
\nn
Using (III.27) and (III.28), we obtain the sought-for differential equation:
$$
{ d\over d \nu }J\ell n \  c(\nu) = - \tr \!
\left({ \int_{ {\bf T}J}^{ }
d s \, J \!
\left({ \nu ; t_b, s }\right) 
Q(s) K \! \left({ \nu ; s, t_a }\right) 
N \! \left({ \nu ;  t_a, t_b }\right) }\right)
 \, . 
\eqno ({\mathop{\rm III.32 }\nolimits} )$$

\medskip
Within focal distance of 
$z_{ {\mathop{\rm  cl }\nolimits} }
\left({ \nu ; t_a }\right)$, the $d$ Jacobi fields
$K_{(\b)} (\nu) $ and the $d$ Jacobi fields 
$J^{(\b)} (\nu)$ are linearly independent, and 
$c(\nu)$ is defined. It gives the rate at which the flow of classical paths
$
\left\lbc{ z_{ {\mathop{\rm  cl }\nolimits}} (\nu) }\right\rbc$
diverges or converges. 

\vfill\eject
\nn
 2.2. {\it Comparison between the two differential equations} (III.32) {\it and }
(III.22).

\smallskip
We shall express $Q_\nu^{- 1}$ in terms of Jacobi fields, as this is always possible
within focal distance of 
$
z_{ {\mathop{\rm  cl }\nolimits}} \left({ t_a }\right)$.
The unique inverse $G_\nu $ of $Q_\nu$ satisfies the equation 
$$Q_\nu G_\nu = \un
\eqno ({\mathop{\rm III.33 }\nolimits} )$$
with 
$$G_\nu \!\left({  t_b, s }\right) = 0 
\eqno ({\mathop{\rm III.34 }\nolimits} )$$
by virtue of (A.40) together with the specialization used in (A.62) when 
applied to the space ${\bf Z}_b$ of paths
vanishing at $t_b$. Therefore (see formula (B.26) in Appendix B)
$$
\eqalign{
G_\nu (t, s) &= \t (s - t) K \!
\left({ \nu ; t, t_a }\right) N \!\left({ \nu ; t_a, t_b }\right)  
J \!\left({ \nu ; t_b, s }\right) \cr
&\ \ \ - \t(t - s) J \!\left({ \nu ; t, t_b }\right) 
\ww N \!\left({ \nu ; t_b, t_a }\right) 
\ww K \!\left({ \nu ; t_a, s }\right) \cr} 
\eqno ({\mathop{\rm III.35 }\nolimits} )$$
where $\ww K$ is the transpose of $K$, and $\ww N$ is $\ww K$'s inverse.

\medskip
Substituting $Q_\nu^{- 1} = G_\nu $ in (III.22) gives:
$${ d\over d \nu }J\ell n \Det  \! \left({ Q_\nu / Q_0 }\right) 
= \tr \!
\left({ \int_{{\bf T}} dt   \,  K \!\left({ \nu ; t, t_a }\right) 
N \!\left({ \nu ; t_a, t_b }\right)  
J \! \left({ \nu ; t_b, t }\right)  Q(t) }\right)
 \, . 
\eqno ({\mathop{\rm III.36 }\nolimits} )$$

\medskip
\nn
Comparison with (III.32) shows that
$$\ell n \Det \left({  Q_\nu / Q_0 }\right)  = - \ell n \, c (\nu) 
\eqno ({\mathop{\rm III.37 }\nolimits} )$$
i.e.\footnote{$^{10}$}{
The credit for the technique used in deriving (III.38) is due to B. Nelson and B. Sheeks [6].
It is made simpler and more general here by not integrating $Q_\nu$ by parts. The first
calculation giving the ratio of functional determinants in terms of finite determinants
can be found in [16, 5].}
$$
\Det \left({ Q_0 / Q_\nu }\right)   =
\det \!
\left({  K \!\left({ \nu ; t_b, t_a }\right)  \cdot  K \!\left({ 0 ; t_b, t_a }\right)^{-1}
}\right)
  \, .
 \eqno ({\mathop{\rm III.38 }\nolimits} )$$

\medskip
\nn
Here $ K \!\left({ 0 ; t, t_a }\right)$ satisfies
$$Q_0  K \!\left({ 0 ; t, t_a }\right) = 0 
\eqno ({\mathop{\rm III.39 }\nolimits} )$$
$$  K \!\left({ 0 ; t_a, t_a }\right) = \un 
\eqno ({\mathop{\rm III.40 }\nolimits} )$$
where
$$Q_0 (z) = - \ds \int_{{\bf T} } dt \, h_{\a  \b }
\ddot z^{\a} (t) z^\b (t)  - h_{\a  \b }
\dot z^\a \! \left({ t_a }\right) z^\b \! \left({ t_a }\right)\, .$$

\medskip
\nn
Hence
 we have
$$ {d^2 \over dt^2} K \! \left({ 0 ; t, t_a }\right) = 0 \ \  , \ \ 
\left.{ 
{d \over dt } K \! \left({ 0 ; t, t_a }\right)
 }\right\vert_{t = t_a}
= 0 
\eqno ({\mathop{\rm III.41 }\nolimits} )$$
therefore
$$K^\a {}_\b \left({ 0 ; t, t_a }\right) = \d^\a{}_\b \, .
\eqno ({\mathop{\rm III.42 }\nolimits} )$$

\medskip
Inserting 
$\Det (Q_0 / Q_\nu) = \det K^\a {}_\b \left({ t_b, t_a }\right)$ 
into (III.17) we obtain:
$$\Psi_{
{\mathop{\rm  WKB}\nolimits}
} \left({ t_b, x_b }\right) =
\left({  \build {\det  }_{ \a, \b}^{ } 
{ \part    z_{ {\mathop{\rm  cl }\nolimits}}^\a   \left({ t_b }\right) \over 
\part    z_{ {\mathop{\rm  cl }\nolimits}}^\b   \left({ t_a }\right)
} 
}\right)_{}^{1/2} \cdot 
\exp \!
\left({ { i  \over \hbar }  \Sa \!\left({ t_b, x_b }\right) }\right)
\cdot \Ta \! \left({ x_{ {\mathop{\rm  cl }\nolimits}} (t_a) }\right) \, .
\eqno ({\mathop{\rm III.43 }\nolimits} )$$

\smallskip
\nn
If $ z_{ {\mathop{\rm  cl }\nolimits}}  \left({ t_b }\right)$
is conjugate to 
$z_{ {\mathop{\rm  cl }\nolimits}}  \left({ t_a}\right)$, in the sense of caustic theory,
one needs [2] to include terms in $\z$ of order higher than $2$ in the calculation of $\Psi$.

\bigskip
\nn
2.3. {\it  End of calculation.}

\smallskip
Equation (III.43) is not the end of the calculation;
$\Psi_{
{\mathop{\rm  WKB}\nolimits}
} $ must be expressed in terms of 
$x_{ {\mathop{\rm  cl }\nolimits}}
: T \ra N$, where
$$x_{ {\mathop{\rm  cl }\nolimits}} (t) =  x \!
\left({ t,  z_{ {\mathop{\rm  cl }\nolimits}} }\right)
= x_b \cdot \Si \left({ t, z_{ {\mathop{\rm  cl }\nolimits}} }\right)
\eqno ({\mathop{\rm III.44 }\nolimits} )$$
but not in terms of 
$z_{ {\mathop{\rm  cl }\nolimits}} : {\bf T} \ra \RR^d$.
Two techniques present themselves: we know the evolution equation satisfied by the wave
function $\Psi$ on $N$, therefore we can construct an action on $N$ and find its critical
points. But we do not need to find 
$x_{ {\mathop{\rm  cl }\nolimits}} $
 corresponding to $z_{ {\mathop{\rm  cl }\nolimits}}$, we only need the expression
corresponding to the prefactor of (III.43). This can be obtained by a simpler technique.

\medskip
Recall the parametrization $P : \Pa_0 \RR^d \ra \Pa_{x_b} N$ given by 
$P(z) = x$. Corresponding to the action functional 
$S $ on $\Pa_0 \RR^d$ we get a functional $\oo S $ on 
$\Pa_{x_b} N $ such that $S = \oo S \circ P$, i.e. 
$\oo S (x) = S(z)$. Making the substitution $z =
z_{ {\mathop{\rm  cl }\nolimits}} + \z$ and expanding up to second order terms in $\z$,
we get
$$\oo S \!\left({ x_{ {\mathop{\rm  cl }\nolimits}} }\right)
= S \!\left({ z_{ {\mathop{\rm  cl }\nolimits}} }\right)
\eqno ({\mathop{\rm III.45 }\nolimits} )$$
$$\oo S' \! \left({ x_{ {\mathop{\rm  cl }\nolimits}} }\right) \cdot 
{ \d x_{ {\mathop{\rm  cl }\nolimits}} \over 
\d z_{ {\mathop{\rm  cl }\nolimits}} } = 
S' \! \left({ z_{ {\mathop{\rm  cl }\nolimits}} }\right)
\eqno ({\mathop{\rm III.46 }\nolimits} )$$
$$
\oo S'' \! \left({ x_{ {\mathop{\rm  cl }\nolimits}} }\right) \cdot 
{ \d x_{ {\mathop{\rm  cl }\nolimits}} \over 
\d z_{ {\mathop{\rm  cl }\nolimits}} } 
{ \d x_{ {\mathop{\rm  cl }\nolimits}} \over 
\d z_{ {\mathop{\rm  cl }\nolimits}} } 
= 
S'' \! \left({ z_{ {\mathop{\rm  cl }\nolimits}} }\right). 
\eqno ({\mathop{\rm III.47 }\nolimits} )$$

\smallskip
\nn
In the beginning of this paragraph we calculated the expansion of
$S \! \left({  z_{ {\mathop{\rm  cl }\nolimits}} + \z }\right)
$ up to second order and obtained the expression
$$S'' \!  \left({ z_{ {\mathop{\rm  cl }\nolimits}} }\right) 
\cdot \z \z =
Q_0 (\z) + Q(\z) = Q_1 (\z)
\eqno ({\mathop{\rm III.48 }\nolimits} )$$
for the second order variation. We also know that the infinite-dimensional determinant 
$\Det (Q_0 /Q_1)$ is equal to the determinant of the Jacobi matrix
$\part z_{ {\mathop{\rm  cl }\nolimits}} 
\left({ t_b  }\right) /
\part z_{ {\mathop{\rm  cl }\nolimits}} 
\left({ t_a  }\right)$.

\medskip
Similarly we write
$$\oo S'' \! 
\left({ x_{ {\mathop{\rm  cl }\nolimits}}  }\right) \cdot 
{ \d x_{ {\mathop{\rm  cl }\nolimits}} \over 
\d z_{ {\mathop{\rm  cl }\nolimits}} } \z 
{ \d x_{ {\mathop{\rm  cl }\nolimits}} \over 
\d z_{ {\mathop{\rm  cl }\nolimits}} } 
\z  \, = \, 
\oo Q_0 (\z) + \oo Q(\z) 
\eqno ({\mathop{\rm III.49 }\nolimits} )$$
and identify $\oo Q_0$ as follows. Let
$$S_0 (z) : = { {}_{1} \over {}^{2} } Q_0 (z)
\eqno ({\mathop{\rm III.50 }\nolimits} )$$
and
$$\oo S_0 (x(z)) : = S_0 (z)
\eqno ({\mathop{\rm III.51 }\nolimits} )$$
then
$$\oo Q_0 (\z) : = 
\oo S''_0 
\left({ x_{ {\mathop{\rm  cl }\nolimits}}   }\right) \cdot 
{ \d x_{ {\mathop{\rm  cl }\nolimits}} \over 
\d z_{ {\mathop{\rm  cl }\nolimits}} } \z
{ \d x_{ {\mathop{\rm  cl }\nolimits}} \over 
\d z_{ {\mathop{\rm  cl }\nolimits}} } 
\z 
\eqno ({\mathop{\rm III.52 }\nolimits} )$$
$$ \det 
{  \part x_{ {\mathop{\rm  cl }\nolimits}}
\left({ t_b  }\right)
 \over \part 
 x_{ {\mathop{\rm  cl }\nolimits}} \left({ t_a  }\right)  } 
= \Det \!
\left({ \oo Q_0 / \! \left({ \oo Q_0 + \oo Q }\right)  }\right)
=\Det \!
\left({  Q_0 / \! \left({  Q_0 +   Q }\right)  }\right)
  \,  .
 \eqno ({\mathop{\rm III.53 }\nolimits} )$$

\smallskip
\nn
Therefore
$$\det 
{  \part x_{ {\mathop{\rm  cl }\nolimits}}
\left({ t_b  }\right)
 \over \part 
 x_{ {\mathop{\rm  cl }\nolimits}} \left({ t_a  }\right)  } 
=  
\det 
{  \part z_{ {\mathop{\rm  cl }\nolimits}}
\left({ t_b  }\right)
 \over \part 
 z_{ {\mathop{\rm  cl }\nolimits}} \left({ t_a  }\right)  } \, .
\eqno ({\mathop{\rm III.54 }\nolimits} )$$

\medskip
In conclusion, we obtain the sought-for semiclassical expansion
$$
\Psi_{
{\mathop{\rm  WKB}\nolimits}
} 
\left({ t_b, x_b }\right)
= \det  \! \left({
{ \part    x_{ {\mathop{\rm  cl }\nolimits}}   \left({ t_b }\right) \over 
\part    x_{ {\mathop{\rm  cl }\nolimits}}   \left({ t_a }\right)
} 
}\right)_{}^{1/2} \cdot \, 
\exp \!
\left({ { i  \over \hbar }  \Sa \left({ t_b, x_b }\right) }\right) \cdot 
\Ta \! \left({ x_{ {\mathop{\rm  cl }\nolimits}} (t_a) }\right)  .
\eqno ({\mathop{\rm III.55 }\nolimits} )$$

\smallskip
\nn
According to equation (III.46), the path 
$ x_{ {\mathop{\rm  cl }\nolimits}}$ in $N$ is a critical point for the action functional
$\oo S : \Pa_{x_b} N \ra \RR$. By construction, we have
$$\oo S' (x (\cdot , z) ) = S' (z)\, .
\eqno ({\mathop{\rm III.56 }\nolimits} )$$

\smallskip
\nn
Using equation (III.5) and (III.6) this can be made more explicit as follows
$$\oo S (x) = { {}_{1} \over {}^{2} }
\int_{ {\bf T} }^{}
dt \, h_{\a \b } \, 
\dot z^\a (t) \dot z^\b (t) - \int_{ {\bf T} }^{}
V(x (t)) dt + \Sa_0 \! \left({ x (t_a) }\right) \, .
\eqno ({\mathop{\rm III.57 }\nolimits} )$$

\smallskip
\nn
In our general setup
, the first integral remains somewhat implicit, but can be used for
practical calculations in the applications given in section IV.

\medskip
The prefactor in (III.55) also gives the volume expansion or contraction of a congruence
of classical paths originating in the neighborhood of 
$\Sa_0 \! \left({ x_{ {\mathop{\rm  cl }\nolimits}} (t_a) }\right)$ with momentum
$\bgt \Sa_0 (m)$, $m \in N$. This determinant depends both on the choice of the initial
wave function and the dynamics of the system.

\medskip
Set $d \om_a$ the volume element on $N$ at 
$x_{ {\mathop{\rm  cl }\nolimits}} (t_a) $ and 
$d \om_b = \det \!
\left({ { \part x_{ {\mathop{\rm  cl }\nolimits}}^\a (t_b) \over
\part x_{ {\mathop{\rm  cl }\nolimits}}^\b (t_a)  } }\right)
d \om_a$; then 
$$
\build {\lim }_{ h = 0 }^{ }
\int_{ C_{t_b } \O }^{ }
\left\vert{  \Psi \! \left({  t_b, x_b }\right) }\right\vert^2 d \om_b
= \int_{ \O }^{ }
\left\vert{ \phi \left({   x_a }\right) }\right\vert^2 d \om_a
$$
where $C_t$ belongs to the group of transformations generated by the classical flow
(see details in [5 p. 299]).

\bigskip
\nn
{\bf 3. Diffusion problems.}

\smallskip
In a diffusion problem, one is interested in computing the functional integral
$$
\eqalign{
\Psi \!\left({  t_b, x_a }\right) =
&
\int_{{\bf Z}  }^{ } \Da_{s, Q_0 }   z \exp \!
\left({- {\pi \over s} Q_0 (z) }\right) \cdot\cr
& \  \ \ 
\exp \!
\left({ {1  \over s\hbar } \int_{t_a }^{t_b } dt \, 
V \!\left({ x_a \cdot \Si (t, z)  }\right)
 }\right)\cdot \phi \left({ x_a \cdot \Si (t_b, z) }\right) \cr} 
\eqno ({\mathop{\rm III.58 }\nolimits} )$$
where ${\bf Z}$ denotes now the space of paths $z$ vanishing at $t = t_a$.
If we choose $\phi$ to be of the form (III.2), namely
$$
\phi \left({ x_a \cdot \Si (t_b, z) }\right)
=
\phi \left({ x (t_b, z) }\right) =
\exp \!
\left({- {1  \over s\hbar } \Sa_0 \left({ x (t_b, z) }\right)  }\right)
\Ta \!\left({ x (t_b, z) }\right) 
\eqno ({\mathop{\rm III.59 }\nolimits} )$$
then computing (III.58), with $s = i$, can be said to be computing the probability
amplitude of a transition from a position $x_a$ to a  momentum of the form
$$p \left({ x_b }\right) =
\bgt \Sa_0  \left({ x_b }\right)  
\eqno ({\mathop{\rm III.60 }\nolimits} )$$
for the end point 
$x_b = x \! \left({ t_b, z }\right)$ of the path $x$.

\medskip
With $\phi$ given by (III.59), one obviously expects the WKB approximation of (III.58) to
be $$\Psi_{
{\mathop{\rm  WKB}\nolimits}
} 
\left({ t_b, x_a }\right)
= \det \! \left({
{ \part    x_{ {\mathop{\rm  cl }\nolimits}}^\a   \left({ t_a }\right) \over 
\part    x_{ {\mathop{\rm  cl }\nolimits}}^\b   \left({ t_b }\right)
} 
}\right)_{}^{1/2}
\exp \! 
\left({ {i  \over \hbar }  \Sa \! \left({ t_b, x_a }\right) }\right) \cdot 
\Ta \!\left({ x_{ {\mathop{\rm  cl }\nolimits}} (t_b) }\right) \, .
\eqno ({\mathop{\rm III.61 }\nolimits} )$$

\smallskip
\nn
Nevertheless, it is gratifying to derive (III.61) by the method followed in paragraph 2.
The only necessary changes are described as follows:
$$
c(\nu) = \det \!
\left({ K\!
\left({ \nu ; t_a , t_b  }\right)
\cdot 
K \!
\left({ 0 ; t_a , t_b  }\right)^{- 1}  }\right)
\eqno ({\mathop{\rm III.20 }\nolimits}^{{\mathop{\rm bis }\nolimits} }  )$$
$$
G_{\nu}^{ {\mathop{\rm ad }\nolimits}  \nu}
(t, s) = \t (t - s) J (\nu ; s, t) =
- \t (s - t) J (\nu ; t, s) 
\eqno ({\mathop{\rm III.28 }\nolimits}^{{\mathop{\rm bis }\nolimits} }  )$$
$$
\eqalign{
G_\nu (t, s) &
= \t (s - t) J \!
\left({\nu ; t,t_a  }\right) \ww N \! \left({\nu ;  t_a, t_b }\right)
\ww K \! \left({\nu ;  t_b, s }\right) \cr
& \ \ \ -
\t  (t - s) K \! \left({ \nu ; t, t_b }\right) N \! \left({\nu ;  t_b, t_a }\right)
J \! \left({\nu ;  t_a, s }\right)\cr} 
\eqno ({\mathop{\rm III.35 }\nolimits}^{{\mathop{\rm bis }\nolimits} }  )$$
\nn
(see formula (B.27) in Appendix B).

\vfill\eject
\vglue 1cm
\nn
\centerline {\titre IV - Illustrations}

\bigskip
\bigskip
The specializations presented in this section illustrate and develop the general
 formulas derived in sections II and
III. They differ by the choice of the manifold $N$ and of the initial wave function
 $\phi$. We treat the case of scalar wave functions, but the case of
a tensor field is easily accommodated by using the Lie derivative of tensor fields.

\medskip
We are basically interested in the Schr\"odinger equation. So the reader should substitute
$i$ to $s$ in the formulas. Moreover, we want to evaluate the value of the wave function
 $\Psi$ at the
{\it final time} $t_b$, and {\it final position } $x_b$. This dictates the choice of the 
space ${\bf Z}_b$
of paths, determined by $z(t_b) = 0$.

\bigskip
\nn
{\bf 1. Point-to-point transitions in a flat space.}

\smallskip
The basic manifold $N$ is the euclidean space $\RR^d$ of dimension $d$, in cartesian 
coordinates (see
paragraph IV.2 for other coordinates). For a vector $x$, with coordinates $x^1, \cdots, x^d$,
 we set its
length to be
$\left\vert{ x }\right\vert =
\left({ \sum_{ \a = 1}^{d } (x^\a)^2}\right)_{}^{1/2}$ as usual.
We consider a particle of mass $m$ moving in the field of a potential $V$. The Lagrangian
 and action are
given as usual by
$$
\eqalignno{
L(x, \dot x) &= 
{m \over 2} \left\vert{ \dot x }\right\vert_{}^{2} - V(x), 
& ({\mathop{\rm IV.1 }\nolimits})\cr
S(x (\cdot)) &= \int_{ {\bf T} }^{ } dt \, 
L(x(t), \dot x(t)) \, . 
& ({\mathop{\rm IV.2 }\nolimits}) \cr }$$

\medskip
In classical mechanics, we solve the equations of motion with suitable boundary conditions:
$$m \! \build { x}_{ }^{  \cdot \cdot }_{
}\!\!{}_{{\mathop{\rm cl }\nolimits}} = - \bgt V \!\left({ x_{{\mathop{\rm cl }\nolimits}}  }\right)
\eqno ({\mathop{\rm IV.3 }\nolimits})$$
$$
x_{{\mathop{\rm cl }\nolimits}} (t_a) = x_a \quad ,
\quad x_{{\mathop{\rm cl }\nolimits}} (t_b) = x_b \, , 
\eqno ({\mathop{\rm IV.4 }\nolimits})$$
where $x_a, x_b$ are points in $\RR^d$. In quantum mechanics, we want to solve the Schr\"odinger
equation
$$i \hbar { \part \Psi \over \part t} 
= - { \hbar^2 \over 2m } \D \Psi + V \Psi
\eqno ({\mathop{\rm IV.5 }\nolimits})$$
with initial condition
$$\Psi \! \left({ t_a, x }\right)
= \phi (x) \, . 
\eqno ({\mathop{\rm IV.6 }\nolimits})$$

\bigskip
\nn
The {\it point-to-point transition amplitudes}  are given by
$$
\left\langle{  t_b ,  x_b \bigm|  t_a , x_a }\right\rangle
=
\Psi \!
\left({ t_b, x_b }\right) \, ,
\eqno ({\mathop{\rm IV.7 }\nolimits})$$
where $\phi (x)$ is  a delta function $\d \!\left({ x - x_a }\right)$.

\bigskip
The original claim of Feynman was that we can solve the Schr\"odinger equation by the following
path integral
$$\Psi \! \left({ t_b, x_b }\right) =
\int_{ \Pa_b }^{ }
\Da x \cdot e^{i S (x) / \hbar} 
\cdot \phi (x (t_a)) \, ,
\eqno ({\mathop{\rm IV.8 }\nolimits})$$
where $\Pa_b$ is the space of all paths $x : {\bf T} \ra \RR^d$ with endpoint at $x_b$, namely 
$x(t_b) = x_b$. The questionable part was the rigorous definition of the integrator 
$\Da x$.

\bigskip
To fit within our general framework, we consider translations acting on $\RR^d$, namely
$$x \cdot 
\s_{(\a) }(r) =
\left({ x^1, \cdots, x^{\a - 1}, x^\a + \l r , x^{\a + 1} , \cdots , x^d }\right) \, . 
\eqno ({\mathop{\rm IV.9 }\nolimits})$$

\bigskip
\nn
The corresponding Lie derivative is given by
$$\La_{X_{(\a )} } f = 
\l {\part f \over \part x^\a }
\eqno ({\mathop{\rm IV.10 }\nolimits})$$
and the parameter $\l$ is chosen equal to $(h/m)^{{1 \over 2}}$.
The general differential equation $dx = X_{(\a)} (x) \cdot d z^\a$ reduces to 
$d x^\a = \l d z^\a$ (for $\a$ in $\{1, \cdots, d \}$). Hence the solution
$$x(t, z) = x_b + \l z(t) 
\eqno ({\mathop{\rm IV.11 }\nolimits})$$
describes the parametrization of the space $\Pa_b$ of paths $x$ with $x(t_b) = x_b$ by the space
${\bf Z}_b$ of paths $z$ with $z(t_b) = 0$. Accordingly, we obtain for the action
$$- {1 \over s \hbar } 
S (x(\cdot, z)) = - 
{\pi \over s} \int_{ {\bf T} }^{ } dt \, 
\left\vert{ \dot z (t) }\right\vert_{}^{2} +
{1 \over s \hbar }  \int_{ {\bf T} }^{ } dt \, V (x(t, z)) \, . 
\eqno ({\mathop{\rm IV.12 }\nolimits})$$

\medskip
\nn
Substituting $z$ for the integration variable $x$ in equation (IV.8), 
we obtain the path integral (remember that $s = i$ in quantum mechanics)
$$\Psi \! \left({ t_b, x_b }\right) =
\int_{ {\bf Z}_b }^{ } \Da z \cdot \exp \!
\left({- { 1\over s \hbar }  S (x(\cdot, z)) }\right)
\cdot \phi (x (t_a, z)) \, . 
\eqno ({\mathop{\rm IV.13 }\nolimits})$$

\medskip
\nn
>From our general results in section II, the function $\Psi$ is a solution of the differential equation
$$
 {\part \Psi \over \part t} =
 {s \over 4 \pi } \sum_{ \a }^{ }
\La_{X_{(\a)} }^2 \Psi + { 1\over s \hbar }  V \Psi \, , 
\eqno ({\mathop{\rm IV.14 }\nolimits})$$
that is 
$$s \hbar  {\part \Psi \over \part t} =
{s^2 \hbar^2 \over 2m} \D \Psi + V \Psi \, . 
\eqno ({\mathop{\rm IV.15 }\nolimits})$$

\medskip
\nn
For $s = i$, this is the sought-for Schr\"odinger equation. The integrator $\Da z$ in the space
${\bf Z}_b$ is invariant under translations and is normalized by\footnote{$^{11}$}
{\nn The scaling factor $\l = (h/m)^{1/2}$ has been chosen in such a way that no physical constant enters
in this normalization. With the notations of paragraph A.3.6, we have the dimensional equations
$$
\left\lbk{ z^\a }\right\rbk = \Ta^{1/2}
\quad , \quad [t] = \Ta \, .$$
}
$$ 
\int_{ {\bf Z}_b }^{ } \Da z \cdot \exp \!
\left({ - { \pi \over s} \int_{ {\bf T} }^{ } dt 
\left\vert{ \dot z (t)  }\right\vert_{}^{2}
 }\right) = 1
\, . 
\eqno ({\mathop{\rm IV.16 }\nolimits})$$

\bigskip
We derive now the {\it semiclassical approximation}\footnote{$^{12}$}{
\nn Since $\l = (h/m)^{1/2}$, the semiclassical expansion will proceed according to the powers of
$h^{1/2}$. }.
We use the initial wave function 
$\phi (x) = \d \left({ x - x_a }\right)$ and reparametrize the paths {\it around the classical
path} $x_{{\mathop{\rm  cl }\nolimits}}$, that is
$$x(t, \z)  = x_{{\mathop{\rm  cl }\nolimits}} (t) + \l \z (t) \, , 
\eqno ({\mathop{\rm IV.17 }\nolimits})$$
with $\z$ in ${\bf Z}_b$.
Since the integrator in ${\bf Z}_b$ is invariant under translations, and 
$x_{{\mathop{\rm  cl }\nolimits}} (t_a) = x_a $, we transform equation (IV.13) into
$$\Psi \! \left({ t_b, x_b }\right) =
\int_{ {\bf Z}_b }^{ } \Da \z \cdot \exp \!
\left({ - { 1 \over s \hbar} 
S \! \left({ x_{{\mathop{\rm  cl }\nolimits}} + \l \z }\right)
 }\right) \d (\l \z (t_a))
\, . 
\eqno ({\mathop{\rm IV.18 }\nolimits})$$

\medskip
\nn
In expanding $S \! \left({ x_{{\mathop{\rm  cl }\nolimits}} + \l \z }\right)$ in powers of $\l$, there
is no term linear in $\l$, since $x_{{\mathop{\rm  cl }\nolimits}}$ is a critical point 
of the action
functional $S$ {\it and} because $\z (t_a) = \z(t_b) = 0$.
Hence we obtain
$$
\eqalignno{
- { 1 \over s \hbar} 
S \! \left({ x_{{\mathop{\rm  cl }\nolimits}} + \l \z }\right)
&=-
{ 1 \over s \hbar} 
S \! \left({ x_{{\mathop{\rm  cl }\nolimits}}  }\right)
- {\pi  \over s } \int_{ {\bf T}}^{ } 
dt \ \vert  \dot \z (t) \vert_{}^{2}
& ({\mathop{\rm IV.19 }\nolimits}) \cr
& \ \ \  + {\pi  \over ms  } \int_{ {\bf T}}^{ }  dt \, 
\bgt_\a \bgt_\b V \! 
\left({ x_{{\mathop{\rm  cl }\nolimits}} (t)  }\right)
\z^\a (t) \z^\b (t) + O(\hbar^{1/2} ) \, . \cr}
$$

\medskip
\nn
The action 
$S \! \left({ x_{{\mathop{\rm  cl }\nolimits}} }\right)$ corresponding to the classical path
$x_{{\mathop{\rm  cl }\nolimits}}$ with endpoints 
$x_{{\mathop{\rm  cl }\nolimits}} (t_a) = x_a$, 
$x_{{\mathop{\rm  cl }\nolimits}}(t_b) = x_b$ is nothing else than the 
{\it classical action function }
$\Sa \! \left({ t_b, x_b ; t_a, x_a }\right)$.

\bigskip 
Omitting terms of order $\hbar^{1/2}$, we obtain the WKB approximation \break
$\Psi_{{\mathop{\rm WKB }\nolimits}} 
\! \left({ t_b, x_b }\right)$ to $\Psi \! \left({ t_b, x_b }\right)$. 
Using formulas (IV.18) and (IV.19), we derive
$$
\Psi_{{\mathop{\rm WKB }\nolimits}} 
\! \left({ t_b, x_b }\right) =
\exp \!
\left({ - {1 \over s \hbar }
\Sa \! \left({ t_b, x_b ; t_a , x_a }\right) 
 }\right) \cdot I \, , 
\eqno ({\mathop{\rm IV.20 }\nolimits})$$
with the integral
$$I =
\int_{ {\bf Z}_b }^{ } \Da \z \cdot \exp \!
\left({ - {\pi  \over s}  Q_1 (\z) }\right)
\d \! \left({ \l \z (t_a) }\right)
\, . 
\eqno ({\mathop{\rm IV.21 }\nolimits})$$

\medskip
\nn
Besides the quadratic form
$$Q_0 (\z) = \int_{ {\bf T} }^{ } dt \, 
\, \vert \dot \z (t) \vert^2
\eqno ({\mathop{\rm IV.22 }\nolimits})$$
corresponding to the free particle, we need the quadratic form
$$Q_V (\z) = - {1 \over m} \int_{ {\bf T} }^{ } dt \, 
h_{\a \b } (t) \, \z^\a (t) \z^\b (t) 
\eqno ({\mathop{\rm IV.23 }\nolimits})$$
with 
$h_{\a \b } (t) = \bgt_\a \bgt_\b V \!
\left({  x_{ {\mathop{\rm cl }\nolimits} } (t)  }\right)$.
In (IV.21), we use the quadratic form $Q_1 = Q_0 + Q_V$.

\bigskip
The easiest method to calculate the functional integral $I$ consists of the following steps.

\medskip
a) {\it Changing the integrator}: according to formulas (A.126), (A.132) and (A.134) in
Appendix A, we obtain $I = I_1 I_2 $ where
$$
\eqalignno{
I_1 &=
\left\vert{ \Det \!
\left({ Q_0 / Q_1 }\right)
 }\right\vert_{}^{1/2}
s^{- {\mathop{\rm Ind  }\nolimits} \! \left({  Q_1 }\right) } \, ,
&({\mathop{\rm IV.24 }\nolimits})\cr
I_2 &=
\int_{ {\bf Z}_b  }^{ }
\Da_{s, Q_1} \z \cdot \exp
\! \left({-  {\pi  \over s } Q_1 (\z)  }\right)
\d \! \left({ \l \z (t_a) }\right) \, .
&({\mathop{\rm IV.25 }\nolimits})\cr} $$
Here ${\mathop{\rm Ind  }\nolimits} \! \left({  Q_1 }\right)$ is the number of negative
directions for the quadratic form $Q_1$.

\medskip
b) {\it Restricting the domain of integration}:
to treat the $\d$ factor in $I_2$, we use the linear change of variables 
$\z \mps \z (t_a) $ from ${\bf Z}_b$ to $\RR^d$. By a method similar to the one used in
paragraph A.3.8, we obtain 
$$I_2 =
\left({ \det s W_1 \! \left({ \d_{t_a}  }\right) }\right)_{}^{- 1/2}
\l^{- d } \, , 
\eqno ({\mathop{\rm IV.26 }\nolimits})$$
where $W_1$ is the quadratic form on ${\bf Z}'_b$ inverse to $Q_1$.
Using the Green's function $G$ given by equation (B.26), we evaluate the matrix
 $W_1 \! \left({ \d_{t_a}  }\right)$ as follows
$$W_1 \! \left({ \d_{t_a}  }\right) =
G \! \left({ t_a, t_a }\right) =
N \! \left({ t_a, t_b }\right)
J \! \left({ t_b, t_a }\right) \, .
\eqno ({\mathop{\rm IV.27 }\nolimits})$$

\medskip
c) {\it Reducing a functional determinant to a finite determinant}:
using the same strategy as in paragraph III.2 we obtain
$$\Det \! \left({ Q_0 / Q_1 }\right) =
\det \!
\left({ K \! \left({ 1 ; t_b, t_a  }\right)^{- 1 }
 K \! \left({ 0 ; t_b, t_a }\right)  }\right) \, . 
\eqno ({\mathop{\rm IV.28 }\nolimits})$$
Here 
$K \! \left({ \nu ; t, t_a }\right)$ is the Jacobi field defined by
$$Q_\nu K \! \left({ \nu ; t, t_a }\right) = 0 \quad , 
\quad 
K \! \left({ \nu ; t_a, t_a }\right)
= \un 
\eqno ({\mathop{\rm IV.29 }\nolimits})$$
for $\nu$ equal to $0$ or $1$.

\bigskip
Finally, collecting the previous equations (IV.20) to (IV.29) and using equations (B.11),
(B.12), (B.7), (B.8), (B.19) and (III.42), we obtain {\it the} WKB {\it approximation to the
point-to-point transition amplitudeJ} (in the case $s = i$):
$$
\left\langle{ t_b, x_b \bigm|  t_a, x_a }\right\rangle_{
 WKB 
}^{}
= c \, 
h^{- d/2}
\left\vert{ \det \part^2 \Sa /
\part x_a^\a \cdot \part x_b^\b }\right\vert_{}^{1/2}
e^{i \Sa / \hbar} \, . 
\eqno ({\mathop{\rm IV.30 }\nolimits})$$
Here 
$\Sa = \Sa \! 
\left({ t_b , x_b ; t_a, x_a }\right)$ is the classical action function and the phase factor 
$c$ is $e^{\pi i (p - q)/4}$ where $p$ $(q)$ is the number of positive (negative) eigenvalues of
the van Vleck-Morette matrix
$\left({ \part^2 \Sa / \part x_a^\a \cdot \part x_b^\b }\right)_{
{1 \leq \a \leq d \atop 1 \leq \b \leq d }}$.

\bigskip
For the higher-order terms in the semiclassical approximation, whether 
$x_{{\mathop{\rm cl}\nolimits}} $ is or is not a degenerate critical point of the action, we
refer the reader to the literature [5] or the references at the end of Appendix B. The reader
can easily transfer these old results using the simpler and more general formalism of this
paper.

\bigskip
\nn
{\bf 2. Polar coordinates.}\footnote{$^{13}$}{Contributed by John La Chapelle.}
\smallskip
Our basic integral is formulated in terms of a transformation 
$\Si(t, z)$ on a manifold $N$, in a form valid for an arbitrary system of coordinates.
Before considering the case of a general Riemannian manifold in paragraph IV.3, we consider
polar coordinates in a plane. The case of cylindrical coordinates in $\RR^3$ is very similar.

\bigskip
``Path integrals in polar coordinates'', the  1964 paper by S.F. Edwards and Y.V.
Gulyaev [17] has been, and still is, at the origin of many investigations. To the best of our
knowledge, all the papers on the subject deal with the discretized version of the path
integral, and propose various path integral prescriptions when the discretized values of the
paths are expressed in coordinates other than cartesian.

\bigskip
The basic manifold is $N = \RR^2 \bl \{ 0 \}$ with coordinates 
$x^1 $ and $x^2$. We also consider the manifold $\ww N =
]0 , + \Inf [ \ts \RR$ with coordinates $r, \t$ (sole restriction $r > 0$) and the covering map 
$\Pi : \ww N \ra N$ taking $(r, \t) $ into 
$ \left({ x^1, x^2 }\right)$ with
$$x^1 = r \cos \t \quad , \quad  x^2 = r \sin \t \, . 
\eqno ({\mathop{\rm IV.31 }\nolimits})$$
Two points of $\ww N$ map onto the same point of $N$ if, and only if, their $\t$-coordinates differ
by an integral multiple of $2 \pi$.

\bigskip
Let $\D = 
{\part^2 \over (\part x^1)^2 } + 
{\part^2 \over (\part x^2)^2 }$ be the Laplacian in cartesian coordinates.
We know how to solve the Schr\"odinger equation
$${\part \Psi \over \part t} =
{s \over 4 \pi } \D \Psi
\eqno ({\mathop{\rm IV.32 }\nolimits})$$
by means of the path integral
$$\Psi \! \left({ t_b, x_b }\right)
= \int_{{\bf Z}_b }^{ }
\Da z \cdot e^{ - \pi Q_0 (z) / s} 
\phi \! \left({  x(t_a, z)  }\right)
\eqno ({\mathop{\rm IV.33 }\nolimits})$$
provided that $x(t, z)$ is a solution of the differential system
$$dx^1 = dz^1 \quad , \quad dx^2 = dz^2
\eqno ({\mathop{\rm IV.34 }\nolimits})$$
with boundary condition 
$x \! \left({ t_b, z }\right) = x_b$.
To solve the Schr\"odinger equation in polar coordinates, we need only to transform the system
(IV.34) in polar coordinates with the help of the transformation equations (IV.31). Hence
$$
\left\{\matrix{
dr 
= \cos \t \cdot d z^1 + \sin \t \cdot dz^2 =
X_{(1)}^1 dz^1 + X_{(2)}^1 dz^2
  \hfill\cr
\noalign{\medskip}
d \t = - 
 {\ds \sin \t  \over \ds r} \cdot dz^1 +
   {\ds \cos \t  \over \ds r} \cdot dz^2 =
X_{(1)}^2 dz^1 + X_{(2)}^2 dz^2 \, . 
  \hfill \cr}
\right. \eqno ({\mathop{\rm IV.35 }\nolimits})
$$
The vector fields $X_{(1)}$ and $X_{(2)}$ can be read off from the above equations:
$$\La_{X_{(1)}}
=
\cos \t \cdot {\part  \over \part r} -
{\sin \t \over r} \cdot {\part  \over \part \t} \,  ,
\eqno ({\mathop{\rm IV.36 }\nolimits})$$
$$\La_{X_{(2)}}
=
\sin \t \cdot {\part  \over \part r} +
{\cos \t \over r} \cdot {\part  \over \part \t} \, .
\eqno ({\mathop{\rm IV.37 }\nolimits})$$
Hence, if 
$(r(t, z), \t(t, z))$ is the solution of the differential system 
(IV.35) such that $r(t_b, z) = r_b$,
$\t(t_b, z) = \t_b$, the path integral
$$\Psi
\! \left({ t_b, r_b, \t_b }\right) : = 
\int_{ {\bf Z}_b }^{ } \Da z
\cdot e^{- \pi Q_0 (z) / s}
\phi (r(t_a, z), \t(t_a, z)) 
\eqno ({\mathop{\rm IV.38 }\nolimits})$$
solves the Schr\"odinger equation (IV.32). As expected, the operator
$$\D = 
\La_{X_{(1)}}^2 + \La_{X_{(2)}}^2
=
{ \part^2\over \part r^2}
+ {1 \over r^2} {\part^2 \over \part \t^2} +
{1 \over r} {\part  \over \part r}
$$
is the Laplacian on $\RR^2 \bl \{ 0 \}$ in polar coordinates.

\bigskip
We want to evaluate the point-to-point  transition amplitudes in polar coordinates, denoted
by $\left\langle{ t_b, r_b, \t_b \bigm| t_a , r_a, \t_a  }\right\rangle$. To obtain them, it suffices to put
$\phi (r, \t) = \d (r - r_a) \d (\t - \t_a)$ in equation (IV.38).
 Solving the system (IV.35) is easy by reverting
to cartesian coordinates, hence
$$
\left\{\matrix{
r(t, z) \cos \t (t, z) = r_b \cos \t_b + z^1 (t)
 \hfill\cr
\noalign{\medskip}
  r(t, z) \sin \t (t, z) = r_b \sin \t_b + z^2 (t) \, . 
  \hfill \cr}
\right. \eqno ({\mathop{\rm IV.39 }\nolimits}) 
$$
For given $x^1$, $x^2$, the equations (IV.31) in $r, \t$ have infinitely many solutions, and we derive easily
$$\d \!
\left({ r \cos \t - r_a \cos \t_a  }\right)
\d \! \left({ r \sin \t - r_a \sin \t_a }\right)
= 
{1 \over r_a}
\sum_{ n \in \ZZ}^{ }  
\d \! \left({ r - r_a  }\right) 
\d \! \left({ \t - \t_a - 2n \pi  }\right) \, .
\eqno ({\mathop{\rm IV.40 }\nolimits}) $$
Substituting $r(t, z)$ to $r$ and $\t (t, z)$ to $\t$ and taking into account equations (IV.39), we obtain the
path integral representation
$$
{1 \over r_a}
\sum_{ n \in \ZZ}^{ }  
\left\langle{t_b, r_b, \t_b \bigm| t_a, r_a , \t_a + 2 n \pi  }\right\rangle
= I \, . 
\eqno ({\mathop{\rm IV.41 }\nolimits}) $$

\medskip
\nn
Here we use the definitions
$$\phi_1 (u) =
\d \! \left({ r_b \cos \t_b - r_a \cos \t_a - u^1 }\right)
\d \! \left({ r_b \sin \t_b - r_a \sin \t_a - u^2 }\right)
\eqno ({\mathop{\rm IV.42 }\nolimits}) $$
for $u = \left({ u^1, u^2 }\right)$ in $\RR^2$, and
$$ I = \int_{ {\bf Z}_b }^{ }
\Da z \cdot e^{- \pi Q_0 (z) / s }
\phi_1 (z(t_a)) \, . 
\eqno ({\mathop{\rm IV.43 }\nolimits}) $$
This integral is nothing else than a point-to-point transition amplitude in cartesian coordinates.
We give now a derivation of the well-known formula for this amplitude, by relying on our methods.

\bigskip
The integrand 
in equation (IV.43) depends only on 
$z^1 (t_a) $, $z^2 (t_a)$, hence as explained in paragraph A.2.2, we introduce the linear map
$L : {\bf Z}_b \ra \RR^2 $ mapping $z$ into $z(t_a)$.
The image of the integrator 
$\Da z \cdot e^{- \pi Q_0 (z) / s }$ on ${\bf Z}_b $ is a Gaussian integrator on $\RR^2$; we
need only the value of the corresponding covariance matrix
$$W^{\a \b } =
\int_{ {\bf Z}_b }^{ }
\Da z \cdot e^{- \pi Q_0 (z) / s } \, 
z^\a (t_a) z^\b (t_a) \, .  
\eqno ({\mathop{\rm IV.44 }\nolimits})$$

\medskip
\nn
According to formulas (A.57) and (A.62), we obtain
$$W^{\a \b } =
{s \over 2 \pi }
G_b^{\a \b }\! \left({ t_a, t_a }\right)
=
{s \over 2 \pi } \d^{\a \b }
\! \left({ t_b - t_a }\right)
\, .
\eqno ({\mathop{\rm IV.45 }\nolimits})
$$

\medskip
\nn
According to formulas (A.27) and (A.32), we transform the integral (IV.43) into
$$I =
\left (s \left({ t_b - t_a }\right)\right)^{-1}
\int_{\RR }^{ } du^1 \, \int_{\RR }^{ } du^2 
\exp \! 
\left({- {\pi  \over s} 
{ \left\vert{ u }\right\vert_{}^{2} \over t_b - t_a }
}\right)
\phi_1 (u) \, . 
\eqno ({\mathop{\rm IV.46 }\nolimits})$$

\medskip
\nn
The function $\phi_1 (u) $ is a $\d$-factor and the final result is
$$
\eqalign{I 
&= \left(s\left({ t_b - t_a }\right)\right)^{-1} 
\exp \! 
\left({ -{ \pi \over s}
{r_a^2 + r_b^2 \over t_b - t_a} }\right)
\exp \! 
\left({ { 2\pi \over s}
{r_a  r_b \cos (\t_b - \t_a) \over t_b - t_a} }\right)
\cr
&
=
\left (s \left({ t_b - t_a}\right)\right) ^{-1}  
\exp \! 
\left({- {\pi  \over s} 
{ \left\vert{ x_b - x_a }\right\vert_{}^{2} \over t_b - t_a }
}\right) \, . \cr} 
\eqno ({\mathop{\rm IV.47 }\nolimits})$$

\bigskip
>From the equations (IV.41) and (IV.47) we cannot derive directly the point-to-point transition amplitudes
in polar coordinates. We defer to paragraph IV.7 a further discussion of this point.

\bigskip
\nn
{\bf 3. Frame bundles over Riemannian manifolds.}\footnote{$^{14}$}{\nn
For a discussion of functional integrals when paths take their values in a Riemannian space, see
[5], [14] and [18]. }

\smallskip
We consider a Riemannian (or pseudo-Riemannian) manifold $M$ of dimension $d$, with metric
$g$, and the orthonormal frame bundle $N = O(M)$ over $M$, with projection 
$\pi : O(M) \ra M$. We want to choose vector fields 
$X_{(1)}, \cdots, X_{(d)}$ on $O(M)$ such that the {\it equation on the bundle}
$O(M)$ 
$${ \part \over \part t_b} 
\Psi \! \left({ t_b, \r_b }\right) =
{s \over  4 \pi} h^{\a \b } X_{(\a)} X_{(\b)}
\Psi \! \left({ t_b, \r_b }\right) \, ,
\eqno ({\mathop{\rm IV.48 }\nolimits})$$
where $\r_b$ belongs to $O(M)$ and the first-order differential operators are given by
$X_{(\a)} = X_{(\a)}^\l (\r_b) 
{\part  \over \part \r_b^\l}$, gives by projection {\it an equation on the base space $M$},
namely:
$$
\eqalign{
{\part  \over \part t_b } \psi \! 
\left({ t_b, x_b  }\right) 
&= { s \over  4 \pi } g^{\l \mu} D_\l D_\mu \psi \! 
\left({ t_b, x_b  }\right)  \cr
& 
= { s \over  4 \pi } \D \psi \! 
\left({ t_b, x_b  }\right)  \, . \cr}
\eqno ({\mathop{\rm IV.49 }\nolimits})$$
Here $D_\l$ is the covariant derivative defined by the Riemannian connection, and 
$\Psi = \psi \circ \pi$; moreover $\D$ is the Laplace-Beltrami operator on $M$.

\bigskip
It has been shown in [13] that the covariant Laplacian $\D$ at a point $x_b$ of $M$ can be lifted
to a sum of products of Lie derivatives $h^{\a \b } X_{(\a)} X_{(\b)}$ at the frame 
$\r_b $ in $O(M)$; the integral curves of the set of vector fields
$\left\lbc{ X_{(\a)} }\right\rbc$ starting from $\r_b$ at time $t_b$ are the horizontal lifts
of a set of geodesics at $x_b$, tangent to the basis 
$\left\lbc{ e_{\a} }\right\rbc$ of $T_{x_b} M$ corresponding to the frame $\r_b$.
The constant matrix $\left({ h_{\a \b} }\right)$ has been chosen with the same signature as
the metric $g$ on $M$, and $g\! \left({ e_\a, e_\b }\right) = h_{\a \b }$.

\bigskip
An explicit construction of $\left\lbc{ X_{(\a)} }\right\rbc$ goes as follows. Let $\r(t) $ be
the horizontal lift of a path $x(t) $ in $M$ defined by the Riemannian connection map\footnote{$^{15}$}{
More explicitly, for a given frame $\r$ at a point $x$ of $M$, $\s (\r)$ is a linear map from 
$T_x M$ to $T_\r O(M) $, mapping $u$ into $\s (\r) \cdot u$ for $u$ in $T_x M$.}
$\s : O (M) \ra L (TM, TO (M))$; it satisfies the differential equation
$$\dot \r (t) = \s (\r (t)) \cdot \dot x (t) \quad, \quad
\r(t_b) = \r_b \, .
\eqno ({\mathop{\rm IV.50 }\nolimits})$$

\medskip
\nn
If we put $\r(t) = 
(x(t), u(t))$ and $\r_b = \left({ x_b, u_b }\right)$, a solution $\r(t)$ of the previous
equation corresponds to the frame $u(t) $ obtained by parallel transport of $u_b$,  along the
path $x$ from $x_b$ to $x(t)$. The frame $u(t)$ is also an admissible map
$$u(t) : \RR^d \ra T_{x(t)} M$$
i.e. $u(t)$ maps a $d$-tuple into a vector whose components in $u(t)$ are the chosen $d$-tuple;
equivalently
$\left({ u(t)^{- 1} \dot x (t)  }\right)_{}^{\a}$
is the $\a$-coordinate of $\dot x (t) $ in the $u(t) $ frame.
Set 
$$\dot z (t) : = u(t)^{-1}
\dot x (t) = \dot z^\a (t) \oo e_\a \, , 
\eqno ({\mathop{\rm IV.51 }\nolimits})$$
where $\left\lbc{ \oo e_\a }\right\rbc$ denotes the canonical basis of the model space $\RR^d$.
Then we can express (IV.50) in the canonical form (II.21)
$$\dot \r (t) = X_{(\a)} \!
\left({ \r(t)   }\right) \dot z^\a (t)
\quad, \quad \r (t_b) = \r_b 
\eqno ({\mathop{\rm IV.52 }\nolimits})$$
where $X_{(\a)}$ is defined by 
$$X_{(\a)} (\r) =
\s(\r) \cdot e_\a = 
\left({ \s (\r) \circ u }\right) \cdot \oo e_\a \, . 
\eqno ({\mathop{\rm IV.53 }\nolimits})$$
Here $x$ is a point of $M$ and $\r = (x, u) $ a frame, where
$u : \RR^d \ra T_x M$ is an admissible map.

\bigskip
If $\dot z(t) = \oo e_\a $, i.e. $\dot z^\b (t) = \d_\a^\b $, the coordinates of $\dot x(t) $ are
constant in the frame $u(t) $ parallel transported along $x(t)$. Therefore $x(t)$ is the
geodesic defined by 
$$x(t_b) = x_b \quad , 
\quad \dot x (t_b) = e_\a \, . 
\eqno ({\mathop{\rm IV.54 }\nolimits})$$
With 
$\dot z (t) = \oo e_\a$, equation (IV.52) reads 
$$\dot \r_{(\a)} (t) =
X_{(\a) } \! \left({ \r_{(\b)} (t) }\right)
\d_\a^\b = 
X_{(\a) } \! \left({ \r_{(\a)} (t) }\right)\, . $$

\medskip
\nn
{\it The horizontal lift $\r_{(\a)} (t)$ of the geodesic} (IV.54) {\it is, as desired, the
integral curve of 
$X_{(\a)}$ going through $\r_b$
at time $t_b$.
}
With $X_{(\a)}$ defined by (IV.53), $\r$ can be expressed in terms of the {\it Cartan
development map
}
$$(\pi \circ \r)(t) = x(t) =
(\Dev z)(t) \, . 
\eqno ({\mathop{\rm IV.55 }\nolimits})$$
The Cartan development is a bijection from a space of pointed paths (paths with a fixed end
point) on $T_{x_b} M$ (identified to $\RR^d$ via the frame $\r_b$) into a space of pointed paths on $M$ -- or
vice versa.  Here $$\Dev : \Pa_0 T_{x_b} M \ra  \Pa_{x_b} M \quad \by \quad z \mps x \, .
\eqno ({\mathop{\rm IV.56 }\nolimits})$$
The path $x$ is said to be the development of $z$, if $\dot x (t)$ parallel transported along $x$ from
$x(t)$ to $x_b$ is equal to $\dot z(t) $ trivially transported to the origin of $T_{x_b} M$,
for every $t \in {\bf T}$.

\bigskip
The path integral solution of (IV.48) is 
$$\Psi
\! \left({ t_b, \r_b  }\right) =
\int_{ {\bf Z_b}J}^{ }
\Da z \cdot \exp
\! \left({ - {\pi  \over s } Q_0 (z)
}\right)  \Phi 
\! \left({ \r (t_a, z) }\right)
\eqno ({\mathop{\rm IV.57 }\nolimits})$$
with $\r (t, z)$ solution of (IV.52). The path integral solution of (IV.49) is 
$$
\psi
\! \left({ t_b, x_b  }\right)
=
\int_{ {\bf Z_b}J}^{ }
\Da z \cdot \exp
\! \left({- {\pi  \over s} Q_0 (z)}\right)
\phi ((\Dev z)(t_a)) \, . 
\eqno ({\mathop{\rm IV.58 }\nolimits})$$
If a scalar potential is desired in (IV.48) or (IV.49), one can proceed along either of the
methods outlined in paragraph II.4.

\bigskip
The semiclassical approximation [5] of $\Psi$, with or without scalar potential, is considerably
more complicated to compute than the semiclassical approximation of $\Psi$ given in section
III. For an initial wave function $\phi$ of type (III.2) the blue-print given in Section III.2
is complete. If one wishes to compute the point-to-point propagator on a Riemannian space $M$,
one chooses the initial wave function on $M$ to be
$$\phi (x) = \d_{x_a } (x) \, . 
\eqno ({\mathop{\rm IV.59 }\nolimits})$$
The detailed  calculation can be found in reference [5, pp. 309-311].
The development map cannot parametrize spaces of paths with two fixed points\footnote{$^{16}$}{
\nn For instance, two geodesics on $S^2$ intersect at two antipodal points; they are the
developments of two halflines with one common origin.}
but (IV.59) with $(\Dev z) (t_a)$ substituted to $x$ restricts the domain of integration appropriately.

\bigskip
The following remarks simplify the calculations. 

\bigskip
\nn
{\it Remark 1.}
According to the formulas
$$
\eqalign{
\int_{{\bf T} }^{ } dt \, 
h_{\a \b } \dot z^\a (t) \dot z^\b (t) 
&= 
\int_{{\bf T} }^{ } dt \, 
g_{\mu \nu } (x(t)) \dot x^\mu (t) \dot x^\nu (t)\cr
& 
= \int_{{\bf T} }^{ } dt \, 
g_{\mu \nu} (x(t)) \dot \r^\mu (t) \dot \r^\nu (t) \, ,\cr}
\eqno ({\mathop{\rm IV.60 }\nolimits})$$
development map and horizontal lift preserve lengths and angles.

\bigskip
\nn
{\it Remark 2.}
If $z$ develops into a classical paths 
$x_{{\mathop{\rm cl }\nolimits}}$, the determinant of the derivative mapping 
$\Dev'(z)$ is unity. For the proof see, for instance, reference [5], p. 308.

\bigskip
\nn
{\bf 4. A multiply connected manifold.}

\smallskip
The domain of integration, a space of pointed paths $\Pa_x N$, is the union of disjoint sets
made of paths in different homotopy classes. We recall in paragraphs a) and b) earlier
calculations of propagators on multiply connected spaces.
 Then, in paragraph c), we explain how these methods fit into our general framework.

\bigskip
a) It was shown [3] in 1971 that, for a system with a multiply connected configuration space
$N$, the propagator $K$ is a linear combination of propagators $K_{(\a)}$
$$
\left\vert{ K^A }\right\vert
=
\left\vert{ \sum_{ g_\a \in \pi_1 (N) }^{ }  
\chi^A (g_\a) K_{(\a)} }\right\vert \, . 
\eqno ({\mathop{\rm IV.61 }\nolimits})$$

\medskip
\nn
Each $K_{(\a)}$ is obtained by summing over paths in the same homotopy class, say $\a$; the set
$\left\lbc{ \chi^A (g_\a }\right\rbc_{_\a}^{}$ forms a representation, labeled $A$, of the
fundamental group $\pi_1 (N)$. Since a homotopy class cannot be identified uniquely with an
element of $\pi_1 (N)$, the propagator is defined modulo an overall phase factor. There are as
many propagators $K^A$ as there are inequivalent representations of $\pi_1 (N)$.

\medskip
The proof of (IV.61) uses two facts:

\smallskip
i) the superposition principle of quantum propagators implies the linear combination of
partial propagators;

\smallskip
ii) the fundamental group based at a point is isomorphic to the fundamental group based at
another point, but not canonically so. Therefore the pairing 
$\left({ g_\a, K_{(\a)} }\right)$ is done by choosing an homotopy mesh (choosing a point for
the fundamental group, and pairing one group element with one homotopy class), then requiring
that the result be independent of the homotopy mesh.

\bigskip
b) Later on [14, see also 18, p. 65] the same result (worked out for a different example) was
obtained from stochastic processes on fibre bundles. 
The basic steps are as follows\footnote{$^{17}$}{We refer the reader to paragraph IV.7 for an explicit example
where we use this strategy.}.

\smallskip
i) A universal covering $\ww N$ is a principal $G$-bundle with projection \break $\Pi : \ww N\ra N$,
where $N = \ww N / G$ and $G$ is a discrete group of automorphisms of $\ww N$ isomorphic to
the fundamental group of $N$. For example, $N = S^1 = \RR / \ZZ$ and $\ww N = \RR$ is a
$\ZZ$-principal bundle over $S^1$.

\smallskip
ii) The wave function for a system with configuration space $N$ is a section of a bundle {\it
weakly} associated to $\ww N$. i.e. a bundle whose typical fibre is associated to a not
necessarily faithful representation of $G$, hence to a representation of a group homomorphic
to $G$, says $G_0$. 
For example if $N = S^1$, a vector bundle over $S^1$ with structure group $U(1)$ is weakly
associated to a $\ZZ$-principal bundle over $S^1$ by a homomorphism 
$h_\a$ of $\ZZ$ into $U(1)$ mapping $n$ into $e^{in\a}$.

\smallskip
iii) There is a unique connection on $\ww N$: the horizontal lift $\ww x$ of a pointed path
$x$ with fixed point $\ww x (t_0) = \ww x_0$ is uniquely defined by the lift $\ww x_0 $ of
$x_0$. This unique connection defines the parallel transport of the wave function $\phi (x
(t))$ back to $x_0$.

\smallskip
Consider for instance the case where $\ww N = \RR$, $N = S^1$ and 
$\Pi (x) = e^{ix}$. Take for $\phi$ a section of a $U(1)$-bundle defined by the homomorphism
$h_\a$ as above. Let $x(t)$ be a map into $S^1$, lifted to a map $\ww x_0 (t)$ into $\RR$ in such a way that 
$x(t) = e^{i \ww x_0 (t)}$. The other liftings are given by
$$\ww x_k (t) = \ww x_0 (t) + 2 \pi k
\eqno ({\mathop{\rm IV.62 }\nolimits})$$
for $k$ in $\ZZ$.
The parallel transport of $\phi 
\! \left({ \ww x_k (t) }\right)
$ to $x_0$ is given by the formulas
$$\tau_{t_0}^t \phi 
\! \left({ \ww x_k (t) }\right) =
\ww h_\a (\ww x_0) \ww h_\a
\! \left({ \ww x_k (t) }\right)_{}^{-1}
\phi (x(t)) 
\eqno ({\mathop{\rm IV.63 }\nolimits})$$
and
$$
\ww h_\a
\! \left({ \ww x_k (t) }\right)_{}^{-1} =
\ww h_\a
\! \left({ \ww x_0 (t) }\right)_{}^{-1}
h_\a (k)^{-1}
= e^{- i k \a }
\ww h_\a
\! \left({ \ww x_0 (t) }\right)_{}^{}\, . 
\eqno ({\mathop{\rm IV.64 }\nolimits})$$
The map $\ww h_\a$ is the bundle map from the $\ZZ$-bundle $\RR $ over $S^1$ to the
$U(1)$-bundle corresponding to the map $h_\a : \ZZ \ra U (1)$.

\smallskip
iv) Summing the wave function $\tau_{t_0}^t \phi \! \left({ \ww x_k (t) }\right)$
over  all $k$'s gives the solution $\Psi (t, x_0)$ of a parabolic
equation with initial value $\Psi (t_0, x_0) = \phi (x_0)$.

\bigskip
In both computations outlined above, one chooses a representation of the fundamental group,
$\left\lbc{ \chi^A (g_\a) }\right\rbc$, or the homomorphism $h : G \ra G_0$.

\bigskip
c) We specialize our basic formula (II.1)
$$
\Psi \! \left({ t_b, x_b }\right)
= \int_{ {\bf Z }_b }^{ }
\Da z 
\cdot e^{- \pi Q_0 (z) / s}
\phi \! \left({ x_b \cdot \Si (t_a, z) }\right) 
\eqno ({\mathop{\rm IV.65 }\nolimits} )$$
to the case where $N$ is multiply-connected.
To calculate the point-to-point transition amplitudes, we select $\phi (x)$ of the form $\d_{x_a} (x)$
with a $\d$-factor centered at a point $x_a$ of $N$.
Denote the evaluation map taking $z$ into $x_b \cdot \Si (t_a, z)$ by 
$\ve : {\bf Z}_b \ra N$.
Since ${\bf Z}_b $ is contractible, we can lift $\ve$ to a map
$\ww \ve : {\bf Z}_b \ra \ww N$ into the universal covering $\ww N $ of $N$ and hence
$\ve = \Pi \circ \ww \ve$. In the path integral (IV.65), the domain of integration is
restricted by the $\d$-factor 
$\d_{x_a} (x) $ to the inverse image 
$\ve^{-1}(x_a)$. It consists of paths such that
$x_b \cdot \Si (t_a, z) = x_a$ and {\it splits as the union of domains} 
${\bf Z}_\G = \ww \ve^{-1} \! \left({ \ww x_\G }\right)$ where $\ww x_\G$ runs over the various
points of $\ww N$ mapping to $x_a$ by $\Pi$. The labels $\G$ correspond to the various
homotopy classes of paths 
$x : {\bf T} \ra N$ such that 
$x(t_a) = x_a$, $x(t_b) = x_b$.
Hence the integral (IV.65) splits into a sum of integrals over the various subdomains
$$\Psi \! \left({ t_b, x_b }\right)
= \sum_{ \G }^{ } \int_{ {\bf Z }_\G }^{ }
\Da z 
\cdot e^{- \pi Q_0 (z) / s} \d_{x_a}
 \! \left({ x_b \cdot \Si (t_a, z) }\right) \, .
\eqno ({\mathop{\rm IV.66 }\nolimits} )$$

\medskip
This equation is the justification of the heuristic idea used in (IV.61) that the building
blocks of $K$ are the propagators 
$K_{(\a)}$ obtained by summing over paths in the same homotopy class.
To explain the coefficients $\chi^A (g_\a)$ we can proceed as follows: using the previous
notations 
$\ww N $, $\Pi$, $G$, consider a homomorphism $\chi^A$ of $G$ into $U(1)$.
The corresponding wave functions are functions $\ww \phi $ on $\ww N$ such that
$\ww \phi (  \ww x  g ) = \chi^A (g) \ww \phi (\ww x)$ for $g$ in $G$ and 
$\ww x $ in $\ww N$. We denote by 
$\ww \Si (t, z) $ the lifting of $\Si (t, z) $ to $\ww N$ and generalize equation (IV. 65)
by
$$\ww \Psi \! \left({ t_b, \ww x_b }\right)
=  \int_{ {\bf Z }_b }^{ }
\Da z 
\cdot e^{- \pi Q_0 (z) / s} \ww \phi
 \! \left({\ww  x_b \cdot \ww \Si (t, z) }\right) \, .
\eqno ({\mathop{\rm IV.67 }\nolimits} )$$
The function $\ww \Psi$ will satisfy the same transformation property as $\ww \phi$ and when $\ww
\phi$ is a $\d$-factor, we can split the integration domain into subdomains ${\bf Z}_\G$
as above. We recover then the results derived in paragraphs IV 4a and b.

\bigskip
\nn
{\bf 5. Gauge fields.}

\smallskip
We begin with the case of an abelian gauge group. 
In physical terms, we consider a particle of mass $m$ and electric charge $e$ moving  under the influence of a
magnetic potential $A$, with components $A_\a (x)$ at the point $x$. We consider generally a
$d$-dimensional space $\RR^d$ in cartesian coordinates $x^\a$ ($\a \in \{1, \cdots, d \}$) and 
metric $\left\vert{ x }\right\vert^2 = \d_{\a \b}
 x^\a x^\b $. The classical Lagrangian is given by
$$ L (x, \dot x) =
{m \over 2} \left\vert{ \dot x }\right\vert^2
+ e \, A_\a (x) \dot x^\a \, ,
\eqno ({\mathop{\rm IV.68 }\nolimits} )$$
hence the action functional
$$S(x) = 
\int_{{\bf T} }^{ } dt \ 
 L(x , \dot x) = {m \over 2} \int_{{\bf T} }^{ }
{ \left\vert{ dx }\right\vert^2 \over dt}  + e 
\int_{{\bf T} }^{ } A_\a d x^\a \, . 
\eqno ({\mathop{\rm IV.69 }\nolimits} )$$

\medskip
\nn
The equation of motion can be derived from this Lagrangian, and can be put into the
Hamiltonian form with the following definitions
$$p_\a = m \dot x^\a + e A_\a \quad , \quad 
H = \left\vert{ p - eA }\right\vert^2 / 2m \, . 
\eqno ({\mathop{\rm IV.70 }\nolimits} )$$
The corresponding Schr\"odinger equation is obtained in the standard way by replacing $p_\a $ by the operator 
${\hbar  \over i} 
{\part  \over \part x^\a} $
in the definition of $H$, and reads as
$$i \hbar { \part \psi \over \part t} 
= {1 \over 2m } 
\sum_{ \a }^{ }
\left({ { \hbar \over i}
{ \part \over \part x^\a} - e A_\a (x)
 }\right)^2 \psi \, . 
\eqno ({\mathop{\rm IV.71 }\nolimits} )$$
Our goal in this paragraph is to fit the well-known path integral solution of this equation
into our general framework.

\bigskip
It has long been recognized [19] that it is desirable to treat $\psi$ as a section of a complex
line bundle (here over $\RR^d$) associated to a principal $U(1)$-bundle via the canonical
representation of $U(1)$ acting on $\CC$ by multiplication. We describe the main steps of this
construction.

\vfill\eject
\nn
5.1. {\it The invariant formalism.}
\smallskip

The {\it base space} is $M = \RR^d$.
The {\it gauge group} $G$ is the set $U(1)$ of complex numbers of modulus one
$g = e^{i \T}$.
We consider a {\it principal bundle } $P$, with $G$ acting from the right via
$(p, g) \mps p \cdot g$, and projection $\Pi : P \ra M$.

\bigskip
The {\it connection } is a differential form $\om$ on $P$ with the transformation rule
$$\om (p \cdot g) = \om (p) + g^{- 1} dg \, . 
\eqno ({\mathop{\rm IV.72 }\nolimits} )$$
With the angular coordinate $\T$ such that 
$g = e^{i \T}$, one obtains $g^{- 1} dg = id \T $ for the invariant differential form on 
$U(1)$, hence the Lie algebra  ${\hbox{\tenfm g}} $ of $U(1) $ is naturally identified with
the set of pure imaginary numbers, and $\om$ is pure imaginary.

\bigskip
For any path $x : {\bf T} \ra M$ and any point $p_b$ in $P$ with 
$\Pi (p_b) = x(t_b)$, the {\it horizontal lifting}
is a curve $\xi : {\bf T} \ra P$ satisfying the following conditions:
$$\Pi \xi (t) = x(t) \quad, \quad \xi (t_b) = p_b 
\eqno ({\mathop{\rm IV.73 }\nolimits} )$$
$$\left\langle{  \om_{\xi (t) }, \dot \xi (t) }\right\rangle
= 0 \, ,
\eqno ({\mathop{\rm IV.74 }\nolimits} )$$
where $\om_{\xi (t) }$ is the value of $\om$ at the point
$\xi (t) $ of $P$, that is a linear form on the tangent space
$T_{\xi (t) } P$.

\bigskip
Let $L$ be the {\it associated line bundle}.
For any point $x$ in $M$, a point $p$ of $P$ with $\Pi (p) = x$ corresponds to an admissible
map
$\hh p : \CC \ra L_x $ where $L_x $ is the fiber of $L$ above the point $x$.
If $\psi$ is a {\it section} of $L$ over $M$, its value at point $x$ is an element $\psi (x) $
of $L_x$, hence $\hh p^{- 1} (\psi (x))$ is a complex number $\Psi (p)$. In this way  (see
e.g. [20], vol. I, p. 404), we identify the section $\psi $ of $L$ to a function $\Psi : P
\ra \CC$ with the transformation rule
$$\Psi (p \cdot g) = g^{- 1} \cdot \Psi (p) 
\eqno ({\mathop{\rm IV.75 }\nolimits} )$$
for $p$ in $P$ and $g$ in $G$.

\bigskip
\nn
5.2. {\it Fixing the gauge.}

\smallskip
A fixing of the gauge corresponds to a section $s : M \ra P$ of the principal bundle. When $s$
is chosen, we may identify $P$ with $M \ts G$ in such a way that $p =
s(x) \cdot g $ in $P$ correspond to $(x, g)$ in $M \ts G$.
A section $\psi$ of the line bundle $L$ corresponds now to  a {\it wave function},
that is to a complex-valued function on $M$ and the corresponding function $\Psi $ on $M \ts
G$ is given by $\Psi (x, g) = g^{- 1} \psi (x) $ and conversely 
$\psi (x) = \Psi (x, 1)$, or intrinsically $\psi = \Psi \circ s$. 

\bigskip
The differential form $\om$ on $P$ gives by pull-back via $s : M \ra P$ a purely imaginary
differential form on $M$, to be written as 
$- { ie\over \hbar} A$ to fit with standard physical dimensions. Hence the differential form $A$
on $M$ can be written as $A_\a (x) dx^\a$ and the functions $A_\a (x)$ are the components of
the magnetic potential. On $M \ts G$ the differential form $\om$ is given by
$$\om = g^{- 1} dg - 
{ ie\over \hbar} A =
i 
\left({ d \T - {e \over \hbar} A_\a dx^\a }\right)
\eqno ({\mathop{\rm IV.76 }\nolimits} )$$
(for $g = e^{i \T}$).

\bigskip
Let $x : {\bf T} \ra M$ be a path.
The horizontal lifting $\xi $ of $x$ is now described by 
$ \xi (t) = 
\left({ x(t) , e^{i \T (t) } }\right)$ and since $\om$ induces a zero form on the image $\xi ({\bf T})
\sbs M \ts G$, we obtain the differential equation
$$\dot \T =  {e \over \hbar}
A_\a (x(t)) \dot x^\a (t)
\eqno ({\mathop{\rm IV.77 }\nolimits} )$$
(see e.g. [18, pp. 64-65]).

\bigskip
Changing the gauge corresponds to choosing another section 
$s_1 : M \ra P$. There exists then a function $R : M \ra \RR$ such that 
$s_1 (x) = s(x) \cdot e^{- i e R (x) / \hbar}$.
In the new gauge, the section of the line bundle $L$ corresponds to a new wave function 
$$\psi_1 (x) = e^{ i e R (x) / \hbar} \psi (x) 
\eqno ({\mathop{\rm IV.78 }\nolimits} )$$
and the new components of the magnetic potential are given by 
$$A_\a^1 (x) = A_\a (x) +
{\part R (x) \over \part x^\a}
\, . 
\eqno ({\mathop{\rm IV.79 }\nolimits} )$$

\bigskip
\nn
5.3. {\it Path integrals.}
\smallskip
We revert to the notations in paragraph IV.1. The constant $\l$ is again $(h / m)^{1/2}$ and
we parametrize the paths $x$ with $x (t_b) = x_b$ by
$$x(t, z) = x_b + \l z (t) 
\eqno ({\mathop{\rm IV.80 }\nolimits} )$$
where $z$ runs over the space ${\bf Z}_b$. The horizontal lift of the previous path is given
by 
$$\xi (t, z) = 
\left({ x(t, z), e^{i \T (t, z) } }\right) \, . 
\eqno ({\mathop{\rm IV.81 }\nolimits} )$$
Taking into account the differential equation (IV.77), we obtain the following differential
system
$$\left\{\matrix{
d x^\a (t) = \l  d z^\a (t) \quad 
{\mathop{\rm for }\nolimits}
\  \a \ {\mathop{\rm in }\nolimits} \ \{ 1, \cdots, d \}
  \hfill\cr
\noalign{\medskip}
 d \T (t) = 
{e \l  \over \hbar} A_\a (x(t)) d z^\a (t) \, . 
  \hfill \cr}
\right.
\eqno ({\mathop{\rm IV.82 }\nolimits} )$$

\medskip
\nn
This system has the canonical form (II.21) where the vector fields 
$X_{(1)}, \cdots , X_{(d)}$, $Y$ are given by
$$\La_{X_{(\a)}} = \l \! 
\left({ {\part  \over \part x^\a} + {e \over \hbar} A_\a (x) {\part \over \part \T } }\right)
\quad, \quad Y = 0 \, .
\eqno ({\mathop{\rm IV.83 }\nolimits} )$$
Notice that for $\Psi (x, \T) = e^{- i \T} \psi (x) $, we have
$$\La_{X_{(\a)}}
\Psi (x, \T) = 
e^{- i \T } \cdot \l D_\a \psi (x) 
\eqno ({\mathop{\rm IV.84 }\nolimits} )$$
with the differential operator (see e.g. [20, p. 405]):
$$D_\a = 
 {\part  \over \part x^\a } 
-  {ie \over \hbar}  A_\a (x) \, . 
\eqno ({\mathop{\rm IV.85 }\nolimits} )$$

\bigskip
Our general partial differential equation
$$
 {\part \Psi \over \part t} 
=
 {s \over 4 \pi} \sum_{\a }^{ }
\La_{X_{(\a)}}^2 \Psi
\eqno ({\mathop{\rm IV.86 }\nolimits} )$$
translates now as
$$
{\part \psi \over \part t} 
=
 {s \over 4 \pi} {h \over m}
\sum_{\a }^{ } D_\a^2 \psi 
\eqno ({\mathop{\rm IV.87 }\nolimits} )$$
for the $\psi$-compoment of $\Psi$.
For $s = i$, this equation coincides with the Schr\"odinger equation (IV.71).

\bigskip
Feynman path integral solution to this equation reads as follows:
$$\psi \!  \left({ t_b, x_b  }\right) =
\int_{\Pa_b  }^{ }
\Da x \cdot e^{i S (x) / \hbar} \phi (x(t_a)) 
\eqno ({\mathop{\rm IV.88 }\nolimits} )$$
where the action $S(x)$ is given by equation (IV.69). Parametrizing the paths $x$ in $\Pa_b$ by the
paths $z$ in ${\bf Z}_b$ (see equation (IV. 80)), we can rewrite the previous path integral
as\footnote{$^{18}$}{
\nn We give the formula in the oscillatory case $s = i$. The reader is invited to work out the
formulas for the case $s = 1$. The definition of $Q_0 (z)$ is given in equation (IV.22).
}
$$
\psi  \! \left({ t_b, x_b  }\right) =
\int_{{\bf Z }_b }^{ }
\Da z \cdot
e^{ \pi i Q_0 (z)}
\exp \!
\left({ { ie\over \hbar  } 
\int_{{\bf T } }^{ }
A_\a (x(t, z)) d x^\a (t, z)
}\right) \phi (x(t_a, z)) \, . 
\eqno ({\mathop{\rm IV.89 }\nolimits} )$$

\medskip
\nn
Replacing $\psi (x) $ by $\Psi (x, \T) = e^{- i \T} \psi (x)$ and similarly $\phi (x) $ by
$\Phi (x, \T) =e^{- i \T} \phi (x)$, we can absorb the phase factor and obtain
$$\Psi \! \left({ t_b, x_b, \T_b }\right) =
\int_{{\bf Z }_b }^{ }
\Da z \cdot
e^{ \pi i Q_0 (z)}
\Phi \! \left({ \left({ x_b , \T_b }\right) \cdot \Si \!
\left({ t_a, z }\right) }\right) \, .
\eqno ({\mathop{\rm IV.90 }\nolimits} )$$
The transformation $\Si (t, z)$ of the bundle space $M \ts G$ takes 
$\left({x_b, \T_b  }\right)$ \break into 
$(x_b +\l z(t), \ \T_b - {e \l  \over \hbar }
\int_{t }^{ t_b} A_\a 
\left({ x_b + \l z (t) }\right) d z^\a (t))$ and corresponds to the integration of the
differential system (IV.82).

\bigskip
\nn
5.4. {\it Various generalizations.}

\smallskip
a) It is easy to incorporate an {\it electric  potential } $V$.
The complete action functional is now
$$S(x) = {m \over 2} \int_{{\bf T } }^{ }
{ \left\vert{ dx }\right\vert^2 \over dt} + e  \int_{{\bf T } }^{ }
 A_\a dx^\a - V dt \, . 
\eqno ({\mathop{\rm IV.91 }\nolimits} )$$
Feynman solution (IV.88) is still valid and can be made explicit as
$$\psi \! \left({ t_b, x_b }\right) =
\int_{{\bf Z }_b }^{ }
\Da z \cdot
e^{ \pi i Q_0 (z)} 
\exp \!
\left({ {ie \over \hbar } \int_{{\bf T } }^{ }  A_\a dx^\a - V dt 
 }\right) \phi\! \left({ x (t_a, z) }\right)
\eqno ({\mathop{\rm IV.92 }\nolimits} )$$
where the line integral 
$\int_{{\bf T } }^{ }  A_\a dx^\a - V dt $ is calculated along the path
$x(\cdot, z)$. The Schr\"odinger equation reads as follows:
$$i \hbar { \part \psi \over \part t}
= {1 \over 2m } \sum_{ \a }^{ } 
\left({ {\hbar \over i } {\part  \over \part x^\a}
- e A_\a (x)  }\right)_{}^{2} \psi + e V \psi \, . 
\eqno ({\mathop{\rm IV.93 }\nolimits} )$$

\bigskip
b) A non-abelian gauge group is for instance 
$G = U(N)$; a more general compact gauge group can always be realized as a closed subgroup of
some unitary group $U(N)$. We mention a few of the required changes. The Lie algebra 
${\hbox{\tenfm g}}$ is the set of matrices of the form $iE$ where $E$ is an $N$-by-$N$
hermitian matrix. The connection form $\om$ on the principal bundle $P$ takes its values in 
${\hbox{\tenfm g}}$ and the transformation rule is now
$$\om(p \cdot g) =
g^{-1} \om (p) g + g^{- 1}  dg
\eqno ({\mathop{\rm IV.94 }\nolimits} )$$
where $\om (p)$, $g$, $dg$ are 
 $N$-by-$N$ matrices and the products are given by matrix multiplication.

\bigskip
For the associated vector bundle $L$, we consider the natural action of $U(N)$ on the complex
vector space $\CC^N$. Hence a section $\psi$ of $L$ corresponds to a function 
$\Psi :P \ra\CC^N $ such that $\Psi (p \cdot g) = g^{- 1} \cdot \Psi (p)$ for $p$ in $P$ and
$g$ in $G = U(N)$.

\bigskip
In a given gauge the form  $\om $ is given by $\om (x, g) = g^{- 1} dg - 
{ ie \over \hbar} A(x) $, where $A(x)  = A_\a (x) dx^\a$ is a hermitian  $N$-by-$N$ matrix of
differential forms on $M$. A gauge transformation is given by the formulas
$$
\eqalignno{
\psi_1 (x) 
&= U(x)^{- 1} \psi (x) 
&({\mathop{\rm IV.95 }\nolimits} )\cr
A_\a^1 (x) 
&= U(x)^{- 1} A_\a (x) U(x) +  {i \hbar  \over e} U(x)^{- 1}
{\part \over \part x^\a } U(x),
&({\mathop{\rm IV.96 }\nolimits} )\cr}$$
where $U(x)$ is a unitary matrix depending on  the point $x$ of $M$. 

\bigskip

The horizontal lift of a path $x : T \ra M$ is of the form
$\xi (t) = (x(t), U(t))$ where the unitary matrix $U(t)$ satisfies the differential equation
$$\dot U (t) =  {i e \over \hbar} U(t) \cdot 
A_\a (x(t)) \dot x^\a (t) \, .
\eqno ({\mathop{\rm IV.97 }\nolimits} )$$
We can solve this equation in the standard way using time-ordered exponentials $\Ta \exp$.

\bigskip
The Schr\"odinger equation is still written in the form (IV.87), but $D_\a$ is now a
matrix of differential operators, namely ${\part  \over \part x^\a }
\cdot \un -  {ie \over \hbar}  A_\a (x)$ where $A_\a (x) $ is an $N$-by-$N$ hermitian matrix.
In the path integral (IV.89) replace the exponential factor by 
$$\Ta \exp \! 
\left({ {i e  \l  \over \hbar  }
\int_{ t_a}^{t_b }
A_\a \! \left({ x_b + \l  z(t)  }\right) dz^\a (t)
 }\right) \, . 
\eqno ({\mathop{\rm IV.98 }\nolimits} )$$

\bigskip
c) We could consider gauge groups over a curved manifold and combine the results of paragraph
IV.3 with those of the present paragraph. 

\bigskip
\nn
{\bf 6. A symplectic manifold.}\footnote{$^{19}$}{\nn Contributed by John LaChapelle.}

\smallskip
Let $N$ be a symplectic manifold $\Ma$ of dimension $d = 2n$. The manifold $\Ma$ represents
the classical phase space of a physical system -- usually the cotangent bundle $T^* Q$ of a
configuration space $Q$ but not necessarily. Paths in $\Ma$ have $n$ initial and $n$ final boundary conditions. These boundary conditions
are consistent with the requirement of quantum uncertainty, and they imply a choice of
polarization\footnote{$^{20}$}{\nn Roughly speaking, a polarization is a foliation of $\Ma$ whose
leaves are Lagrangian submanifolds of dimension $n$. }.
 Hence, it is possible to cover $\Ma$ with a family of open subsets $
\left\lbc{ U_i }\right\rbc$ such that each $U_i$ is diffeomorphic to a product of two
transverse Lagrangian submanifolds $L_i \ts L'_i$. For simplicity consider a symplectic
manifold which admits global transverse Lagrangian submanifolds $L$ and $L'$,
and identify $\Ma$ with $L \ts L'$. Set $x(t) $ in the form 
$\left({ x_L (t), x_{L'} (t)}\right)$ with $x_L (t)$ in $L$ and $x_{L'}(t)$ in $L'$, in such a way
 that
$x_L (t_b) = x_b$ and  that 
$x_{L'}(t_a) = x_a$.

\bigskip
We require the group generated by the transformations $\Si ({\bf T}, z)$ to be a subgroup of the group of
symplectomorphisms which leave the polarization invariant. This implies that the set
$ \left\lbc{ X_{(\a)} }\right\rbc$ is of the form
$\left\lbc{ X_{(a)},  X_{(a')} }\right\rbc$ such that
$ X_{(a)} \! \left({ x_{L'} (t) }\right) = 0$,
$X_{(a')} \! \left({ x_{L} (t) }\right) = 0$, and 
$\left\lbk{ X_{(a)}, X_{(a')} }\right\rbk = 0$. Here 
$a \in \left\lbc{ 1, \cdots, k }\right\rbc$,
$a' \in \left\lbc{k +  1, \cdots, d }\right\rbc$, and $k$  is a fixed integer between $1$ and
$d$,  possibly, but not necessarily, equal to $d/2$.
Consequently, a path satisfies the differential equations
$$
\left\{\matrix{
dx_L (t)
=  X_{(a)} \! \left({ x_{L} (t) }\right) dz^a + 
Y \! \left({ x_{L} (t) }\right) dt
  \hfill\cr
\noalign{\medskip}
 dx_{L'} (t)
=  X_{(a')} \! \left({ x_{L'} (t) }\right) dz^{a'} + 
Y \! \left({ x_{L'} (t) }\right) dt  \, ,
  \hfill \cr}
\right.
\eqno ({\mathop{\rm IV.99 }\nolimits} )$$
and the general formula (II.1) becomes
$$
\left.\matrix{ 
\left({ U_{t_b, t_a} \phi }\right)
\! \left({ x_b, x_a }\right) := 
\int_{{\bf Z}_{L} }^{ } \int_{{\bf Z}_{L'} }^{ }
\Da z_L \, \Da z_{L'} \cdot 
e^{
 - \pi Q_0 \! \left({ z_L, z_{L'} }\right) /s
}
\cr
\noalign{\medskip}
\ts \phi
\! \left({ x_b \cdot \Si
\! \left({ t_a, z_L }\right), x_a \cdot \Si
\! \left({ t_b, z_{L'} }\right)
 }\right) \, , 
  \cr}\right.
\eqno ({\mathop{\rm IV.100 }\nolimits} )
$$
where now $h^{\a \b}
= \pmatrix{
h^{ab} & 0 \cr
0 & h^{a'b'} \cr
}$. Here ${\bf Z}_L$ is the space of paths 
$z_L : {\bf T} \ra \RR^k$ such that $z_L (t_b) = 0$, and 
${\bf Z}_{L'}$ is the space of paths $z_{L'} : {\bf T} \ra \RR^{d-k}$
such that $z_{L'}(t_a) = 0$. Each functional integral separately satisfies a partial differential equation:
$$
\left\{\matrix{
 {\ds\part \Psi_L \over \ds \part t_a} 
=
 {\ds s \over \ds 4 \pi } 
h^{ab} \La_{X_{(a)}} \La_{X_{(b)}} \Psi_L + \La_{Y} \Psi_L
  \hfill\cr
\noalign{\medskip}
  {\ds\part \Psi_{L'} \over \ds \part t_b} 
=
 {\ds s \over \ds 4 \pi } 
h^{a'b'} \La_{X_{(a')}} \La_{X_{(b')}} \Psi_{L'} + \La_{Y} \Psi_{L'} \, , 
  \hfill \cr}
\right. \eqno ({\mathop{\rm IV.101 }\nolimits} )
$$
where $\Psi_L := 
\left.{ \Psi }\right\vert_{L }$
and 
$\Psi_{L'} := 
\left.{ \Psi }\right\vert_{L' }$.

\vfill\eject
\nn
6.1. {\it The case of a cotangent bundle.}
\smallskip
Choosing an initial function $\phi $ is choosing a transition amplitude.
We consider the case where $\Ma$ is the cotangent bundle of a flat configuration space $Q$,
hence $\Ma = Q \ts P$ where $Q$ and $P$ are finite-dimensional vector spaces in duality.
In order to define {\it position-to-position transition amplitudes}, we choose
 the initial function of the form
$\phi (q, p) = \d \! \left({q  - q_a }\right)
$.
Equation (IV.100) yields
$$\Ka \! \left({  q_b, t_b ; q_a, t_a 
 }\right) =
\int_{{\bf Z}_Q }^{ } \int_{{\bf Z}_P }^{ }
\Da z_Q \, \Da z_P \cdot 
\exp \! 
\left({ - {\pi  \over s} Q_0  \! \left({ z_Q, z_P }\right) }\right)
h \! \left({ z_Q, z_P  }\right)
\eqno ({\mathop{\rm IV.102 }\nolimits} )$$
where the integrand is given by 
$$
h \! \left({ z_Q, z_P  }\right)
=
\d \! 
\left({ q_b \cdot \Si \! \left({ t_a, z_Q }\right) - q_a }\right)
\, . \eqno ({\mathop{\rm IV.103 }\nolimits} )$$

\medskip
In order for the transition amplitudes to be consistent with the initial wave function, we require
$\build { \lim }_{ t_a \ra t_b }^{ }
\Ka \! \left({ q_b, t_b ; q_a, t_a }\right)$ to be equal to $\d \! 
\left({ q_b - q_a }\right)$. But the integrand 
$h \! \left({ z_Q, z_P }\right)$ tends to 
$\d \! \left({ q_b \cdot \Si \! \left({ t_b, z_Q }\right) - q_a}\right)
= \d \! \left({ q_b - q_a }\right)$ when $t_a$ tends to $t_b$, a limit independent of 
$z_Q$, $z_P$. Hence, provided we can interchange limit and integration, we get
$$\build { \lim }_{ t_a \ra t_b }^{ }
\Ka \! \left({ q_b, t_b ; q_a, t_a }\right) =
\int_{{\bf Z}_Q }^{ } \int_{{\bf Z}_P }^{ }
\Da z_Q \, \Da z_P \cdot
\exp \! 
\left({ - {\pi  \over s} Q_0  \! \left({ z_Q, z_P }\right) }\right)
\d \! \left({ q_b - q_a }\right) \, .
$$
Hence this is equal to $\d \! \left({ q_b - q_a }\right)$ by the normalization of our integrator.

\medskip
We handle the other cases in a similar way.

\smallskip
a) {\it Momentum-to-position amplitude}
$\Ka \! \left({ q_b, t_b ; p_a, t_a }\right)$: use the initial function $
 h^{-n/2} e^{i q_b \cdot p / \hbar }$
and the integrand
$$ h^{-n/2}
\exp \! 
\left({ { i \over \hbar} q_b \cdot 
\left({ p_a \cdot \Si \! \left({ t_b, z_P }\right) }\right)
 }\right)
\, . \eqno ({\mathop{\rm IV.104 }\nolimits}_{{\mathop{\rm a }\nolimits}} )
$$

\smallskip
b) {\it Position-to-momentum amplitude}
$\Ka \! \left({ p_a, t_a ; q_b, t_b }\right)$: use the initial function $
h^{-n/2} e^{- i q \cdot p_a / \hbar }$
and the integrand
$$ h^{-n/2}
\exp \! 
\left({ - { i \over \hbar} 
\left({ q_b \cdot \Si \! \left({ t_a, z_Q }\right) }\right) \cdot p_a
 }\right)
\, . \eqno ({\mathop{\rm IV.104 }\nolimits}_{{\mathop{\rm b }\nolimits}} )$$

\smallskip
c) {\it Momentum-to-momentum amplitude}
$\Ka \! \left({ p_b, t_b ; p_a, t_a }\right)$: use the initial function $\d \! \left({ p - p_b }\right)$
and the integrand
$$\d \! 
\left({ p_a \cdot \Si \! \left({ t_b, z_P }\right) - p_b}\right) 
\, . \eqno ({\mathop{\rm IV.104 }\nolimits}_{{\mathop{\rm c }\nolimits}} )
$$

\medskip
\nn
It follows easily from these definitions that
$$
\eqalignno{
\Ka \! \left({ p_b, t_b ; p_a, t_a }\right) 
&=
h^{- n/2 }
\int_{Q }^{ } dq_b \  e^{- i q_b \cdot p_b / \hbar }  \Ka \! \left({ q_b, t_b ; p_a, t_a }\right)
& ({\mathop{\rm IV.105 }\nolimits} )  \cr
\Ka \! \left({ q_b, t_b ; q_a, t_a }\right) 
&=
h^{- n/2 }
\int_{P }^{ } dp_b \  e^{ i q_b \cdot p_b / \hbar } \Ka \! \left({ p_b, t_b ; q_a, t_a }\right)  \, . &
({\mathop{\rm IV.106 }\nolimits} ) \cr} $$

\bigskip
\nn
6.2. {\it Coherent states.}

\smallskip
More general transition amplitudes, which cannot be interpreted as posi-tion-to-momentum
transitions, are possible by choosing more complicated initial wave functions, different
polarizations, and/or by having non-trivial phase spaces. For instance, coherent state
transitions can be calculated -- given a suitable characterization of coherent states.

\medskip

Choose a K\"ahler polarization on a non-trivial symplectic manifold, so that $\Ma$ is a
K\"ahler manifold $M$ of complex dimension $n$. It is convenient (though not necessary) to take the phase space
to be the product manifold $ M \ts \oo M$ with coordinates $(\z, \oo z)$ and complex dimension $2n$.
Following the work of Berezin [21] and Bar-Moshe and Marinov [22], we use generalized coherent
state wave functions
$\phi_{\z '} (\z) = \exp  K(\z, \oo \z') $
where 
$K(\z, \oo \z')$ is the K\"ahler potential.
Consider the function $\phi : M \ts \oo M \ra \CC  $ given by 
$\phi (\z, \oo \z') = 
e^{K \! \left({ \z, \oo \z_a }\right) }
$
and take  $z_M : {\bf T} \ra \CC^n$ such that $z_M (t_b) = 0$ and 
$z_{\oo M} : {\bf T} \ra \CC^n$ such that $z_{\oo M}(t_a) = 0$.
Then
$$ \Ka \!  
\left({ \z_b, t_b ; \oo \z_a, t_a  }\right): = \int_{ {\bf Z}_M }^{ }  \int_{ {\bf Z}_{\oo M} }^{ }
\Da z_M \, \Da z_{\oo M} \cdot  
e^{
 - \pi   Q_0  \left({ z_M, z_{\oo M } }\right) /s
}
e^{K
\! \left({  \z_b \cdot \Si  \left({ t_a, z_{ M} }\right), \oo \z_a  }\right)
}  
$$
 are generalized coherent state transition amplitudes.
Note that
$\Ka 
 \!  
\left({\oo \z_b, t_b ;  \z_a, t_a  }\right) $ \break $=
\oo { \Ka  \!  
\left({ \z_b, t_b ; \oo \z_a, t_a  }\right) } $,
and 
$\build {\lim }_{ t_a \ra t_b  }^{ }
\Ka 
 \!  
\left({ \z_b, t_b ; \oo \z_a, t_a  }\right)
=
\exp  K \! \left({ \z_b, \oo \z_a }\right) 
$ provided the limit can be taken inside the integral.

\bigskip
\nn
{\bf 7. A simple model of the Bohm-Aharonov effect.}

\smallskip
We illustrate  on an elementary example various techniques expounded in this section.
Our basic manifold $N$ is the euclidean plane $\RR^2$ with the origin removed. As in paragraph IV.2, we denote
the cartesian coordinates by $x^1, x^2$. The universal covering of $N$ is the set $\ww N$ of pairs of real
numbers $r, \t$ with $r > 0$, and the covering map $\Pi : \ww N \ra N $ is described by equation (IV.31), namely
$x^1 = r \cos \t$, $x^2 = r \sin \t$.

\medskip
We introduce a {\it magnetic potential} $A$ with components 
$$A_1 = - {Fx^2 \over 2 \pi \left\vert{ x }\right\vert^2} 
\quad , \quad  A_2 = 
{Fx^1 \over 2 \pi \left\vert{ x }\right\vert^2} \, .
\eqno ({\mathop{\rm IV.107 }\nolimits} )$$
The magnetic field has one component $B$ perpendicular to our plane, given by
$$B = {\part A_1 \over \part x^2} - {\part A_2 \over \part x^1} \, ,
\eqno ({\mathop{\rm IV.108 }\nolimits} )$$
and an obvious calculation gives $B = 0 $ {\it outside of the origin.}
But the circulation $\oint  A_1 dx^1 + A_2 dx^2 
$ on any loop around the origin is equal to $F$, hence
$$ B = F \d (x) \, .
\eqno ({\mathop{\rm IV.109 }\nolimits} )$$
Physically, we have a wire perpendicular to the plane through its
origin, carrying a magnetic flux equal to $F$.

\medskip
Denote by $m$ the mass and by $e$ the electric charge of a particle. For a path $x : {\bf T} \ra N$, with
components
$x^1 (t)$, $x^2(t)$ at time $t$, the action is given by 
$$S(x) = S_0 (x) + S_M (x) \, ,
\eqno ({\mathop{\rm IV.110 }\nolimits} )$$
with the {\it kinetic action}
$$S_0 (x) = {m \over 2}
\int_{ {\bf T}}^{ }
 { \left\vert{  dx }\right\vert_{}^{2 } \over dt}
\eqno ({\mathop{\rm IV.111 }\nolimits} )$$
and {\it the magnetic action}
$$S_M (x) = {eF \over 2 \pi } \int_{ {\bf T}}^{ }
{ x^1 dx^2 - x^2 d x^1 \over (x^1)^2 + (x^2)^2}
\, . 
\eqno ({\mathop{\rm IV.112 }\nolimits} )$$
We fix a point $x_b$ in $N$, and denote by ${\bf X}_b$ the set of paths $x : {\bf T} \ra N$ for which the
kinetic action $S_0 (x)$ is finite (that is the $L^{2,1}$ paths) with the boundary condition 
$x(t_b) = x_b$. On ${\bf X}_b$, we denote by $\Da x$ the translation-invariant integrator normalized by 
$$\int_{{\bf X}_b}
\Da x \cdot \exp \! 
\left({ { i \over \hbar}  S_0 (x) }\right) = 1 \, . 
\eqno ({\mathop{\rm IV.113 }\nolimits} )$$

\medskip
The free position-to-position transition amplitudes are defined by
$$
\left\langle{ t_b, x_b \bigm| t_a , x_a }\right\rangle_{0}^{}
= \int_{ {\bf X}_b}^{ }
\Da x \cdot \exp \! 
\left({ { i \over \hbar}  S_0 (x) }\right)
\d \!
\left({ x(t_a) - x_a }\right) \, .
\eqno ({\mathop{\rm IV.114 }\nolimits} )$$
According to the calculation made at the end of paragraph IV.2, we get\footnote{$^{21}$}{We express our formulas
directly in terms of the path $x$, and have no need for the scaling $x = \l z$ introduced earlier.} 
$$\left\langle{ t_b, x_b
\bigm| t_a , x_a }\right\rangle_{0}^{} = {- mi \over h(t_b - t_a)} 
\exp \!
\left({ { i \over \hbar}  S_0 (x_{{\mathop{\rm cl }\nolimits}}) }\right)
\, , \eqno ({\mathop{\rm IV.115 }\nolimits} )$$
where $x_{{\mathop{\rm cl }\nolimits}}$ is the classical path of a free particle, namely:
$$x_{{\mathop{\rm cl }\nolimits}} (t) =
{ x_a (t_b - t) + x_b (t - t_a) \over t_b - t_a} 
\, . 
\eqno ({\mathop{\rm IV.116 }\nolimits} )$$
Hence, we obtain explicitely
$$S_0 (x_{{\mathop{\rm cl }\nolimits}}) = 
 {m \over 2} 
{ \left\vert{  x_b - x_a}\right\vert_{}^{2} \over t_b - t_a} \, , 
\eqno ({\mathop{\rm IV.117 }\nolimits} )$$
that is, the WKB approximation (IV.30) is exact in our case.

\medskip
We consider now the transition amplitudes in the given magnetic potential:
$$\left\langle{ t_b, x_b
\bigm| t_a , x_a }\right\rangle_{F}^{}
=
\int_{ {\bf X}_b}^{ }
\Da x \cdot \exp \! 
\left({ { i \over \hbar}  S (x) }\right)
\d \!
\left({ x(t_a) - x_a }\right) \, .
\eqno ({\mathop{\rm IV.118 }\nolimits} )$$
According to (IV.87), we consider the Schr\"odinger equation
$$
{ \part \psi_F \over \part t} =
{ i \hbar \over 2m}
 \left({ D_1^2 + D_2^2 }\right)
\psi_F 
\eqno ({\mathop{\rm IV.119 }\nolimits} )$$
with the differential operators 
$D_\a = {\part  \over \part x^\a} - { ie \over \hbar} A_\a$.
Explicitly, we obtain
$${\part \psi_F \over \part t}
=
{ i \hbar \over 2m} L \psi_F \, ,
\eqno ({\mathop{\rm IV.120 }\nolimits} )$$
with
$$L =
\left({ {\part \over \part x^1} }\right)_{}^{2} +
\left({ {\part \over \part x^2} }\right)_{}^{2} -
 {2 ci \over \left\vert{ x }\right\vert^2 }
\left({ x^1{ \part \over \part x^2} - x^2{ \part \over \part x^1} }\right)
- {c^2 \over \left\vert{ x }\right\vert^2 } \, . 
\eqno ({\mathop{\rm IV.121 }\nolimits} )$$
The constant $c$ is equal to $eF / h$; it is dimensionless. The solution to the equation (IV.119) is given by
$$\psi_F \! \left({ t_b, x_b  }\right) =
\int_{ N }^{ } dx_a 
\left\langle{ t_b, x_b \bigm| t_a, x_a }\right\rangle_F \psi_F \! 
\left({  t_a, x_a}\right) \, ,
\eqno ({\mathop{\rm IV.122 }\nolimits} )$$
a formula equivalent to our more familiar one
$$
\psi_F \!\left({ t_b, x_b  }\right) = 
\int_{ {\bf X}_b}^{ }
\Da x \cdot \exp \! 
\left({ { i \over \hbar}  S (x) }\right)
\phi \!
\left({ x(t_a)  }\right)
\eqno ({\mathop{\rm IV.123 }\nolimits} )$$
if we take into account the initial wave function $\phi(x) = 
\psi_F \!\left({ t_a, x }\right)$ at time $t_a$.

\medskip
We lift now everything to the universal covering $\ww N$. By means of the formulas (IV.31), the wave function is
now a function $\ww \psi_F (t, r, \t)$ with the restriction
$$\ww \psi_F (t, r, \t + 2 \pi) = \ww \psi_F (t, r, \t)
\eqno ({\mathop{\rm IV.124 }\nolimits} )$$
that is, a function on ${\bf T} \ts \ww N$ invariant under the group $\ZZ$ acting on $\ww N$ by 
$(r , \t)  \cdot n = (r , \t + 2 \pi n)$. The Schr\"odinger equation (IV.120) keeps its form with
 $L$ changed into the new
operator
$$
\ww L = 
 {\part^2 \over \part r^2} +  {1 \over r^2}   {\part^2\over \part \t^2} +
{1 \over r}   {\part \over \part r}
-  {2ci \over r^2}  {\part  \over \part \t} -  {c^2 \over r^2} \, .
\eqno ({\mathop{\rm IV.125 }\nolimits} )$$
This operator can be written as $\ww L = e^{c i \t} \, \D e^{- ci \t}$ where
$\D $ is the Laplacian in polar coordinates, namely:
$$
\D = {\part^2 \over \part r^2} +  {1 \over r^2}   {\part^2\over \part \t^2} +
{1 \over r}   {\part \over \part r} \, .
\eqno ({\mathop{\rm IV.126 }\nolimits} )$$
Hence the new wave function
$$\ww \psi_0 (t, r, \t) = e^{- ci \t }
\ww \psi_F (t, r, \t)
\eqno ({\mathop{\rm IV.127 }\nolimits} )$$
satisfies the equation
$$
 {\part \ww \psi_0 \over \part t} =
 {i \hbar \over 2m}  \D \ww \psi_0 \, ,
\eqno ({\mathop{\rm IV.128 }\nolimits} )$$
that is the Schr\"odinger equation for a free particle, written in polar coordinates. As we shall see in a
moment, equation (IV.127) expresses a {\it gauge transformation} which apparently removes the magnetic
potential. But the periodicity condition (IV.124) transforms into
$$\ww \psi_0 (t, r , \t + 2 \pi) = e^{- 2 \pi i c}
\ww \psi_0 (t, r, \t) \, ,
\eqno ({\mathop{\rm IV.129 }\nolimits} )$$
hence the wave function $\ww \psi_0$ on $\ww N$ is not the lifting of a function $\psi_0 $ on $N$.
Rather, the linear representation $n \mps e^{ 2 \pi i n c }$ of the group $\ZZ$ into $U(1)$ defines a line
bundle on $N = \RR^2 \bl \{ 0 \} $ and $\ww \psi_0$ {\it corresponds to a section $\psi$ of this line
bundle.}
Otherwise stated, despite the fact that the magnetic field $B = \part A_1 / \part x^2 - \part A_2 / \part x^1$
is identically zero on $N$, there is no single-valued function $R$ on $N$ such that $A_\a = - \part R / \part
x^\a$: {\it the space $N$ is multiply connected.}

\bigskip
To explain gauge transformations, let us introduce the principal bundles 
$P = N \ts U(1)$ and $\ww P = \ww N \ts U(1)$. According to the general formulas the connection form $\om_F$ on
$P$ is given by
$$\om_F =
i 
\left({d \T -   {e \over \hbar } A}\right)
= i \left({d \T - c { x^1 dx^2 - x^2 d x^1\over 
\left({ x^1 }\right)^2 + \left({ x^2 }\right)^2
} }\right)
\, ,
\eqno ({\mathop{\rm IV.130 }\nolimits} )$$
where $\T$ is the angular coordinate on $U(1)$ (taken modulo $2\pi$). When there is no magnetic field, it reduces
to 
$\om_0 = id\T$. We cannot transform $\om_F$ into $\om_0$
 by a gauge transformation, but if we lift these differential forms to $\ww N$, we
obtain
$$\ww \om_F = i(d \T - cd\t) \quad , \quad \ww \om_0 = id \T \, .
\eqno ({\mathop{\rm IV.131 }\nolimits} )$$
The transformation $\ww U$ of $\ww P$ into $\ww P$ taking $(r, \t, \T)$ into $(r, \t, \T + c\t)$ is an
automorphism of $U(1)$-bundle, and 
$\ww U^* \ww \om_F = \ww \om_0$, but $\ww U$ is not the lifting of an automorphism $U$ of $P$, 
{\it except when $c$ is an integer.}

\bigskip
To the wave function $\psi_F (x) $ on $N$ we associate the function
$$\Psi_F (x, \T) = e^{- i \T} \psi_F (x)
\eqno ({\mathop{\rm IV.132 }\nolimits} )$$
on $P$. The Schr\"odinger equation (IV.119) translates into our standard form
$$
 {\part \Psi_F \over \part t} 
=  {i \over 4 \pi} h^{\a \b}
\La_{X_{(\a)}^F} \La_{X_{(\b)}^F} \Psi_F
\eqno ({\mathop{\rm IV.133 }\nolimits} )$$
provided we take $h_{\a \b}$ equal to $m \d_{\a \b}/h$, with the vector fields
$$
\left\{\matrix{
\La_{X_{(1)}^F} = {\ds \part  \over \ds\part x^1} - {\ds cx^2 \over \ds\left\vert{ x }\right\vert^2 }
{ \ds\part \over \ds\part \T} 
  \hfill\cr
\noalign{\medskip}
\La_{X_{(2)}^F} = {\ds\part  \over \ds\part x^2} + {\ds cx^1 \over \ds\left\vert{ x }\right\vert^2 }
{\ds \part \over \ds\part \T}   
\, .
  \hfill \cr}
\right.
\eqno ({\mathop{\rm IV.134 }\nolimits} )$$
We lift now everything to $\ww P$, so $\psi_F$ lifts to $\ww \psi_F, \ \Psi_F $ to $\ww \Psi_F$ and the lifted
vector fields are given by\footnote{$^{22}$}
{For $F = 0$, they reduce to the vector fields given in equations (IV.36) and (IV.37).}
$$
\left\{\matrix{
\La_{\ww X_{(1)}^F} = \cos \t \cdot {\ds\part  \over \ds\part r} - {\ds\sin \t \over\ds r } \cdot 
\left({ {\ds \part \over \ds\part \t} + c {\ds \part \over \ds\part \T} }\right) 
  \hfill\cr
\noalign{\medskip}
\La_{\ww X_{(2)}^F} = \sin \t \cdot {\ds\part  \over \ds\part r} + {\ds\cos \t \over \ds r } \cdot 
\left({ {\ds \part \over \ds\part \t} + c {\ds \part \over \ds\part \T} }\right)   
\, .
  \hfill \cr}
\right.
\eqno ({\mathop{\rm IV.135 }\nolimits} )$$
The Schr\"odinger equation (IV.133) remains the same with $\Psi_F$ replaced by $\ww \Psi_F$ and
$X_{(\a)}^F $ by $\ww X_{(\a)}^F$.

\medskip
We describe the effect of the gauge transformation $\ww U$ given by 
$\ww U (r, \t, \T) = (r, \t, \T + c \t)$.
We noticed that 
$\ww U^* \ww \om_F $ is equal to $\ww \om_0$, and similarly, we get 
$\ww \Psi_F \circ \ww U = \ww \Psi_0$ where $\ww \Psi_0$ is defined by analogy to $\ww \Psi_F$, namely
$$\ww \Psi_0 (r, \t, \T) =
e^{- i \T }
\ww \psi_0 (r, \t) \,.
\eqno ({\mathop{\rm IV.136 }\nolimits} )$$
Dually, $\ww U$ transforms the vector field $\ww X_{(\a)}^0$ into 
$\ww X_{(\a)}^F$ for $\a = 1$ or $2$.
Hence $\ww \Psi_F$ is a solution of the Schr\"odinger equation (IV.133) if, and only if,
$\ww \Psi_0$ is a solution of the corresponding equation for $F = 0$.

\medskip
We conclude this paragraph by a discussion of path integrals.
We parame-trize paths in ${\bf X}_b$ by paths in ${\bf Z}_b$, using the correspondence
$$x(t) = x_b + z(t) \, . 
\eqno ({\mathop{\rm IV.137 }\nolimits} )$$
Fix $r_b$ and $\t_b$ in such a way that
$$x_b = 
\left({ r_b \cos \t_b, r_b \sin \t_b }\right)
\, . 
\eqno ({\mathop{\rm IV.138 }\nolimits} )
$$
Then $x$ can be lifted in a unique way to a path $\ww x : {\bf T} \ra \ww N$ of the form 
$\ww x (t) = (r(t), \t (t))$ such that 
$\ww x (t_b) = 
\left({ r_b, \t_b }\right)$,
that is 
$r(t_b) = r_b $ and $\t (t_b) = \t_b$.
The initial point $x_a = x(t_a)$ is lifted to $\ww x (t_a) =
\left({ r_a, \t_a }\right)$ where 
$r(t_a) = r_a$, $\t(t_a) = \t_a$.
Moreover, the magnetic action of the original path is given by 
$$S_M (x) / \hbar = c \cdot 
\left({ \t_b - \t_a }\right)
\, . 
\eqno ({\mathop{\rm IV.139 }\nolimits} )$$

\smallskip
\nn
Consider the space ${\bf X}_{a, b}$ of paths $x : {\bf T} \ra N$ such that 
$x(t_a) = x_a$, $x(t_b) = x_b$. {\it Then two such paths are in the same homotopy class if, and only if,
the lifted paths in $\ww N$ correspond to the same determination of the angular coordinate
$\t_a$ of $x_a$.
An equivalent condition is that they have the same magnetic action.}

\medskip
Using the connection of the principal bundle $\ww P$ over $\ww N$, we can define the horizontal lifting of
$\ww x$; it is the unique path of the form 
$y(t) = (r(t), \t(t), \T(t))$ on which $\ww \om_F$ induces a zero form, that is 
$\T(t) - c \t(t)$ is a constant in time. For the phase factors, we get
$$e^{i \T (t_b)} =
e^{i \T (t_a)} e^{i S_M (x) / \hbar}
\, . 
\eqno ({\mathop{\rm IV.140 }\nolimits} )$$
Moreover, the horizontal lifting is a solution of our standard differential equation $dy =
X_{(\a)} (y) \cdot dz^\a$.

\medskip
>From all this, it follows that the transition amplitudes as defined in the Feynman way (IV.118) agree with our
standard expression (see for instance (IV.38)). The result obtained in paragraph IV.2 (see formulas (IV.41) and
(IV.47)) can be restated as
$$
\left\langle{ t_b, x_b \bigm| t_a, x_a  }\right\rangle_0
= { 1 \over r_a}
\sum_{ n\in \ZZ }^{ } 
\left\langle{ t_b, r_b,\t_b \bigm| t_a, r_a, \t_a + 2 n \pi}\right\rangle
\, . 
\eqno ({\mathop{\rm IV.141 }\nolimits} )$$
A similar reasoning yields a more general formula:
$$
\left\langle{ t_b, x_b\bigm| t_a, x_a  }\right\rangle_F
= { 1 \over r_a}
\sum_{ n\in \ZZ }^{ }  e^{ic (\t_b - \t_a - 2n \pi) }
\left\langle{ t_b, r_b,\t_b \bigm| t_a, r_a, \t_a + 2 n \pi}\right\rangle
\, . 
\eqno ({\mathop{\rm IV.142 }\nolimits} )$$
We can invert this formula\footnote{$^{23}$}{Notice that the transition probabilities are functions of $F$ which
admit a period $h/e$, a well-known physical effect.}
and obtain
$$\left\langle{ t_b, r_b, \t_b \bigm| t_a, r_a, \t_a  }\right\rangle
= r_a \int_{ 0}^{ 1} dc\,  e^{- ic \left({ \t_b- \t_a }\right) }
\left\langle{ t_b, x_b\bigm| t_a, x_a }\right\rangle_{hc/e} \, .
\eqno ({\mathop{\rm IV.143 }\nolimits} )$$
Hence, {\it from the knowledge of the transition amplitude 
$\left\langle{ t_b, x_b\bigm| t_a, x_a }\right\rangle_{F}$
as a function of the magnetic flux $F$, we can infer the value of the transition amplitudes in polar
coordinates.  }

\vfill \eject
\centerline {{\bf Acknowledgements.} }

\bigskip
Our collaboration has been made possible by the financial support of the NATO collaborative
Research Grant \# 910101, the University (of Texas) Research Institute grant R-270, the Jane
and Roland Blumberg Centennial Professorship in Physics, and the Institut des Hautes Etudes
Scientifiques (I.H.E.S.) at Bures-sur-Yvette.

\medskip
The Center for Relativity of the University of Texas, the Ecole Normale Sup\'erieure, Paris
and the Institut des Hautes Etudes Scientifiques, Bures-sur-Yvette have given us ideal
working conditions.

\medskip
The manuscript which went through  several revisions has been diligently prepared by
Marie-Claude Vergne at I.H.E.S. and Debbie Hajji at the Center for Relativity.

\vfill\eject
\vglue 1cm
\centerline{J\bf Appendix A } 
\medskip
\centerline {\titre   Functional integration}
\bigskip
\bigskip

In this appendix, we develop the basic properties of our integrators. In the oscillating case $(s
= i)$, our theory is, up to some inessential changes in notation, the same as the one expounded by
Albeverio and H\o egh-Krohn in [9]. The introduction of a parameter $s$ equal to $1$ or $i$
enables us to treat in a unified way the oscillating integrators of Fresnel type
$e^{\pi i Q(x) } \Da x$ and the Gaussian integrators of type 
$e^{- \pi Q(x)} \Da x$ (with $Q(x) > 0$ for $x \ne 0$ in the latter case).
We content ourselves by giving here the {\it basic formulas} and the {\it computational tools}.
In the Gaussian case, we would have to justify these formal manipulations, since, as it is
well-known, the $L^{2, 1}$ {\it functions are a set of measure $0$ for the Wiener measure}, and 
{\it we
are considering functions on this set of measure $0$.} Our claims can be fully vindicated, but we
defer the complete justification to another publication.

\bigskip
\nn
{\bf 1. Gaussian integrators.}
\smallskip
\nn
1.1. {\it The setup.}
\smallskip

We denote by  ${\bf X}$ a real, separable, Banach space, by ${\bf X'}$ its dual and by 
$\left\langle{ x', x }\right\rangle$ (or sometimes
$\left\langle{ x', x }\right\rangle_{{\bf X} }$) the duality between 
${\bf X}$ and ${\bf X'}$. We suppose given a continuous linear map
$D : {\bf X} \ra {\bf X'}$ with the following properties:

\smallskip
-- (symmetry) $\left\langle{ Dx, y }\right\rangle =
\left\langle{ Dy, x }\right\rangle$ for $x, y $ in ${\bf X}$;

\smallskip
-- (invertibility) there exists a continuous linear map
$G : {\bf X'} \ra {\bf X}$ inverse of $D$, that is $DG = \un$ and $GD = \un$. 

\smallskip
\nn
Out of these data one constructs two quadratic forms, $Q$ on ${\bf X}$ and $W$ on ${\bf X'}$, by
the rules
$$Q(x)  = 
\left\langle{ Dx, x }\right\rangle \quad, \quad
W(x') = \left\langle{ x', Gx' }\right\rangle. 
\eqno ({\mathop{\rm A.1 }\nolimits})$$

\smallskip
\nn
They are related to each other as follows:
$$Q(x)  = 
W (Dx) \quad, \quad
W(x') = Q(Gx')
\eqno ({\mathop{\rm A.2 }\nolimits})$$
for $x$ in ${\bf X}$
and $x'$ in ${\bf X'}$.

\bigskip
We denote by $s$ a parameter equal to $1$ or $i$. The function 
$e^{- \pi s W}
$ on ${\bf X'}$ is continuous. It is bounded in the following cases:

\smallskip
-- $s = i$;

\smallskip
-- $s = 1$ and $W(x') > 0$ for 
$x' \ne 0$ in ${\bf X'}$. Equivalently, by (A.2), $s = 1$ and $Q(x) > 0$
 for $x \ne 0$ in ${\bf X}$.

\bigskip
\nn
1.2. {\it The oscillatory case $(s = i)$.}
\smallskip
The integrator $\Da x$ is characterized by the following integration formula:
$$
\int_{ {\bf X} }^{ }
\Da x \cdot \exp \!
\left({ \pi i Q(x) - 2 \pi i 
\left\langle{ x', x }\right\rangle
 }\right) =
\exp (- \pi i W (x'))
\eqno ({\mathop{\rm A.3 }\nolimits})$$
for every $x'$ in ${\bf X'}$. This relation should be interpreted as follows.

\bigskip
Since the metric
space ${\bf X'}$ is complete and separable, we know the notion of a complex bounded Borel
measure\footnote{$^{25}$}{That is: a $\s$-additive functional from the $\s$-algebra of Borel
subsets of ${\bf X'}$ into the complex numbers (see [23] for instance). We use a simplified
terminology by refering to $\mu$ as a ``measure'' on ${\bf X'}$.} $\mu$ on ${\bf X'}$.
The {\it Fourier-Stieltjes transform} $\Fa \mu$ of $\mu$ is given the customary definition:
$$\left({ \Fa \mu }\right)
(x) = \int_{{\bf X'}}^{ } d \mu (x') \,  e^{- 2 \pi i 
\left\langle{ x', x }\right\rangle } .
\eqno ({\mathop{\rm A.4 }\nolimits})$$ 
It is a continuous bounded function on ${\bf X}$.

\bigskip
We denote by\footnote{$^{26}$}{The initial $\Fa$ stands for Fresnel or Feynman according to the
worshipping habits of the reader.} $\Fa({\bf X}) $ the set of functions 
on ${\bf X}$ of the form $\Phi_\mu =
e^{\pi i Q} \cdot \Fa \mu$, where $\mu$ runs over the  measures on
${\bf X'}$. Since the map $\mu \mps \Fa \mu$ is injective, $\Fa({\bf X})$ is a Banach space, the
norm of the function $\Phi_\mu$ being taken equal to the total variation\footnote{$^{27}$}{According
to the standard definition, this is the l.u.b. of the set of numbers
$\ds \sum_{i =1 }^{ p} \left\vert{ \mu (A_i)  }\right\vert$ where 
$\left({ A_1, \cdots, A_p }\right)$ runs over the set of all partitions of ${\bf X'}$ into a
finite number of Borel subsets.
} 
 $\Var (\mu)$ of $\mu$.
On this Banach space, one defines a continuous linear form, denoted as an integral, by
$$\int_{ {\bf X}  }^{ } \Da x \cdot  \Phi_\mu (x) =
\int_{ {\bf X'}  }^{ }
d \mu  (x') \, e^{- \pi i W (x') }.
\eqno ({\mathop{\rm A.5 }\nolimits})$$

\smallskip
\nn
Hence by definition, we have
$$
\int_{ {\bf X}  }^{ } \Da x 
\int_{ {\bf X'}  }^{ }
d \mu (x') \exp \!
\left({  \pi i Q(x) - 2 \pi i 
\left\langle{ x', x }\right\rangle
}\right)
=
\int_{ {\bf X'}  }^{ }
d \mu  (x') \exp (- \pi i W(x')) \, .
\eqno ({\mathop{\rm A.6 }\nolimits})$$

\smallskip
\nn
Formally, this is the equation obtained by integrating equation (A.3) w.r.t. $d \mu (x') $ and
then interchanging the integrations:
$$\int_{ {\bf X'}  }^{ }
d \mu  (x') \int_{ {\bf X}  }^{ } \Da x  =
\int_{ {\bf X}  }^{ } \Da x  
\int_{ {\bf X'}  }^{ }
d \mu  (x') \, .$$

\bigskip
The space $\Fa ({\bf X})$ of Feynman-integrable functionals on ${\bf X}$
 is invariant under translations by elements of ${\bf X}$, and so is the integral, namely:
$$\int_{ {\bf X}  }^{ } \Da x  \cdot F (x) =
\int_{ {\bf X}  }^{ } \Da x \cdot 
F \! \left({ x + x_0 }\right)
\eqno ({\mathop{\rm A.7 }\nolimits})$$ 
or, in shorthand notation, 
$\Da x = \Da \! \left({ x + x_0 }\right)$, for any fixed element $x_0$ of ${\bf X}$.

\bigskip
According to equations (A.1) and (A.2), formula (A.3) can be rewritten as
$$
\int_{ {\bf X}  }^{ } \Da x  \cdot e^{\pi i Q(x - Gx')}
= 1. 
\eqno ({\mathop{\rm A.8 }\nolimits})$$ 

\smallskip
\nn
Hence {\it assuming  the invariance under translations
 of the integral}, the normalization of $\Da x$ is achieved by
$$
\int_{ {\bf X}  }^{ } \Da x \cdot 
 e^{\pi i Q(x )}
= 1. 
\eqno ({\mathop{\rm A.9 }\nolimits})$$ 

\bigskip
\nn
1.3. {\it The positive case $(s = 1)$.}
\smallskip
Henceforth, we assume that $s = 1$ and that $Q$ is positive-definite (that is $Q(x) > 0$ for $x
\ne 0$ in ${\bf X}$).
We take as basic integration formula
$$\int_{ {\bf X}  }^{ } \Da x 
\cdot  \exp \!
\left({ - \pi  Q(x) - 2 \pi i 
\left\langle{ x', x }\right\rangle
}\right) = \exp (- \pi W (x')) \, .
\eqno ({\mathop{\rm A.10 }\nolimits})$$ 
We can interpret this relation as above.
The conclusion is the integration formula
$$
\int_{ {\bf X}  }^{ } \Da x 
\int_{ {\bf X'}  }^{ }
d \mu (x') \exp \!
\left({  - \pi  Q(x) - 2 \pi i 
\left\langle{ x', x }\right\rangle
}\right)
=
\int_{ {\bf X'}  }^{ }
d \mu  (x') \exp (- \pi  W(x')) \, .
\eqno ({\mathop{\rm A.11 }\nolimits})$$ 

\smallskip
\nn
But the space $\Fa ({\bf X})$ is no longer invariant under translations. 
To restore this invariance, we have to replace Fourier-Stieltjes transforms by
{\it Laplace-Stieltjes transforms}. 
For that purpose, we have to consider the complex dual space
${\bf X}'_\CC$
consisting of the continuous real-linear maps
$x' : {\bf X}\ra \CC$ and measures $\mu$ on 
${\bf X}'_\CC$.
We leave the details to a forthcoming publication.

Formulas (A.3) and (A.10) are obtained as the specializations of formula (I.7) for $s = i$ and $s
= 1$ respectively.

\bigskip
\nn
{\bf 2. Linear changes of variables.}
\smallskip
\nn
2.1. {\it The case of Lebesgue integrals.}
\smallskip
Let ${\bf Y}$ be another real, separable, Banach space, and let 
$L : {\bf X} \ra {\bf Y}$ be a Borel-measurable map.
If $\om_{\bf X}$ is a measure on ${\bf X}$, its image under $L$ is the measure
$\om_{\bf Y}$ on ${\bf Y} $ defined by 
$\om_{\bf Y} (B) = \om_{\bf X} \!
\left({ L^{- 1} (B) }\right)$ for any Borel subset 
$B$ of ${\bf Y}$. In functional terms, this definition is tantamount to
$$
\int_{ {\bf X}  }^{ }  d \om_{\bf X} (x) \, 
g(L(x)) =
\int_{ {\bf Y}  }^{ }  d \om_{\bf Y} (y) \, g(y)
\eqno ({\mathop{\rm A.12 }\nolimits})$$ 
for any bounded (or non-negative) Borel-measurable function $g$ on ${\bf Y}$.

\bigskip
Assume now that $L$ is linear and continuous.
We consider the transpose $\ww L$ of $L$, that is the map
$\ww L : {\bf Y'} \ra {\bf X'}$ such that
$$
\langle \ww L y', x \rangle_{\bf X}
=
\langle y', Lx \rangle_{\bf Y}
\eqno ({\mathop{\rm A.13 }\nolimits})$$
for $x$ in ${\bf X}$ and $y'$ in ${\bf Y'}$.
The measures $\om_{\bf X}$ on ${\bf X}$ and 
$\om_{\bf Y}$ on ${\bf Y}$ being as above, introduce their Fourier-Stieltjes transforms 
$\Fa \om_{\bf X}$ and $\Fa \om_{\bf Y}$ respectively.
Hence $\Fa \om_{\bf X}$ is a bounded continuous function on 
${\bf X'}$, and similarly for 
$\Fa \om_{\bf Y}$ on ${\bf Y'}$.
By definition, we have
$$
\eqalignno{
\Fa \om_{\bf X} (x')
&
= \int_{ {\bf X}  }^{ } d \om_{\bf X} (x) 
\exp \!
\left({  - 2 \pi i 
\left\langle{ x', x }\right\rangle_{{\bf X}}
}\right)&({\mathop{\rm A.14 }\nolimits})\cr
\Fa \om_{\bf Y} (y')
&
= \int_{ {\bf Y}  }^{ } d \om_{\bf Y} (y) 
\exp \!
\left({  - 2 \pi i 
\left\langle{ y', y }\right\rangle_{{\bf Y}}
}\right) \, . &({\mathop{\rm A.15 }\nolimits})\cr}$$

\bigskip
By specializing 
$g(y) = \exp \!
\left({  - 2 \pi i 
\left\langle{ y', y }\right\rangle_{{\bf Y} } }\right)$
into formula (A.12) and taking into account formulas (A.13) to (A.15), we get
$\Fa \om_{\bf Y} (y') =
\Fa \om_{\bf X} (\ww L y')$ for every $y'$ in ${\bf Y'}$. Hence the composition formula:
$$\Fa \om_{\bf Y}  =
\Fa \om_{\bf X} \circ \ww L \, .
\eqno ({\mathop{\rm A.16 }\nolimits})$$ 
[The figure can be found in the Journal of Mathematical Physics {\bf 36} p.~2284]
\vglue 5cm 

\bigskip
The case of a translation is similar.
Assume now that $x_0$ is a given vector in ${\bf X} $
and denote by $T$ the translation taking $x$ into 
$x + x_0$ in ${\bf X}$.
If $\om$ is any measure on ${\bf X}$
and $\om_{x_0}$ its image under the translation $T$, a suitable specialization of formula (A.12)
gives the Fourier transform of $\om_{x_0}$, namely
$$\left({ \Fa \om_{x_0} }\right) (x') =
\exp \! \left({  - 2 \pi i 
\left\langle{ x', x_0 }\right\rangle
}\right)
\cdot \Fa \om (x') \, . 
\eqno ({\mathop{\rm A.17 }\nolimits})$$ 

\smallskip
\nn
2.2. {\it Infinite-dimensional integrators.}
\smallskip
We go back to the setup of paragraph A.1.1. We denote by $\Da x$ the integrator characterized by
$$\int_{ {\bf X}  }^{ } \Da x \cdot 
\exp \! \left({- { \pi \over s } Q(x)  }\right) \cdot 
\exp \! \left({  - 2 \pi i 
\left\langle{ x', x }\right\rangle
}\right)
=
\exp (- \pi s W (x')) \, .
\eqno ({\mathop{\rm A.18 }\nolimits})$$

\medskip
\nn
Formally, this means that the integrator $\Da \om $ defined by 
$$\Da \om (x) = \exp \!
\left({- { \pi \over s } Q(x)  }\right) 
\cdot \Da x
\eqno ({\mathop{\rm A.19 }\nolimits})$$ 
has a Fourier transform equal to 
$\exp (- \pi sW)$, namely:
$$\int_{ {\bf X}  }^{ } \Da \om (x) 
\exp \! \left({  - 2 \pi i 
\left\langle{ x', x }\right\rangle
}\right) = \exp 
\left({  -  \pi s 
W(x')
}\right). 
\eqno ({\mathop{\rm A.20 }\nolimits})$$ 

\smallskip
\nn
This can be interpreted as follows: the integrator $\Da \om$ is a continuous linear form on the
Banach space of Fourier-Stieltjes transforms $\Fa \mu$, given by
$$
\int_{ {\bf X}  }^{ } \Da \om (x) \Fa \mu (x) =
\int_{ {\bf X'}  }^{ } d \mu (x') \exp
\left({  -  \pi s 
W(x')
}\right). 
\eqno ({\mathop{\rm A.21 }\nolimits})$$  
This is another form of the Parseval relation.

\bigskip
We proceed to study the image of the integrator $\Da \om$ under a linear map.
We use the following notations:

\smallskip
-- $Q_{\bf X}$ is a quadratic form on ${\bf X}$, with inverse $W_{{\bf X}}$ 
on ${\bf X'}$;

\smallskip
-- 
$Q_{\bf Y}$ is a quadratic form on ${\bf Y}$, with inverse $W_{{\bf Y}}$ 
on ${\bf Y'}$;

\smallskip
-- $L$ is a continuous linear map from ${\bf X}$ into ${\bf Y}$.

\medskip
\nn
We assume that $L$ is surjective, hence $\ww L$ is injective, and that the quadratic forms
$W_{{\bf X}}$ and $W_{{\bf Y}}$ are related by
$$W_{{\bf Y}} = W_{{\bf X}} \circ \ww L \, .
\eqno ({\mathop{\rm A.22 }\nolimits})$$  
 
\bigskip
Consider now the integrators\footnote{$^{28}$}{
In the space ${\bf X} $, the integrator $\Da x$ depends on $Q_{{\bf X} } $
and $s$, and should be written more explicitly as 
$\Da_{s, Q_{{\bf X} } } x$. Similarly for $\Da y$ in space ${\bf Y}$.}
$$\eqalignno{
\Da \om_{{\bf X}} (x)
&= \exp \! 
\left({ - {\pi  \over s}  Q_{{\bf X} } (x) }\right)\cdot \Da x
&({\mathop{\rm A.23 }\nolimits})\cr
\Da \om_{{\bf Y}} (y)
&= \exp \! 
\left({ - {\pi  \over s}  Q_{{\bf Y} } (y) }\right) \cdot \Da y \, . 
&({\mathop{\rm A.24 }\nolimits})\cr
}
$$

\smallskip
\nn
The Fourier transforms are given respectively by 
$$\Fa \om_{{\bf X}} 
= \exp \! \left({ - \pi s  W_{{\bf X} }  }\right) \quad ,
\quad
\Fa \om_{{\bf Y}} 
= \exp \! \left({ - \pi s  W_{{\bf Y} }  }\right)
\eqno ({\mathop{\rm A.25 }\nolimits})$$  
and according to formula (A.22), we obtain
$$
\Fa \om_{{\bf Y}} 
= \Fa \om_{{\bf X}} 
\circ \ww L \, . 
\eqno ({\mathop{\rm A.26 }\nolimits})$$  

\smallskip
\nn
This is the same as formula (A.16), hence 
$\om_{{\bf Y}} $ is the image of $\om_{{\bf X}} $ under the linear mapping $L$.
Explicitly stated, we obtain the integration formula:
$$\int_{ {\bf X} }^{ } \Da x \cdot \exp \! 
\left({ - {\pi  \over s}  Q_{{\bf X} } (x) }\right)
\cdot g (Lx) =
\int_{ {\bf Y} }^{ } \Da y \cdot \exp \! 
\left({ - {\pi  \over s}  Q_{{\bf Y} } (y) }\right)
\cdot g (y) \, . 
\eqno ({\mathop{\rm A.27 }\nolimits})$$  

\smallskip
\nn
This holds if $g$ is a Fourier-Stieltjes transform
$\Fa \nu$ for some measure $\nu$
on
${\bf Y'}$.

\bigskip
The relationship between the quadratic forms 
$Q_{{\bf X}}$ on ${\bf X}$ and $Q_{{\bf Y}}$ on ${\bf Y}$
is expressed by the relation (A.22) between their inverses. Let
$D_{{\bf X}} : {\bf X} \ra {\bf X'}
$ be the continuous linear map such that 
$Q_{{\bf X}} (x) = 
\left\langle{ D_{{\bf X}} x, x }\right\rangle_{ {\bf X} }$
and 
$\left\langle{ D_{{\bf X}} x_1, x_2 }\right\rangle_{ {\bf X} }=
\left\langle{ D_{{\bf X}} x_2, x_1 }\right\rangle_{ {\bf X} }$,
 and define similarly 
 $ D_{{\bf Y}}$. Between the inverses $G_{\bf X}$ of $D_{{\bf X}}$ and 
$G_{{\bf Y}}$ of $D_{{\bf Y}}$ there holds the relation:
$$G_{{\bf Y}} = L \circ G_{{\bf X}} \circ \ww L \, . 
\eqno ({\mathop{\rm A.28 }\nolimits})$$  

\medskip
\nn
We mention a few particular cases:

\smallskip
-- if $L$ is invertible, with inverse $L^{- 1}$, then
$Q_{{\bf Y}} = Q_{{\bf X}} \circ L^{- 1}$;

\smallskip
-- if $Q_{{\bf X}}$ is positive-definite, then
$Q_{{\bf Y}} (y)$ is the infimum of $Q_{{\bf X}}$ over the set of elements $x$ in 
${\bf X}$ such that  $Lx = y$;

\smallskip
--
for any given $y$ in ${\bf Y}$ the equation
$$D_{{\bf X}} x =
\left({ \ww L \circ D_{{\bf Y}} }\right) (y)
\eqno ({\mathop{\rm A.29 }\nolimits})$$  
has a unique solution $x = x(y)$ in ${\bf X}$ and we obtain
$$Q_{{\bf Y}} (y) = Q_{{\bf X}} (x(y)) \, .
\eqno ({\mathop{\rm A.30 }\nolimits})$$

\bigskip
\nn
{\bf 3. Examples and applications.}

\smallskip
\nn
3.1. {\it The finite-dimensional case.} 
\smallskip
We record here the basic formulas. 
Assume that ${\bf X}$ is of finite dimension $d$; after choosing a linear frame, we represent a
vector $x$ by a column matrix $(x^\a)$ where
$\a \in 
\left\lbc{ 1, \cdots, d }\right\rbc$.
The elements of ${\bf X'}$ are represented by row matrices and the duality is given by
$\left\langle{ x', x }\right\rangle
= x'_\a x^\a$ (Einstein's summation convention).
The volume element is given by 
$d x = dx^1 \cdots dx^d$. 

\bigskip
The quadratic forms $Q$ and $W$ correspond to  symmetric matrices, namely:
$$Q(x) = h_{\a \b} x^\a x^\b \quad , \quad 
W(x') = h^{\a \b} x'_\a x'_\b
\eqno ({\mathop{\rm A.31 }\nolimits})$$  
with 
$h_{\a \b} h^{\b \g } = \d_\a^\g$.
The integrators are given by 
$$
\eqalignno{
\Da x
& = \left\vert{  \det h_{\a \b } }\right\vert^{1/2}
dx^1 \cdots dx^d 
& ({\mathop{\rm A.32 }\nolimits}) \cr
\Da \om (x)
& = 
\exp \! 
\left({ - \pi h_{\a \b} x^\a x^\b }\right) \cdot \Da x
& ({\mathop{\rm A.33 }\nolimits}) \cr}$$
when $s = 1$ and the matrix 
$\left({  h_{\a \b}  }\right)$ is positive-definite.
In the oscillating case, we have to multiply 
$\Da x$ by $e^{\pi i (q - p) /4}$ where the symmetric matrix
$\left({  h_{\a \b}  }\right)$ has $p$ positive and $q$ negative eigenvalues.

\bigskip
The quadratic form $W$ gives the covariance. More precisely, one obtains
$$\int_{ {\bf X} }^{ } \Da \om (x) 
\left\langle{ x', x }\right\rangle^2 =
{ s \over 2 \pi  } W(x')  
\eqno ({\mathop{\rm A.34 }\nolimits})$$
hence
$$\int_{ {\bf X} }^{ } \Da \om (x) 
x^\a x^\b =
{ s \over 2 \pi  } h^{\a \b }  \, .
\eqno ({\mathop{\rm A.35 }\nolimits})$$

\bigskip
\nn
3.2. {\it Image under a linear form.}
\smallskip
Let $x'$ in the dual ${\bf X'}$ of  ${\bf X}$. By specializing the results of paragraph A.2.2 to
the linear map $L : x \mps \left\langle{ x', x }\right\rangle$ from 
${\bf X}$ into $\RR$, we obtain the following  result. We identify $\RR$ with its dual, hence
$\ww L $ takes a number $\l$ to $\l x'$ in ${\bf X'}$. With $Q_{{\bf X}}$ equal to $Q$, hence
$W_{{\bf X}}$ to $W$ , we obtain
$$W_\RR (\l) = \l^2 W (x'), 
\eqno ({\mathop{\rm A.36 }\nolimits})$$
hence
$$Q_\RR (u) = u^2  / W (x') . 
\eqno ({\mathop{\rm A.37 }\nolimits})$$

\medskip
\nn
>From formula (A.27), we deduce the following integration formula:
$$\int_{{\bf X} }^{ }
\Da x \cdot e^{- \pi Q (x) / s } 
g 
\left({ \left\langle{ x', x }\right\rangle }\right)
 =
C \int_{ \RR }^{ } du \, e^{- \pi u^2  / s W (x')} 
g(u) 
\eqno ({\mathop{\rm A.38 }\nolimits})$$
where the normalization constant $C$ is 
$1 / (s W (x'))^{1/2}$ (principal branch of the square root). More explicitly:

\smallskip
-- for $s = 1$, hence $W(x') > 0$, then 
$ C = {\ds  1\over \ds\sqrt{ W (x')  }  } $;

\smallskip
-- for $s = i$, and $W(x') > 0$, then 
$ C = { \ds e^{- \pi i / 4} \over \ds  \sqrt{ W (x')  }  } $;

\smallskip
--  for $s = i$, and $W(x') < 0$, then 
$ C = {\ds  e^{\pi i / 4}  \over \ds \sqrt{ |W (x')|  }  } $.

\medskip
\nn
If we take in particular $g(u) = u^2 $, we get
$$\int_{{\bf X} }^{ }
\Da x \cdot e^{- \pi Q (x) / s } 
 \left\langle{ x', x }\right\rangle^2 =
 { s \over 2 \pi  } W(x') \, . 
\eqno ({\mathop{\rm A.39 }\nolimits})$$

\medskip
\nn
By polarization, we obtain the more general formula
$$
\int_{{\bf X} }^{ }
\Da x \cdot e^{- \pi Q (x) / s } 
 \left\langle{ x', x }\right\rangle
\left\langle{ y', x }\right\rangle
= 
{s  \over 2 \pi }
\left\langle{ x', Gy' }\right\rangle
\eqno ({\mathop{\rm A.40 }\nolimits})$$
for $x'$ and $y'$ in ${\bf X'}$, thus giving the covariance of our integrator. 

\bigskip
\nn
{\it Remark.}
In the application IV.1 to point-to-point transitions, we need to restrict the domain of integration (see also
paragraph A.3.8). This can be achieved by inserting a delta function
$\d \left({ \left\langle{ x'_0, x }\right\rangle }\right)$ and then using equation 
(A.38): we restrict the domain of integration from the space ${\bf X}$ to the hyperplane
${\bf X}_0$ with equation $\left\langle{ x'_0, x }\right\rangle  = 0$.

\bigskip
In general, any integral of the form
$$
\int_{{\bf X} }^{ }
\Da x \cdot e^{- \pi Q (x) / s } 
f \!
\left({  \left\langle{ x'_1, x }\right\rangle , \cdots , 
\left\langle{ x'_n, x }\right\rangle }\right)
$$
with fixed elements 
$x'_1, \cdots , x'_n $ in ${\bf X}'$ can be reduced to an $n$-dimensional integral.
For instance, combining formulas (A.38) and (A.40), the Green's function 
$G_{a,b} (t, u)$ on ${\bf Z}_{a, b}$ given by formulas (A.58) and (B.28), can be evaluated using an integral
$$
\int_{{\bf Z}_b }^{ }
\Da_b z \cdot \exp 
\left({ - { \pi \over s }
\int_{{\bf T} }^{ }dt \, z(t)^2
 }\right)
\d \! \left({ z (t_a) }\right) z(t) z(u)$$
and reducing it to a 3-dimensional integral.

\vfill\eject
\nn
3.3. {\it Space of paths of finite action.}
\smallskip
The time interval ${\bf T} = 
\left\lbk{ t_a,t_b }\right\rbk$ being given, we denote by $L^{2,1}$
(or more accurately $L^{2,1} ({\bf T})$) the space of real-valued functions 
$z(\cdot)$ on ${\bf T}$ with square-integrable derivative
$\dot z (\cdot)$ and we define the quadratic form $Q_0$ on $L^{2,1}$ by
$$Q_0 (z) =
\int_{ {\bf T} }^{ } dt \, 
\dot z (t)^2 \, .
\eqno ({\mathop{\rm A.41 }\nolimits})$$

\medskip
\nn
The definition of $L^{2,1}$ can be rephrased as follows:
the function $z$ belongs to $L^{2,1}$ if and only if there exists a function $\dot z$ in 
$L^2 ({\bf T})$ such that
$$z (t') - z(t) =
\int_{ t }^{ t' } du\, \dot z (u)  
\eqno ({\mathop{\rm A.42 }\nolimits})$$
whenever $t, t' $ are epochs in ${\bf T}$ such that $t < t'$.
The function $\dot z$ is then defined up to a null set; indeed by Lebesgue's derivation theorem,
one gets
$$\dot z(t) = 
\build { \lim }_{ \tau = 0 }^{ }
\left({ z(t + \tau) - z(t) }\right)/ \tau
\eqno ({\mathop{\rm A.43 }\nolimits})$$
for almost all $t$ in ${\bf T}$. By Cauchy-Schwarz inequality, one deduces from (A.42) the
inequality
$$
\left\vert{ z(t') - z(t) }\right\vert_{}^{2}
\leq Q_0 (z) \cdot 
\left\vert{ t' - t }\right\vert
\eqno ({\mathop{\rm A.44 }\nolimits})$$
for $t$, $t'$ in ${\bf T}$. Hence any function $z$ in $L^{2,1}$ satisfies a Lipschitz condition
of order $1/2$, and {\it a fortiori} it is a continuous function. 

\bigskip
The quantity $Q_0 (z)$ in equation (A.41) can be calculated using a {\it discretization of time}.
Consider a subdivision $\Ta$ of the time interval ${\bf T}$ by epochs
$$t_a \leq  t_0 < t_1 < \dots < t_{N-1} < t_N \leq  t_b\, .$$

\medskip
\nn
Set $\D t_i = t_i - t_{i-1}$ for 
$0 \leq i \leq N + 1$ with the convention $t_{-1} = t_a$, $t_{N+1}=t_b$  and denote by 
$\d (\Ta) $ the largest among the increments $\D t_i$, that is the {\it mesh} of the subdivision
$\Ta$. For any function $z : {\bf T} \ra \RR$ set
$$z_i = z(t_i) \quad , \quad \D z_i = z_i - z_{i-1}
\eqno ({\mathop{\rm A.45 }\nolimits})$$
for $0 \leq i \leq N + 1 $.
The {\it quadratic variation of} $z$ w.r.t. the subdivision $\Ta$ is defined as
$$Q_\Ta (z) = 
\sum_{ i = 1}^{ N } 
\left({ \D z_i }\right)_{}^{2}
/
\D t_i \, .
\eqno ({\mathop{\rm A.46 }\nolimits})$$

\medskip
\nn
Then the function $z$ belongs to $ L^{2,1}$ if and only if the set of quadratic variations 
$Q_{\Ta} (z)$ is bounded when $\Ta$ runs over all subdivisions of $ {\bf T} $. Then
$$Q_0 (z) = {\mathop{\rm l.}\nolimits} 
\build { {\mathop{\rm u}\nolimits}  }_{ \Ta}^{ }
 {\mathop{\rm .b. }\nolimits}\ Q_\Ta (z) \, .
\eqno ({\mathop{\rm A.47 }\nolimits})$$

\medskip
\nn
More precisely, for any sequence of subdivisions $\Ta (n)$ (for $n \in \NN$) whose mesh $\d (\Ta
(n))$ tends to $0$, one gets
$$Q_0(z) = 
\build { \lim }_{ n = \Inf  }^{ }
Q_{\Ta (n)}  (z) \, .
\eqno ({\mathop{\rm A.48 }\nolimits})$$

\medskip
\nn
It is therefore justified to write  $Q_0 (z) $ in the form 
$\ds \int_{{\bf T} }^{ }
{(dz)^2 \over  dt} .$

\bigskip
\nn
3.4. {\it Green's functions.}
\smallskip
Fix an element $t_0$ in ${\bf T}$ and denote by 
${\bf Z}_0$ the space of functions $z$ in $L^{2,1}$ such that 
$z (t_0) = 0$. This is a (real) Hilbert space, with scalar product
$$
\left\langle{ z_1 \bigm|  z_2 }\right\rangle
= \int_{{\bf T} }^{ } dt\, 
\dot z_1 (t) \dot z_2 (t)  ,
\eqno ({\mathop{\rm A.49 }\nolimits})$$
hence $Q_0 (z)$ is equal to 
$\left\langle{ z \bigm|  z }\right\rangle$.
>From (A.44), one deduces the inequality
$$\left\vert{ z(t) }\right\vert^2 \leq 
\left\vert{ t - t_0 }\right\vert
\cdot
\left\langle{ z \bigm|  z }\right\rangle,
\eqno ({\mathop{\rm A.50 }\nolimits})$$
and therefore there is an element $\d_t$ in the dual 
${\bf Z}'_0 $ of ${\bf Z}_0$ such that $z(t) =
\left\langle{ \d_t, z }\right\rangle$
(that is: $\d_t $ is a Dirac ``function'' centered at $t$).
According to the general theory, one introduces
a linear continuous and invertible map
$G_0 : {\bf Z}'_0 \ra {\bf Z}_0$ corresponding to the quadratic form $Q_0$ on ${\bf Z}_0$. 
The 
{\it Green's function} is defined as follows:
$$G_0(t, u) = 
\left\langle{ \d_t , G_0 \d_u }\right\rangle.
\eqno ({\mathop{\rm A.51 }\nolimits})$$
That is, the function $t \mps G_0(t, u)$ belongs to ${\bf Z}_0$ and is equal to 
$G_0 \d_u$. By definition, one gets
$$z(u) = 
\left\langle{ G_0 \d_u, z }\right\rangle
\eqno ({\mathop{\rm A.52 }\nolimits})$$
for every $z$ in ${\bf Z}$. More explicitly, these conditions on $G_0(t, u)$ can be expressed
as follows:
$$ \quad \qquad \qquad \qquad\left\{\matrix{
  G_0 \!\left({ t_0 , u }\right) = 0
   \hfill \cr
\noalign{\medskip}
\ds \int_{ {\bf T} }^{ } dt \, {\ds \part  \over \ds \part t} G_0(t, u) {\ds \part  \over \ds \part
t} z(t)  = z(u). 
  \hfill \cr}
\right. \qquad    \quad    \quad \qquad \qquad
\left.\matrix{
   ({\mathop{\rm A.53 }\nolimits})\cr
\noalign{\medskip}
 ({\mathop{\rm A.54}\nolimits})
  \cr}
\right. 
$$

\medskip
\nn
The unique solution to these equations is given by
$$
G_0(t, u) = 
\left\{\matrix{
\inf \left({ t - t_0, u - t_0 }\right)\hfill &
{\mathop{\rm for }\nolimits}
\ t \geq t_0 , \ u \geq t_0 , 
  \hfill\cr
\noalign{\medskip}
\inf \left({  t_0 - t ,  t_0 - u  }\right)\hfill &
{\mathop{\rm for }\nolimits}
\ t \leq t_0 , \ u \leq t_0 , 
   \hfill \cr
\noalign{\medskip}
0\hfill &
{\mathop{\rm otherwise }\nolimits}
\, .
  \hfill \cr}
\right. \eqno ({\mathop{\rm A.55 }\nolimits})$$

\medskip
We consider some special cases:

\smallskip
-- space ${\bf Z}_a$ of functions $z$ in $L^{2,1}$ with 
$z \!\left({ t_a }\right) = 0$;

\smallskip
-- space ${\bf Z}_b$ of functions $z$ in  $L^{2,1}$ with 
$z \!\left({ t_b }\right) = 0$;

\smallskip
-- space ${\bf Z}_{a, b} = {\bf Z}_a \cap {\bf Z}_b $ with boundary conditions 
$z \!\left({ t_a }\right) = z \!\left({ t_b }\right) = 0.$

\medskip
\nn
The corresponding Green's functions are as follows:
$$\eqalignno{
G_a (t, u) 
&=
\t(t - u) 
\left({ u - t_a }\right) + \t (u - t) \left({ t - t_a }\right)
& ({\mathop{\rm A.56 }\nolimits}) \cr
G_b (t, u) 
&=
\t(t - u) 
\left({ t_b - t }\right) + \t (u - t) \left({ t_b - u }\right)
& ({\mathop{\rm A.57 }\nolimits}) \cr
G_{a, b} (t, u) 
&=
 \t (t - u) \left({ t - t_b }\right) \left({ t_b - t_a }\right)^{-1 } 
\left({ t_a - u }\right)  
 \cr
& \ \  \  - \t (u - t) \left({ t - t_a }\right) \left({ t_a - t_b }\right)^{-1 } 
\left({ t_b - u }\right)  
 \, . 
& ({\mathop{\rm A.58 }\nolimits})\cr
}$$

\bigskip
\nn
{\it Remark. }
Equation (A.58) can be obtained by integrating
$$
\int_{{\bf Z}_b }^{ }
\Da \om (z) \, \d \! \left({ z(t_a) }\right)
\, z(t) \, z(u) \, .
$$

\bigskip
\nn
3.5. {\it Vector-valued functions.}
\smallskip
As in paragraph II.1.1, we consider vector functions
$z = \left({ z^1, \cdots, z^d }\right)
$ with components in $L^{2,1}$.
We introduce a real symmetric matrix
$\left({ h_{\a \b} }\right)$ with an inverse denoted by 
$\left({ h^{\a \b} }\right)$. The basic quadratic form is given by
$$Q_0 (z) = 
\int_{ {\bf T} }^{ } dt\,  h_{\a\b }\, 
\dot z^\a (t) \, \dot z^\b (t) \, . 
\eqno ({\mathop{\rm A.59 }\nolimits})$$

\medskip
\nn
The auxiliary condition 
$z (t_0) = 0$ defines the space ${\bf Z}_0$. Similarly, the spaces
${\bf Z}_a$, ${\bf Z}_b$ and ${\bf Z}_{a,b}$ are described by the respective boundary conditions:

\smallskip
-- $z(t_a) = 0 $ for ${\bf Z}_a$;

\smallskip
-- $z(t_b) = 0 $ for ${\bf Z}_b$;

\smallskip
-- $z(t_a) = z(t_b) = 0 $ for ${\bf Z}_{a,b}$.

\medskip
\nn
In each case, the linear form
$\d_t^\a$ taking $z$ into $z^\a (t)$ belongs to the dual of the corresponding of paths, and the
Green's function is characterized by 
$$G_0^{\a  \b} (t, u) =
\left\langle{ \d_t^\a, G_0 \d_u^\b  }\right\rangle \, . 
\eqno ({\mathop{\rm A.60 }\nolimits})$$

\medskip
\nn
Explicitly, we obtain
$$ G_0^{\a  \b} (t, u)  =
h^{\a \b } G_0(t, u)
\eqno ({\mathop{\rm A.61 }\nolimits})$$
where the Green's function $G_0$ refers to the boundary condition of the relevant space of paths
 (see formulas (A.55) to (A.58)).
The Green's function can also be obtained from formula (A.40) by specialization:
$$\int_{ {\bf Z} }^{ }
\Da z \cdot e^{- \pi Q_0 (z) / s}
\, z^\a (t) \, z^\b (u) =
{s \over 2 \pi } G_0^{\a \b }\,  (t, u) \, .
\eqno ({\mathop{\rm A.62 }\nolimits})$$
Here $\Da z $ stands for $\Da_{s, Q_0} z$.

\bigskip
\nn
{\it Remark}.
In general, the Green's functions $G$ satisfying a second-order differential equation $DG = \un$ are uniquely
determined by $d$ conditions at $t = t_a$ and $d$ conditions at $t = t_b$. These conditions are obvious from equation
(A.62) for the space ${\bf Z}_{a,b}$, namely $G^{\a \b } (t, u) = 0 $ if $t$ (or $u$) is one of the end points $t_a,
t_b$.

\medskip

In the case of the space ${\bf Z}_a$, the boundary conditions are
$$
\left.{ G }\right\vert_{t = t_a} = 0 \quad , \quad 
\left.{ \part G / \part t }\right\vert_{t = t_b} = 0 \, . $$
According to equation (A.62), these conditions can be expressed as 
$$
\eqalign{
\int_{ {\bf Z}_a }^{ } \Da \om (z) \cdot z^\a (t_a) z^\b (u) &= 0, \cr
\int_{ {\bf Z}_a }^{ } \Da \om (z) \cdot \dot z^\a (t_b) z^\b (u) &= 0\, .\cr}$$
The interpretation is as follows: any path $z$ in ${\bf Z}_a$ satisfies
$z^\a (t_a) = 0$; the time derivative 
$\dot z^\a (t_b) $ at $t_b$ is totally unspecified, nevertheless it vanishes in a statistical sense being
uncorrelated to the position and velocity at any other time.

\bigskip
On the space ${\bf Z}'_0$ dual to ${\bf Z}_0$, we have defined the quadratic form
$W_0$ inverse to $Q_0$. The Greens's function can be expressed as follows:
$$G_0^{\a \b} (t, u) =
{ {}_{1} \over {}^{2} }
\left\lbk{ W_0 \!
\left({ \d_t^\a + \d_u^\b }\right) -
W_0\!
\left({ \d_t^\a  }\right) -
W_0\! 
\left({  \d_u^\b }\right)
 }\right\rbk \, .
\eqno ({\mathop{\rm A.63 }\nolimits}) $$

\medskip
\nn
Conversely, given an element $z'$ of ${\bf Z'}$ represented by 
$$
\left\langle{ z', z }\right\rangle =
\int_{ {\bf T}J}^{ }dt \,  z'_\a (t) z^\a (t) , 
\eqno ({\mathop{\rm A.64 }\nolimits})$$
the quadratic form is given by
$$W_0(z') = \int_{ {\bf T}J}^{ } dt \int_{ {\bf T}J}^{ }
dt \, du \,
G_0^{\a \b} (t, u) 
z'_\a (t) \, z'_\b (u) 
\, .
\eqno ({\mathop{\rm A.65 }\nolimits})$$

\bigskip
\nn
3.6. {\it Scaling the paths.}

\smallskip
Here is a concise {\it dimensional analysis} of our quantities.

\nn
{\it Conventions} : $\La \ (\Ta)$ stands for the dimension of length (time) and $[X]$ for the
dimensional content of a quantity $X$.

\medskip
Since $Q_0 (z)$ appears in 
$e^{- \pi Q_0 (z) /s}$, it has to be a pure number, hence we can use formula (A.59) to deduce
$\left\lbk{ h_{\a \b } }\right\rbk
= \La^{- 2} \Ta$. From formulas (A.55) to (A.58) we infer the dimension of 
$G_0 (t, u)$ to be $\Ta$. 
>From (A.62) we infer that
$\left\lbk{G_0^{\a \b } }\right\rbk = \La^2$, hence 
$\left\lbk{h^{\a \b } }\right\rbk = \La^2 \Ta^{-1}$ from
(A.61). This is in accordance with the matrix relation 
$h_{\a \b } h^{\b \g }
= \d_\a^{\ \g}$. Notice also that for a particle of mass $m$ in a flat space, we have
$$
h_{\a \b } = m \d_{\a \b } /h 
\quad , \quad
h^{\a \b } = h \d_{\a \b } / m 
\eqno ({\mathop{\rm A.66 }\nolimits})$$
(where $h$ is the Planck constant and $\d_{\a \b}$ the Kronecker delta), and that the dimension of
$h/m $ is $\La^2 \Ta^{-1}$. We summarize our findings :

\bigskip
\centerline{TABLE 1}
$$\vbox{
\offinterlineskip
\halign
{\tv#&#  &\tv#&#  &#  &#  &#  &#  &#  &\tv#  \cr
\noalign{\hrule}
\tvi
\cc{Quantity}
&&\cc {$z^\a$} &\cc {$t$}
 &\cc {$h_{\a \b }$} &\cc  {$G_0 (t, u)$} &{$ G_0^{\a \b} (t, u)$} &
\cc {$h^{\a \b }$}  &&\cr
\tvi 
\cc{Dimension}
&&\cc {$\La$} &\cc {$\Ta$}
 &\cc {$\La^{-2} \Ta$} &\cc  {$\Ta$} &\cc {$ \La^2$} &
\cc {$\La^2 \Ta^{-1}$}  &&\cr
\noalign{\hrule} }}$$

\medskip
We can confirm these results by {\it scaling } our paths.
Let $\L > 0$ be a numerical constant, and denote by $\L  z$ the path with components
$\L z^\a (t) $ at time $t$. The basic integrator being written as 
$\Da \om (z) = e^{- \pi Q_0 (z) / s } \Da z$, we scale it into an integrator 
$\Da \om^\L (z) $ according to the formula
$$
\int_{ {\bf Z}_0 }^{ } \Da \om^\L (z) \, F(z) =
\int_{ {\bf Z}_0 }^{ } \Da \om (z) \, F (\L z) \, . 
\eqno ({\mathop{\rm A.67 }\nolimits})$$

\medskip
\nn
By suitably specializing $F(z)$, we obtain the new Green's function:
$$
\int_{ {\bf Z} }^{ } \Da \om^\L (z) \, 
z^\a (t) z^\b (u) =
\L^2 G_0^{\a \b} (t, u)\, ,
\eqno ({\mathop{\rm A.68 }\nolimits})$$
in accordance with 
$\left\lbk{ G_0^{\a \b } }\right\rbk = \La^2$. Similarly any $n$-point function scales as
$\La^n$. 
The map taking $z$ into $\L z$ is a linear map from ${\bf Z}_0 $ to ${\bf Z}_0 $.
By the general theory in paragraph A.2.2, $\Da \om^\L $ is another Gaussian integrator. From the
formulas (A.65) and (A.68), it corresponds to the quadratic form
$\L^2 W_0 (z')$ on ${\bf Z}'_0$, with an inverse quadratic form on ${\bf Z}_0$ given by
$$\L^{-2} Q_0 (z) = 
\int_{ {\bf T} }^{ } dt \, \L^{-2} \cdot
h_{\a \b} \, \dot z^\a (t) \dot z^\b (t) \, .
\eqno ({\mathop{\rm A.69 }\nolimits})$$

\medskip
\nn
This relation confirms $\left\lbk{  h_{\a \b} }\right\rbk = \La^{- 2} $.
A similar analysis applies to the scaling of time.

\bigskip
We can rewrite formula (A.67) in the following form:
$$
\int_{ {\bf Z}_0 }^{ }
\Da^\L z \cdot e^{- \pi \L^{-2} Q_0 (z) / s  }
H \!\left({ \L^{- 1} z }\right)
= \int_{ {\bf Z}_0 }^{ } \Da z \cdot e^{- \pi Q_0 (z) / s } 
H (z) \, ,
\eqno ({\mathop{\rm A.70 }\nolimits})$$
with a new integrator 
$\Da^\L z$  which is invariant under translations.
It is justified to summarize the previous formula by 
$\Da^\L z = \Da \left({ \L^{- 1} z}\right) $, hence
$\Da \om^\L (z) = \Da \om \!\left({ \L^{- 1} z}\right) $ according to formula (A.67).
In the standard heuristic derivations, one writes $\Da z$ in the form
$C \cdot \prod_{ t, \a  }^{ }
d z^\a (t) $. According to our normalization
$$
\int_{ {\bf Z}_0 }^{ } \Da z \cdot e^{- \pi Q_0 (z) / s } 
= 1, 
\eqno ({\mathop{\rm A.71 }\nolimits})$$
we should write
$$
\Da z =
{ \ds \prod_{t, \a  }^{ } d z^\a (t)  \over
\ds \int_{ }^{ } \ds \prod_{t, \a  }^{ } d z^\a (t) \cdot
\exp \left({ - \pi Q_0 (z) / s }\right)
 } \, , 
\eqno ({\mathop{\rm A.72 }\nolimits})$$
and similarly
$$
\Da^\L z =
{ \ds \prod_{t, \a  }^{ } d z^\a (t)  \over 
\ds \int_{ }^{ } \ds \prod_{t, \a  }^{ } d z^\a (t) \cdot
\exp \left({ - \pi Q_0 \left({ \L^{- 1} z }\right)  / s }\right)
} \, . 
\eqno ({\mathop{\rm A.73 }\nolimits})$$
Replacing $z$ by $\L^{- 1} z$ amounts to replacing 
$d z^\a (t) $ by $\L^{-1} dz^\a (t)$, hence the volume element 
$\prod_{t, \a  }^{ } d z^\a (t) $ is multiplied by $\L^{- N}$ where 
$N$ is the (infinite) number of degrees of freedom $t, \a$.
If we calculate $\Da \left({ \L^{- 1} z }\right)$ in accordance with (A.72), both numerator and 
denominator acquire a factor $\L^{- N}$ which drops out, and the correct formula
$\Da^\L z = \Da \left({ \L^{- 1} z }\right)$ is obtained from the heuristic formulas
(A.72) and (A.73).

\bigskip
The heuristic constant $\L^{- N}$ is equal to $\Inf$, 1 or 0 according to the three cases
$0 < \L < 1$, $\L = 1$, $\L > 1$. This is reflected in the rigorous theory by the following fact :
{\it for $\L \ne 1$  no functional $F (z)$, except the constant $0$, is such that both $F(z)$ and 
$F (\L z)$ can be simultaneously integrated } w.r.t. $\Da z$ In a finite-dimensional space $\RR^N$, the
volume element $dx^1 \cdots dx^N$ is the only one, up to a multiplicative constant, which is
invariant under translations. In our infinite-dimensional setup there exist many
translation-invariant integrators, but {\it they act in different functional sectors}.

\bigskip
\nn
3.7. {\it White noise representations.}

\smallskip
We consider now the Hilbert space $L^2 ({\bf T}) $ (or $L^2$) consisting of the square-integrable functions
$\xi : {\bf T} \ra \RR$, with the standard quadratic form
$$H(\xi) = \int_{ {\bf T} }^{ } dt \, \xi (t)^2 \, . 
\eqno ({\mathop{\rm A.74 }\nolimits})$$
As usual we identify 
$L^2 ({\bf T})$ with its own dual by means of the following scalar product
$$
\left\langle{ \xi', \xi  }\right\rangle
=
\int_{ {\bf T} }^{ } dt \, \xi' (t)  \xi (t) \, . 
\eqno ({\mathop{\rm A.75 }\nolimits})$$
Let us denote for a while the space
$L^2 ({\bf T})$ by ${\bf H}$. Under our general conventions, both
$D : {\bf H} \ra {\bf H' } $ and 
$G : {\bf H'} \ra {\bf H } $ are identity operators, hence the quadratic form 
${\bf H}$ is its own inverse.
The corresponding covariance is given by the kernel of the identity operator, namely
$\d (t - u)$. Hence denoting by $\Da \xi$ the basic integrator, we get from equation (A.40):
$$
\int_{ L^2 ({\bf T})  }^{ }
\Da \xi \cdot \exp 
\left({ 
- { \pi \over s}
\int_{ {\bf T} }^{ } dt \, \xi (t)^2 
 }\right)
\xi (t) \xi (u) =
{s \over 2 \pi}
\d (t - u) \, . 
\eqno ({\mathop{\rm A.76 }\nolimits})$$

\bigskip
Let 
$L_0^{2,1}({\bf T})$ (or $L_0^{2,1}$) be the subspace of functions $z$ in 
$L^{2,1} ({\bf T})$ such that $z(t_0) = 0$.
The spaces $L^2 ({\bf T})$ and $L_0^{2,1}({\bf T})$ are isomorphic under a correspondance 
$\xi \lra z$ where
$\xi = \dot z $ and conversely 
$z(t) = \int_{ t_0}^{ t} dt' \, \xi (t')$. This can be expressed by the formulas
$$z(t) =
\int_{ {\bf T} }^{ } dt'\, \T (t', t) \xi (t') \, , 
\eqno ({\mathop{\rm A.77 }\nolimits})$$
$$
\T (t', t) =
\left\{\matrix{
 1\hfill &  {\mathop{\rm for }\nolimits} \hfill & t_0 \leq t' \leq t,  \hfill\cr
\noalign{\medskip}
   - 1\hfill &  {\mathop{\rm for }\nolimits} \hfill & t \leq t' \leq t_0,  \hfill\cr
\noalign{\medskip}
 0\hfill &  {\mathop{\rm otherwise .}\nolimits} \hfill &   \hfill\cr}
\right. \eqno ({\mathop{\rm A.78 }\nolimits})$$

\medskip
\nn
The following transformation rule holds for functional integrals:
$$\int_{ L_0^{2,1} }^{ }
\Da z \cdot e^{- \pi Q_0 (z) / s} F(z) =
\int_{ L^{2} }^{ }
\Da \xi \cdot e^{- \pi H (\xi) / s}
\Phi (\xi) \, , 
\eqno ({\mathop{\rm A.79 }\nolimits})$$
where $\Phi (\xi)$ is equal to $F(z)$ if $\xi$ corresponds to $z$ by (A.77).
The covariance in 
$L_0^{2,1}$ is obtained by specializing $F(z)$ to $z(t) z(u) $ in (A.79).
Using (A.76), we obtain
$$G_0 (t, u) =
\int_{ {\bf T} }^{ } dt' \, \T (t', t) \T (t', u)
\eqno ({\mathop{\rm A.80 }\nolimits})$$
in agreement with formula (A.55).

\bigskip
To make contact with the heuristic definitions, we introduce the 
``coordinates'' 
$X_t = \xi (t) \sqrt{ dt } $ for the function 
$\xi$, hence
$$H(\xi) = 
\sum_{ t }^{ } X_t^2 \quad , 
\quad \Da \xi = \prod_{ t }^{ }
\left({ dX_t / \sqrt{ s}  }\right) \, .
\eqno ({\mathop{\rm A.81 }\nolimits})$$
We urge the reader to substantiate these claims by resorting to a subdivision 
$\Ta$ of the time interval ${\bf T}$, like in paragraph A.3.3.

\bigskip
Let ${\bf X}$ be any closed vector subspace of
$L^2 ({\bf T})$, and let $\Pi$ be the orthogonal projector from
$L^2 ({\bf T} )$ onto ${\bf X}$, represented as an integral operator with kernel $\Pi (t, u)$.
The quadratic form $H_{\bf X}$ on ${\bf X}$ is obtained by restriction of $H$, namely:
$$H_{\bf X} (\xi) = 
\int_{ {\bf T} }^{ } dt \, \xi (t)^2
\eqno ({\mathop{\rm A.82 }\nolimits})$$
for $\xi$ in ${\bf X}$.
Denote by 
$\Da_{{\bf X}} \xi$ the corresponding integrator on ${\bf X}$.
Then the covariance is expressed as follows:
$$\int_{ {\bf X} }^{ }
\Da_{{\bf X} } \xi \cdot \exp 
\left({ -{ \pi \over s}
\int_{ {\bf T} }^{ }
dt \, \xi (t)^2
 }\right)
\cdot \xi (t) \xi (u) =
{ s \over 2 \pi } \Pi (t, u) \, .
\eqno ({\mathop{\rm A.83 }\nolimits})$$

\medskip
\nn
For instance, if we denote by $L_0^2 ({\bf T})$ the subspace of $L^2 ({\bf T})$ defined by the
condition $\int_{ {\bf T}}^{ } dt \, \xi (t) = 0$, the orthogonal projector is given by
$$(\Pi \xi) (u) = \xi (u) - 
{1 \over t_b - t_a} \int_{ {\bf T}}^{ } dt \, \xi (t),
\eqno ({\mathop{\rm A.84 }\nolimits})$$
hence
$$\Pi (t, u) = \d(t - u) - 
{1 \over t_b - t_a} \, . 
\eqno ({\mathop{\rm A.85 }\nolimits})$$
The map $\xi \mps z$ where $z(t) =
\int_{t_a }^{t } dt' \, \xi (t')$ takes isomorphically the space
$L_0^2 ({\bf T})$ onto the space denoted by 
${\bf Z}_{a,b}$ in paragraph A.3.4. It follows that the Green's function 
$G_{a,b}$ corresponding to the space ${\bf Z}_{a,b}$ is given by
$$G_{a,b}(t, u) = \int_{t_a }^{t } dv \int_{t_a }^{u} dv' \, \Pi (v, v') \, . 
\eqno ({\mathop{\rm A.86 }\nolimits})$$
Using (A.85), we obtain easily (A.58).

\bigskip
Explicit formulas for transformations, mapping a given quadratic form on
$L^2 ({\bf T}) $ into any given quadratic form on 
$L^{2,1} ({\bf T})$, can be found in 
[5, p. 274]. The reader is urged to extend the results of this paragraph to vector-valued
functions.

\bigskip
\nn
3.8. {\it Fixing the endpoints.}
\smallskip
We denote by 
$\Da_a z$, $\Da_b z$ and $\Da_{a, b} z$ respectively the integrators on the spaces
${\bf Z}_a$, ${\bf Z}_b$ and ${\bf Z}_{a, b}$ (see paragraph A.3.5.). These integrators are
related by the following formula:
$$
\int_{ {\bf Z}_a }^{ } \Da_a z \cdot e^{- \pi Q_0 (z) / s }
F(z) \d \!\left({ z \!\left({ t_b }\right) }\right)
= C^{- 1 }
\int_{ {\bf Z}_{a,b} }^{ } \Da_{a,b} \, z \cdot e^{- \pi Q_0 (z) / s }
F(z) \, ,
\eqno ({\mathop{\rm A.87 }\nolimits})$$
where the constant $C$ is given by
$$\eqalign{
C &= 
\left({ \det s G_a^{\a, \b } \left({ t_b, t_b  }\right) }\right)_{}^{1/2}\cr
&= \left({ \det sh^{\a \b } }\right)_{}^{1/2}
\left({ t_b - t_a }\right)_{}^{d/2}
\, . \cr} 
\eqno ({\mathop{\rm A.88 }\nolimits})$$
A similar formula (where the roles of $t_a$ and $t_b$ are exchanged) can be found in 
[5, pp. 284, 357].

\bigskip
We give a new derivation of this formula to illustrate our methods.
We simplify the notations by putting $\Da \om_a (z)$ equal to 
$\Da_a z \cdot e^{- \pi Q_0 (z) / s }$, and defining 
$\Da \om_{a, b} (z)$ similarly. According to our general strategy, we need only to prove formula
(A.87) in the special case
$$F(z) = e^{- 2 \pi i \left\langle{  z', z}\right\rangle }
\eqno ({\mathop{\rm A.89 }\nolimits})$$
where $z'$ is an element of ${\bf Z}'_a$. That is, we have to establish the following equality:
$$\int_{ {\bf Z}_a }^{ } \Da \om_a (z) \cdot e^{- 2\pi i
\left\langle{  z', z}\right\rangle  }
\d \! \left({ z \! \left({ t_b }\right) }\right)
= C^{- 1} 
e^{- \pi s W_{a, b } (z') } \, ,  
\eqno ({\mathop{\rm A.90 }\nolimits})$$
where $ W_{a, b }$ is the quadratic form on 
${\bf Z}'_{a,b}$ inverse of the restriction of $Q_0$ to the subspace 
${\bf Z}_{a,b}$ of ${\bf Z}_{a}$.

\bigskip
To prove (A.90), we use the well-known formula 
$$\d(x) =
\int_{\RR  }^{ }
e^{- 2 \pi i u x} \, du \, .
\eqno ({\mathop{\rm A.91 }\nolimits})$$
We shall consider in detail the scalar case $d = 1$ and leave the general case to the reader.
The left-hand side $L$ of (A.90) can be rewritten as
$$L = \int_{ \RR }^{ } du 
\int_{ {\bf Z}_a }^{ } \Da \om_a (z) \, 
\exp 
\left({ - 2 \pi i 
\left\langle{ z' + u \d_{t_b}, z }\right\rangle
 }\right), 
\eqno ({\mathop{\rm A.92 }\nolimits})$$
hence
$$L = \int_{ \RR }^{ } du  \, \exp 
\left({ - \pi s W_a 
\left({ z' + u \d_{t_b} }\right)
 }\right)
\eqno ({\mathop{\rm A.93 }\nolimits})$$
by (A.18). We develop now the exponent
$$W_a  \left({ z' + u \d_{t_b }  }\right)
= u^2 G_a \! \left({ t_b, t_b }\right)
+ 2 u 
\int_{ {\bf T} }^{ } dt \, z'(t) G_a 
\! \left({ t_b, t  }\right)
+ W_a (z') 
\eqno ({\mathop{\rm A.94 }\nolimits})$$
and use the integration formula
$$
\int_{ \RR }^{ } du \, 
e^{- \pi s  
\left({ au^2 + 2 bu + c  }\right) }=
{1 \over \sqrt{ as }  }
\, \exp 
\left({ - \pi s 
\left({ c - b^2/a }\right)
}\right)
\eqno ({\mathop{\rm A.95 }\nolimits})$$
to obtain 
$$L = 
{1  \over \sqrt{s G_a \! \left({ t_b, t_b  }\right) } } \cdot 
\eqno ({\mathop{\rm A.96 }\nolimits})
$$
$$\exp 
\left({ - \pi s \! \left({  W_a (z')  
- 
G_a \! \left({ t_b, t_b  }\right)^{- 1 } 
\int_{ {\bf T} }^{ } dt 
\int_{ {\bf T} }^{ } du \,
z'(t) z'(u) 
G_a \!\left({ t_b, t  }\right)
G_a \! \left({ t_b, u  }\right)  }\right) }\right) \, . 
$$

\medskip
\nn
This can be transformed in the desired right-hand side of (A.90) with 
$C = \sqrt{s G_a \!\left({ t_b, t_b  }\right) }$ provided we establish the following identity:
$$G_{a, b} (t, u) = 
G_a(t, u) - G_a \!\left({ t, t_b }\right) 
G_a \!\left({ t_b, t_b }\right)^{- 1 }
G_a \!\left({ t_b, u }\right) \, .
\eqno ({\mathop{\rm A.97 }\nolimits})$$
We leave it as an exercice (see formulas (A.56) and (A.58)).

\bigskip
We can prove a more general formula by a similar reasoning. We fix points 
$z_a$ and $z_b$ and consider $L^{2,1}$ paths $z(\cdot)$ with 
$z(t_a) = z_a$. We obtain the affine space 
$z_a + {\bf Z}_a$, and transport to it the integrator 
$\Da_a z$. Then, with a suitable constant $C_{a, b}$ one gets
$$
\int_{ z_a + {\bf Z}_a }^{ }
\Da_a z \cdot 
e^{- \pi Q_0 (z) / s } F(z) 
\d \left({   z(t_b) - z_b}\right)$$
$$= C_{a, b}^{-1}
\int_{ {\bf Z}_{a, b} }^{ } \Da_{a,b}  \z \cdot 
e^{- \pi Q_0 (\z) / s } F
\left({ z_{{\mathop{\rm cl }\nolimits}}  + \z }\right)
\, .
\eqno ({\mathop{\rm A.98 }\nolimits})$$
Here $ z_{{\mathop{\rm cl }\nolimits}}$ is the affine-linear map
$ z_{{\mathop{\rm cl }\nolimits}} (t) = \l t + \mu $ with endpoints 
$z(t_a) = z_a $, $z(t_b) = z_b$, and the path $\z$ in 
${\bf Z}_{a,b}$, that is 
$\z (t_a) = \z (t_b) = 0$, is a {\it quantum fluctuation.}

\bigskip
\nn
{\bf 4. Infinite-dimensional determinants.}
\smallskip
For this section, we refer the reader to [24], [25, chapter 5] and [26].

\bigskip
\nn
4.1. {\it The case of operators.}
\smallskip
We return to our Banach space ${\bf X}$.
We assume that there exists on ${\bf X}$ an invertible\footnote{$^{29}$}{
A positive-definite continuous quadratic form is not necessarily invertible.
For instance, let ${\bf X}$ be the space $\ell^2$ of sequences
$\left({ x_1, x_2, \cdots }\right)$ of real numbers with
$\sum_{n = 1 }^{ \Inf } \left({ x_n }\right)^2 < \Inf
$ and define the norm by
$\left\Vert{ x }\right\Vert_{}^{2}=
\sum_{ n=1}^{ \Inf }
\left({ x_n }\right)^2
$.
We can identify ${\bf X}$ with its dual ${\bf X'}$, the scalar product being given by 
$\sum_{ n=1}^{ \Inf } x'_n x_n$.
The quadratic form 
$Q(x) = \sum_{ n=1}^{ \Inf } \left({ x_n / n }\right)^2$ corresponds to the map
$D : {\bf X } \ra {\bf X' }$ taking 
$\left({  x_1, x_2, \cdots}\right)$ into
$\left({ x_1/1, x_2/2, \cdots }\right)$. The inverse of $D$ does not exist as  a map from
$\ell^2$ into $\ell^2$ since the sequence $1, 2, 3, ...$ is unbounded.
 }
positive-definite quadratic form $H$ (for instance
$H(z) = \sum_{ \a }^{ }
\int_{ {\bf T } }^{ }
dt\, \dot z^{\a } (t)^2$ on the space ${\bf Z}_0$). 
We use the norm defined on ${\bf X}$ by
$\left\Vert{ x }\right\Vert
= H(x)^{1/2}$, and derive from it the dual norm
$\left\Vert{ x' }\right\Vert$ on ${\bf X'}$ as usual.

\bigskip
According to Grothendieck, an operator $T$ in ${\bf X}$ is called {\it nuclear} if it admits a
representation in the form
$$Tx = \sum_{ n \geq 0 }^{ }
\left\langle{ x'_n, x }\right\rangle
\cdot x_n
\eqno ({\mathop{\rm A.99 }\nolimits})$$
with elements 
$x_n$ in ${\bf X}$ and $x'_n$ in ${\bf X'}$ such that
$\sum_{ n \geq 0 }^{ }
\left\Vert{ x_n }\right\Vert
\cdot
\left\Vert{ x'_n }\right\Vert$ be finite.
The g.l.b. of all such sums 
$\sum_{ n  }^{ }
\left\Vert{ x_n }\right\Vert
\cdot
\left\Vert{ x'_n }\right\Vert$ is called the {\it nuclear norm} of $T$,
to be denoted by
$\left\Vert{ T }\right\Vert_1$. The nuclear operators in ${\bf X}$ form a Banach space, denoted by
$\La^1 ({\bf X})$, with norm 
$\left\Vert{ \cdot }\right\Vert_1$.
On $\La^1 ({\bf X})$, there exists a continuous linear form, the {\it trace}, such that
$$\Tr (T) = \sum_{ n \geq 0 }^{ }
\left\langle{ x'_n, x_n }\right\rangle 
\eqno ({\mathop{\rm A.100 }\nolimits})$$
for an operator $T$ given by (A.99).

\bigskip
We introduce a power series in $\l$, namely:
$$
\sum_{p \geq 0  }^{ } 
\s_p (T) \l^p : = 
\exp 
\left({ \l \Tr (T) -
{\l^2 \over 2}  \Tr (T^2) + {\l^3 \over 3} \Tr (T^3) - \cdots
 }\right) \, . 
\eqno ({\mathop{\rm A.101 }\nolimits})$$
By using Hadamard's inequality about determinants, we obtain the {\it basic estimate}:
$$\left\vert{ \s_p (T) }\right\vert
\leq p^{p/2}
\left\Vert{ T }\right\Vert_{1}^{p}
/ p!\, .
\eqno ({\mathop{\rm A.102 }\nolimits})$$

\medskip
\nn
It follows that the power series $\sum_{ p \geq 0 }^{ }  \s_p (T) \l^p $ has an infinite radius of
convergence. We can therefore define the determinant as follows:
$$\Det(1 + T):= \sum_{ p \geq 0 }^{ } \s_p (T)
\eqno ({\mathop{\rm A.103 }\nolimits})$$
for any nuclear operator $T$. By homogeneity, we obtain more generally
$$\Det (1 + \l T) = 
\sum_{ p \geq 0 }^{ } \s_p (T) \l^p \, .
\eqno ({\mathop{\rm A.104 }\nolimits})$$
The fundamental property of determinants is, as expected, the
{\it multiplicative rule}:
$$\Det \left({ U_1 \circ U_2 }\right)
= \Det \left({ U_1 }\right)
\Det \left({ U_2 }\right) \, , 
\eqno ({\mathop{\rm A.105 }\nolimits})$$
where $U_i$ is of the form $1 + T_i$, with $T_i$ nuclear
(for $i = 1, 2)$. From this, and the relation $\s_1 (T) = \Tr (T)$, we get a
{\it variation formula} (for $U - 1$ and $\d U$ nuclear):
$$
{ \Det (U + \d U) \over \Det (U) } 
= 1 + \Tr 
\left({ U^{- 1} \cdot \d U}\right)
+ O
\left({ \left\Vert{ \d U }\right\Vert_{1}^{2} }\right) \, .
\eqno ({\mathop{\rm A.106 }\nolimits})$$
Otherwise stated, if $U(\nu)$ is an operator of the form 
$1 + T(\nu)$, where $T(\nu)$ is nuclear, depending smoothly on the parameter $\nu$,
we get the {\it derivation formula}:
$${ d \over d \nu } 
\ell n \Det (U (\nu)) =
\Tr 
\left({ U (\nu)^{-1} {d \over d \nu} U (\nu) }\right) \, .
\eqno ({\mathop{\rm A.107 }\nolimits})$$

\bigskip
\nn
{\it Remark.}
For any other norm 
$\left\Vert{ \cdot }\right\Vert^1$ defining the topology of ${\bf X}$, we have an estimate
$$C^{- 1}
\left\Vert{ x }\right\Vert
\leq \left\Vert{ x }\right\Vert^1 \leq C 
\left\Vert{ x }\right\Vert \, , 
\eqno ({\mathop{\rm A.108 }\nolimits})$$
with a  finite numerical constant $C > 0$.
It follows easily that the previous definitions are independent of the choice of the particular
norm
$\left\Vert{ x }\right\Vert = H (x)^{1/2}$ in ${\bf X}$.

\bigskip
\nn
4.2. {\it Explicit formulas.}
\smallskip
Introduce a basis
$\left({  e_n}\right)_{n \geq 1}^{}$ of ${\bf X}$ orthonormal for the quadratic form 
$H$, hence $H ( \sum_{ n }^{ } t_n e_n) = \sum_{ n }^{ } t_n^2$.
An operator $T$ in ${\bf X}$ has a matrix $(t_{mn})$ such that
$$T e_n = \sum_{ m }^{ } e_m \cdot t_{mn} \, . 
\eqno ({\mathop{\rm A.109 }\nolimits})$$

\bigskip
Assume that $T$ is nuclear. 
Then the series $\sum_{ n }^{ }  t_{nn}$ of diagonal terms in the matrix converges absolutely and 
the trace
$\Tr (T)$ is equal to $\sum_{ n }^{ }  t_{nn}$, as it should be.
Furthermore $\s_p (T)$ is the sum of the series made of the principal minors of order $p$:
$$\s_p(T) = 
\sum_{ i_1 < \cdots < i_p}^{ }
\det \left({ t_{i_\a, i_\b} }\right)_{
{1 \leq \a \leq p \atop 1 \leq \b \leq p} }
\, . 
\eqno ({\mathop{\rm A.110 }\nolimits})$$

\medskip
\nn
For the operator $U = 1 + T$, with the matrix with elements 
$u_{mn} = \d_{mn} + t_{mn}$, we obtain the determinant as a 
{\it limit of finite-size determinants}:
$$\Det(U) = 
\build { \lim }_{ N = \Inf }^{ }
\det 
\left({ u_{mn} }\right)_{ {1 \leq m \leq N \atop 1 \leq n \leq N}  }^{}
\, . 
\eqno ({\mathop{\rm A.111 }\nolimits})$$

\bigskip
As a special case, suppose that the basic vectors $e_n$ are 
{\it eigenvectors} for $T$, namely
$$T e_n = \l_n e_n \, . 
\eqno ({\mathop{\rm A.112 }\nolimits})$$
Then we get
$$\Tr (T) = \sum_{ n }^{ } \l_n \quad , 
\quad \Det (1 + T) = 
\prod_{ n}^{ }
\left({ 1 + \l_n }\right)$$
where both the series and the infinite product converge absolutely.

\bigskip
The nuclear norm 
$\left\Vert{ T }\right\Vert_1$ can also be computed as follows: there exists an orthonormal basis
$(e_n)$ such that the vectors $T e_n$ are mutually orthogonal (for the quadratic form $H$) and
then
$\left\Vert{ T }\right\Vert_1 =
\sum_{ n }^{ } 
\left\Vert{ T e_n }\right\Vert$.

\bigskip
\nn
{\it Remark.}
Let $T$ be a continuous linear operator in ${\bf X}$. Assume that the series of diagonal terms 
$\sum_{ n }^{ } t_{nn}$ converges absolutely for {\it every} orthonormal basis.
Then $T$ is nuclear. When $T$ is symmetric and positive, it is enough to assume that this
statement holds for {\it one} given orthonormal basis, and then it holds for all. There are
counterexamples when $T$ is not symmetric and positive. 

\bigskip
\nn
4.3. {\it The case of quadratic forms.}
\smallskip
Contrary to a widespread misbelief, {\it there is no such thing like the determinant of a
quadratic form.}
Consider for instance a quadratic form $Q$ on some finite-dimensional space with coordinates
$x^1, \cdots , x^d$, namely
$$Q(x) = h_{\a \b } x^\a x^\b \, . 
\eqno ({\mathop{\rm A.113 }\nolimits})$$
If we introduce a new system of coordinates 
$\oo x^1, \cdots, \oo x^d$ such that 
$x^\a = u_\l^\a \oo x^\l $, then we obtain
$$Q(x) = \oo h_{\l \mu} \oo x^\l \oo x^\mu
\eqno ({\mathop{\rm A.114 }\nolimits})$$
with a new matrix
$$\oo h_{\l \mu } = u_\l^\a u_\mu^\b h_{\a \b }
\, . 
\eqno ({\mathop{\rm A.115 }\nolimits})$$
The determinants 
$D = \det \left({ h_{\a \b } }\right)$ and 
$\oo D = \det \left({ \oo h_{\l \mu } }\right)$ are connected by the 
{\it scaling relation}:
$$\oo D = D \cdot 
\left({ \det (u_\l^\a) }\right)_{}^{2} \, . \eqno ({\mathop{\rm A.116 }\nolimits})$$

\medskip
\nn
Hence, what makes sense is the {\it ratio of determinants}
$$\det 
\left({ h_{\a \b }^{(1)} }\right)
/
\det 
\left({ h_{\a \b }^{(0)} }\right)
\eqno ({\mathop{\rm A.117 }\nolimits})$$
associated to two quadratic forms
$$Q_0 (x) =
h_{\a \b }^{(0)} x^\a x^\b \quad , \quad 
Q_1 (x) = h_{\a \b }^{(1)} x^\a x^\b
\eqno ({\mathop{\rm A.118 }\nolimits})$$
on the same space. We shall denote it by 
$\det (Q_1 / Q_0)$.

\bigskip
To obtain an intrinsic definition, let us consider two continuous quadratic forms $Q_0$ and $Q_1$
on  a Banach space ${\bf X}$ and assume that $Q_0$ is invertible.
We associate to $Q_0$ and $Q_1$ two continuous linear maps
$D_i : {\bf X} \ra {\bf X'}$ such that
$$
\left\langle{ D_i x_1, x_2 }\right\rangle =
\left\langle{ D_i x_2, x_1 }\right\rangle
\quad , \quad
Q_i (x) =
\left\langle{ D_i x, x  }\right\rangle \, ,
\eqno ({\mathop{\rm A.119 }\nolimits})$$
for $x$, $x_1$ and $x_2$ in ${\bf X}$.
By assumption, $D_0$ is invertible, hence there exists a unique continuous 
operator $U$ in ${\bf X}$ such that 
$D_1 = D_0 \circ U$.
We denote $U$ by\footnote{$^{30}$}{A better notation should perhaps be $Q_0 \bl Q_1$. }
$Q_1/Q_0$. In case the determinant of $U$ is defined,
namely when $U - 1$ is nuclear, we denote it  by 
$\Det \left({ Q_1 / Q_0 }\right)$. Of course, when ${\bf X}$ is finite-dimensional, this
definition agrees with the previous one and
$$\Det \left({ Q_1 / Q_0 }\right) =
\det \left({ Q_1 / Q_0 }\right)
\eqno ({\mathop{\rm A.120 }\nolimits})$$
in this case.

\bigskip
A procedure to calculate this determinant is as follows. Let $V$ be a finite-dimensional subspace
of ${\bf X}$. By restricting $Q_0$ and $Q_1$ to $V$ we obtain two quadratic forms 
$Q_{0, V}$ and $Q_{1, V}$ on $V$. 
Assume now that $V$ runs through an increasing sequence of subspaces, whose union is dense in ${\bf X}$,
and  that $Q_{0, V}$ be invertible for every $V$.
Then 
$$\Det \left({ Q_1 / Q_0 }\right) =
\build {\lim }_{ V}^{ }
\det 
\left({ Q_{1,V} / Q_{0,V} }\right) \, . 
\eqno ({\mathop{\rm A.121 }\nolimits})$$

\bigskip
\nn
4.4. {\it A pencil of quadratic forms.}

\smallskip
Here are our assumptions:

\smallskip
-- {\it there exists an invertible positive-definite quadratic form on ${\bf X}$};

\smallskip
-- {\it $Q_0$ and $Q_1$ are continuous quadratic forms on ${\bf X}$};
\smallskip

-- {\it $Q_0$ is invertible};

\smallskip
-- {\it if we set $Q = Q_1 - Q_0$, the operator $T = Q / Q_0$ on ${\bf X}$ is nuclear.}

\smallskip
The last condition is an intrinsic property of the quadratic form $Q$, namely the existence of a
representation like
$$Q(x) =
\sum_{ n \geq 1 }^{ } 
\a_n 
\left\langle{ x'_n, x }\right\rangle^2, 
\eqno ({\mathop{\rm A.122 }\nolimits})$$
where the series
$\sum_{ n }^{ } \a_n$ converges absolutely, and the sequence of numbers
$\left({ 
\left\langle{ x'_n, x }\right\rangle
 }\right)_{n \geq 1}^{} $ is bounded for every vector $x$ in ${\bf X}$.
We express this property by saying that the {\it quadratic form $Q$ is nuclear.}

\bigskip
We say  that $\l$ is an {\it eigenvalue } of a quadratic form $Q'$ 
w.r.t. $Q_0$ if it is an eigenvalue of the operator $Q' / Q_0$.
This is tantamount to saying that the quadratic form $Q' - \l Q_0 $ 
is not invertible.

\bigskip
We interpolate between the quadratic forms $Q_0$ and 
$Q_1 = Q_0 + Q$ by putting $Q_\nu = Q_0 + \nu Q$ where $\nu$ is a real or complex parameter.
The determinant 
$$\D (\nu) = \Det (Q_\nu / Q_0)
\eqno ({\mathop{\rm A.123 }\nolimits})$$
is defined, being equal to $\Det (1 + \nu T)$.
This is an entire function of the complex variable $\nu$, hence it vanishes for a discrete
set of values of $\nu$ (possibly empty). Furthermore, $\l$ is an eigenvalue of $Q$ w.r.t. $Q_0$ if and only if the
quadratic form $Q_{- 1/\l}$ is non invertible, that is if and only if $\D (- 1 /\l) = 0$.

\medskip

As a consequence, the function of real variable $\nu \mps \D (\nu)$ has, at most, a finite number
of zeroes $\nu_{(1)}, \cdots , \nu_{(p)}$ in the interval $[0, 1]$ with 
$0 < \nu_{(1)} < \cdots < \nu_{(p)} \leq 1$. According to the relation (A.107), between two such
zeroes, the following differential equation holds:
$${ d \over d \nu } \ell n \Det (Q_\nu / Q_0) = \Tr (Q / Q_\nu) 
\eqno ({\mathop{\rm A.124 }\nolimits})$$
where $Q_\nu = Q_0 + \nu Q$.

\bigskip
\nn
4.5. {\it Some integration formulas.}
\smallskip
With the previous notations, consider any $\nu$ in $[0, 1]$ distinct from
$\nu_{(1)}, \cdots ,$ $ \nu_{(p)}$. Hence the quadratic form $Q_\nu$ on ${\bf X}$ is invertible,
with an inverse quadratic form $W_\nu$ on ${\bf X}'$. There exists an integrator 
$\Da_{\nu} \, x$ on ${\bf X}$ characterized by the formula
$$
\int_{ {\bf X} }^{ } \Da_\nu x\cdot \exp 
\left({ - { \pi \over s} Q_\nu (x) - 2 \pi i 
\left\langle{ x', x }\right\rangle
}\right) = \exp 
\left({ - \pi s \, W_\nu (x') }\right) \, . 
\eqno ({\mathop{\rm A.125 }\nolimits})$$
Our claim is that $\Da_\nu \, x$ {\it is proportional to } $\Da_0 x$, namely that there exists a
constant $I(\nu)$ such that $\Da_\nu \, x = I (\nu) \Da_0 x$. 
More explicitly, we assert the formula
$$
\int_{ {\bf X} }^{ } \Da_\nu x \cdot F (x) = I (\nu) 
\int_{ {\bf X} }^{ } \Da_0 x \cdot F(x) \, , 
\eqno ({\mathop{\rm A.126 }\nolimits})$$
and that a functional $F(\cdot)$ on ${\bf X}$ is Feynman-integrable for 
$\Da_\nu x$ if and only if it is Feynman-integrable for $\Da_0 x$.
The constant $I(\nu)$ can be obtained by putting 
$F(x) =
\exp \left({- { \pi \over s}  Q_\nu (x) }\right)$ into equation (A.126), hence
$$I(\nu)^{- 1 } = \int_{ {\bf X} }^{ }  
\Da_0 x \cdot \exp 
\left({ - { \pi \over s} 
\left({ Q_0 (x) + \nu Q (x)  }\right)
}\right) \, . 
\eqno ({\mathop{\rm A.127 }\nolimits})$$

\bigskip
\nn
{\bf First case $s = 1$:}
\smallskip
According to our conventions, the quadratic forms $Q_0$ and $Q_1$ are positive-definite and
invertible. Since the quadratic form $Q = Q_1 - Q_0$ is nuclear, it follows from the spectral
theory that $Q_0$ and $Q$ can be simultaneously diagonalized. Hence there exists a basis
$\left({ e_n }\right)_{n\geq 1}^{}$ for ${\bf X}$ such that
$$
\eqalignno{
Q_0 \left({ \sum_{ n=1}^{\Inf } t_n e_n  }\right)
&= \sum_{ n=1}^{\Inf } \left({ t_n }\right)^2
& ({\mathop{\rm A.128 }\nolimits}) \cr
Q \left({ \sum_{ n=1}^{\Inf } t_n e_n  }\right)
&= \sum_{ n=1}^{\Inf } \l_n \left({ t_n }\right)^2
& ({\mathop{\rm A.129 }\nolimits})\cr}$$
with real constants $\l_n$ such that 
$\sum_{ n=1}^{\Inf } \left\vert{ \l_n }\right\vert < \Inf $.
Since $Q_1 = Q_0 + Q$ is positive-definite, we have $1 + \l_n > 0$,
hence
$$Q_\nu  \left({ \sum_{ n=1}^{\Inf } t_n e_n  }\right)
= \sum_{ n=1}^{\Inf }
\left({ (1 - \nu) + \nu (1 + \l_n) }\right)
\cdot \left({ t_n }\right)^2
\eqno ({\mathop{\rm A.130 }\nolimits})$$
is again positive-definite and invertible for 
$\nu $ in $[0, 1]$.
Hence the determinant of 
$Q_\nu / Q_0$ is defined and
$$\Det (Q_\nu / Q_0) =
\prod_{n=1 }^{ \Inf}
\left({ 1 + \nu \l_n }\right) > 0 \, .
\eqno ({\mathop{\rm A.131 }\nolimits})$$

\medskip

The main result is given by the following formula:
$$I(\nu) = \Det (Q_\nu / Q_0)^{{1 \over 2}} 
\, . 
\eqno ({\mathop{\rm A.132 }\nolimits})$$
The proof is obtained without difficulty using equations (A.127) to (A.131) and the 
{\it approximation formula}:
$$\int_{ {\bf X} }^{ } \Da_0 x \cdot F (x) =
\build { \lim }_{ N = \Inf }^{ }
\int_{ \RR^N }^{ } d^N t \, 
F 
\left({ t_1 e_1 + \cdots + t_N e_n   }\right)  \, . 
\eqno ({\mathop{\rm A.133 }\nolimits})$$

\bigskip
\nn
{\bf Second case $s = i$:}
\smallskip
Here the basic formula is given by 
$$I (\nu) =
\left\vert{ \Det
\left({ Q_\nu / Q_0  }\right)
 }\right\vert^{1/2} 
i^{
{\mathop{\rm Ind  }\nolimits}
\left({ Q_\nu / Q_0  }\right)
}
\eqno ({\mathop{\rm A.134 }\nolimits})$$
where the {\it index} ${\mathop{\rm Ind  }\nolimits}
\left({ Q_\nu / Q_0  }\right)$ is the number of negative eigenvalues of 
$Q_\nu$ w.r.t. $Q_0$. To simplify the statements, we shall assume that $Q_0$ is
positive-definite. When $\nu$ runs over the interval $[0, 1]$, this index remains constant except
when $\nu$ is crossing a singular value $\nu_{(k)}$, where it experiences a jump.
Hence the index ${\mathop{\rm Ind  }\nolimits}
\left({ Q_\nu / Q_0  }\right)$ is a sum of local contributions from the exceptional values
$\nu_{(1)} , \cdots, \nu_{(p)}$, a phenomenon reminiscent of caustics.

\medskip
The easiest proof of formula (A.134) is obtained by following a strategy initiated by Nelson and
Sheeks in [6]. It works in both cases $s = 1$ and $s = i$.
Starting from equation (A.127), we obtain by derivation
$$
\eqalign{
{d  \over d \nu } I(\nu)^{-1}
&= \int_{{\bf X} }^{ } \Da_0 x \cdot \exp 
\left({ - {\pi  \over s}
\left({ Q_0 (x) + \nu Q (x)  }\right)
 }\right)
\left({ - {\pi  \over s} Q (x)  }\right) \cr
&= I (\nu)^{- 1 } 
\int_{{\bf X} }^{ } \Da_\nu x \cdot \exp 
\left({ - {\pi  \over s}
 Q_\nu (x)  
 }\right)
\left({ - {\pi  \over s} Q (x)  }\right) , \cr}
\eqno ({\mathop{\rm A.135 }\nolimits})$$
that is
$$
{d  \over d \nu } \, \ell n \, I(\nu)
= {\pi  \over s} \int_{ {\bf X} }^{ } 
\Da \om_\nu (x) \, Q(x) \, .
\eqno ({\mathop{\rm A.136 }\nolimits})$$

\medskip
\nn
Expanding $Q(x) $, according to formula (A.122) and using formula (A.39) for the covariance, we
obtain
$${d  \over d \nu } \, \ell n \, I (\nu) = 
{\pi  \over s} \cdot 
\sum_{n=1 }^{\Inf  }
\a_n \cdot {s  \over 2 \pi }
W_\nu (x'_n) 
\eqno ({\mathop{\rm A.137 }\nolimits})$$
which can be transformed easily into
$$
{d \over d \nu } \, \ell n \, I (\nu) = 
{ {}_{1} \over {}^{2} } \Tr \left({ Q / Q_\nu }\right)
\, . \eqno ({\mathop{\rm A.138 }\nolimits})$$

\medskip
\nn
According to equation (A.124), we conclude
$${d  \over d \nu } \, \ell n \, I (\nu)
= { {}_{1} \over {}^{2} } {d  \over d \nu } \, \ell n \, 
\Det \left({ Q_\nu / Q_0 }\right) \, . 
\eqno ({\mathop{\rm A.139 }\nolimits})$$

\medskip
\nn
It remains to study the shift in phase when $\nu$ goes through an exceptional value
$\nu_{(k)}$.

\bigskip
\nn
{\it Remark.}
To simplify matters, assume again  that $Q_0$ is positive-definite. Let $L$ be an invertible operator in
${\bf X}$, of the form $1 + T$ where $T$ is nuclear.
Putting $Q_1 = Q_0 \circ L$, it can be shown that the quadratic form 
$Q = Q_1 - Q_0$ is nuclear and that 
$\Det \left({ Q_1 / Q_0 }\right) $ is equal to 
$\Det (L)^2$. 
The quadratic from $Q_1$ is positive-definite, hence
${\mathop{\rm Ind }\nolimits}
\left({ Q_1 / Q_0 }\right)  = 0$.
>From our main formulas (A.132) and (A.134), we derive
$$\Da_{s, Q_0} (Lx) =
\left\vert{ \Det (L)  }\right\vert
\cdot \Da_{s, Q_0}  x \, .
\eqno ({\mathop{\rm A.140 }\nolimits})$$

\bigskip

Coming back to the general situation, where $Q_\nu = Q_0 + \nu Q$, the reader will find in [5, p. 277-279] examples
of operators $L(\nu)$ such that $Q_\nu = Q_0 \circ L (\nu)$. In such a case, we obtain 
$\Da_\nu x = \left\vert{ \Det L(\nu) }\right\vert
\cdot \Da_0 x $.

\bigskip
Let us summarize the main results obtained in this paragraph:

\medskip
{\it Let $Q_0$ and $Q_1$ be continuous quadratic forms on the Banach space ${\bf X}$.
 We assume that $Q_0$ and $Q_1$
are invertible and that $Q_1-Q_0$ is nuclear. 

\smallskip
(A) Assuming that both $Q_0$ and $Q_1$  be positive-definite, we obtain
$$\int_{{\bf X}}{} \Da_{Q_0} x \cdot e^{- \pi Q_1 (x) } =
\Det \! \left({Q_1/Q_0 }\right)_{}^{-1/2} \, .
\eqno ({\mathop{\rm A.141 }\nolimits})$$

\smallskip
(B) For the the oscillatory integral, we obtain
$$\int_{{\bf X}}{} \Da_{i, Q_0} x \cdot e^{ \pi i Q_1 (x) } =
\left\vert{ \Det \! \left({Q_1/Q_0 }\right) }\right\vert_{}^{-1/2}
i^{-{\mathop{\rm Ind  }\nolimits} \! \left({Q_1/Q_0 }\right)  } 
\eqno ({\mathop{\rm A.142 }\nolimits})$$
where the index ${\mathop{\rm Ind  }\nolimits} \! \left({Q_1/Q_0 }\right)$ counts 
the number of negative eigenvalues
of $Q_1 $ w.r.t. $Q_0$.

\bigskip
Equations } (A.141) {\it and } (A.142) {\it  justify the basic formula } (A.126)
{\it
$$ 
\int_{{\bf X}}{} \Da_{\nu} x \cdot F(x) =
I(\nu) 
\int_{{\bf X}}{} \Da_{0} x \cdot F(x)
$$
with $I(\nu) $ given by } (A.132) {\it or}  (A.134).

\vfill\eject
\vglue 1cm
\centerline{J\bf Appendix B } 
\medskip
\centerline {\titre   Functional determinants of Jacobi operators}
\bigskip
\bigskip

Many functional determinants of Jacobi operators have been computed in works on
semiclassical expansions. They are conveniently expressed in terms of finite-dimensional
determinants of Jacobi matrices. We give here the backbone of the method, and
applications to Jacobi operators along a classical path $x_{{\mathop{\rm cl }\nolimits}}$ 
characterized by $d$ initial
conditions (position or momentum at $t_a$) and $d$ final conditions (position or momentum
at $t_b$). Details, proofs, applications and generalizations of the equations presented
here are scattered in the literature and we give a few selected references at the end of
this appendix. A general presentation on the properties of Jacobi operators will be found
in the Ph.D dissertation of John La Chapelle [University of Texas (Austin) Ph. D. expected in 1995].

\bigskip
\nn
{\bf 1. Jacobi fields and Jacobi matrices.}
\smallskip
A {\it Jacobi field} is a vector field along a classical path obtained by variation through
classical paths. A $d \ts d$ Jacobi matrix is built up from Jacobi fields; each column
consists of the components of a Jacobi field. We give explicit constructions when the action
functional is
$S(x) = 
\int_{ {\bf T} }^{ }
dt \, L(x(t), \dot x (t)) $.
 
\medskip
Let $\{ x_{{\mathop{\rm cl }\nolimits}} (\mu) \} $ be a $2d$-parameter
family of classical paths (critical paths of the functional $S$)  with values in a
manifold  $M^d$ 
$$x_{{\mathop{\rm cl }\nolimits}}(\mu): {\bf T} \ra M^d \quad , \quad {\bf T} = [0, T]\, .$$

\smallskip
\nn
We introduce the notations 
$$
\eqalignno{
 x_{{\mathop{\rm cl }\nolimits}} (t ; \mu )&: = (x_{{\mathop{\rm cl }\nolimits}} (\mu)) (t) 
& ({\mathop{\rm B.1 }\nolimits})\cr
\dot x_{{\mathop{\rm cl }\nolimits}} (t ; \mu )&: = \part x_{{\mathop{\rm cl }\nolimits}} (t ; \mu
) / \part t  & ({\mathop{\rm B.2 }\nolimits}) \cr
x'_{{\mathop{\rm cl }\nolimits}, \a }(t ; \mu )&: =
 \part x_{{\mathop{\rm cl }\nolimits}} (t ; \mu ) / \part  \mu^\a.
& ({\mathop{\rm B.3 }\nolimits}) \cr} $$

\medskip
\nn
We choose $\mu = (\mu^1, \cdots , \mu^{2d})$ to stand for $2d$ initial
conditions that characterize the classical path $x_{{\mathop{\rm cl }\nolimits}} 
(\mu)$, assumed to be unique for the
time being; namely,
$$S' (x_{{\mathop{\rm cl }\nolimits}} (\mu)) = 0 \quad
\eqno ({\mathop{\rm B.4 }\nolimits})$$
has a unique solution for a given set $\mu$ of parameters.

\medskip
By varying successively the $2d$ initial conditions, one obtains $2d$ Jacobi fields
$\{ \part x_{{\mathop{\rm cl }\nolimits}} / \part \mu^\a \} $ for 
$\a \in \{ 1, \cdots, 2 d \} $. It is easy to show that
$$S'' (x_{{\mathop{\rm cl }\nolimits}} (\mu)) \cdot 
{\part x_{{\mathop{\rm cl }\nolimits}} \over \part \mu^\a } = 0 \, .
\eqno ({\mathop{\rm B.5 }\nolimits})$$

\medskip
\nn
We assume for the time being that the quadratic form 
$S'' (x_{{\mathop{\rm cl }\nolimits}} (\mu))\cdot \xi \xi  $ is not degenerate for variations $\xi
$ which respect the boundary conditions specified by $\mu$.  Therefore the $2d$ Jacobi fields are
linearly independent. 

\medskip
There are several convenient basis for the $2d$-dimensional space of Jacobi
fields. Choose as initial condition
$$ \mu = (x_a, p_a)
\eqno ({\mathop{\rm B.6a }\nolimits})$$
where
$$ x_a: = x (t_a)
\eqno ({\mathop{\rm B.6b }\nolimits})$$
$$ p_a: = \part L / \!\!
\left.{ \part \dot x (t) }\right\vert_{t = t_a} \, . 
\eqno ({\mathop{\rm B.6c }\nolimits})$$
The corresponding Jacobi fields are
$$j^{\bu \b } (t) = \part x_{{\mathop{\rm cl }\nolimits}}^\bu 
\left({ t; x_a, p_a }\right)  / \part p_{a ,\b}
\eqno ({\mathop{\rm B.7 }\nolimits})$$
$$k^\bu{}_\b (t) = \part x_{{\mathop{\rm cl }\nolimits}}^\bu 
\left({ t; x_a, p_a }\right) / \part x_a^\b \, .
\eqno ({\mathop{\rm B.8 }\nolimits})$$
Having introduced $p_a$ as one of the initial conditions, the reader suspects that we shall
also use\footnote{$^{31}$}{The classical momentum 
$p_{ {\mathop{\rm  cl }\nolimits}, \bu }  \left({ t; x_a, p_a }\right)  $
is obtained by evaluating $\part L/ \part \dot x (t) $ at the point
$x (t) = x_{ {\mathop{\rm  cl }\nolimits} } 
\left({ t; x_a, p_a }\right) $.
}
$$\ww k_\bu{}^\b (t) = 
\part p_{{\mathop{\rm cl }\nolimits}, \bu } 
\left({ t; x_a, p_a }\right)  / \part p_{a ,\b}
\eqno ({\mathop{\rm B.9 }\nolimits})$$
$$\ell_\bu{}_\b (t) = 
\part p_{{\mathop{\rm cl }\nolimits}, \bu  } 
\left({ t; x_a, p_a }\right) / \part x_a^\b \, .
\eqno ({\mathop{\rm B.10 }\nolimits})$$

\medskip
\nn
Let the {\it Jacobi matrices} $J, K, \ww K, L$ be made of the Jacobi fields $j, k, \ww k, \ell$
respectively
$$J^{\a \b} (t, t_a) = j^{\a \b } (t)
\eqno ({\mathop{\rm B.11 }\nolimits})$$
$$
K^\a {}_\b (t, t_a) = k^\a {}_\b (t)
\eqno ({\mathop{\rm B.12 }\nolimits})$$
$$\ww K_\a {}^\b (t, t_a) = \ww k_\a {}^\b (t)
\eqno ({\mathop{\rm B.13 }\nolimits})$$
$$L_{\a \b } (t, t_a) = \ell_{\a \b } (t) \, .
\eqno ({\mathop{\rm B.14 }\nolimits})$$

\medskip
\nn
The $2d \ts 2d $ matrix constructed from the four $d \ts d$ blocks
$J, K, \ww K, L$ is a solution of the equation defined by the Jacobi operator in phase
space. The properties of the Jacobi matrices are numerous. We note only their properties at
$t = t_a$, which can be read off from their definitions :
$$J (t_a, t_a) = 0 \quad , \quad L (t_a, t_a) = 0 
\eqno ({\mathop{\rm B.15 }\nolimits})$$
$$K(t_a, t_a) = \un \quad , \quad \ww K (t_a, t_a) = \un \, .
\eqno ({\mathop{\rm B.16 }\nolimits})$$

\medskip
\nn
The matrices 
$K$ and $\ww K$ are indeed the transposed of each other
$$\ww K_\a {}^\b (t, t_a) = K^\b {}_\a (t_a, t) .
\eqno ({\mathop{\rm B.17 }\nolimits})$$

\medskip
\nn
Let $M, N, \ww N, P$ be the matrix inverses of $J, K, \ww K, L$ respectively.
Hence we have
$$J^{\a \b } (t, t_a) M_{\b \g } (t_a, t) = \d_\g^\a
\eqno ({\mathop{\rm B.18 }\nolimits})$$
and three similar equations expressed in condensed form as follows
$$KN = \un, \quad \ww K \ww N = \un, \quad LP = \un \, .
\eqno ({\mathop{\rm B.19 }\nolimits})$$
The matrices $M, N, \ww N, P$ are the hessians of the corresponding action functions (Van
Vleck matrices):
$$M_{\b \a } (t_b, t_a) =
{\part^2 \Sa  \over \part x_{{\mathop{\rm cl }\nolimits}}^\b (t_b) 
\part x_{{\mathop{\rm cl }\nolimits}}^\a (t_a)} ,
\eqno ({\mathop{\rm B.20 }\nolimits})$$
and similarly $N$ and $P$ are the hessians of 
$\Sa (x_{{\mathop{\rm cl }\nolimits}} (t_b), p_{{\mathop{\rm cl }\nolimits}} (t_a)) $ and 
$\Sa (p_{{\mathop{\rm cl }\nolimits}} (t_b), p_{{\mathop{\rm cl }\nolimits}} (t_a))$ respectively.

\bigskip
\nn
{\bf 2. Jacobi operators and their Greens's functions.}
\smallskip
The Green's functions of the Jacobi operators in phase space can be expressed in terms of
the Jacobi matrices $J, K, \ww K, L$; they include the Green's functions of the Jacobi
operators in configuration space. For the sake of brevity, but at the cost of elegance, we
consider here only the Green's functions of the Jacobi operators in configuration space.

\medskip
In the previous paragraph, $\ $ Jacobi fields were obtained by variation \break through classical
paths. In this paragraph, we consider a one-parameter family of paths $ x (\l)  $
satisfying $d$ fixed boundary conditions at $t_a$ and $d$ fixed boundary conditions at 
$t_b$, 
referred in brief as ``$a$'' and ``$b$''. The corresponding path space is
 denoted by $\Pa_{a, b} M$, hence
$$x(\l) \in \Pa_{a, b} M \sbs \Pa M \, . 
\eqno ({\mathop{\rm B.21 }\nolimits})$$

\medskip
\nn
We assume that for $\l = 0$, the path 
$x_{{\mathop{\rm cl }\nolimits}} = x(0)$ is the unique critical point of the action functional 
$S$ restricted to $\Pa_{a, b} M$. We set
$$\dot x (\l, t): = \part x (\l , t) / \part t \ , \quad  x' (\l, t): = \part x (\l, t) / \part 
\l \ ,  \quad x' (0, t) = \xi (t) .$$

\medskip
\nn
The second variation of the action functional gives the Jacobi operator 
at $x_{{\mathop{\rm cl }\nolimits}}$:
$$\eqalign{
S'' (x_{{\mathop{\rm cl }\nolimits}}) \cdot \xi \xi
& : = {d^2 \over d \l^2 } 
\left.{  S (x (\l)) }\right\vert_{\l = 0} 
\cr
& = 
\int_{{\bf T}}^{ } dt 
\left({ L_{1 \a, 1 \b} \xi^\a (t) + L_{2 \a, 1 \b } \dot \xi^\a (t) }\right)
\xi^\b (t) \cr
&\ \ \ + \int_{{\bf T} }^{ } dt 
\left({ L_{1 \a, 2 \b} \xi^\a (t) + 
L_{2 \a, 2 \b} \dot \xi^\a (t) }\right)
\dot \xi^\b (t) \cr} \eqno ({\mathop{\rm B.22 }\nolimits}) $$
where
$L_{2 \a, 2 \b}$ is $  \part^2 L / \part \dot x^\a (t) \part \dot x^\b (t)$ evaluated at
 $x_{{\mathop{\rm cl }\nolimits}}$,
etc...
\smallskip

\medskip
If we integrate (B.22) by parts, we obtain the Jacobi operator in its differential garb
(second order differential operator on the space of vector fields $\xi$ 
along $x_{{\mathop{\rm cl }\nolimits}}$
together with boundary terms at $t_a$ and $t_b$). It is often simpler to work with the
quadratic form 
$S'' (x_{{\mathop{\rm cl }\nolimits}}) \cdot \xi \xi$ as written in (B.22). We call functional
Jacobi operator
 the kernel corresponding to the quadratic form 
$S'' (x_{{\mathop{\rm cl }\nolimits}}) \cdot \xi \xi$, namely:
$$ {1 \over 2}
 {\d^2 S'' (x_{{\mathop{\rm cl }\nolimits}}) \cdot \xi \xi \over \d \xi^\a (s) \d \xi^\b (t)}
=: \Ja_{\a \b } (x_{{\mathop{\rm cl }\nolimits}}, s, t ) \, .
\eqno ({\mathop{\rm B.23 }\nolimits})$$

\medskip
\nn
Its inverse will be called its {\it Green's function }
$G^{\bu \bu } \left({ x_{{\mathop{\rm cl }\nolimits}}, t, u }\right) $, namely
$$\int_{\bf T }  dt \, \Ja_{\a \b}
(x_{{\mathop{\rm cl }\nolimits}}, s, t)  G^{\b \g}
(x_{{\mathop{\rm cl }\nolimits}}, t, u) =
\d_\a^\g \d (s - u ) \, .
\eqno ({\mathop{\rm B.24 }\nolimits})$$

\medskip
\nn
Provided the quadratic form 
$S'' (x_{{\mathop{\rm cl }\nolimits}})\cdot  \xi \xi$ is non degenerate, the functional Jacobi
operator has a unique inverse. We say $S'' (x_{{\mathop{\rm cl }\nolimits}})$ is  {\it non
degenerate} if for any $\xi \ne 0$ in $T_{x_{{\mathop{\rm cl }\nolimits}} } \Pa_{a, b} M$
there exists $\eta$ in this vector space with 
 $$S'' (x_{{\mathop{\rm cl }\nolimits}}) \cdot \xi \eta \ne 0 \, .
\eqno ({\mathop{\rm B.25 }\nolimits})$$

\medskip
\nn
This equation says that there are no (nonzero) Jacobi field in the tangent space
 $T_{x_{{\mathop{\rm cl }\nolimits}}}
\Pa_{a, b} M$ to the space of paths with boundary conditions $a$ and $b$.

\bigskip
We list below the Green's functions of Jacobi operators at classical paths with different
boundary conditions. A more abstract formula could encode all cases. 
Explicit formulas may be more
useful for applications.

\medskip
i) 
{\it The classical path is characterized by }
$ p_{{\mathop{\rm cl }\nolimits}} (t_a) = p_a, \ x_{{\mathop{\rm cl }\nolimits}}
 (t_b) = x_b$ 
$$
\eqalign{ G(t, s) &= \t (s - t) K (t, t_a) N (t_a, t_b) J (t_b , s)\cr
&\ \ \ - 
\t (t- s) J (t, t_b) \ww N (t_b, t_a) \ww K (t_a, s) \, .\cr}
\eqno ({\mathop{\rm B.26 }\nolimits})$$

\medskip
ii)
{\it The classical path is characterized by } 
 $x_{{\mathop{\rm cl }\nolimits}} (t_a) = x_a$, $p_{{\mathop{\rm cl }\nolimits}} (t_b) =p_b$
$$\eqalign{ G(t, s) &=
\t (s - t) J (t, t_a) \ww N (t_a, t_b) \ww K (t_b, s) \cr
&\ \ \ - \t (t - s) K (t, t_b) N (t_b, t_a) J (t_a, s) \, .\cr}
\eqno ({\mathop{\rm B.27 }\nolimits})$$

\medskip
iii)
{\it The classical path is characterized by }
$x_{{\mathop{\rm cl }\nolimits}} (t_a) = x_a$, $x_{{\mathop{\rm cl }\nolimits}} (t_b) = x_b$
$$\eqalign{ G(t, s) &= \t (s - t) J (t, t_a) M (t_a, t_b) J (t_b, s) \cr
&\ \ \  - \t (t - s) J
(t, t_b) M (t_b, t_a) J (t_a, s) \, .\cr}
\eqno ({\mathop{\rm B.28 }\nolimits})$$

\medskip
iv)
{\it The classical path is characterized by 
} $p_{{\mathop{\rm cl }\nolimits}} (t_a) = p_a$, 
$p_{{\mathop{\rm cl }\nolimits}} (t_b) = p_b$
$$\eqalign{ G(t, s) &=
\t (s - t) K (t, t_a) P (t_a, t_b) \ww K (t_b, s) \cr
&\ \ \ -
\t (t - s) K (t, t_b) \ww P (t_b, t_a) \ww K (t_a, s) \, .\cr}
\eqno ({\mathop{\rm B.29 }\nolimits})$$

\bigskip
\nn
{\bf 3. Semiclassical expansions.}
\smallskip
We have given examples of WKB approximations in two cases:

i) The classical path of reference is characterized by initial momentum, final position (paragraph
III.2).

\medskip
ii) The classical path is characterized by initial position, final position (paragraph IV.1).

\medskip
Two other cases are often needed:

iii) The classical path is characterized by initial position and final momentum. The transposition
from case i) is straightforward.

\medskip
iv) The classical path is characterized by initial momentum and final momentum.
The critical points of the action functional are degenerate when the classical system is
constrained by a conservation law. We refer the reader to the reference [15] for this case.

\medskip
The computation of semiclassical expansions along the lines of paragraphs III.2 and IV.1, but
with general action functionals $S$, brings in a second variation
$$S'' (x_{{\mathop{\rm cl }\nolimits}}) \cdot \xi \xi = Q_0 (\xi) + Q (\xi)
 \eqno ({\mathop{\rm B.30 }\nolimits})$$
where
$$Q_0 (\xi) =  
\int_{ {\bf T} }^{ } dt \,  L_{2 \a, 2 \b }
\dot \xi^\a (t) \dot \xi^\b (t) \, .
\eqno ({\mathop{\rm B.31 }\nolimits})$$

\medskip
\nn
Provided the Legendre matrix
$$L_{2 \a, 2 \b } 
(x_{{\mathop{\rm cl }\nolimits}} (t), \dot x_{{\mathop{\rm cl }\nolimits}} (t)) =
\part^2 L / \part \dot x_{{\mathop{\rm cl }\nolimits}}^\a \part 
\dot x_{{\mathop{\rm cl }\nolimits}}^\b
\eqno ({\mathop{\rm B.32 }\nolimits})$$
is invertible, the Gaussian integrator defined by the quadratic form(B.30) can be handled
 by the same
techniques as the Gaussian integrator defined by the quadratic form 
$$\int_{{\bf T }} dt \,  h_{\a \b} \dot \xi^\a (t) \dot \xi^\b (t) \, .
\eqno ({\mathop{\rm B.33 }\nolimits})$$

\medskip
The contribution of the second variation 
$S'' (x_{{\mathop{\rm cl }\nolimits}}) \cdot \xi \xi $ to the semiclassical expansion of
$\Psi(t_b, x_b)$ is
$$\int \Da_{ s, Q_0}^{ } \xi\, \,   \exp 
\left({ - { \pi \over s } (Q_0 (\xi) + Q(\xi) }\right)=
 \Det\left({ Q_0 / (Q_0 + Q) }\right)_{}^{1/2}.
\eqno ({\mathop{\rm B.34 }\nolimits})$$
Then one uses the Green's functions (B.26-29) to identify the ratio of Jacobi matrices whose
determinant is equal to $\Det (Q_0 / (Q_0 + Q))$. The results are identical to the results obtained
by discretizing the functional determinants, and reported in [16].

\bigskip
\nn
{\bf 4. References for Appendix B.}
\smallskip
The most complete {\it summary}, to date, of Jacobi fields in phase space, including
degenerate critical points of the action, can be found in Appendix A of 

\smallskip
- C. DeWitt-Morette ``Feynman path integrals'' Acta Physica Austriaca, Suppl. XXVI,
101-170 (1984).

\medskip
For the {\it proofs}  of the results presented in the present  paper, and for properties of Jacobi
matrices used in the present  paper, see [2, 5, 15, 16] and 

- C. DeWitt-Morette and T.-R. Zhang ``Feynman-Kac formula in phase space with application to
coherent state transitions'' Phys. Rev. {\bf D 28}, 2517-2525 (1983).

\medskip
\nn
More on degenerate critical points can be found in 

- C. DeWitt-Morette, B. Nelson and T.-R. Zhang ``Caustic problems in quantum mechanics
with applications to scattering theory'' Phys. Rev. {\bf D 28}, 2526-2546 (1983).

- C. DeWitt-Morette and B.L. Nelson ``Glories and other degenerate points of the action''
 Phys. Rev.
{\bf D 29}, 1663-1668 (1984).

\medskip
A short introduction to the Hamiltonian techniques underlying the above results can be found in 

- P. Cartier ``Some fundamental techniques in the theory of integrable systems''
in {\it Lectures on Integrable Systems} (O. Babelon, P. Cartier and Y. Kosmann-Schwarzbach edit.),
World Scientific, Singapore (1994).

\vfill\eject
\vglue 1cm
\centerline{J\bf Appendix C } 
\medskip
\centerline {\titre   A new class of ordinary differential equations}
\bigskip
\bigskip

The purpose of this Appendix is to extend to the $L^{2,1}$ case the familiar theorems about the existence and
uniqueness of solutions of differential equations, and to describe the parametrization of paths in a curved
space by means of paths in a flat space.

\bigskip
\nn
{\bf 1. Solutions of differential equations: t$\!$he classical case.}

\smallskip
We shall follow the usual strategy, as described in any standard textbook, for instance Bourbaki  [27]. More
precisely, consider a domain $\O$ in the Euclidean space $\RR^n$ and a vector field in $\O$ associating to
every epoch $t$ in the time interval ${\bf T}  = 
\left\lbk{ t_a, t_b }\right\rbk
$ and to every point $x$ in $\O$ a velocity vector $v(t, x) $ in $\RR^n$. 
We assume that $v(t, x)$ is a continuous function of $t$ and $x$.
The two basic remarks are as follows:

\smallskip
a) {\it any trajectory of the vector field $v$ can be prolongated as long as it does not reach the boundary of
$\O$;

}

\smallskip
b)J$\,${\it by the mean value theorem, if the absolute velocity $|v|$ is bounded by a constant $V$ along a given
trajectory leading from $x_a$ at time $t_a$ to $x_b$ at time $t_b$, then the mean velocity
$ { \left\vert{x_b - x_a  }\right\vert \over t_b - t_a } $ is bounded by $V$. 

}

\smallskip
\nn
>From these remarks follows the following existence theorem (Peano):

\smallskip
{\it Assume that $\O$ is a closed ball centered at $x_a$ of radius $L$, and that 
$$\left\vert{ v (t, x) }\right\vert
< L / T
\eqno ({\mathop{\rm C.1 }\nolimits})$$
holds uniformly for $t$ in  ${\bf T}$ and $x$ in $\O$, where $T = t_b - t_a$ is the length of ${\bf T}$.
Then there exists a solution $x : {\bf T} \ra \O$ of the differential equation $\dot x = v (t, x)$ with the
initial condition $x(t_a) = x_a$. }

\medskip
For the proof, we construct first an approximate solution by the Euler method. We select epochs 
$t_1, \cdots , t_{N-1}$ such that $t_a < t_1 < \cdots < t_{N-1} < t_b $ and set
$t_0 = t_a$ , $t_N = t_b$. Then we define inductively points $x_0, x_1, \cdots, x_N$ in $\O$ by $x_0 = x_a$ and
$$x_i = x_{i-1} + \left({ t_i - t_{i-1} }\right) v \! \left({ t_{i-1}, x_{i-1} }\right)
\eqno ({\mathop{\rm C.2 }\nolimits})$$
for $1 \leq i \leq N$; the estimate (C.1) guarantees that the points $x_0, x_1,\cdots, x_N$ are in $\O$.
We then interpolate linearly in each subinterval 
$\left\lbk{ t_{i-1}, t_i }\right\rbk$, and generate a function 
$x_\Ta : {\bf T} \ra \O$ depending on the subdivision 
$\Ta = \left({ t_0 < t_1 < \cdots < t_N }\right)$ of ${\bf T}$.
By the mean value theorem and the estimate (C.1) we obtain
$$\left\vert{ x_\Ta (t) - x_\Ta (t') }\right\vert
\leq V \left\vert{ t - t' }\right\vert
\eqno ({\mathop{\rm C.3 }\nolimits})$$
for $t, t'$ in ${\bf T}$, where $V = L/ T$.  We then use Ascoli's theorem: it asserts 
the existence of a sequence of subdivisions $\Ta (n) $ of ${\bf T}$, whose mesh $\D(n) $ tends to $0$, such that
$x_{\Ta (n)} (t)$ tends to a limit $x(t)$ uniformly for $t$ in ${\bf T}$. Then, one shows that 
$x_{\Ta (n)} $ satisfies an approximate integral equation
$$
\left\vert{ x_{\Ta (n)}(t) - x_a - \int_{ t_a }^{  t}
ds \ v (s, x(s))
 }\right\vert \leq \ve_n \, ,
\eqno ({\mathop{\rm C.4 }\nolimits})$$
where $\build { \lim }_{ n = \Inf  }^{ } \ve_n = 0$.
By uniform convergence, the limit function $x(\cdot) $ satisfies the integral equation
$$x(t) = x_a + \int_{ t_a}^{ t}
ds \  v(s, x(s)) \, ,
\eqno ({\mathop{\rm C.5 }\nolimits})$$
fully equivalent to the differential equation
$$\left\{\matrix{
\dot x (t) = v(t, x (t))
 \hfill\cr
\noalign{\medskip}
x(t_a) = x_a \, .   
 \hfill \cr}
\right.
\eqno ({\mathop{\rm C.6 }\nolimits})$$

\medskip
Both the {\it uniqueness } of a solution to the previous equation and the {\it continuous dependence of the
unique solution on the initial position} $x_a$ are established by an analysis of the {\it stability} of
trajectories. Suppose given another trajectory $y : {\bf T} \ra \O$ with initial position $y_a = y (t_a)$; such
a trajectory exists if $y_a$ is sufficiently close to $x_a$. Hence we have the integral relation
$$y(t) = y_a + \int_{t_a }^{ t} ds \, v(s, y (s)) \, .
\eqno ({\mathop{\rm C.7 }\nolimits})$$
Furthermore, assume that the velocity field satisfies a {\it Lipschitz condition}:
$$\left\vert{ v(t, x) - v(t, x') }\right\vert
\leq k \left\vert{ x - x' }\right\vert
\eqno ({\mathop{\rm C.8 }\nolimits})$$
for $x, x'$ in $\O$ and $t$ in ${\bf T}$, with a fixed constant $k > 0$. Denote by $\d(t) $ the distance
between $x(t) $ and $y(t)$. From (C.5), (C.7) and (C.8) one derives
$$\left\vert{ \d(t) - \d(t_a)  }\right\vert
\leq k \int_{t_a }^{ t} ds \, \d (s) \, .
\eqno ({\mathop{\rm C.9 }\nolimits})$$

\nn
By using an iteration procedure reminiscent of Picard's method, one derives  the following estimate:
$$\d (t) \leq \d (t_a) e^{k (t - t_a)}
\eqno ({\mathop{\rm C.10 }\nolimits})$$
(``{\it  Gronwall's lemma}''). 
If $x_a = y_a$, then $\d (t_a) = 0$, hence $\d(t) = 0$ for all $t$ and $x(t) = y(t)$: uniqueness. Furthermore
if $\d(t_a) = \left\vert{ y_a - x_a }\right\vert$ tends to $0$, then $y(t)$ tends
uniformly to $x(t)$ for  $t$ in the finite interval ${\bf T}$.

\medskip
The results obtained so far are local.
To get a global existence theorem, assume that $N$ is a compact manifold, of dimension $n$, and that $X$ is a
time-dependent vector field on $N$, namely $X(t, x) $ belongs to $T_x N$ for $t $ in ${\bf T}$, and $x$ in $N$.
Assume furthermore that $X(t, x) $ is continuous in $t, x$, with continuous first derivatives in $x$.
Choose a point $x_a$ in $N$ and a coordinate system around $x_a$. Then $X$ is expressed in this coordinate
system by a function $v : {\bf T} \ts \O \ra \RR^n$, where $\O$ is a  closed ball centered at $x_a$, with radius
$L$. Since a continuous function is bounded on any compact set, there exists a constant $\tau > 0$ such that
$| v(t, x) | \leq L / 2 \tau$ uniformly for $t$ in ${\bf T}$ and $x$ in $\O$.
>From the local existence theorem, it follows that for any point $x_0$ in the open ball centered at $x_a$ with
radius $L/2$, and any $t_0$ in ${\bf T}$, there exists a unique trajectory 
$x : \left\lbk{ t_0 , t_0 + \tau }\right\rbk \ra N$ such that 
$x(t_0) = x_0$. By compacity of $N$, we can cover $N$ by a finite number of such balls.
Defining $\s $ as the minimum of the life-times $\tau$, we conclude that the existence and uniqueness of a
trajectory $x : \left\lbk{ t_0 , t_0 + \s }\right\rbk \ra N$ with $x(t_0) = x_0$ hold for every point $x_0$ in
$N$.

\medskip
Consider now the time interval ${\bf T} = \left\lbk{ t_a , t_b }\right\rbk $ and subdivide it into subintervals
$\left\lbk{ t_i , t_{i+1} }\right\rbk $ of lengths $\leq \s$ (for $i$ in 
$\{ 0, 1, \cdots, N \}$)
where $t_0 = t_a$, $t_{N+1} = t_b$. So given any point $x_a$ in $N$, there exists a trajectory $x_0(t)$
 (for $t_0
\leq t \leq t_1$) such that $x_0(t_0) = x_a$.
Set $x_1 = x_0(t_1)$ and consider the trajectory $x_1(t) $ (for $t_1 \leq t \leq t_2$) with initial point $x_1$.
By repeating this procedure, we construct step by step the required trajectory $x : {\bf T}\ra N$ such that 
$x(t_a) = x_a$.

\medskip
Proceeding backwards rather than forwards, we conclude that given any point $x_b$ in $N$, there exists a unique
solution $x(t) = x(t ; x_b)$ of the differential equation 
$\dot x (t) = X (t, x(t)) $ with final position $x(t_b) = x_b$.
For fixed $t$, we define the transformation $\Si (t) $ in $N$ taking $x_b$ to $x(t ; x_b)$, hence
$$x (t; x_b) = x_b \cdot \Si (t)\, .
\eqno ({\mathop{\rm C.11 }\nolimits})$$

\bigskip
\nn
{\bf 2. Solutions of differential equations: the $L^{2,1}$ case.}

\smallskip
We sketch the necessary modifications. A more detailed account will appear elsewhere.

\smallskip
Consider an $L^{2,1}$ path $x : {\bf T} \ra \O$; its action is defined by
$$A = \int_{ {\bf T} }^{ }
{|dx|^2\over dt} \, .
\eqno ({\mathop{\rm C.12 }\nolimits})$$
The mean value theorem is replaced by the following estimate
$$
{ \left\vert{ x_b - x_a }\right\vert^2 \over t_b - t_a }
\leq A \, , 
\eqno ({\mathop{\rm C.13 }\nolimits})$$
which follows easily from Cauchy-Schwarz inequality.

\medskip
We make the following assumptions about the velocity field $v(t, x)$.

\smallskip
(A) {\it The function $v(t, x)$ is jointly-measurable in $(t, x)$, and continuous in $x$ for a given $t$, and
furthermore there exists a numerical $L^2$ function \break $V : {\bf T} \ra [0, + \Inf[$ such that
$$
\left\vert{ v(t, x) }\right\vert \leq V(t)
\eqno ({\mathop{\rm C.14 }\nolimits})$$
holds for $x$ in $\O$ and $t$ in ${\bf T}$.
}

\smallskip
(B) {\it
Denoting by $A$ the action $\int_{ {\bf T}}^{ } dt \, V(t)^2$, the radius $L$ of the ball $\O$ centered at
$x_a$, and the length $T$ of the time interval ${\bf T}$ obey the estimate
$$A < L^2 / T \, . \eqno ({\mathop{\rm C.15 }\nolimits})$$
}
\smallskip
Under these conditions, one proves an existence theorem for the differential equation (C.6) with an $L^{2,1}$
solution $x : {\bf T} \ra \O$. As a first step, one replaces Euler's approximation (C.2) by the following
approximate solution
$$x_\Ta (t) = x_{\Ta} (t_{i-1}) +
\int_{t_i }^{t }ds \, v \! 
\left({ s, x_\Ta (t_{i-1}) }\right)
\eqno ({\mathop{\rm C.16 }\nolimits})$$
for $t$ in the subinterval $\left\lbk{ t_{i-1}, t_i }\right\rbk$.
The inequalities (C.13) to (C.15) guarantee that this trajectory $x_\Ta (t)$ remains in the closed ball $\O$.
Furthermore, the velocity $\dot x_\Ta (t)$ of this trajectory satisfies the estimate
$\left\vert{ \dot x_\Ta (t) }\right\vert
\leq V(t)$ for $t$ in ${\bf T}$, and by (C.13), we obtain the estimate
$$\left\vert{ x_\Ta (t) - x_\Ta (t') }\right\vert^2
\leq A \left\vert{ t - t' }\right\vert \, .
\eqno ({\mathop{\rm C.17 }\nolimits})$$
We can then invoke Ascoli's theorem, and find a sequence of approximate solutions $x_{\Ta (n)}$ converging
uniformly on ${\bf T}$ towards an $L^{2,1}$ function $x$. Each approximation $x_{\Ta (n)}$ satisfies an
approximate integral equation, and in the limit the integral equation (C.5) is obtained using Lebesgue's
dominated convergence.

\medskip
For the uniqueness, we need a Lipschitz condition of the type
$$
 \left\vert{ v(t, x) - v(t, x') }\right\vert
\leq k(t) \left\vert{ x - x' }\right\vert
\eqno ({\mathop{\rm C.18 }\nolimits})$$
where the integral $\int_{ {\bf T} }^{ } dt \, k(t)^2$ is finite. 
We use again a variant of Gronwall's lemma.

\medskip
The global results are obtained for a compact manifold $N$ by the reasonings used at the end of paragraph C.1.
Only minor modifications are needed.

\bigskip
\nn
{\bf 3. Parametrization of paths.}
\smallskip
We consider a manifold $N$ and $d$ vector fields 
$X_{(1)}, \cdots , X_{(d)}$ on $N$.
We assume that they are of class $C^1$ and linearly independent at each point of $N$.
For $x$ in $N$ denote by $H_x$ the vector subspace of $T_x N$ generated by 
$X_{(1)}(x), \cdots , X_{(d)}(x)$. The collection of these vector spaces is a subbundle\footnote{$^{32}$}{The
letter $H$ stands for ``horizontal''.} $H$ of the tangent bundle $TN$ to $N$. We consider also a
real  symmetric invertible
matrix
$\left({ h_{\a \b} }\right)$ of size $d \ts d$, and define a field of quadratic forms $h_x$ on $H_x$ by
$$h_x 
\left({ X_{(\a)} (x), X_{(\b)} (x) }\right)
= h_{\a \b} \, .
\eqno ({\mathop{\rm C.19 }\nolimits})$$

\medskip
Fix a point $x_b$ in $N$.
We denote by $\Pa_{x_b}^H N$ the set of $L^{2,1}$ paths 
$x : {\bf T} \ra N$ which satisfy the following conditions:

\smallskip
(A) {\it The endpoint $x(t_b) $ is equal to $x_b$.}

\smallskip
(B) {\it For each epoch $t$ in ${\bf T}$, the velocity vector $\dot x (t)$ lies in the subspace 
$H_{x(t)}$ of $T_{x(t)} N$.

}

\medskip
\nn
We define the {\it action}
 of such a path by 
$$A(x) = 
\int_{{\bf T} }^{ } dt \, h_{x(t)}
\! \left({ \dot x(t), \dot x(t) }\right)\, .
\eqno ({\mathop{\rm C.20 }\nolimits})$$
Since the vectors $X_{(\a)} (x(t))$ (for $1 \leq \a \leq d$) form a basis of 
$H_{x(t)}$, and from the hypothesis that the path $x$ is of class $L^{2,1}$, we infer that there exist functions
$\dot z_\a $ in $L^2 ({\bf T})$ such that
$$\dot x(t) = X_{(\a)} (x(t)) \dot z^\a (t) \, . 
\eqno ({\mathop{\rm C.21 }\nolimits})$$
Furthermore, the function $\dot z^\a$ is the derivative of a function $z^\a$ in $L^{2,1}({\bf T})$ normalized
by 
$z^\a (t_b) = 0$. The vector function $z = \left({ z^1, \cdots z^d }\right)$ is an element of the space denoted
by ${\bf Z}_b$ in paragraph A.3.8. {\it This construction associates to a path $x$ in $\Pa_{x_b}^H N$ a path $z$
in ${\bf Z}_b$ with conservation of action:} 
$$A(x) = Q_0 (z) \, , 
\eqno ({\mathop{\rm C.22 }\nolimits})$$
where as usual $Q_0 (z)$ is equal to 
$\int_{{\bf T} }^{ } dt \, h_{\a \b}
\dot z^\a (t) \dot z^\b (t) $. 

\medskip
Assume now that $N$ is compact. By using the theory of $L^{2, 1}$ differential equations sketched in paragraph
C.2, it can be shown that we can invert the transformation $x \mps z$. Given any $z$ in ${\bf Z}_b$, the
differential equation (C.21) has a unique solution $x$ in $\Pa_{x_b}^H N $, and we obtain {\it a
parametrization
 $$P : {\bf Z}_b \ra \Pa_{x_b}^H N$$
of a space of paths in a curved space $N$ by a space of paths in a flat space $\RR^d$.}
For $z$ in ${\bf Z}_b$, we denote by $x(t, z)$ the solution of the differential equation (C.21) with endpoint 
$x(t_b, z) = x_b$. We can also introduce a global transformation $\Si (t, z) : N \ra N$ taking $x_b$ into $x(t,
z)$.
If necessary, we include 
$t_b$ in the notation and denote this transformation by
$\Si ({\bf T} ; z)$ or
 $\Si \! \left({  t_b, t_a ; z}\right)$ in the case $t = t_a$. The chain rule
$$\Si \! \left({  t_b, t_a ; z_{ba}J}\right)
=
\Si \! \left({  t_b, t_c ; z_{bc}}\right)
\cdot
\Si \! \left({  t_c, t_a ; z_{ca}}\right)
\eqno ({\mathop{\rm C.23 }\nolimits})$$
is a consequence of the uniqueness of the solution of the differential equation (C.21).
Here $z_{ba}$ is $z$, $t_c$ is an intermediate epoch, 
and the paths $z_{bc} : \left\lbk{  t_c, t_b }\right\rbk \ra \RR^d$
and $z_{ca} : \left\lbk{  t_a, t_c }\right\rbk \ra \RR^d$
are given by 
$$z_{bc}(t) = z(t) \quad , \quad
z_{ca} (t) = z(t) - z(t_c) \, . 
\eqno ({\mathop{\rm C.24 }\nolimits})$$

\medskip
\nn
{\it Remark.}
The more general differential equation
$$\dot x (t) = X_{(\a)} 
(x(t)) \dot z^\a (t)  + Y (x(t)) 
\eqno ({\mathop{\rm C.25 }\nolimits})$$
can be handled in a similar way.
In this case, we replace $\Pa_{x_b}^H N$ by the space
$
\Pa_{x_b}^{H, Y} N$ of $L^{2,1}$ paths 
$x : {\bf T} \ra N$ such that
$x (t_b) = x_b $ and $\dot x(t) - Y(t) $ belong to 
$H_{x (t)} $ for every $t$ in ${\bf T}$.

\medskip
These constructions are related to the {\it Cartan development map}. Take for $M$ a compact Riemannian manifold,
and let $N$ be the corresponding bundle of orthonormal frames; it is a compact manifold.
Fix $x_b$ in $M$, and a frame $\r_b = 
\left({ x_b, u_b }\right)$ at $x_b$.
Then, by the Riemannian connection, there is defined a  ``horizontal'' subspace 
$T_{\r_b}^H N$ of the tangent space $T_{\r_b} N$,
and the projection $\Pi : N \ra M$ induces an identification of $T_{\r_b}^H N$
with $T_{x_b} M$.
Since $u_b$ defines an orthonormal basis of $T_{x_b} M$, we obtain a basis
$X_{(1)} (\r_b) , \cdots , X_{(d)} (\r_b)$
of $T_{\r_b}^H N$. This construction is valid for every point in $N$, hence we define vector fields
$X_{(1)} , \cdots , X_{(d)}$ on
 $N$.
By using the previous construction, we get a parametrization
$$P : {\bf Z}_b \ra \Pa_{\r_b}^H N \, . 
$$
But the paths in $\Pa_{\r_b}^H N$ are nothing else than the horizontal liftings of the 
$L^{2,1}$ paths in $M$.
More precisely denote by $\Pa_{x_b} M$ the set of $L^{2,1}$
paths $x : {\bf T} \ra M$, such that $x(t_b) = x_b$. Then, the projection 
$\Pi : N \ra M$ induces a mapping 
$\ww x \mps \Pi \circ \ww x$ of $\Pa_{\r_b}^H N$ into $\Pa_{x_b} M$, and this map is a bijection.
To conclude, we get a diagram
$${\bf Z}_b \ra \Pa_{\r_b}^H N \ra \Pa_{x_b} M$$
and by composition a parametrization of 
$\Pa_{x_b} M$ by ${\bf Z}_b$.
This is the Cartan development map for $L^{2,1}$ paths.
The standard theory works for $C^1$ paths.

\vfill \eject
\centerline {{\bf REFERENCES} }

\bigskip
\bigskip

\item{[1]}
C. Morette 
``On the definition and approximation of Feynman's path integral''
Phys. Rev. {\bf 81}, 848-852 (1951).
 
\medskip

\item{[2]}
C. DeWitt-Morette, B. Nelson, and T.R. Zhang
``Caustic problems in quantum mechanics with applications to scattering theory''
Phys. Rev. D {\bf 28}, 2526-2546 (1983).
 
\medskip
\item{[3]}
M.G.G. Laidlaw and C. Morette-DeWitt 
``Feynman functional integrals for systems of indistinguishable particles'' 
Phys. Rev. D {\bf  3}, 1375-1378 (1971).

\medskip
\item{[4]}
R.H. Cameron and W.T. Martin
``Transformation of Wiener integrals under a general class of linear transformations''
Trans. Amer. Math. Soc. {\bf 58}, 184-219 (1945).

\medskip
\item{[5]}
C. DeWitt-Morette, A. Maheshwari, and B. Nelson 
``Path Integration in Non-Relativistic Quantum Mechanics''
Phys. Rep. {\bf 50}, 266-372 (1979).
This article includes a summary of earlier articles.

\medskip
\item{[6]}
B. Nelson and B. Sheeks 
``Fredholm determinants associated with Wiener integrals''
J. Math. Phys. {\bf 22},
2132-2136 (1981).

\medskip
\item{[7]}
K.D. Elworthy
``{\it Stochastic Differential Equations on Manifolds}''
Cambridge University Press, Cambridge U.K. (1982).

\medskip
\item{[8]}
C. Morette-DeWitt 
``Feynman's Path Integral, definition without limiting procedure''
Commun. Math. Phys. {\bf 28}, 47-67 (1972), and ``Feynman Path Integrals, I. Linear and
affine techniques, II. The Feynman-Green function'' 
Commun. Math. Phys. {\bf 37}, 63-81 (1974).

\medskip
\item{[9]}
S.A. Albeverio and R.J. H\o egh-Krohn
{\it ``Mathematical Theory of Feynman Path Integrals''J}
Springer Verlag Lecture Notes in Mathematics 523 (1976).

\medskip
\item{[10]}
P. Kr\'ee 
``Introduction aux th\'eories des distributions en dimension infinie''
Bull. Soc. Math. France {\bf 46}, 143-162 (1976), and references therein, in particular, {\it
Seminar P. Lelong} Springer Verlag Lecture Notes in Mathematics {\bf 410} and {\bf 474}
(1972-1974).

\medskip
\item{[11]}
P. Cartier and C. DeWitt-Morette
``Int\'egration fonctionnelle; \'el\'ements d'axiomatique''
 C.R. Acad. Sci. Paris, t.316, S\'erie II,  733-738 (1993).

\medskip
\item{[12]}
A. Young and C. DeWitt-Morette ``Time substitutions in stochastic Processes as a Tool in Path
Integration'' Ann. of Phys. {\bf 69}, 140-166 (1986).
 
\medskip
\item{[13]}
Y. Choquet-Bruhat and C. DeWitt-Morette
``Supplement to {\it Analysis, Manifolds and Physics}''
Armadillo preprint, Center for Relativity, University of Texas, 
Austin TX 78712.

\medskip
\item{[14]}
C. DeWitt-Morette, K.D. Elworthy, B.L. Nelson, and G.S. Sammelman
``A stochastic scheme for constructing solutions of the Schr\"odinger equation''
Ann. Inst. H. Poincar\'e A {\bf 32}, 327-341 (1980).

\medskip
\item{[15]}
C. DeWitt-Morette and T.-R. Zhang 
``Path integrals and conservation laws'' 
Phys. Rev. D {\bf  28}, 2503-2516 (1983).

\medskip
\item{[16]}
C. DeWitt-Morette
``The semiclassical expansion''
Ann. of Phys. {\bf 97}, 367-399 (1976).
(Correct a misprint p. 385, l. 5, the reference is [5].)

\medskip
\item{[17]}
S.F. Edwards and Y.V. Gulyaev 
``Path integrals in polar coordinates''
Proc. Roy. Soc. London A{\bf 279}, 224 (1964).

\medskip
\item{[18]}
C. DeWitt-Morette 
``Quantum mechanics in curved spacetimes; stochastic processes on frame bundles''
pp. 49-87 in {\it Quantum Mechanics in Curved Space-Time} 
Eds J. Audretsch and V. de Sabbata, Plenum Press, New York (1990).

\medskip
\item{[19]}
W. Greub and H.-R. Petry 
``Minimal coupling and complex line bundles''
J. Math. Phys. {\bf 16}, 1347-1351 (1975).

\medskip
\item{[20]}
Y. Choquet-Bruhat and C. DeWitt-Morette {\it ``Analysis, Manifolds and Physics Part I:
Basics, Part II: 92 applications''}, North Holland, Amsterdam (1989).

\medskip
\item{[21]}
F.A. Berezin ``Quantization'' Math. USSR, Izvestija {\bf 8}, 1109-1165 (1974).

\medskip
\item{[22]}
D. Bar-Moshe and M.S. Marinov
``Berezin quantization and unitary representations of Lie groups'',
preprint, to appear in Berezin Memorial volume (1994).

\medskip
\item{[23]}
K.R. Parthasarathy
{\it `` Probability measures on metric spaces''}
Academic Press, New York (1967).

\medskip
\item{[24]}
N. Bourbaki
{\it ``Int\'egration, chapitre 9''}
Masson, Paris (1982).

\medskip
\item{[25]}
N. Bourbaki
{\it ``Espaces vectoriels topologiques''}
Masson, Paris (1981).

\medskip
\item{[26]}
P. Cartier ``A course on determinants''
in {\it Conformal invariance and string theory} Eds P. Dita and V. Georgescu, Academic Press, New York (1989).

\medskip
\item{[27]}
N. Bourbaki {\it ``Fonctions d'une variable r\'eelle''}
Masson, Paris (1982).

\vfill\eject
\vglue 	1cm
\centerline {Figure captions}
\vglue 	2cm
{\it Figure} A.1: ``Linear change of variables''

\bye

------------------------------------------------------------------------
Mijan Huq                          mijan@hoffmann.ph.utexas.edu
Center for Relativity,             Phone 512 471 4700
Department of Physics,
University of Texas at Austin,
Austin, TX 78712-1081

URL: http://godel.ph.utexas.edu/Members/mijan/welcome.html
------------------------------------------------------------------------
"A conservative is a statesman who is enamored of existing
 evils, as distinguished from the liberal, who wishes to 
 replace them with others."
          Ambrose Bierce (1842-1914) in Devil's Dictionary

From mijan@hoffmann.ph.utexas.edu Thu Feb  8 13:41:24 1996
Return-Path: <mijan@hoffmann.ph.utexas.edu>
Received: from hoffmann.ph.utexas.edu by helmholtz.ph.utexas.edu (4.1/SMI-4.1)
	id AA24488; Thu, 8 Feb 96 13:41:14 CST
From: mijan@hoffmann.ph.utexas.edu (Mijan Huq)
Posted-Date: Thu, 8 Feb 1996 13:41:12 -0600 (CST)
Message-Id: <9602081941.AA08642@hoffmann.ph.utexas.edu>
Received: by hoffmann.ph.utexas.edu (931110.SGI/5.51)
	id AA08642; Thu, 8 Feb 96 13:41:13 -0600
Subject: Re: one more time
To: debbie@helmholtz.ph.utexas.edu (Debbie Hajji)
Date: Thu, 8 Feb 1996 13:41:12 -0600 (CST)
In-Reply-To: <9602081939.AA24485@helmholtz.ph.utexas.edu> from "Debbie Hajji" at Feb 8, 96 01:39:54 pm
X-Mailer: ELM [version 2.4 PL21]
Content-Type: text
Content-Length: 267782    
Status: R

\\
Title: A new perspective on Functional Integration
Authors: Pierre Cartier and C\'ecile DeWitt-Morette
Comments: 102 pages
Journal-ref: Journal of Mathematical Physics, vol 36, 2137-2340 (1995)
\\
The core of this article is a general theorem with a large number of specializations. Given a
manifold $N$ and a finite number of one-parameter groups of point transformations on $N$ with
generators $Y, X_{(1)}, \cdots, X_{(d)} $, we obtain, via
functional integration over spaces of pointed paths on $N$ (paths with one fixed point), a
one-parameter group of functional operators acting on tensor or spinor fields on $N$. The
generator of this group is a quadratic form in the Lie derivatives $\La_{X_{(\a)}}$ in the 
$X_{(\a)}$-direction plus a term linear in $\La_Y$. 

\smallskip
The basic functional integral is over $L^{2,1}$ paths $x: {\bf T} \ra N$ (continuous paths
with square integrable first derivative). Although the integrator is invariant under time
translation, the integral is powerful enough to be used for systems which are not time
translation invariant. We give seven non trivial applications of the basic formula, and we
compute its semiclassical expansion.

\smallskip
The methods of proof are rigorous and combine Albeverio H\o egh-Krohn oscillatory integrals
with Elworthy's parametrization of paths in a curved space. Unlike other approaches we solve
 Schr\"odinger type equations directly, rather than solving first diffusion equations
and then using analytic continuation.

\\
\def\today{\ifcase\month\or January\or February\or March\or
April\or May\or June\or July\or August\or September\or
October\or November\or December\fi \space\number\day,
\number\year}

\def\a{\alpha}
\def\b{\beta}
\def\d{\delta}

\def\g{\gamma}

\def\l{\lambda}
\def\om{\omega}
\def\r{\rho}
\def\s{\sigma}
\def\t{\theta}

\def\ve{\varepsilon}
\def\vp{\varphi}

\def\z{\zeta}

\def\D{\Delta}
\def\G{\Gamma}
\def\L{\Lambda}
\def\O{\Omega}
\def\Si{\Sigma}
\def\T{\Theta}

\let\nn=\noindent

\font\tenbb=msym10
\font\sevenbb=msym8
\font\fivebb=msym5
\newfam\bbfam
\textfont\bbfam=\tenbb \scriptfont\bbfam=\sevenbb
\scriptscriptfont\bbfam=\fivebb
\def\bb{\fam\bbfam}

\def\CC{{\bb C}}

\def\NN{{\bb N}}

\def\RR{{\bb R}}

\def\ZZ{{\bb Z}}

\font\titre=cmbx12
\font\tenfm=eufm10

\def\part{\partial}
\def\Inf{\infty}
\def\bl{\backslash}

\def\bu{\bullet}

\def\bgt{\nabla}
\def\mps{\mapsto}
\def\ts{\times}

\def\sbs{\subset}

\def\ra{\rightarrow}

\def\lra{\leftrightarrow}

\def\lbc{\lbrace}
\def\rbc{\rbrace}
\def\lbk{\lbrack}
\def\rbk{\rbrack}

\def\oo{\overline}
\def\ww{\widetilde}
\def\hh{\widehat}

\def\ds{\displaystyle }

\def\Da{{\cal D}}

\def\Fa{{\cal F}}

\def\Ja{{\cal J}}
\def\Ka{{\cal K}}
\def\La{{\cal L}}
\def\Ma{{\cal M}}

\def\Pa{{\cal P}}

\def\Sa{{\cal S}}
\def\Ta{{\cal T}}
\def\Ua{{\cal U}}

\def\and{\mathop{\rm and}\nolimits}

\def\by{\mathop{\rm by}\nolimits}

\def\cos{\mathop{\rm cos}\nolimits}

\def\det{\mathop{\rm det}\nolimits}

\def\diag{\mathop{\rm diag}\nolimits}

\def\Det{\mathop{\rm Det}\nolimits}

\def\Dev{\mathop{\rm Dev}\nolimits}

\def\Eucl{\mathop{\rm Eucl}\nolimits}
\def\Eucl.{\mathop{\rm Eucl.}\nolimits}

\def\exp{\mathop{\rm exp}\nolimits}

\def\or{\mathop{\rm or}\nolimits}

\def\sin{\mathop{\rm sin}\nolimits}

\def\supp{\mathop{\rm supp}\nolimits}
\def\supp.{\mathop{\rm supp.}\nolimits}

\def\space{\mathop{\rm space}\nolimits}

\def\Tr{\mathop{\rm Tr}\nolimits}
\def\tr{\mathop{\rm tr}\nolimits}

\def\Var{\mathop{\rm Var}\nolimits}

\catcode`\@=11
\def\displaylinesno #1{\displ@y\halign{
\hbox to\displaywidth{$\@lign\hfil\displaystyle##\hfil$}&
\llap{$##$}\crcr#1\crcr}}

\def\ldisplaylinesno #1{\displ@y\halign{
\hbox to\displaywidth{$\@lign\hfil\displaystyle##\hfil$}&
\kern-\displaywidth\rlap{$##$}
\tabskip\displaywidth\crcr#1\crcr}}
\catcode`\@=12

\def\buildrel#1\over#2{\mathrel{
\mathop{\kern 0pt#2}\limits^{#1}}}

\def\build#1_#2^#3{\mathrel{
\mathop{\kern 0pt#1}\limits_{#2}^{#3}}}

\def\hfl#1#2{\smash{\mathop{\hbox to 6mm{\rightarrowfill}}
\limits^{\scriptstyle#1}_{\scriptstyle#2}}}

\def\hfll#1#2{\smash{\mathop{\hbox to 6mm{\leftarrowfill}}
\limits^{\scriptstyle#1}_{\scriptstyle#2}}}

\def\up#1{\raise 1ex\hbox{\sevenrm#1}}

\def\per{|\!\raise -4pt\hbox{$-$}}
 \def\cqfd{\unskip\kern
6pt\penalty 500 \raise
-2pt\hbox{\vrule\vbox to10pt{\hrule
width 4pt \vfill\hrule}\vrule}\par}

\def\trait{\hbox to 12mm{\hrulefill}}
\def\2{{\mathop{\rm 
I }\nolimits}\!{\mathop{\rm  
I}\nolimits}}

\def\1{{\mathop{\rm I }\nolimits}}

\def\og{\leavevmode\raise.3ex\hbox{$
\scriptscriptstyle\langle\!\langle$}}
\def\fg{\leavevmode\raise.3ex\hbox{$
\scriptscriptstyle \,\rangle\!\rangle$}}

\def\picture #1 by #2 (#3)
{\vcenter{\vskip #2
\special{picture #3}
\hrule width #1 height 0pt depth 0pt
\vfil}}

\def\scaledpicture #1 by #2 (#3 scaled #4){{
\dimen0=#1 \dimen1=#2
\divide\dimen0 by 1000 \multiply\dimen0 by #4
\divide\dimen1 by 1000 \multiply\dimen1 by #4
\picture \dimen0 by \dimen1 (#3 scaled #4)}}

\def\un{{\rm 1\mkern-4mu l}}
\def\[{{[\mkern-3mu [}}
\def\]{{]\mkern-3mu ]}}

\def\tvi{\vrule height 12pt depth 5pt width 0pt}

\def\tv{\tvi\vrule}

\def\cc#1{\hfill\kern .7em#1\kern .7em\hfill}

\def\TeX{T\kern-.1667em\lower.5ex\hbox{E}\kern-.125em X}

\def\ins{{ \raise -2mm\hbox{$<$}  
\atop \raise 2mm\hbox{$\sim$}J}
}

\def\sus{{ \raise -2mm\hbox{$>$}  
\atop \raise 2mm\hbox{$\sim$}J}
}

\def\lta{\hbox{\raise.5ex\hbox{$<$}
\kern-1.1em\lower.5ex\hbox{$\sim$}}}
\def\gta{\hbox{\raise.5ex\hbox{$>$}
\kern-1.1em\lower.5ex\hbox{$\sim$}}}

\magnification=1200
\overfullrule=0mm
\baselineskip=12pt
 
\hsize=118mm 
\hoffset=7mm
\vsize=185mm
\voffset=10mm

\vglue 15mm

\centerline{\titre A new perspective on}
\smallskip
\centerline{\titre  Functional Integration}

\vglue 5mm
{\baselineskip=11pt
\centerline{\bf Pierre Cartier}

\centerline{ Ecole Normale Sup\'erieure }

\centerline{45 rue d'Ulm}

\centerline{F-75230 Paris C\'edex 05}

\centerline{France} 
\smallskip
\centerline{\it and}
\smallskip

\centerline{\bf C\'ecile DeWitt-Morette}

\centerline{Center for Relativity}

\centerline{and Department of Physics}

\centerline{The University of Texas at Austin}

\centerline{Austin, Texas 78712-1081 USA}
\bigskip

\centerline{Journal of Mathematical Physics, {\bf 36}, 2137-2340 (1995)}
\par }

\vglue 20mm

\nn
\centerline {\bf Abstract}

\bigskip

The core of this article is a general theorem with a large number of specializations. Given a
manifold $N$ and a finite number of one-parameter groups of point transformations on $N$ with
generators $Y, X_{(1)}, \cdots, X_{(d)} $, we obtain, via
functional integration over spaces of pointed paths on $N$ (paths with one fixed point), a
one-parameter group of functional operators acting on tensor or spinor fields on $N$. The
generator of this group is a quadratic form in the Lie derivatives $\La_{X_{(\a)}}$ in the 
$X_{(\a)}$-direction plus a term linear in $\La_Y$. 

\smallskip
The basic functional integral is over $L^{2,1}$ paths $x: {\bf T} \ra N$ (continuous paths
with square integrable first derivative). Although the integrator is invariant under time
translation, the integral is powerful enough to be used for systems which are not time
translation invariant. We give seven non trivial applications of the basic formula, and we
compute its semiclassical expansion.

\smallskip
The methods of proof are rigorous and combine Albeverio H\o egh-Krohn oscillatory integrals
with Elworthy's parametrization of paths in a curved space. Unlike other approaches we solve
 Schr\"odinger type equations directly, rather than solving first diffusion equations
and then using analytic continuation.

\vfill\eject

\nn
\centerline {\titre I - Introduction}

\bigskip
\bigskip

We have studied many applications of functional integration looking for its
{\it substantifique moelle}\footnote{$^{1}$ }{
F. Rabelais.
Literally ``bone marrow as a producer of substance''
in {\it Gargantua}, Prologue de l'auteur.}. Little by little, several ideas have taken shape 
concerning the domains of integration, the integrators, and the integrands.

\bigskip
\nn
{\bf 1. The domain of integration is a function space.} 

\smallskip
Working with an infinite-dimensional
space is easier than working with the limit for large $n$ of the product of $n$ copies of a
finite-dimensional space. For example, a space of pointed paths (paths with one fixed
point) is contractible even when the paths take their values in a non contractible space.
Over the years the advantages -- often the necessity -- of working with spaces of paths rather
than with the discretized version of the paths have become increasingly apparent:

\medskip
-- Semiclassical approximations, even in the presence of caustics, are obtained by expanding
a functional on the space of paths around a dominating contribution [1,2].

\medskip
-- If the paths take their values in a multiply-connected space the topology of the
space of paths plays a central role [3].

\medskip
-- A change of variable of integration, regarded as a map on the domain of integration
gives, in two lines,  the Cameron-Martin formula [4] obtained from discretized paths via a
lengthy derivation, and generalizes its applications [5].

\medskip
-- Computing functional determinants using properties of linear maps on Banach spaces is
simpler than computing limits of finite-dimensional determinants [6].

\medskip
--  A major progress in the definition and computation of functional integrals was
achieved by Elworthy [7] when he parametrized the space of paths on a Riemannian manifold
$M$, with fixed initial point $a$, by the space of paths on $T_a M$ starting at the
origin of $T_a M$.

\medskip
The importance of the domain of integration, noted in the above examples, is even more
striking in the formulation presented here. A key point is a generalization of 
 Elworthy's
idea.
 Consider a finite-dimensional manifold $N$, and denote by ${\bf T}$ some finite time interval.
 Here
we parametrize a space $\Pa_{x_0} N$ of pointed paths $x: {\bf T} \ra N$ by a space $\Pa_0
\RR^d$ of pointed paths $z: {\bf T} \ra \RR^d$. A general construction for maps
$$\Pa_0 \RR^d \ra \Pa_{x_0} N \qquad  {\mathop{\rm by }\nolimits}
\qquad z \mps x$$
is given in Section II by solving suitable  first-order {\it differential} equations, not by
constructing stochastic processes. It is stated in terms of vector fields $Y$, $X_{(1)}, \cdots , X_{(d)}$ on $N$. By
specializing $N $ and the vector fields, {\it one basic functional integral yields, with rigor and no Ansatz, a great
variety of functional integrals which are solutions of complex problems.}
A number of them are pre\-sented in Section IV.

\medskip

In general, to define the domain of integration, one needs to specify:

-- the analytic nature of the paths\footnote{\nn$^2$}{ \nn A path $z: {\bf
T} \ra \RR^{\bf d}$ is said to be $L^{2, 1}$ if 
$\int_{\bf T} dt \, |\dot z (t)|^2  < \Inf$ where 
$\dot z (t) = dz (t) / dt$.
 }  (continuous, $L^{2}, L^{2,1}, \cdots $);

-- the domains of the paths  and their ranges;

-- the behaviour of the paths  at the boundaries of their domains.

\bigskip
\nn
{\bf 2. Integrators cannot be expected to be universal.}
\smallskip

The naive approach to the definition of an integral in an infinite number of variables is to
take a limit $d = \Inf$ in a $d$-dimensional integral. Because of scaling problems this
procedure is well known to abort. For instance, if we wish to evaluate (for $a > 0$) the
integral $$I_{\Inf}: =  \int_{ \RR^\Inf }^{ }  d^\Inf x \,
\exp 
\left({- { \pi \over 	a} |x|^2 }\right) ,
\eqno  ({\mathop{\rm I.1 }\nolimits})$$
where 
$|x|^2: = \ds \sum_{ \a = 1 }^{ \Inf }
\left({ x^\a }\right)_{}^{2}$,
we may first evaluate the corresponding integral
$$I_d: = 
\int_{ \RR^d }^{ } d x \,
\exp 
\left({- { \pi \over 	a} |x|^2 }\right) 
\eqno  ({\mathop{\rm I.2 }\nolimits})$$
for finite $d$ and then set $d = \Inf$.
Since $I_d = a^{d/2}$ we get
$$I_\Inf =
\left\{\matrix{ 
0 \hfill& {\mathop{\rm if }\nolimits} & 0 < a < 1 \hfill\cr
1 \hfill& {\mathop{\rm if }\nolimits} & a = 1 \hfill\cr
\Inf \hfill& {\mathop{\rm if }\nolimits} & 1 < a \hfill\cr}\right.
\eqno  ({\mathop{\rm I.3 }\nolimits})$$
and this fails to be continuous in the parameter $a$, as should be reasonably desired.

\medskip
A way out of this difficulty is to introduce, for each value of the scaling parameter $a >
0$, an integrator $\Da_a x$ in the $d$-dimensional space $\RR^d$, namely
$$\Da_a x = a^{- d/2} d x^1 \cdots dx^d
\eqno  ({\mathop{\rm I.4 }\nolimits})$$
and to remark that it is characterized by the following integration formula
$$
\int_{ \RR^d }^{ }  \Da_a x \cdot
\exp 
\left({- { \pi \over 	a} |x|^2  - 2 \pi i 
\left\langle{ x', x }\right\rangle
}\right)=
\exp 
\left({  - \pi a |x'|^2 }\right) .
\eqno  ({\mathop{\rm I.5 }\nolimits})$$

\medskip
\nn
Here $x'$ runs over the  space $\RR_d$ dual to $\RR^d$ and the scalar product is
given\footnote{$^3$}{\nn Here and in the rest of
this paper we use the Einstein summation convention over repeated indices.  } by 
$\left\langle{ x', x }\right\rangle = x'_\a x^\a$.

\medskip

This formula is dimension-independent and hence suitable for the generalization from $\RR^d$
to a (real) Banach space ${\bf X}$.
Let ${\bf X'}$ be its dual and consider two continuous quadratic forms $Q$ on 
${\bf X}$ and $W$ on ${\bf X'}$.
Assume that $Q$ {\it and  $W$ are inverse to each other}
in the following sense. There exist continuous linear maps
$$D: {\bf X} \ra {\bf X'} \ ,
\qquad
G: {\bf X'} \ra {\bf X}$$
such that\footnote{$^4$}{\nn
\nn In standard applications, $D$ is a differential operator and $G$ the corresponding Green
operator, taking into account the boundary conditions of the domain of $D$. }
$$\left\{\matrix{
DG = GD = \un
  \hfill\cr
\noalign{\medskip}
\left\langle{ Dx, y }\right\rangle =
\left\langle{ Dy, x }\right\rangle 
   \hfill \cr
\noalign{\medskip}
Q(x) = 
\left\langle{ Dx, x }\right\rangle \  ,
\  
W(x') =
\left\langle{ x', Gx' }\right\rangle \, .
  \hfill \cr}
\right. 
\eqno  ({\mathop{\rm I.6 }\nolimits})$$
Here $\left\langle{ x', x }\right\rangle $ denotes the duality between 
 ${\bf X}$ (elements $x$) and its dual $
{\bf X'}$ (elements $x'$).

\medskip
Then we define the integrator 
$\Da_{s, Q} x$ (also denoted $\Da x$ for simplicity) by the following requirement 
$$\int_{ {\bf X} }^{ } \Da_{s, Q} x \cdot 
\exp \!
\left({- { \pi \over 	s} Q(x)  - 2 \pi i 
\left\langle{ x', x }\right\rangle
}\right)  =
\exp 
 \! \left({ - \pi s W (x') }\right) \, .
\eqno  ({\mathop{\rm I.7 }\nolimits})$$

\medskip
\nn
Here $x'$ runs over ${\bf X'}$ and there are two cases:

\smallskip
i) $s = 1$ and  $Q$ is positive definite, namely $Q(x) > 0$ for $x \ne 0 $;

ii) $s = i$ and there is  no restriction on $Q$ except it be real.

\medskip
\nn
{\it Introduction of the parameter $s$ enables us to treat in a unified way the diffusion and
the Schr\"odinger equations}\footnote{$^5$}{ \nn
It is not true, as is still often stated, that the case $s = i$ has no mathematical
foundation. See for example references [8, 9, 10, 11].  }.
We can introduce a suitable space $\Fa ({\bf X})$ of functionals on ${\bf X}$ integrable  by $\Da_{s,
Q}$ and  a norm on $\Fa ({\bf X}) $, and then compute integrals of the type
$$I = \int_{ {\bf X} }^{ } \Da_{s, Q} x \cdot F(x) 
\eqno  ({\mathop{\rm I.8 }\nolimits})$$
for $F $ in $\Fa ({\bf X})$. In both cases,
$s = 1$ and $s = i$,  we shall call $\Da_{s, Q} x$ a {\it Gaussian integrator}.

\medskip

 Gaussian integrators have the following properties:
$$\Da (x + x_0) = \Da x, \qquad
x_0
\ \ 
{\mathop{\rm  a \ fixed \ element \ of }\nolimits}
\ \ {\bf X} \eqno  ({\mathop{\rm I.9 }\nolimits})$$
$$ \Da (L x) = | \Det L | \cdot \Da x, \qquad L: {\bf X} \ra {\bf X}
\eqno  ({\mathop{\rm I.10 }\nolimits})$$
where $L$ is in a suitable class of linear changes of variables of integration,
including the obvious case  where $Q(Lx) = Q(x) $ and 
$| \Det L | = 1$ (cf. formula (A.140) in Appendix A).

\medskip
Gaussian integrators are not the only possible integrators. 
In reference [11] we  have developed an axiomatic for functional integrals on a Banach space
$\Phi$ expressed in terms of integrators 
$\Da_{\T, Z}$ defined by 
$$\int_\Phi \Da_{\T, Z} \vp \cdot \T (\vp, J)
= Z(J) 
\eqno  ({\mathop{\rm I.11 }\nolimits})$$
for $\vp$ in $\Phi$, $J$ in the dual $\Phi'$ of $\Phi$,
where $\T $ and $Z$ are two given continuous bounded functionals
$$\T: \Phi \ts \Phi' \ra \CC \quad ,
\quad
Z: \Phi' \ra \CC.$$

\smallskip
\nn
In quantum field theory, we interpret $\vp$ as a field and $J$ as a source,
$Z(J)$ is then the Schwinger generating functional for the $n$-point functions.

\bigskip
\nn
{\bf 3. Integrands and integrators. }

\smallskip
Splitting the quantity inside the integral sign into ``integrator'' and
``integrand'' belongs to the art of integration, but rules of thumb apply:

\smallskip
-- When the functional integral has its origin in physics try not to break up the action into,
say, kinetic and potential contributions. On the other hand, do not hesitate to work with a
potential which is a functional of a path rather than a function of its {\it value} (e.g. in
equation (III.1), $V$ is a functional of $z$, not a function of $z(t)$).

\smallskip

-- Look for a possible change of variable of integration;  this may suggest
 a practical choice for the integrand.

\smallskip
-- Gaussian integrators have a wealth of simple, powerful properties;
 look for exponentials of quadratic forms,  this suggests a practical choice for
the integrator. \medskip

\bigskip
\nn
{\bf 4. A basic functional integral. }
\smallskip

The core of this article is a theorem which provides the mathematical
underpinning for a great variety of functional integrals. 
It consists of two parts: the definition of a functional integral, and the partial
differential equation satisfied by 
the {\it value} of the functional integral, as a function of a set of parameters. Given a
manifold $N$  consider the  $L^{2,1}$ paths over a finite time interval ${\bf T}$ with values
in $N$,  $x: {\bf T} \ra  N$, and with a fixed point 
$x_0 \in N$; for instance, if 
${\bf T} = \left\lbk{ t_a ,  t_b  }\right\rbk$, the fixed point
can be either 
$x(t_a) $ or $x(t_b)$. Given $d + 1$ vector fields $Y$ and 
$ X_{(\a)}  $ on $N$,  define a map $P$ from a space $\Pa_0 \RR^d $ of $L^{2, 1}$
paths into a space $\Pa_{x_0} N$ of $L^{2, 1} $ paths, 
$$P: \Pa_0 \RR^d \ra \Pa_{x_0} N,$$
by $P(z) = x $, where
$$
\left\{\matrix{
 dx (t, z)= 
 X_{(\a)}
(x (t, z)) dz^\a +
Y
(x(t, z)) dt 
\hfill\cr
\noalign{\medskip}
   x(t_0, z)  = x_0 .\hfill \cr}
\right. \eqno ({\mathop{\rm I.12 }\nolimits} )
$$
Here $t_0$ and $x_0$ are fixed, with $t_0$ in ${\bf T}$ and $x_0$ in $N$, and the paths
$z : {\bf T} \ra \RR^d$  satisfy $z \!\left({ t_0 }\right) = 0$. 
In general, the vector fields do not commute, that is: 
$$\left\lbk{  X_{(\a)} ,  X_{(\b)}  }\right\rbk
\ne 0 \quad , \quad \left\lbk{  Y ,  X_{(\a)}  }\right\rbk
\ne 0\, ;$$
therefore the solution of (I.12) is of the form
$$x(t, z) = x_0 \cdot \Si(t, z)
\eqno ({\mathop{\rm I.13 }\nolimits} )$$
where $x$ is a {\it function of}  $x_0$ {\it and  $t$, and} {\it a functional of}  $z$.
Here $\Si (t, z)$ is a transformation in $N$, depending on $t$ and $z$ as stated.
Only when 
$\left\lbk{  X_{(\a)} ,  X_{(\b)}  }\right\rbk = \left\lbk{  Y ,  X_{(\a)}  }\right\rbk = 0$ can one express
$\Si(t, z) $ as a function of $t$ and $z(t)$.

\medskip
In Appendix C the familiar theorems of differential equations are extended to cover the properties of (I.12).
\medskip

Given (I.12) one can express an integral over $\Pa_{x_0} N $ as an integral over
$\Pa_0 \RR^d$, which is an integral that one can manipulate and compute. The partial
differential equation satisfied by the integral on 
$\Pa_{x_0} N$ is expressed in terms of Lie derivatives along the vector fields $Y$ and
$  X_{(\a)}  $. Conversely, given a parabolic 
partial differential equation, one can construct in many cases the path integral representation of its solutions
(see for instance in Section IV, the lift of a covariant Laplacian at a point of  a Riemannian manifold in terms
of  Lie derivatives at a point of its frame bundle [13]). 

\medskip

The basic equations are given in the first
paragraph of Section II, followed by a summary of notations used throughout the paper.

\bigskip
\nn
{\bf 5. Examples.}

\smallskip
In Section III we compute semiclassical approximations of the basic functional integral
(II.1), and in Section IV we specialize the basic integral by making particular choices of the
manifold $N$.
 We treat
in detail the following examples:
\smallskip
-- $N = \RR^d$, in cartesian or polar coordinates. 
\smallskip
-- $N$ is a frame bundle $O(M) $ over a Riemannian manifold $M$.  We give
the explicit functional integral representing  the solution $\Psi$ of the Schr\"{o}dinger
equation on $M$, with initial wave function $\phi$.

\smallskip

-- $N$ is  a multiply-connected space.

\smallskip
-- $N$ is a $U(1)$-bundle; the basic integral solves the Schr\"odinger equation for a
particle in an electromagnetic field.

\smallskip
-- $N$ is a symplectic manifold; coherent-state transitions can be
obtained from the basic integral.

\smallskip
\nn
We conclude this section by an analysis of the Bohm-Aharonov effect, where all the previous techniques are
brought to bear.

\bigskip
\nn
{\bf 6. Techniques.}
\smallskip

In the course of computing our basic functional integral  in various situations, we have
used the properties of linear changes\footnote{\nn$^6$}{  \nn
Linear changes of variable of integration in an infinite-dimensional space are sufficiently
powerful and varied for the purposes of this paper. In a later publication we shall present
nonlinear changes, simplifying and generalizing earlier works such as [12].  } of
variable of integration and properties of functional determinants; they are given in two
appendices.  The transformation $\Si$ in formula (I.13) is related to
the Cartan development map; the properties of this map are discussed in the
third appendix.

\vfill\eject
\vglue 1cm
\nn
\centerline {\titre II - A general theorem}

\bigskip
\bigskip

The primary goal of this section is to define functional integrals over $L^{2,1}$ pointed paths
(paths with a fixed point) taking their values in a manifold $N$, by reducing them to functional
integrals over paths taking their values in a flat space $\RR^d$. The main definition is given as
follows
$$\left({ U_T \phi }\right)
\left({ x_0 }\right) : = \int_{ {\bf Z }_0 }^{ }
\Da_{s, Q_0 } z\cdot \exp 
\left({ - {\pi  \over s } Q_0 (z)  }\right)
\phi \left({ x_0 \cdot \Si (T, z) }\right) .
\eqno ({\mathop{\rm II.1 }\nolimits} )$$
All the notations are given in paragraph 1. The previous integral, being also denoted $\Psi \! \left({ T, x_0
}\right)$, is a solution of the {\it generalized Schr\"odinger equation }
$$
{ \part \Psi  \over \part T} =
{ s \over 4 \pi } h^{\a \b } 
\La_{X_{(\a) } }
\La_{X_{(\b) } } \Psi + \La_{Y }\Psi
\eqno ({\mathop{\rm II.2 }\nolimits} )$$
with initial condition $\Psi \! \left({ 0, x_0 }\right)
= \phi \left({ x_0 }\right)$.

\medskip
We give also a general construction of {\it time-ordered products} in the form of a functional
integral generalizing equation (II.1), namely
$$
\left({ U_T^F \phi  }\right)
\left({ x_0  }\right)
= \int_{ {\bf Z }_0J}^{ } \Da_{s, Q_0 } z \cdot 
\exp \left({ - {\pi  \over s } Q_0 (z)  }\right)
F ({\bf T} , z) \phi \left({ x_0 \cdot \Si (T, z)  }\right).
\eqno ({\mathop{\rm II.3 }\nolimits} )$$
In paragraph 2, we construct the simplest functional integral of type (II.1) and prove that it
satisfies equation (II.2). In paragraph 3, we study the general case.

\medskip
\nn
{\bf 1. The setup, and a summary of notations.}
\smallskip
\nn 1.1. {\it A manifold and vector fields.}

\smallskip
-- A finite-dimensional manifold $N$.

\smallskip
-- One-parameter groups acting on $N$, denoted $\s_{(\a)} (r) $; 
here $\a$ takes the values  $0, 1, 2, \cdots, d$ and $r$ is a real parameter; the transform of a
point $x$ in $N$ under $\s_{(\a) }  (r) $ is denoted by 
$x \cdot \s_{(\a) }  (r)$, and
$$\s_{(\a) }  (r) \circ \s_{(\a) } (s) = \s_{(\a) }
(r + s) \,  . \eqno ({\mathop{\rm II.4 }\nolimits} )$$

\smallskip
-- The generator of 
$\left\lbc{ \s_{(\a) } (r) }\right\rbc$ is the vector field 
$X_{(\a) } $ in $N$ such that 
$${ d \over dr } \left({ x \cdot  \s_{(\a) } (r)  }\right)
= X_{(\a) } \left({ x \cdot \s_{(\a) }  (r) }\right) 
\eqno ({\mathop{\rm II.5 }\nolimits} )$$
and in particular
$$X_{(\a) } (x) = 
\left.{ { d \over dr } \left({ x \cdot \s_{(\a) }  (r) }\right)  }\right\vert_{r = 0}
\eqno ({\mathop{\rm II.6 }\nolimits} )$$
for any point $x$ in $N$.
We do not assume that the vector fields 
$X_{(\a)}$ commute, hence
$\left\lbk{ X_{(\a)}, X_{(\b)} }\right\rbk \ne 0$
in general. We often write $Y$ for $X_{(0)}$ emphasizing its special role.

\smallskip

-- $\La_X$ denotes the Lie derivative w.r.t. the vector field $X$.

\bigskip
\nn 1.2. {\it Pointed paths on the flat space $\RR^d$.}

\smallskip

-- ${\bf T} $ is a time interval of length $T$, hence
$$ {\bf T } 
= \left\lbk{ t_a, t_b }\right\rbk \ ,
\qquad T = t_b - t_a \, .
\eqno ({\mathop{\rm II.7 }\nolimits} )$$ 

\smallskip

-- $t_0$ is a chosen element of ${\bf T}$; the standard choices are $t_0 = t_a$ or $t_0 = t_b$.

\smallskip

-- ${\bf Z}_0$ (or ${\bf Z}_{0, {\bf T}}$ if we need to specify ${\bf T}$)
consists of the real vector-valued functions $z = 
\left({ z^1 , \cdots,  z^d }\right)$ whose components 
$t \mps z^\a (t) $  are continuous functions with square-integrable derivatives $\dot z^\a$. We
assume the normalization  $z^\a \!\left({ t_0 }\right) = 0$.

\smallskip

-- ${\bf Z}_0'$ is the space dual to ${\bf Z}_0$. Its elements are interpreted as vector-valued
distributions
$z' = \left({ z'_1 , \cdots, z'_d }\right)$, each component being the derivative of an
$L^2$-function. The duality is given by
$$\left\langle{ z', z }\right\rangle =
\int_{ {\bf T }J}^{ } dt \, z'_\a (t) z^\a (t) 
\eqno ({\mathop{\rm II.8 }\nolimits} )$$ 
(summation over $\a$ and integration over $t$).

\smallskip
--
$s$ is a parameter equal to $1$ or $i$; its square root 
$\sqrt{ s }$ is normalized as follows
$$\sqrt{ s } = 
\left\{\matrix{
 1 \hfill& {\mathop{\rm  if }\nolimits} \hfill& s = 1
   \hfill \cr
\noalign{\medskip}
e^{\pi i / 4} \hfill& {\mathop{\rm  if }\nolimits} \hfill& s = i \, .
  \hfill \cr}
\right. \eqno ({\mathop{\rm II.9 }\nolimits} )$$ 

\smallskip
--
$h = \left({  h_{\a \b } }\right)$ is a constant invertible symmetric real matrix of size $d \ts d$. 
We denote by $\left({  h^{\a \b } }\right)$ the inverse matrix and we assume that $h$ is positive
definite in case $s$ is equal to $1$. By a suitable linear change of coordinates, we may take $h$
into a diagonal form
$$h = \diag (1, \cdots, 1, - 1, \cdots, - 1)
\eqno ({\mathop{\rm II.10 }\nolimits} )$$ 
with $p$ elements $+ 1$, $q$ elements $- 1$ and $p + q = d$. 
Hence $p = d$, $q = 0$ in the case $s = 1$.

\smallskip

-- We introduce a quadratic form $Q_0 $ on ${\bf Z}_0$ as follows
$$Q_0 (z) = \int_{ {\bf T }J}^{ }
dt \, h_{\a \b } \, 
\dot z^\a (t) \dot z^\b (t) \, .
\eqno ({\mathop{\rm II.11 }\nolimits} )$$
The corresponding kernel is given by 
$$D_{\a \b } (u, r) =  
{ {}_{1} \over {}^{2} } 
{\d^2 Q_0 (z)  \over \d z^\a (u) \d z^\b (r) } \, , 
\eqno ({\mathop{\rm II.12 }\nolimits} )$$ 
hence the representation 
$$Q_0 (z) = 
\int_{ {\bf T }   }^{ }du
\int_{ {\bf T }J}^{ }
 dr\, D_{\a \b } (u, r) \, z^\a(u) \, z^\b (r) \, .
\eqno ({\mathop{\rm II.13 }\nolimits} )$$

\smallskip

-- On ${\bf Z}'_0$ we consider a quadratic form $W_0$, with kernel 
$G^{\a \b } (u, r)$. As above, we have the relations 
$$
G^{\a \b } (u, r) = 
{ {}_{1} \over {}^{2} } 
{\d^2 W_0 (z')  \over \d z'_\a (u) \d z'_\b (r) } \, ,
\eqno ({\mathop{\rm II.14 }\nolimits} )$$ 
$$
W_0 (z') = \int_{ {\bf T}  }^{ } du
\int_{ {\bf T }J}^{ }
 dr\, G^{\a \b } (u, r) \, z'_\a (u)  z'_\b (r) \, .
\eqno ({\mathop{\rm II.15 }\nolimits} )$$

\smallskip

-- We assume that the quadratic forms
$Q_0 $ on ${\bf Z}_0$ and $W_0$ on ${\bf Z}'_0$ are inverse to each other in the sense of relation
(I.6). In terms of kernels, this is expressed as follows
$$\int_{ {\bf T }J}^{ } dt \, D_{\a \b }
\left({  s_1, t }\right)
G^{\b \g }
\left({ t, s_2  }\right)
= \d_\a^\g \, 
\d \left({ s_1 - s_2  }\right).
\eqno ({\mathop{\rm II.16 }\nolimits} )$$
Here are explicit formulas for the kernels:
$$ D_{\a \b } (u, r) = 
\int_{ {\bf T }J}^{ } dt \, h_{\a \b } \, 
\d' (t - u) \d' (t - r) = - h_{\a \b } \, 
\d'' (u - r) \, , 
\eqno ({\mathop{\rm II.17 }\nolimits} )$$
that is $D : {\bf Z }_0 \ra {\bf Z }'_0 $ is the differential operator with matrix 
$\left({ - h_{\a \b } { d^2 \over dt^2} }\right) $.
Hence, 
$G^{\a \b }(u, r)$ is the corresponding Green's function  taking into account the boundary condition 
$z^\a \!\left({ t_0  }\right) = 0$, that is
$$G^{\a \b } (u, r) 
=
\left\{\matrix{ h^{\a \b }
\inf 
\left({ u - t_0, r - t_0 }\right)
\hfill & {\mathop{\rm for \  }\nolimits} u \geq t_0, r \geq t_0 ,
  \hfill\cr
\noalign{\medskip}
h^{\a \b } \inf \left({ t_0 - u , t_0 - r  }\right) \hfill &{\mathop{\rm for \  }\nolimits}
 u \leq t_0, r \leq t_0 ,
   \hfill \cr
\noalign{\medskip}
0 \hfill &
{\mathop{\rm otherwise. }\nolimits}
 \hfill \cr}
\right. \eqno ({\mathop{\rm II.18 }\nolimits} )$$

\medskip
\nn
{\it Remark.}
We shall refrain from integrating by parts in formulas like (II.11) in order not to have to make
explicit statements about boundary conditions.

\bigskip
\nn
1.3. {\it  Functional integrals on $\Pa_0 \RR^d $. }

\smallskip
 
-- The space of paths ${\bf Z}_0$ shall also be denoted by 
$\Pa_0 \RR^d$ to remind us of the flat space $\RR^d$ where the paths lie and of the fixed point $0$
(origin in $\RR^d$) of the paths $z$. We shall consider a variety of pairs $(Q, W)$ consisting of
a quadratic form $Q$ on a space ${\bf Z}$ and a quadratic form $W$ on its dual ${\bf Z}'$ satisfying the
analogues of relations (II.12), (II.14) and (II.16).

\smallskip

-- For such a pair $(Q, W)$ we have a translation invariant integrator 
$\Da_{s, Q} z$ on ${\bf Z}$ characterized by 
$$
\int_{ {\bf Z }J}^{ }
\Da_{s, Q} z \cdot 
\exp 
\left({ -{ \pi  \over s}  Q(z) - 2 \pi i \left\langle{ z', z }\right\rangle }\right)
=
\exp 
(-  \pi s W (z')) 
\eqno ({\mathop{\rm II.19 }\nolimits} )$$ 
for $z'$ in ${\bf Z'}$.
The normalizations are chosen so that Gaussian integrators and Fourier transforms do not include
powers of $\pi$ depending on the dimension of their domain of definition.

\medskip
We also write
$$\Da \om_{s, Q} (z) = \Da_{s, Q} z \cdot 
\exp 
\left({ - { \pi  \over s}  Q(z)  }\right)
 \eqno ({\mathop{\rm II.20 }\nolimits} )$$ 
and refer to both 
$\Da_{s, Q}$ and $\Da \om_{s, Q} $ as 
``Gaussian integrators''. When working with the basic pair
$ \left({ Q_0, W_0  }\right) $ we simply write 
$\Da_s$ and 
$\Da \om_s$. In the context of an application where $s$ has been chosen once for all, we omit it
in the notations.

\bigskip
\nn
1.4. {\it  Functional integrals on $\Pa_{x_0} N $. }

\smallskip

-- We fix a point $x_0$ in $N$ and consider continuous paths 
$x : {\bf T} \ra N$ with the fixed point $x \!\left({ t_0 }\right) = x_0$ and square-integrable
velocity\footnote{$^7$}{More precisely, for every smooth function $f$ in $C^\Inf (N)$, we assume
that the continuous function $t \mps f(x(t))$ on ${\bf T}$ is the primitive of a function  in 
$L^2 ({\bf T })$. }. The set of all such paths is denoted by
$\Pa_{x_0} N$.

\smallskip

-- The time interval ${\bf T}$ being given, consider an element $z$ of ${\bf Z}_0$. As we shall see
in Appendix C, the differential equation 
$$
\left\{\matrix{
  dx(t) = X_{(\a)}
( x(t) ) dz^\a (t) + Y (x(t)) dt 
   \hfill \cr
\noalign{\medskip}
x \!\left({ t_0 }\right) = x_0
  \hfill \cr}
\right. 
\eqno ({\mathop{\rm II.21 }\nolimits} )$$
admits a unique solution $x(\cdot)$ in 
$\Pa_{x_0} N$.

\smallskip

-- The previous construction defines a parametrization $P$ of the space
$\Pa_{x_0} N$ of pointed paths on $N$ by the space $\Pa_0 \RR^d$ of pointed paths in $\RR^d$, namely
$$P : \Pa_0 \RR^d \ra \Pa_{x_0} N$$
by taking $z$ into $x$. If necessary we shall denote by $x(t, z)$ the solution of the differential
equation (II.21) for given $z$ in $\Pa_0 \RR^d$. Hence $x(t, z)$ is {\it a function of $t$ and a
functional of $z$}.

\smallskip

--  Assume now that ${\bf T} = \left\lbk{ 0, T }\right\rbk$ and $t_0 = 0$.
 With the previous definitions, define 
$\Si (T, z)$ as the transformation taking a point $x_0$ in $N$ into the point
$x(T, z)$. 

\smallskip

-- Take a good\footnote{$^8$}{Any function in $C^\Inf (N)$ with compact support will do.}
 function $\phi$ on $N$. Define a functional $\ww \phi$ on ${\bf Z}_0$ by 
$$\ww \phi (z) = \phi
\left({  x_0 \cdot \Si (T, z)}\right) .
\eqno ({\mathop{\rm II.22 }\nolimits} )$$
It is integrable under the integrator $\Da \om_s$ on ${\bf Z}_0$.
By integrating we get
$$
\eqalign{ I \! \left({  \phi, T,  x_0}\right)
&= \int_{ {\bf Z}_0  }^{ }
\Da \om_s (z)  \ww \phi (z) \cr
& = \int_{ {\bf Z}_0  }^{ }
\Da_s z \cdot \exp 
\left({ - {\pi  \over s }  Q_0 (z) }\right)
\cdot \phi 
\left({ x_0 \cdot \Si (T, z)  }\right) .
\cr}
\eqno ({\mathop{\rm II.23 }\nolimits} )$$
We could replace $Q_0$ by another quadratic form on ${\bf Z}_0$.

\smallskip

-- 
The functional operator $U_T$ associates to the function $\phi$ on $N $ the function 
$x_0 \mps I \!\left({  \phi, T, x_0 }\right)$ on $N$.

\bigskip
\nn
1.5. {\it  Variational techniques.}
\smallskip

A frequently used technique consists in introducing one-parameter variations in the space of paths,
or in the space of functionals on the space of paths. We work in particular with the following
variations.

\smallskip
(i)  Fix a time interval 
${\bf T} = \left\lbk{ t_a, t_b  }\right\rbk$ and consider the space 
$\Pa \RR^d$ of paths $x : {\bf T} \ra \RR^d$, say of class\footnote{$^9$}{That is, continuous with
continuous first-order derivatives.} $C^1$. 
Assuming an action functional $S : \Pa \RR^d \ra \RR$ (for instance, the time integral of a
Lagrangian 
$L(x, \dot x, t)$), the critical points  of $S$ will form a $2d$-dimensional manifold 
$\Ua^{2d}$, the so-called {\it space of } (classical) {\it motions}, denoted by 
$x_{{\mathop{\rm cl }\nolimits} }$. We can parametrize them by a set of parameters $\mu =
\left({ \mu^1, \cdots, \mu^{2d}  }\right)$; for a given $k$ between $1$ and $2d$, the derivative 
$\part x_{{\mathop{\rm cl }\nolimits} }
 (\mu) / \part \mu^k$ defines a variation through classical paths, and we get $2d$ of
them.

\smallskip

(ii)
If $t_0$ is a given epoch in ${\bf T}$, the pointed paths are defined by the condition 
$x \!\left({ t_0 }\right) = 0$. They form the space $\Pa_0 \RR^d$. In paragraph III.1 we shall use a
one-parameter family of pointed paths $x (\l)$ in $\Pa_0 \RR^d$ (for $\l$ running over $[0, 1]$).

\smallskip

(iii) Introducing $d$ independent boundary conditions at $t_a$, and similarly at $t_b$, we define the
subspace
$\Pa_{a, b} \RR^d $ of $\Pa \RR^d$.
In paragraph B.2, we consider a one-parameter family of paths 
$x(\l) $ in $\Pa_{a, b} \RR^d$, with $\l$ in $[0, 1]$, such that $x(0)$ belongs to 
$\Ua^{2d} \cap \Pa_{a, b} \RR^d$. The presence and nature of the caustics is conveniently analyzed in
terms of this intersection.

\smallskip

(iv) In paragraph  III.2  we consider a one-parameter family of
action functionals $S(\nu)$ defined on a space $\Pa_0 \RR^d $ of pointed paths. 

\bigskip
\nn
{\bf 2. The one-dimensional case.}
\smallskip

We begin by the simple case $d = 1$. Here $X = X_{(1)}$ is a vector field on $N$ generating the one-parameter
group of transformation $\s(r)$ on $N$. Hence $\s(r)$ obeys the differential equation
$$d \! \left({ x_0 \cdot \s (r)  }\right) =
X \! \left({ x_0 \cdot \s (r)  }\right) dr 
\eqno ({\mathop{\rm II.24 }\nolimits})$$
for a fixed $x_0$ in $N$. We set $Y  = 0$ and consider the time interval
${\bf T}  = [0, T]$. The differential equation (II.21) reads now as 
$$
\left\{\matrix{
dx (t) = X (x(t)) \cdot dz (t)
  \hfill\cr
\noalign{\medskip}
x(0) = x_0     \, . 
  \hfill \cr}
\right. \eqno ({\mathop{\rm II.25 }\nolimits})$$
Substituting $r = z(t)$ in (II.24) we see that the solution to the previous equation is given by $x(t) =
x_0 \cdot \s (z(t))$. Hence the transformation $\Si (T, z)$ is simply $\s (z(T))$.

\bigskip
The path space ${\bf Z}_0$ consists of the $L^{2,1}$ functions 
$z : [0, T] \ra \RR$ such that $z(0) = 0$. It is endowed with the quadratic form
$$Q_0 (z) = \int_{0 }^{T } dt \, \dot z (t)^2 \, . 
\eqno ({\mathop{\rm II.26 }\nolimits})$$
The corresponding integrator is, for simplicity, denoted by $\Da_s z$, hence our basic path integral
specializes to
$$\Psi (T, x) =
\int_{ {\bf Z}_0 }^{ }
\Da_s z \cdot \exp 
\left({- {\pi  \over s } \int_{ 0}^{ T}  dt \, \dot z (t)^2}\right)
\phi (x \cdot \s (z(T))) \, . 
\eqno ({\mathop{\rm II.27 }\nolimits})$$
Our goal is to show that we have solved the differential equation
$$
{\part  \over \part T}
\Psi (T, x) = {s \over 4 \pi }
\La_X^2 \Psi (T, x) \, .
\eqno ({\mathop{\rm II.28 }\nolimits})$$

\bigskip
Fix a point $x$ in $N$ and define a function of a real variable 
$h(r) = $ \break  $\phi (x \cdot \s (r))$. From the group property 
$\s(r_1) \circ \s(r_2) = \s(r_1 + r_2)$, we get
$$\phi (x \cdot \s (r) \cdot \s (z(T))) =
h (r + z (T)) \, . 
\eqno ({\mathop{\rm II.29 }\nolimits})$$
Then define
$$H  (T, r) : = \int_{ {\bf Z}_0 }^{ } \Da_s z \cdot 
e^{- \pi Q_0 (z) / s }
h (r + z (T)) 
 \, . 
\eqno ({\mathop{\rm II.30 }\nolimits})$$

\medskip
\nn
Since the integrand $h(r + z(T))$ depends on the path $z$ through $z(T)$, a  linear change of variable $z \mps
z(T)$ transforms immediately this functional integral over ${\bf Z}_0$ into an ordinary integral over $\RR$
(see formula (A.38) in Appendix A). It is easier to use directly the properties of Fourier transforms of
Gaussian integrators underpinning (A.38). Denoting by $\hh h$ the Fourier transform of $h$, we obtain
$$h(r + z(T)) = \int_{ \RR }^{ }
d \r \, \hh h(\r) \, e^{2 \pi i \r (r + z (T))}
 \, . 
\eqno ({\mathop{\rm II.31 }\nolimits})$$
Therefore, after changing the order of integration, we get
$$H(T, r) = \int_{ \RR }^{ }
d \r \, \hh h(\r) \, e^{2 \pi i \r r}
\int_{ {\bf Z}_0 }^{ } \Da_s z \cdot \exp \!
\left({- {\pi  \over s } Q_0 (z) + 2 \pi i \r z (T)  }\right)
 \, . 
\eqno ({\mathop{\rm II.32 }\nolimits})$$
By the definition (II.19) and the duality 
$\r z (T) = 
\left\langle{ \r \d_T, z }\right\rangle$,
we get
$$\int_{ {\bf Z}_0 }^{ } \Da_s z \cdot \exp \! 
\left({- {\pi  \over s } Q_0 (z) + 2 \pi i \r z (T)  }\right)
=
\exp \! \left({ - \pi s W_0 (\r \d_T) }\right)
\, . 
\eqno ({\mathop{\rm II.33}\nolimits})$$

\medskip
\nn
By (II.15) and (II.18) we obtain
$$ W_0(\r \d_T) = \r^2 \, G (T, T) = \r^2 T
\, . 
\eqno ({\mathop{\rm II.34}\nolimits})$$
Collecting equations (II.32) to (II.34) we conclude
$$H (T, r) = \int_{ \RR }^{ }
d \r \, \hh h (\r) e^{2 \pi i \r r -  \pi s \r^2 T}
\eqno ({\mathop{\rm II.35}\nolimits})$$
and by derivation under this integral sign, we conclude
$$
{\part H \over \part T} 
= { s \over 4 \pi }  
{\part^2 H \over \part r^2} 
\, . 
\eqno ({\mathop{\rm II.36}\nolimits})$$

\bigskip
The vector field $X$ is the generator of the one-parameter group of transformations
$\s(r)$.
Hence, for every integer $m \geq 0$, we get
$$
\left({ {\part  \over \part r } }\right)_{}^{m} 
h (r + z(T)) =
\left({ \La_X^m \phi }\right)
(x \cdot \s (r + z (T))) .
\eqno ({\mathop{\rm II.37 }\nolimits})$$

\smallskip
\nn
By differentiating under the integral sign in formula (II.30), we deduce
$$\left({ {\part  \over \part r } }\right)_{}^{m} 
H(T, r) =
\int_{ {\bf Z}_0 }^{ }
\Da_{s} z \cdot \exp \! 
\left({ -  {\pi \over s}  Q_0 (z)  }\right) 
 \La_X^m \phi
(x \cdot \s (r + z (T))) \, .
\eqno ({\mathop{\rm II.38 }\nolimits})$$
Using definition (II.27), we get therefore
$$
\left.{ \left({ {\part  \over \part r } }\right)_{}^{m} 
H(T, r)  }\right\vert_{r = 0 } = 
 \La_X^m \Psi
(T, x)
\eqno ({\mathop{\rm II.39 }\nolimits})$$
by applying $m$ times the differential operator 
$\La_X$ acting on functions of $x$. In particular, for $m = 0$, we get
$$
\left.{ H(T, r)  }\right\vert_{r = 0}
=
\Psi (T, x)
\eqno ({\mathop{\rm II.40 }\nolimits})$$
and both sides can be derivated with respect to $T$, giving
$$
\left.{  {\part  \over \part T }  H (T, r)  }\right\vert_{r = 0}
= {\part  \over \part T }
\Psi (T, x) .
\eqno ({\mathop{\rm II.41 }\nolimits})$$
Setting $r = 0$ in the relation 
$\ds {\part H \over \part T } =
{s  \over 4 \pi }
{\part^2  \over \part r^2 } H$ and using relation 
(II.39) for $m = 2$, and relation (II.41), we conclude the proof of the differential
equation 
$${\part  \over \part T }
\Psi (T, x) =
{s \over 4 \pi }
\La_X^2 \Psi (T, x) .$$

\bigskip
\nn
{\it Remark. }
The previous proof can be recast in the following symbolical way [5].
We begin with a consequence of equations (II.33) and (II.34), namely
$$\int_{ {\bf Z}_0 }^{ }
\Da_{s} z \cdot \exp \!
\left({ -  {\pi \over s}  Q_0 (z) }\right)  \exp  ( 2 \pi i \r z (T) )
= \exp ( - \pi s \r^2 T )\, ,
\eqno ({\mathop{\rm II.42 }\nolimits})$$
one of many characterizations of our integrator 
$\Da_s z$. Make the formal substitution
$$\r = 
{1 \over 2 \pi i}
\La_X,$$
and get
$$\int_{ {\bf Z}_0 }^{ }
\Da_{s} z \cdot \exp \!
\left({ -  {\pi \over s}  Q_0 (z) }\right) 
\exp (z (T) \La_X) =
\exp \! \left({   {sT  \over 4 \pi }   \La_X^2 }\right) .
\eqno ({\mathop{\rm II.43 }\nolimits})$$

\smallskip
\nn
Apply this operator identity to the function $\phi (x)$ and notice that the operator 
$\exp (r \La_X)$ transforms $\phi (x) $ into 
$\phi (x \cdot \s (r))$.
Hence we get
$$\int_{ {\bf Z}_0 }^{ }
\Da_{s} z \cdot \exp \!
\left({ -  {\pi \over s}  Q_0 (z) }\right) 
\phi (x \cdot \s (z (T)) =
\left({  \exp  {sT  \over 4 \pi }   \La_X^2 }\right)
\phi (x) 
\eqno ({\mathop{\rm II.44 }\nolimits})$$
that is 
$$\Psi (T, x) = 
\exp \!
\left({   {sT  \over 4 \pi }   \La_X^2 }\right)
\phi (x) .
\eqno ({\mathop{\rm II.45 }\nolimits})$$
This integrated form is equivalent to the differential equation (II.28) together with the
initial condition 
$$\Psi (0, x) = \phi (x) \, . 
\eqno ({\mathop{\rm II.46 }\nolimits})$$

\bigskip
\nn
{\bf 3. The general case.}

\smallskip
\nn
3.1. {\it The group property.}
\smallskip
The functional operator $U_T$ defined by formula (II.1) satisfies the group property
$$U_T \circ U_{T'} = U_{T + T'}
\eqno ({\mathop{\rm II.47 }\nolimits})$$
for $T > 0$, $T' > 0$. The proof rests on three facts.

\medskip
a) {\it Group property for the point-transformations} $\Si (T, z)$:

The relevant function space ${\bf Z}_{0, T}$ consists of the paths 
$z : [0, T] \ra \RR^d$ with $L^2$ derivative.
For $z$ in ${\bf Z}_{0, T}$ and $z'$ in ${\bf Z}_{0, T'}$, we define a new path $z \ts z'$ by the rule
$$(z \ts z') (t) = 
\left\{\matrix{
z(t) \hfill & {\mathop{\rm for  }\nolimits} &
0 \leq t \leq T
  \hfill\cr
\noalign{\medskip}
   z(T) + z'(t- T) \hfill & {\mathop{\rm for  }\nolimits} &
T \leq t \leq T + T' \, .
  \hfill \cr}
\right. 
\eqno ({\mathop{\rm II.48 }\nolimits})$$
It is obvious that $z \ts z'$ is an element of ${\bf Z}_{0, T + T'}$.
Furthermore, the map $(z, z') \mps z \ts z'$ is an isomorphism of 
${\bf Z}_{0, T} \ts {\bf Z}_{0, T'}$
onto ${\bf Z}_{0, T+T'}$ and,  by the uniqueness of the solution of the differential equation (II.21), we obtain
$$x_0 \cdot \Si (T + T', z \ts z') = x_0 \cdot \Si (T, z) \cdot \Si (T', z') 
\eqno ({\mathop{\rm II.49 }\nolimits})$$
(see also formula (C.23) in Appendix C).

\medskip
b) {\it Quadratic forms}:

We denote by $Q_{0, T}$ the basic quadratic form on ${\bf Z}_{0, T}$ given by equation (II.11), namely:
$$Q_{0,T}(z) = \int_{0 }^{ T} dt \, 
h_{\a \b} \dot z^\a (t) \dot z^\b (t) \, .
\eqno ({\mathop{\rm II.50 }\nolimits})$$
Using similar definitions for $Q_{0, T'}$ and $Q_{0, T + T'}$, one obtains immediately
$$Q_{0, T + T'} (z \ts z') =
Q_{0,T} (z) + Q_{0,T'} (z')
\eqno ({\mathop{\rm II.51 }\nolimits})$$
and by exponentiating
$$\exp \! 
\left({ - { \pi \over s } Q_{0,T+T'} (z \ts z') }\right)
=
\exp \! 
\left({ - { \pi \over s } Q_{0,T} (z ) }\right)
\exp \! 
\left({ - { \pi \over s } Q_{0,T'} (z' ) }\right) \, . 
\eqno ({\mathop{\rm II.52 }\nolimits})$$

\medskip
We can identify 
${\bf Z}_{0, T+T'}$ to ${\bf Z}_{0, T} \ts {\bf Z}_{0, T'}$. This identification entails an identification of
the dual space ${\bf Z}'_{0, T+T'}$ to the product space 
${\bf Z}'_{0, T} \ts {\bf Z}'_{0, T'}$.
We use the notations $\z$ for elements of ${\bf Z}_{0, T}$, $\z'$ in ${\bf Z}_{0, T'}$
and denote by $\z \ts \z' $ the corresponding element in 
${\bf Z}'_{0, T+T'}$.
Using simple algebra, one derives the identity
$$W_{0, T + T'} (\z \ts \z') 
= W_{0, T} (\z) + W_{0, T'} (\z') 
\eqno ({\mathop{\rm II.53 }\nolimits})$$
from (II.51).
Here $W_{0,T}$ denotes the quadratic form on ${\bf Z}'_{0, T}$ inverse to 
$Q_{0,T}$, etc.

\medskip
c) {\it Integrators}:

The following form of Fubini's theorem holds
$$
\int_{ {\bf Z}_{0, T+T'} }^{ }
\Da_s u \cdot F(u) = 
\int_{ {\bf Z}_{0, T} }^{ } \Da_s z
\int_{ {\bf Z}_{0, T'} }^{ } \Da_s z' 
\cdot F(z \ts z') \, . 
\eqno ({\mathop{\rm II.54 }\nolimits})$$
According to the general method explained in Appendix A, 
it is enough to check this formula for a function of the form
$$F (z \ts z') =
\exp \!
\left({ - 2 \pi i \!
\left({ \left\langle{ \z, z }\right\rangle + \left\langle{ \z', z' }\right\rangle  }\right)
 }\right)
\, . \eqno ({\mathop{\rm II.55 }\nolimits})$$
But then our contention follows from the characterization (II.19) of the integrator, and the relations (II.51)
and (II.53).

\medskip
We combine this formula with equation (II.52) and obtain another form of Fubini's theorem
$$\int_{ {\bf Z}_{0, T+T'} }^{ }
\Da \om_s (u) \, G(u) =
\int_{ {\bf Z}_{0, T} }^{ }
\Da \om_s  (z)  
\int_{ {\bf Z}_{0, T'} }^{ }
\Da \om_s (z') \, G (z \ts z') \, . 
\eqno ({\mathop{\rm II.56 }\nolimits})$$

\medskip
We conclude the proof of equation (II.47). Indeed
$$
\eqalign{
\left({ U_T \! \left({ U_{T'} \phi }\right) }\right)
\left({x_0 }\right)
&= \int_{ {\bf Z}_{0, T} }^{ }
\Da \om_s (z) \left({ U_{T'} \phi }\right)
\left({x_0 \cdot \Si (T, z) }\right)\cr
&=
\int_{ {\bf Z}_{0, T} }^{ }
\Da \om_s (z)
\int_{ {\bf Z}_{0, T'} }^{ }
\Da \om_s (z')
\phi \! \left({ x_0 \cdot \Si (T, z) \cdot \Si (T', z') }\right)\cr
&=
\int_{ {\bf Z}_{0, T} }^{ }
\Da \om_s (z)
\int_{ {\bf Z}_{0, T'} }^{ }
\Da \om_s (z')
\phi \! \left({ x_0 \cdot \Si (T+ T' , z \ts z') }\right)\cr
&=
\int_{ {\bf Z}_{0, T + T'} }^{ }
\Da \om_s (u)
\phi \! \left({ x_0 \cdot \Si (T+ T' ,u) }\right)\cr
&
=
\left({ U_{T + T'} \phi }\right)   \left({x_0 }\right) \, .\cr}$$

\vfill\eject
\nn
3.2. {\it The differential equation.}

\smallskip
>From the group property (II.47) it follows that we need to establish the partial differential equation (II.2)
at the time $T = 0$, and that the general case will follow. That is, we want to prove
$$
\left({ U_T \phi }\right)
\left({ x_0 }\right) =
\phi \! \left({ x_0 }\right) + 
T \! \left({ {s \over 4 \pi} h^{\a \b} \La_{X_{(\a)}} \La_{X_{(\b)}} \phi 
\! \left({ x_0 }\right)
+ \La_Y \phi \! \left({ x_0 }\right) }\right)
 + o(T) 
\, . 
\eqno ({\mathop{\rm II.57 }\nolimits})$$
This relation implies also the initial condition
$$\build {\lim }_{ T = 0 }^{ }
\! \left({ U_T \phi }\right) \! \left({ x_0 }\right) =
\phi \! \left({ x_0 }\right)
\, . 
\eqno ({\mathop{\rm II.58 }\nolimits})$$

\medskip
For the proof, we rely on the scaling properties of paths as described in paragraph A.3.6. For $z$ in 
${\bf Z}_{0,1}$ we define the scaled path $z_T $ in ${\bf Z}_{0,T}$ by
$$z_T (t) = T^{ 1/2} z (t/T) 
\, . 
\eqno ({\mathop{\rm II.59 }\nolimits})$$
By scaling the differential equation (II.21), we get
$$dx (t/T) =
T^{1/2} X_{(\a)} (x (t/T)) 
dz_T^\a (t) + TY 
( x (t/T)) d (t/T) 
\, . 
\eqno ({\mathop{\rm II.60 }\nolimits})$$
Hence the transformation $\Si (T, z_T)$ in $N$ defined using the vector fields $X_{(\a)}$ and $Y$ is the same
as the transformation $\Si (1, z)$ using the vector fields $T^{1/2} X_{(\a)}$ and $TY$.
Moreover we can transfer the path integration from the variable space 
${\bf Z}_{0,T}$ to the fixed space ${\bf Z}_{0,1}$, that is 
$$\left({ U_T \phi }\right) \! \left({ x_0 }\right)  =
\int_{ {\bf Z}_{0,1} }^{ } \Da \om_s (z) \, \phi \! 
\left({ x_0 \cdot \Si \! \left({ T, z_T }\right) }\right)
\, . 
\eqno ({\mathop{\rm II.61 }\nolimits})$$
We use then a limited expansion of $\phi$ around $\phi (x_0)$, namely
$$
\eqalign{
\phi \! 
\left({ x_0 \cdot \Si \! \left({ T, z_T }\right) }\right)
&
= \phi (x_0) + T^{1/2} \La_{X_{(\a)}} \phi (x_0) z^\a (1) \cr
&
\quad +  {T \over 2}  \La_{X_{(\a)}}   \La_{X_{(\b)}} \phi (x_0) 
z^\a (1) z^\b (1) \cr
& 
\quad + T \La_Y \phi (x_0) + O \! \left({ T^{3/2} }\right)\, .\cr}
\eqno ({\mathop{\rm II.62 }\nolimits})$$

\smallskip
\nn
Recall the integration formulas
$$\int_{ {\bf Z}_{0,1} }^{ } \Da \om_s (z) \, z^\a (1) = 0
\eqno ({\mathop{\rm II.63 }\nolimits})$$
 and
$$
\int_{ {\bf Z}_{0,1} }^{ } \Da \om_s (z) \, z^\a (1) z^\b (1)  = sh^{\a \b }/2 \pi \, .
\eqno ({\mathop{\rm II.64 }\nolimits})$$
Collecting equations (II.61) to (II.64), we conclude the proof of (II.57).

\medskip
The previous calculation can be extended to give the complete Taylor expansion of 
$\left({ U_T \phi }\right) (x_0) $ around $T = 0$.
Since the vector fields $X_{(\a)}$ and $Y$ do not commute, we have to use a time-ordered exponential to express
the solution of the differential equation (II.21).

\bigskip
\nn
3.3. {\it About the construction of } $\Si (T, z)$.
\smallskip
For simplicity, assume $Y = 0 $ and $h_{\a \b} = \d_{\a \b }$.
Consider again a path $z : [0, T] \ra \RR^d $ of class $L^{2,1}$ with 
$z(0) = 0$.
For each $t$ in $[0, T]$, denote by $z^1 (t) , \cdots, z^d (t)$ the components of the vector $z(t)$ and define
the transformation 
${\bf s} (z(t)) $ on $N$ by
$$x_0 \cdot {\bf s} (z(t)) =
x_0 \cdot \s_{(1)} (z^1 (t)) \cdot \ldots \cdot \s_{(d)} (z^d (t)) \, . 
\eqno ({\mathop{\rm II.65 }\nolimits})$$
Here $\left\lbc{ \s_{(\a)} (r) }\right\rbc$ denotes the one-parameter group of transformations in $N$ with
generator $X_{(\a)}$.
In general, the vector fields $X_{(\a)}$ do not commute, hence the transformations 
${\bf s}(z(t))$ do not form a group.
In the general case, we can use a {\it multistep method} to solve the differential equation (II.21), hence
$$x_0 \cdot \Si (T, z) =
\build { \lim }_{n = \Inf }^{ }
x_0 \cdot {\bf s} (z(T/n)) \cdot {\bf s} (z (2T/n) - z(T/n))
\cdot \ldots \cdot {\bf s} (z(T) - z(T - T/n)) \, . 
\eqno ({\mathop{\rm II.66 }\nolimits})$$
Putting this into the integral (II.1), we obtain after some calculations the following variant of the Lie,
Trotter, Kato, Nelson formula
$$U_T \phi = \build { \lim }_{n = \Inf }^{ }
\left({ U_{T/n}^{(1)} \cdots U_{T/n}^{(d)}}\right)^n \phi \, . 
\eqno ({\mathop{\rm II.67 }\nolimits})$$

\nn
Here $U_T^{(\a)}$ for $\a$ in $\{ 1, \cdots, d \}$ corresponds to the one-dimensional case studied in paragraph
II.2, hence
$$ U_T^{(\a)} = \exp \! 
\left({ { sT\over 4 \pi } \La_{X_{(\a)}}^2 }\right)\, . 
\eqno ({\mathop{\rm II.68 }\nolimits})$$
These relations are in agreement with
$$ U_T = \exp \! 
\left({ { sT\over 4 \pi } \sum_{ \a }^{ }\La_{X_{(\a)}}^2 }\right)\, . 
\eqno ({\mathop{\rm II.69 }\nolimits})$$

\bigskip
\nn
{\bf 4. Some generalizations.}
\smallskip
\nn
4.1. {\it Including a potential.}
\smallskip

Consider the Schr\"odinger equation with potential
$${ \part \Psi\over \part t}
= {s \over 4 \pi } h^{\a \b }\, \La_{X_{(\a)}} \La_{X_{(\b)}}\Psi + V \Psi \, .
\eqno ({\mathop{\rm II.70 }\nolimits})$$

\nn
A path integral solution is obtained as follows
$$\Psi \! \left({ T, x_0 }\right)
=
\int_{ {\bf Z}_{0,T} }^{ } \Da_sz \cdot \exp \!  
\left({ -{\pi \over s} Q_0 (z) + \int_{{\bf T} }^{ }
dt \, V\! \left({ x_0 \cdot \Si (t, z) }\right)
 }\right) \phi \! \left({ x_0 \cdot \Si (T, z) }\right)
\,
.
\eqno ({\mathop{\rm II.71 }\nolimits})$$

\nn
The proof can be obtained by a slight generalization of the arguments presented in paragraph II.3.
We can also use the following trick:
we add one variable $\T$, considering the manifold $N \ts \RR$. The vector fields $X_{(\a)}$
in
 $N$ give vector fields, also denoted by $X_{\a}$, in 
$N \ts \RR$, which have a zero component on the factor $\RR$.
Moreover 
$Y = V(x) \part / \part \T$.
The function $\Psi (T, x) $ satisfies the equation (II.70) if, and only if, the function $\Psi (T, x)
\exp \T$ satisfies the equation (II.2) on $N \ts \RR$.
The differential system (II.21) is now written as
$$
\left\{\matrix{
dx = X_{(\a) } (x) d z^\a
  \hfill\cr
\noalign{\medskip}
    d \T = V(x) dt \, . 
\hfill \cr}
\right. 
\eqno ({\mathop{\rm II.72 }\nolimits})$$
Hence the transformation $\Si (T, z) $ takes $
\left({ x_0, \T_0 }\right)$ into 
$$\left({  x_0 \cdot \Si (T, z) , \ 
\T_0 + \int_{ {\bf T} }^{ } dt \  V 
(x_0 \cdot \Si (t, z))
}\right)\, ,$$
hence
$$
\phi \! \left({ (x_0, 0) \cdot \Si (T, z) }\right)
= \phi \! \left({ x_0 \cdot \Si (T, z) }\right)
\exp \! \left({ \int_{ {\bf T}}^{ } dt \, 
V \! \left({ x_0 \cdot \Si (t, z) }\right)
 }\right)
\, .\eqno ({\mathop{\rm II.73 }\nolimits})$$
Equation (II.71) follows easily from these remarks.
Notice that the exponent in this formula does not have a simple dynamical interpretation in general -- but
see our illustrations in section IV.

\bigskip
\nn
4.2. {\it Time-ordered product.}

\smallskip
Suppose that $F$ is a suitable functional on the space ${\bf Z}_{0,T}$.
We generalize equation (II.1) as follows
$$\left({ U_T^F \phi }\right)(x_0) =
\int_{ {\bf Z}_{0, T} }^{ } \Da \om_s (z) \, F(z) \, \phi  
\! \left({ x_0 \cdot \Si (T, z) }\right)
\, . \eqno ({\mathop{\rm II.74 }\nolimits})$$
By imitating the proof of (II.47), we get
$$U_{T+T'}^G = U_T^F U_{T'}^{F'}
\eqno ({\mathop{\rm II.75 }\nolimits})$$
for functionals $F$ on ${\bf Z}_{0, T}$ and $F'$ on ${\bf Z}_{0, T'}$, where the functional 
$G$ on 
${\bf Z}_{0, T+T'}$ is defined by
$$G( z \ts z') = F(z) F' (z') \,.
\eqno ({\mathop{\rm II.76 }\nolimits})$$

\vfill\eject
\vglue 1cm
\nn
\centerline {\titre III - Semiclassical Expansions}

\bigskip
\bigskip
We shall compute the semiclassical approximation of the expression in the basic equation
(II.1). The scalar potential $V$ is included via coupled equations rather than via an
additional variable, so that the equations will take readily their familiar forms; i.e. we
compute the wave function:
$$\eqalign{
\Psi \!
\left({ t_b, x_b  }\right)
= \int_{ {\bf Z }  }^{ }
&\Da_{s, Q_0 } z \cdot \exp \! 
\left({ - { \pi \over s } Q_0 (z)  }\right)  . \cr
&\exp \!
\left({ { 1 \over s\hbar  }
\int_{ t_a}^{ t_b } dt \, V \!
\left({ x_b \cdot \Si (t, z) }\right)
 }\right) \cdot 
\phi \left({ x_b \cdot \Si (t_a, z) }\right).\cr} 
\eqno ({\mathop{\rm III.1 }\nolimits} )
$$

\medskip
\nn
For ready use of our results in quantum physics, we compute $\Psi \!\left({ t_b, x_b  }\right)$.
Therefore in the general set up, we set $t_0 = t_b$, $x_0 = x_b$, hence ${\bf Z}$ is the {\it
space  ${\bf Z}_b$ of paths vanishing at} $t_b$ and $x(t_b, z) = x_b$. In paragraph 3, we give the
corresponding results for $\Psi \!\left({ t_b, x_a  }\right)$.

\medskip
We choose an initial wave function given by
$$\phi (x) = \exp \!
\left({-  { 1 \over s\hbar  } \Sa_0 (x)  }\right) \cdot 
\Ta (x), 
\eqno ({\mathop{\rm III.2 }\nolimits} )
$$
where $\Ta$ is a smooth function on $N$ with compact support, and 
$\Sa_0$ is an arbitrary, but reasonable function on $N$. 

\medskip
The initial wave function given by (III.2) generalizes plane waves on $\RR^d$, but is,
obviously, not a momentum eigenstate. We can nevertheless call the semiclassical
expansion of $\Psi \! \left({ t_b, x_b  }\right)$ with this initial wave function a
``{\it momentum-to-position  transition amplitude}'' for three reasons:

\smallskip
(i) In the limit $h = 0$, assuming $s = i$, the current density corresponding to the
initial wave function $\phi$ is
$$
\build { \lim }_{h = 0 }^{ }
{ \hbar \over 2 i  } 
\left({ \phi^* \bgt \phi - (\bgt \phi)^* \phi  }\right)
= 
\left\vert{ \Ta  }\right\vert^2 p\ , 
\quad 
{\mathop{\rm where  }\nolimits}
\ \ p(x) = \bgt \Sa_0 (x) \, .
\eqno ({\mathop{\rm III.3 }\nolimits} )
$$

\medskip
\nn
Consequently, $\Psi \! \left({ t_b, x_b  }\right)$ is the amplitude corresponding to the
transition from momentum $\bgt \Sa_0 (x)$ to position $x_b$.

\smallskip
(ii) We shall expand (III.1) around a ``classical path'' $z_{
{\mathop{\rm  cl }\nolimits}
 }$ characterized by initial momentum and final position.

\smallskip
(iii) The leading terms (henceforth labeled WKB) of semiclassical approximations combine
as if they were transitions between eigenstates, 
e.g.
$$
\eqalign{
\left\langle{  {\mathop{\rm position }\nolimits}
\bigm| {\mathop{\rm  position}\nolimits}  }\right\rangle_{
{\mathop{\rm  WKB}\nolimits}}
&= \int d \ {\mathop{\rm momentum  }\nolimits}
\left\langle{  {\mathop{\rm position }\nolimits}
\bigm| {\mathop{\rm  momentum }\nolimits}  }\right\rangle_{
{\mathop{\rm  WKB}\nolimits}}\cr
&\ \ \ \ \ \  \left\langle{  {\mathop{\rm momentum }\nolimits}
\bigm| {\mathop{\rm  position}\nolimits}  }\right\rangle_{
{\mathop{\rm  WKB}\nolimits}} \cr} $$
[see reference 15, \S 7 in the Appendix].

\medskip
Another useful initial wave function is
$$\phi (x) = \d_{x_a} (x) \, .
\eqno ({\mathop{\rm III.4 }\nolimits} )
$$
The corresponding expression 
 $\psi \!\left({ t_b, x_b  }\right)$
gives a {\it position-to-position transition} amplitude.
In general, this case is more complicated than the previous one, because now the initial
wave function restricts the domain of integration ${\bf Z}$.
Moreover the paths $z \in {\bf Z}$ such that $x \!
\left({t_a, z }\right) = x_a$ do not usually have a common origin.
Even at the very best, when $N = \RR^d$, and 
$z \left({t_a }\right) \simeq x_a - x_b$, one needs to be careful because the domain of
integration, say ${\bf Z}_{a, b} \sbs {\bf Z}$, is not a  vector space but an affine
space. Therefore, we shall compute position-to-position transitions  in section IV where
we specialize the basic equation (III.1). We shall treat in detail the case $N = \RR^d$
and give the necessary indications and references when $N = O (M)$, a frame bundle over a
Riemannian space $M^d$.

\medskip
Semiclassical expansions are best analyzed in the broader context of spaces of paths with
no requirement on their boundary values. For instance let $\Pa \RR^d$ be a space of paths
$z$ with no requirement on $z \left({t_a }\right) $, nor 
$z \left({t_b }\right) $, and let 
$\Pa_{a, b} \RR^d$ be the subspace of $\Pa \RR^d$ with $d$ requirements at 
$t_a$ and $d$ requirements at $t_b$. Let 
$\left\lbc{ S(\nu, z) }\right\rbc_{\nu}^{}$ be a one-parameter family of actions, and let
$U^{2d}(\nu)$ be the 2$d$-dimensional space of motions made of the critical points of
$S(\nu, z)$. A classical paths 
$z_{ {\mathop{\rm  cl }\nolimits} } \in \Pa_{a, b} \RR^d$ for the action $S(\nu, z)$ is at
the intersection of 
$\Pa_{a, b} \RR^d$ with $U^{2d} (\nu)$. If the intersection is transversal, no Jacobi field
of  $z_{ {\mathop{\rm  cl }\nolimits} }$ is in the tangent space to 
$\Pa_{a, b } \RR^d$. In this paper we assume that such is the case. Otherwise there are
caustics and we refer the reader to the literature [e.g. 2].

\bigskip
\nn
{\bf 1. General strategy.}
\smallskip
We introduce a Lagrangian
$$L(t, z) = {{}_1 \over {}^2} 
h_{\a \b } \dot z^\a (t) \dot z^\b (t) - V(x (t, z))
\eqno ({\mathop{\rm III.5 }\nolimits} )
$$
where as usual $x(t, z)$ is the solution of the differential equation (II.21) with
boundary condition $x \! \left({ t_b, z }\right) = x_b$. It is a function of $\dot z (t)$
and a functional of $z$. From this Lagrangian we deduce the {\it action functional}
$$S(z) = 
\int_{ {\bf T}J}^{ }
dt \, L(t, z) + \Sa_0 (x (t_a, z))
\eqno ({\mathop{\rm III.6 }\nolimits} )
$$
on the space ${\bf Z}_b = \Pa_0 \RR^d$ of paths $z$ obeying the boundary condition 
$z \!\left({ t_b }\right) = 0$.
Let 
$z_{ {\mathop{\rm  cl }\nolimits}} $
in ${\bf Z}_b$ be a critical point of the action functional $S$.

\medskip
As it is customary in the calculus of variations, we take a one-parameter variation
$$z(\l) = z_{ {\mathop{\rm  cl }\nolimits}}  + \l \z
\eqno ({\mathop{\rm III.7 }\nolimits} )
$$
with $\z$ in ${\bf Z}_b$
and the equation
$$
\left.{ {d\over d \l  } S (z (\l)) }\right\vert_{\l = 0} = 0
\eqno ({\mathop{\rm III.8 }\nolimits} )
$$
has to be satisfied for all $\z$. That yields, after integrating by parts, a functional
differential equation for 
$ z_{ {\mathop{\rm  cl }\nolimits} } (t_a)$.

\bigskip
Under the affine change of variable from $z $ to $\z$ given
by (III.7), we obtain
$$Q_0 (z) = Q_0 \!
\left({ z_{ {\mathop{\rm  cl }\nolimits}} + \l \z  }\right)
= Q_0 \!
\left({ z_{ {\mathop{\rm  cl }\nolimits}}  }\right)
+ 2 \l 
Q_0 \!\left({ z_{ {\mathop{\rm  cl }\nolimits}} , \z  }\right)
+ \l^2 Q_0 (\z) 
\eqno ({\mathop{\rm III.9 }\nolimits} )
$$
and
$$\Da_{s, Q_0} z = 
\Da_{s, \l^2 Q_0} \z .
\eqno ({\mathop{\rm III.10 }\nolimits} )
$$

\medskip
\nn
The expansion of $x(\cdot, z(\l))$ around 
$x(\cdot, z_{ {\mathop{\rm  cl }\nolimits}})$ reads 
$$P
\left({  z_{ {\mathop{\rm  cl }\nolimits}} + \l \z }\right)
=
P
\left({  z_{ {\mathop{\rm  cl }\nolimits}} }\right)
+ \l P' 
\left({  z_{ {\mathop{\rm  cl }\nolimits}} }\right)
\cdot \z +
{1\over 2}\l^2 P'' 
\left({  z_{ {\mathop{\rm  cl }\nolimits}}  }\right)
\cdot \z \z + O(\l^3).
\eqno ({\mathop{\rm III.11 }\nolimits} )$$

\medskip
\nn
Here 
$P' \left({  z_{ {\mathop{\rm  cl }\nolimits}} }\right)
$ and $P'' \left({  z_{ {\mathop{\rm  cl }\nolimits}} }\right)$
are the first and the second derivative mappings of $P$ at
$  z_{ {\mathop{\rm  cl }\nolimits}} $, where 
$P : {\bf Z}_b \ra \Pa_{x_b} N$ takes $z$ into 
$x (\cdot, z)$. We abbreviate $x(t,z_{{\mathop{\rm cl}\nolimits}})$ to $x_{ {\mathop{\rm cl }\nolimits}}(t)$. They are of the form
$$
\eqalign{
\left({ P'  \left({  z_{ {\mathop{\rm  cl }\nolimits}} }\right) \cdot \z  }\right)^\a (t) 
&= - \int_{t }^{t_b  } ds \, 
{ \d x_{ {\mathop{\rm  cl }\nolimits} }^\a  (t) \over
\d z_{ {\mathop{\rm  cl }\nolimits} }^\b  (s)  } \, 
\z^\b (s)\cr
&= - \int_{t }^{t_b  } ds \, 
k^\a {}_\b (t, s) \, \z^\b (s) = : \xi^\b (t)\cr}
\eqno ({\mathop{\rm III.12 }\nolimits} )$$
$$\eqalignno{
\left({ P''  \left({  z_{ {\mathop{\rm  cl }\nolimits}} }\right) \cdot \z \z  }\right)^\a
(t)  &= - \int_{t }^{t_b  } ds 
\left({  - \int_{t }^{t_b  } du }\right)
{ \d^2 x_{ {\mathop{\rm  cl }\nolimits} }^\a  (t) \over
\d z_{ {\mathop{\rm  cl }\nolimits} }^\b  (s)  
\d z_{ {\mathop{\rm  cl }\nolimits} }^\g  (u)
}
\z^\b (s) \z^\g (u) \cr
&=  \int_{t }^{t_b  } ds  \int_{t }^{t_b  } du \, k^\a 
 {}_{\b \g}  (t, s, u ) \z^\b (s) \z^\g (u) .
&({\mathop{\rm III.13 }\nolimits} )
 \cr} 
$$

\bigskip
We could obtain the short time propagator by expanding the quantity
$\phi \left({  x_b \cdot \Si (t, z) }\right)$ in equation (III.1) around 
$\phi \left({ x_b }\right)$ with $x_b \in N$.
Here, we shall expand 
$\phi \left({  x_b \cdot \Si (t, z) }\right)$ around
$\phi \left({  x_b \cdot \Si (t, z_{ {\mathop{\rm  cl }\nolimits} } ) }\right)
$
with $z_{ {\mathop{\rm  cl }\nolimits} } \in {\bf Z}_b$.
When we choose for $z_{ {\mathop{\rm  cl }\nolimits} }$  a critical point of the action, we
obtain the WKB approximation.
These are two different expansions, and in general, the short time propagator is
different from the WKB approximation [14]. But in some simple cases, they are equal [1].

\bigskip
\nn
{\it Remark.}
A simple and powerful argument of Stephen A. Fulling clarifies this issue: The Schr\"odinger operator 
$\exp (- it H / \hbar)$ can be written in terms of dimensionless terms 
${tH \over \hbar } = - {1 \over 2}  A \D + B V $, where 
$A = {\hbar t \over m} $ and $B = {\l t \over \hbar} $ with $\l$ a coupling constant. There are 
4 {\it different} possible expansions of physical interest:

-- expansion in $B$, equivalently expansion in $\l$;

-- expansion in $A$, equivalently expansion in $m^{- 1}$;

-- expansion in $A \cdot B$, equivalently expansion in $t$;

-- expansion in $A/B$, equivalently expansion in $\hbar$.

\bigskip
\nn
{\bf 2. Momentum-to-position transitions.}
\smallskip
In this paragraph we take $s = i$ and write 
$\Da_{s, Q_0}$ simply as $\Da_{Q_0}$. We expand the integrand in 
(III.1) by taking $z$ in the form
$z_{ {\mathop{\rm  cl }\nolimits} } + \z$ where 
$z_{ {\mathop{\rm  cl }\nolimits} }
$ is a critical point of the action functional $S$, 
 with boundary condition
$z_{ {\mathop{\rm  cl }\nolimits} }
\left({ t_b }\right) = 0$.
The initial wave function is given by (III.2).

\bigskip
The terms independent of $\z$ combine to make the action function 
$$\Sa \! \left({t_b, x_b   }\right) =
{{}_1 \over {}^2}
Q_0 \!
\left({ z_{ {\mathop{\rm  cl }\nolimits}}  }\right)
- \int_{t_a  }^{t_b  } dt \, 
V \!\left({ x_{ {\mathop{\rm  cl }\nolimits}} (t)  }\right)
+ \Sa_0 \!
\left({ x_{ {\mathop{\rm  cl }\nolimits}} (t_a)  }\right)
\eqno ({\mathop{\rm III.14 }\nolimits} )$$
where $x_{ {\mathop{\rm  cl }\nolimits}} (t) $ is by definition 
$x
\left({  t, z_{ {\mathop{\rm  cl }\nolimits}} .
}\right) $.
The terms linear in $\z$
$$
\eqalign{ \int_{\bf T }  dt & 
\left({  h_{\a \b } 
\dot z_{ {\mathop{\rm  cl }\nolimits} }^\a (t) 
\dot z_{ {\mathop{\rm  cl }\nolimits} }^\b (t) 
- \bgt_\a 
V \left({ x_{ {\mathop{\rm  cl }\nolimits}} (t) }\right)
\left({ P' \left({ z_{ {\mathop{\rm  cl }\nolimits}}   }\right) \cdot \z  }\right)^\a
(t)
}\right)
\cr 
& \ \ \ 
+ \bgt_\a \,  \Sa_0 
\left({ x_{ {\mathop{\rm  cl }\nolimits}} (t_a) }\right)
\left({ P' \left({ z_{ {\mathop{\rm  cl }\nolimits}} }\right) \cdot  \z  }\right)^\a 
(t_a)
\cr} \eqno ({\mathop{\rm III.15 }\nolimits} ) $$
vanish since $z_{ {\mathop{\rm  cl }\nolimits}}$ is a critical point of $S$ such that
$z_{ {\mathop{\rm  cl }\nolimits}} (t_b) = 0$. The  quadratic terms  can 
be written in the form
$\left({ Q_0 + Q }\right) (\z)$. Thanks to the equations  (A.127) and (A.132), the
corresponding integration is:
$$
\int_{{\bf Z} }^{ }
\Da_{Q_0 } \z \cdot
\exp \! \left({ -{\pi  \over s } \left({ Q_0 + Q }\right) (\z) }\right)
=
\Det \left({   Q_0 / (Q_0 + Q)  }\right)_{}^{1/2} \, . 
\eqno ({\mathop{\rm III.16 }\nolimits} )$$

\medskip
\nn
Thus the dominating term of the semiclassical expansion of 
$\Psi \!\left({ t_b, x_b  }\right)$ is:
$$\Psi_{ {\mathop{\rm  WKB }\nolimits}
 } \left({ t_b, x_b  }\right) =
\exp \!
\left({ {i \over \hbar} \Sa \! \left({ t_b, x_b  }\right) }\right) \cdot 
\Det \left({   Q_0 / (Q_0 + Q)  }\right)_{}^{1/2} \cdot
\Ta \!
\left({ x_{ {\mathop{\rm  cl }\nolimits}} (t_a) }\right) \, . 
\eqno ({\mathop{\rm III.17 }\nolimits} )$$

\medskip
We proceed to calculate 
$\Det (Q_0 / Q_\nu)$, with 
$Q_\nu = Q_0 + \nu Q$, in the case $\nu = 1$.
Physical arguments as well as previous calculations (see [6] and the appendix of [15])
suggest that $\Det (Q_0/Q_1)$ is equal to 
$\det K^\a {}_\b \left({ t_b, t_a  }\right)$, where the Jacobi matrix
$ K^\a {}_\b \left({ t_b, t_a  }\right)$ is defined by:
$$ K^\a {}_\b \left({ t_b, t_a  }\right)
= { \part z_{ {\mathop{\rm  cl }\nolimits}}^\a (t_b) \over
\part z_{ {\mathop{\rm  cl }\nolimits}}^\b (t_a) }\, .
\eqno ({\mathop{\rm III.18 }\nolimits} )$$

\medskip
\nn
Therefore we shall construct a one-parameter family of actions 
$S(\nu, z)$ and a one-parameter family of Jacobian matrices:
$$K^\a{}_\b 
\left({ \nu ; t_b, t_a }\right)
: = { \part z_{{\mathop{\rm  cl }\nolimits} }^\a  
\left({ \nu ; t_b  }\right)\over 
\part z_{{\mathop{\rm  cl }\nolimits} }^\b  
\left({ \nu ; t_a }\right) } \, .
\eqno ({\mathop{\rm III.19 }\nolimits} )$$

\medskip
\nn
We shall show that:
$$ \det \! \left({  K \!\left({ \nu ; t_b, t_a }\right)
\cdot  K\! \left({ 0 ; t_b, t_a }\right)^{- 1}  }\right) 
= : c (\nu) 
\eqno ({\mathop{\rm III.20 }\nolimits} )$$
satisfies the same boundary condition as 
$\Det \left({ Q_0 / Q_\nu }\right) $, namely
$$c(0) = \Det \left({ Q_0 / Q_0 }\right)  = 1, 
\eqno ({\mathop{\rm III.21 }\nolimits} )$$
and the same differential equation as 
$\Det \left({ Q_0 / Q_\nu  }\right) $, namely,
$${d  \over d \nu}  \ell n \, 
\Det \left({ Q_\nu / Q_0  }\right)   =
\Tr \! \left({ Q_\nu^{- 1}{ d \over d\nu } Q_\nu}\right)
= \Tr \! \left({ Q_\nu^{- 1} Q  }\right)
\eqno ({\mathop{\rm III.22 }\nolimits} )$$
(See equation (A.124) in Appendix A).

\bigskip
\nn
2.1. {\it  Differential equation satisfied by $c(\nu)$.}

\smallskip
We choose the family $S(\nu, z)$ such that its second variation evaluated at
 $z_{ {\mathop{\rm  cl }\nolimits} }$ be
$$S'' \left({ \nu, z_{ {\mathop{\rm  cl }\nolimits} } }\right)
 \cdot \z \z =
Q_0 (\z) + \nu Q (\z) = : Q_\nu (\z) \, .
\eqno ({\mathop{\rm III.23 }\nolimits} )$$

\medskip
\nn
It follows from its definition (III.18) that the $\b$-column 
$K_\b \left({ \nu ; t, t_a }\right)$ of the Jacobi matrix
$K \! \left({ \nu ; t, t_a }\right)$ is the following Jacobi field along the classical path
$z_{ {\mathop{\rm  cl }\nolimits} } (\nu, t)$ of the system governed by the action $S(\nu, z)$
$$Q_\nu K_{(\b)}  \! \left({ \nu ; t, t_a }\right) = 0
\eqno ({\mathop{\rm III.24 }\nolimits} )$$
$$
 K^{\a}_{\ (\b)} \! \left({ \nu ; t_a, t_a }\right) = \d_\b^\a \, .
\eqno ({\mathop{\rm III.25 }\nolimits} )$$

\medskip
\nn
Relation (III.24) encodes a functional differential equation for $K$ and an implicit
equation for 
$d K \! \left({ \nu ; t, t_a }\right) / dt$ at $t = t_a$.
It also gives, after differentiating w.r.t. $\nu$
$$Q_\nu {d \over d \nu } K(\nu) = - Q K (\nu)
\eqno ({\mathop{\rm III.26 }\nolimits} )$$
which can be solved with the {\it retarded Green's function}
$G_\nu^{ {\mathop{\rm  ret}\nolimits} }  $ of $Q_\nu$, namely:
$${d \over d \nu }
 K \!\left({ \nu ; t, t_a }\right)
= - \int_{ {\bf T } }^{ }
G_\nu^{ {\mathop{\rm  ret}\nolimits} } (t, s) Q(s) 
K \!\left({ \nu ; s, t_a }\right) ds \, .
\eqno ({\mathop{\rm III.27 }\nolimits} )$$

\medskip
The retarded Green's function can be expressed in terms of Jacobi fields, namely:
$$G_\nu^{ {\mathop{\rm  ret}\nolimits} } (t, s)
= \t (t - s) J (\nu ; t, s)
\eqno ({\mathop{\rm III.28 }\nolimits} )$$
where $\t$ is the step function, equal to $1$ for $t > s$ and $0$ for $t < s$; moreover
$J(\nu ; t, s)$ is a Jacobi matrix, {\it i.e.}  each column is a Jacobi field. The
$\a$-component of the $\b$-Jacobi field
$J^{(\b)} \left({ \nu ; t, t_a }\right)$ is
$$J^{\a (\b)}  \left({ \nu ; t, t_a }\right)
= {\part 
z_{ {\mathop{\rm  cl }\nolimits} }^\a (\nu, t) 
\over \part 
p_{ {\mathop{\rm  cl }\nolimits}, \b } \left({ \nu , t_a }\right)
}
\eqno ({\mathop{\rm III.29 }\nolimits} )$$
where
$
p_{ {\mathop{\rm  cl }\nolimits}  } (\nu, t) : =
\d \Sa (\nu, z) /
\d z_{ {\mathop{\rm  cl }\nolimits}} (\nu, t)$, here
$p_{ {\mathop{\rm  cl }\nolimits}, \a } (\nu, t) =
h_{\a \b } \dot z_{ {\mathop{\rm  cl }\nolimits}}^{\b}
(\nu, t).$

\medskip

>From (III.20), we get the equation for $c(\nu)$:
$${ d\over d \nu }J\ell n \, c (\nu) =
\tr \!
\left({ N  \! \left({ \nu ; t_a, t_b }\right) { d\over d \nu }J 
K  \! \left({ \nu ; t_b, t_a }\right) 
}\right)\eqno ({\mathop{\rm III.30 }\nolimits} )$$
where the matrix $N$ is the inverse of $K$
$$N \! \left({ \nu ; t_a, t_b }\right)  
K \! \left({ \nu ; t_b, t_a }\right) 
= \un . 
\eqno ({\mathop{\rm III.31 }\nolimits} )$$

\medskip
\nn
Using (III.27) and (III.28), we obtain the sought-for differential equation:
$$
{ d\over d \nu }J\ell n \  c(\nu) = - \tr \!
\left({ \int_{ {\bf T}J}^{ }
d s \, J \!
\left({ \nu ; t_b, s }\right) 
Q(s) K \! \left({ \nu ; s, t_a }\right) 
N \! \left({ \nu ;  t_a, t_b }\right) }\right)
 \, . 
\eqno ({\mathop{\rm III.32 }\nolimits} )$$

\medskip
Within focal distance of 
$z_{ {\mathop{\rm  cl }\nolimits} }
\left({ \nu ; t_a }\right)$, the $d$ Jacobi fields
$K_{(\b)} (\nu) $ and the $d$ Jacobi fields 
$J^{(\b)} (\nu)$ are linearly independent, and 
$c(\nu)$ is defined. It gives the rate at which the flow of classical paths
$
\left\lbc{ z_{ {\mathop{\rm  cl }\nolimits}} (\nu) }\right\rbc$
diverges or converges. 

\vfill\eject
\nn
 2.2. {\it Comparison between the two differential equations} (III.32) {\it and }
(III.22).

\smallskip
We shall express $Q_\nu^{- 1}$ in terms of Jacobi fields, as this is always possible
within focal distance of 
$
z_{ {\mathop{\rm  cl }\nolimits}} \left({ t_a }\right)$.
The unique inverse $G_\nu $ of $Q_\nu$ satisfies the equation 
$$Q_\nu G_\nu = \un
\eqno ({\mathop{\rm III.33 }\nolimits} )$$
with 
$$G_\nu \!\left({  t_b, s }\right) = 0 
\eqno ({\mathop{\rm III.34 }\nolimits} )$$
by virtue of (A.40) together with the specialization used in (A.62) when 
applied to the space ${\bf Z}_b$ of paths
vanishing at $t_b$. Therefore (see formula (B.26) in Appendix B)
$$
\eqalign{
G_\nu (t, s) &= \t (s - t) K \!
\left({ \nu ; t, t_a }\right) N \!\left({ \nu ; t_a, t_b }\right)  
J \!\left({ \nu ; t_b, s }\right) \cr
&\ \ \ - \t(t - s) J \!\left({ \nu ; t, t_b }\right) 
\ww N \!\left({ \nu ; t_b, t_a }\right) 
\ww K \!\left({ \nu ; t_a, s }\right) \cr} 
\eqno ({\mathop{\rm III.35 }\nolimits} )$$
where $\ww K$ is the transpose of $K$, and $\ww N$ is $\ww K$'s inverse.

\medskip
Substituting $Q_\nu^{- 1} = G_\nu $ in (III.22) gives:
$${ d\over d \nu }J\ell n \Det  \! \left({ Q_\nu / Q_0 }\right) 
= \tr \!
\left({ \int_{{\bf T}} dt   \,  K \!\left({ \nu ; t, t_a }\right) 
N \!\left({ \nu ; t_a, t_b }\right)  
J \! \left({ \nu ; t_b, t }\right)  Q(t) }\right)
 \, . 
\eqno ({\mathop{\rm III.36 }\nolimits} )$$

\medskip
\nn
Comparison with (III.32) shows that
$$\ell n \Det \left({  Q_\nu / Q_0 }\right)  = - \ell n \, c (\nu) 
\eqno ({\mathop{\rm III.37 }\nolimits} )$$
i.e.\footnote{$^{10}$}{
The credit for the technique used in deriving (III.38) is due to B. Nelson and B. Sheeks [6].
It is made simpler and more general here by not integrating $Q_\nu$ by parts. The first
calculation giving the ratio of functional determinants in terms of finite determinants
can be found in [16, 5].}
$$
\Det \left({ Q_0 / Q_\nu }\right)   =
\det \!
\left({  K \!\left({ \nu ; t_b, t_a }\right)  \cdot  K \!\left({ 0 ; t_b, t_a }\right)^{-1}
}\right)
  \, .
 \eqno ({\mathop{\rm III.38 }\nolimits} )$$

\medskip
\nn
Here $ K \!\left({ 0 ; t, t_a }\right)$ satisfies
$$Q_0  K \!\left({ 0 ; t, t_a }\right) = 0 
\eqno ({\mathop{\rm III.39 }\nolimits} )$$
$$  K \!\left({ 0 ; t_a, t_a }\right) = \un 
\eqno ({\mathop{\rm III.40 }\nolimits} )$$
where
$$Q_0 (z) = - \ds \int_{{\bf T} } dt \, h_{\a  \b }
\ddot z^{\a} (t) z^\b (t)  - h_{\a  \b }
\dot z^\a \! \left({ t_a }\right) z^\b \! \left({ t_a }\right)\, .$$

\medskip
\nn
Hence
 we have
$$ {d^2 \over dt^2} K \! \left({ 0 ; t, t_a }\right) = 0 \ \  , \ \ 
\left.{ 
{d \over dt } K \! \left({ 0 ; t, t_a }\right)
 }\right\vert_{t = t_a}
= 0 
\eqno ({\mathop{\rm III.41 }\nolimits} )$$
therefore
$$K^\a {}_\b \left({ 0 ; t, t_a }\right) = \d^\a{}_\b \, .
\eqno ({\mathop{\rm III.42 }\nolimits} )$$

\medskip
Inserting 
$\Det (Q_0 / Q_\nu) = \det K^\a {}_\b \left({ t_b, t_a }\right)$ 
into (III.17) we obtain:
$$\Psi_{
{\mathop{\rm  WKB}\nolimits}
} \left({ t_b, x_b }\right) =
\left({  \build {\det  }_{ \a, \b}^{ } 
{ \part    z_{ {\mathop{\rm  cl }\nolimits}}^\a   \left({ t_b }\right) \over 
\part    z_{ {\mathop{\rm  cl }\nolimits}}^\b   \left({ t_a }\right)
} 
}\right)_{}^{1/2} \cdot 
\exp \!
\left({ { i  \over \hbar }  \Sa \!\left({ t_b, x_b }\right) }\right)
\cdot \Ta \! \left({ x_{ {\mathop{\rm  cl }\nolimits}} (t_a) }\right) \, .
\eqno ({\mathop{\rm III.43 }\nolimits} )$$

\smallskip
\nn
If $ z_{ {\mathop{\rm  cl }\nolimits}}  \left({ t_b }\right)$
is conjugate to 
$z_{ {\mathop{\rm  cl }\nolimits}}  \left({ t_a}\right)$, in the sense of caustic theory,
one needs [2] to include terms in $\z$ of order higher than $2$ in the calculation of $\Psi$.

\bigskip
\nn
2.3. {\it  End of calculation.}

\smallskip
Equation (III.43) is not the end of the calculation;
$\Psi_{
{\mathop{\rm  WKB}\nolimits}
} $ must be expressed in terms of 
$x_{ {\mathop{\rm  cl }\nolimits}}
: T \ra N$, where
$$x_{ {\mathop{\rm  cl }\nolimits}} (t) =  x \!
\left({ t,  z_{ {\mathop{\rm  cl }\nolimits}} }\right)
= x_b \cdot \Si \left({ t, z_{ {\mathop{\rm  cl }\nolimits}} }\right)
\eqno ({\mathop{\rm III.44 }\nolimits} )$$
but not in terms of 
$z_{ {\mathop{\rm  cl }\nolimits}} : {\bf T} \ra \RR^d$.
Two techniques present themselves: we know the evolution equation satisfied by the wave
function $\Psi$ on $N$, therefore we can construct an action on $N$ and find its critical
points. But we do not need to find 
$x_{ {\mathop{\rm  cl }\nolimits}} $
 corresponding to $z_{ {\mathop{\rm  cl }\nolimits}}$, we only need the expression
corresponding to the prefactor of (III.43). This can be obtained by a simpler technique.

\medskip
Recall the parametrization $P : \Pa_0 \RR^d \ra \Pa_{x_b} N$ given by 
$P(z) = x$. Corresponding to the action functional 
$S $ on $\Pa_0 \RR^d$ we get a functional $\oo S $ on 
$\Pa_{x_b} N $ such that $S = \oo S \circ P$, i.e. 
$\oo S (x) = S(z)$. Making the substitution $z =
z_{ {\mathop{\rm  cl }\nolimits}} + \z$ and expanding up to second order terms in $\z$,
we get
$$\oo S \!\left({ x_{ {\mathop{\rm  cl }\nolimits}} }\right)
= S \!\left({ z_{ {\mathop{\rm  cl }\nolimits}} }\right)
\eqno ({\mathop{\rm III.45 }\nolimits} )$$
$$\oo S' \! \left({ x_{ {\mathop{\rm  cl }\nolimits}} }\right) \cdot 
{ \d x_{ {\mathop{\rm  cl }\nolimits}} \over 
\d z_{ {\mathop{\rm  cl }\nolimits}} } = 
S' \! \left({ z_{ {\mathop{\rm  cl }\nolimits}} }\right)
\eqno ({\mathop{\rm III.46 }\nolimits} )$$
$$
\oo S'' \! \left({ x_{ {\mathop{\rm  cl }\nolimits}} }\right) \cdot 
{ \d x_{ {\mathop{\rm  cl }\nolimits}} \over 
\d z_{ {\mathop{\rm  cl }\nolimits}} } 
{ \d x_{ {\mathop{\rm  cl }\nolimits}} \over 
\d z_{ {\mathop{\rm  cl }\nolimits}} } 
= 
S'' \! \left({ z_{ {\mathop{\rm  cl }\nolimits}} }\right). 
\eqno ({\mathop{\rm III.47 }\nolimits} )$$

\smallskip
\nn
In the beginning of this paragraph we calculated the expansion of
$S \! \left({  z_{ {\mathop{\rm  cl }\nolimits}} + \z }\right)
$ up to second order and obtained the expression
$$S'' \!  \left({ z_{ {\mathop{\rm  cl }\nolimits}} }\right) 
\cdot \z \z =
Q_0 (\z) + Q(\z) = Q_1 (\z)
\eqno ({\mathop{\rm III.48 }\nolimits} )$$
for the second order variation. We also know that the infinite-dimensional determinant 
$\Det (Q_0 /Q_1)$ is equal to the determinant of the Jacobi matrix
$\part z_{ {\mathop{\rm  cl }\nolimits}} 
\left({ t_b  }\right) /
\part z_{ {\mathop{\rm  cl }\nolimits}} 
\left({ t_a  }\right)$.

\medskip
Similarly we write
$$\oo S'' \! 
\left({ x_{ {\mathop{\rm  cl }\nolimits}}  }\right) \cdot 
{ \d x_{ {\mathop{\rm  cl }\nolimits}} \over 
\d z_{ {\mathop{\rm  cl }\nolimits}} } \z 
{ \d x_{ {\mathop{\rm  cl }\nolimits}} \over 
\d z_{ {\mathop{\rm  cl }\nolimits}} } 
\z  \, = \, 
\oo Q_0 (\z) + \oo Q(\z) 
\eqno ({\mathop{\rm III.49 }\nolimits} )$$
and identify $\oo Q_0$ as follows. Let
$$S_0 (z) : = { {}_{1} \over {}^{2} } Q_0 (z)
\eqno ({\mathop{\rm III.50 }\nolimits} )$$
and
$$\oo S_0 (x(z)) : = S_0 (z)
\eqno ({\mathop{\rm III.51 }\nolimits} )$$
then
$$\oo Q_0 (\z) : = 
\oo S''_0 
\left({ x_{ {\mathop{\rm  cl }\nolimits}}   }\right) \cdot 
{ \d x_{ {\mathop{\rm  cl }\nolimits}} \over 
\d z_{ {\mathop{\rm  cl }\nolimits}} } \z
{ \d x_{ {\mathop{\rm  cl }\nolimits}} \over 
\d z_{ {\mathop{\rm  cl }\nolimits}} } 
\z 
\eqno ({\mathop{\rm III.52 }\nolimits} )$$
$$ \det 
{  \part x_{ {\mathop{\rm  cl }\nolimits}}
\left({ t_b  }\right)
 \over \part 
 x_{ {\mathop{\rm  cl }\nolimits}} \left({ t_a  }\right)  } 
= \Det \!
\left({ \oo Q_0 / \! \left({ \oo Q_0 + \oo Q }\right)  }\right)
=\Det \!
\left({  Q_0 / \! \left({  Q_0 +   Q }\right)  }\right)
  \,  .
 \eqno ({\mathop{\rm III.53 }\nolimits} )$$

\smallskip
\nn
Therefore
$$\det 
{  \part x_{ {\mathop{\rm  cl }\nolimits}}
\left({ t_b  }\right)
 \over \part 
 x_{ {\mathop{\rm  cl }\nolimits}} \left({ t_a  }\right)  } 
=  
\det 
{  \part z_{ {\mathop{\rm  cl }\nolimits}}
\left({ t_b  }\right)
 \over \part 
 z_{ {\mathop{\rm  cl }\nolimits}} \left({ t_a  }\right)  } \, .
\eqno ({\mathop{\rm III.54 }\nolimits} )$$

\medskip
In conclusion, we obtain the sought-for semiclassical expansion
$$
\Psi_{
{\mathop{\rm  WKB}\nolimits}
} 
\left({ t_b, x_b }\right)
= \det  \! \left({
{ \part    x_{ {\mathop{\rm  cl }\nolimits}}   \left({ t_b }\right) \over 
\part    x_{ {\mathop{\rm  cl }\nolimits}}   \left({ t_a }\right)
} 
}\right)_{}^{1/2} \cdot \, 
\exp \!
\left({ { i  \over \hbar }  \Sa \left({ t_b, x_b }\right) }\right) \cdot 
\Ta \! \left({ x_{ {\mathop{\rm  cl }\nolimits}} (t_a) }\right)  .
\eqno ({\mathop{\rm III.55 }\nolimits} )$$

\smallskip
\nn
According to equation (III.46), the path 
$ x_{ {\mathop{\rm  cl }\nolimits}}$ in $N$ is a critical point for the action functional
$\oo S : \Pa_{x_b} N \ra \RR$. By construction, we have
$$\oo S' (x (\cdot , z) ) = S' (z)\, .
\eqno ({\mathop{\rm III.56 }\nolimits} )$$

\smallskip
\nn
Using equation (III.5) and (III.6) this can be made more explicit as follows
$$\oo S (x) = { {}_{1} \over {}^{2} }
\int_{ {\bf T} }^{}
dt \, h_{\a \b } \, 
\dot z^\a (t) \dot z^\b (t) - \int_{ {\bf T} }^{}
V(x (t)) dt + \Sa_0 \! \left({ x (t_a) }\right) \, .
\eqno ({\mathop{\rm III.57 }\nolimits} )$$

\smallskip
\nn
In our general setup
, the first integral remains somewhat implicit, but can be used for
practical calculations in the applications given in section IV.

\medskip
The prefactor in (III.55) also gives the volume expansion or contraction of a congruence
of classical paths originating in the neighborhood of 
$\Sa_0 \! \left({ x_{ {\mathop{\rm  cl }\nolimits}} (t_a) }\right)$ with momentum
$\bgt \Sa_0 (m)$, $m \in N$. This determinant depends both on the choice of the initial
wave function and the dynamics of the system.

\medskip
Set $d \om_a$ the volume element on $N$ at 
$x_{ {\mathop{\rm  cl }\nolimits}} (t_a) $ and 
$d \om_b = \det \!
\left({ { \part x_{ {\mathop{\rm  cl }\nolimits}}^\a (t_b) \over
\part x_{ {\mathop{\rm  cl }\nolimits}}^\b (t_a)  } }\right)
d \om_a$; then 
$$
\build {\lim }_{ h = 0 }^{ }
\int_{ C_{t_b } \O }^{ }
\left\vert{  \Psi \! \left({  t_b, x_b }\right) }\right\vert^2 d \om_b
= \int_{ \O }^{ }
\left\vert{ \phi \left({   x_a }\right) }\right\vert^2 d \om_a
$$
where $C_t$ belongs to the group of transformations generated by the classical flow
(see details in [5 p. 299]).

\bigskip
\nn
{\bf 3. Diffusion problems.}

\smallskip
In a diffusion problem, one is interested in computing the functional integral
$$
\eqalign{
\Psi \!\left({  t_b, x_a }\right) =
&
\int_{{\bf Z}  }^{ } \Da_{s, Q_0 }   z \exp \!
\left({- {\pi \over s} Q_0 (z) }\right) \cdot\cr
& \  \ \ 
\exp \!
\left({ {1  \over s\hbar } \int_{t_a }^{t_b } dt \, 
V \!\left({ x_a \cdot \Si (t, z)  }\right)
 }\right)\cdot \phi \left({ x_a \cdot \Si (t_b, z) }\right) \cr} 
\eqno ({\mathop{\rm III.58 }\nolimits} )$$
where ${\bf Z}$ denotes now the space of paths $z$ vanishing at $t = t_a$.
If we choose $\phi$ to be of the form (III.2), namely
$$
\phi \left({ x_a \cdot \Si (t_b, z) }\right)
=
\phi \left({ x (t_b, z) }\right) =
\exp \!
\left({- {1  \over s\hbar } \Sa_0 \left({ x (t_b, z) }\right)  }\right)
\Ta \!\left({ x (t_b, z) }\right) 
\eqno ({\mathop{\rm III.59 }\nolimits} )$$
then computing (III.58), with $s = i$, can be said to be computing the probability
amplitude of a transition from a position $x_a$ to a  momentum of the form
$$p \left({ x_b }\right) =
\bgt \Sa_0  \left({ x_b }\right)  
\eqno ({\mathop{\rm III.60 }\nolimits} )$$
for the end point 
$x_b = x \! \left({ t_b, z }\right)$ of the path $x$.

\medskip
With $\phi$ given by (III.59), one obviously expects the WKB approximation of (III.58) to
be $$\Psi_{
{\mathop{\rm  WKB}\nolimits}
} 
\left({ t_b, x_a }\right)
= \det \! \left({
{ \part    x_{ {\mathop{\rm  cl }\nolimits}}^\a   \left({ t_a }\right) \over 
\part    x_{ {\mathop{\rm  cl }\nolimits}}^\b   \left({ t_b }\right)
} 
}\right)_{}^{1/2}
\exp \! 
\left({ {i  \over \hbar }  \Sa \! \left({ t_b, x_a }\right) }\right) \cdot 
\Ta \!\left({ x_{ {\mathop{\rm  cl }\nolimits}} (t_b) }\right) \, .
\eqno ({\mathop{\rm III.61 }\nolimits} )$$

\smallskip
\nn
Nevertheless, it is gratifying to derive (III.61) by the method followed in paragraph 2.
The only necessary changes are described as follows:
$$
c(\nu) = \det \!
\left({ K\!
\left({ \nu ; t_a , t_b  }\right)
\cdot 
K \!
\left({ 0 ; t_a , t_b  }\right)^{- 1}  }\right)
\eqno ({\mathop{\rm III.20 }\nolimits}^{{\mathop{\rm bis }\nolimits} }  )$$
$$
G_{\nu}^{ {\mathop{\rm ad }\nolimits}  \nu}
(t, s) = \t (t - s) J (\nu ; s, t) =
- \t (s - t) J (\nu ; t, s) 
\eqno ({\mathop{\rm III.28 }\nolimits}^{{\mathop{\rm bis }\nolimits} }  )$$
$$
\eqalign{
G_\nu (t, s) &
= \t (s - t) J \!
\left({\nu ; t,t_a  }\right) \ww N \! \left({\nu ;  t_a, t_b }\right)
\ww K \! \left({\nu ;  t_b, s }\right) \cr
& \ \ \ -
\t  (t - s) K \! \left({ \nu ; t, t_b }\right) N \! \left({\nu ;  t_b, t_a }\right)
J \! \left({\nu ;  t_a, s }\right)\cr} 
\eqno ({\mathop{\rm III.35 }\nolimits}^{{\mathop{\rm bis }\nolimits} }  )$$
\nn
(see formula (B.27) in Appendix B).

\vfill\eject
\vglue 1cm
\nn
\centerline {\titre IV - Illustrations}

\bigskip
\bigskip
The specializations presented in this section illustrate and develop the general
 formulas derived in sections II and
III. They differ by the choice of the manifold $N$ and of the initial wave function
 $\phi$. We treat the case of scalar wave functions, but the case of
a tensor field is easily accommodated by using the Lie derivative of tensor fields.

\medskip
We are basically interested in the Schr\"odinger equation. So the reader should substitute
$i$ to $s$ in the formulas. Moreover, we want to evaluate the value of the wave function
 $\Psi$ at the
{\it final time} $t_b$, and {\it final position } $x_b$. This dictates the choice of the 
space ${\bf Z}_b$
of paths, determined by $z(t_b) = 0$.

\bigskip
\nn
{\bf 1. Point-to-point transitions in a flat space.}

\smallskip
The basic manifold $N$ is the euclidean space $\RR^d$ of dimension $d$, in cartesian 
coordinates (see
paragraph IV.2 for other coordinates). For a vector $x$, with coordinates $x^1, \cdots, x^d$,
 we set its
length to be
$\left\vert{ x }\right\vert =
\left({ \sum_{ \a = 1}^{d } (x^\a)^2}\right)_{}^{1/2}$ as usual.
We consider a particle of mass $m$ moving in the field of a potential $V$. The Lagrangian
 and action are
given as usual by
$$
\eqalignno{
L(x, \dot x) &= 
{m \over 2} \left\vert{ \dot x }\right\vert_{}^{2} - V(x), 
& ({\mathop{\rm IV.1 }\nolimits})\cr
S(x (\cdot)) &= \int_{ {\bf T} }^{ } dt \, 
L(x(t), \dot x(t)) \, . 
& ({\mathop{\rm IV.2 }\nolimits}) \cr }$$

\medskip
In classical mechanics, we solve the equations of motion with suitable boundary conditions:
$$m \! \build { x}_{ }^{  \cdot \cdot }_{
}\!\!{}_{{\mathop{\rm cl }\nolimits}} = - \bgt V \!\left({ x_{{\mathop{\rm cl }\nolimits}}  }\right)
\eqno ({\mathop{\rm IV.3 }\nolimits})$$
$$
x_{{\mathop{\rm cl }\nolimits}} (t_a) = x_a \quad ,
\quad x_{{\mathop{\rm cl }\nolimits}} (t_b) = x_b \, , 
\eqno ({\mathop{\rm IV.4 }\nolimits})$$
where $x_a, x_b$ are points in $\RR^d$. In quantum mechanics, we want to solve the Schr\"odinger
equation
$$i \hbar { \part \Psi \over \part t} 
= - { \hbar^2 \over 2m } \D \Psi + V \Psi
\eqno ({\mathop{\rm IV.5 }\nolimits})$$
with initial condition
$$\Psi \! \left({ t_a, x }\right)
= \phi (x) \, . 
\eqno ({\mathop{\rm IV.6 }\nolimits})$$

\bigskip
\nn
The {\it point-to-point transition amplitudes}  are given by
$$
\left\langle{  t_b ,  x_b \bigm|  t_a , x_a }\right\rangle
=
\Psi \!
\left({ t_b, x_b }\right) \, ,
\eqno ({\mathop{\rm IV.7 }\nolimits})$$
where $\phi (x)$ is  a delta function $\d \!\left({ x - x_a }\right)$.

\bigskip
The original claim of Feynman was that we can solve the Schr\"odinger equation by the following
path integral
$$\Psi \! \left({ t_b, x_b }\right) =
\int_{ \Pa_b }^{ }
\Da x \cdot e^{i S (x) / \hbar} 
\cdot \phi (x (t_a)) \, ,
\eqno ({\mathop{\rm IV.8 }\nolimits})$$
where $\Pa_b$ is the space of all paths $x : {\bf T} \ra \RR^d$ with endpoint at $x_b$, namely 
$x(t_b) = x_b$. The questionable part was the rigorous definition of the integrator 
$\Da x$.

\bigskip
To fit within our general framework, we consider translations acting on $\RR^d$, namely
$$x \cdot 
\s_{(\a) }(r) =
\left({ x^1, \cdots, x^{\a - 1}, x^\a + \l r , x^{\a + 1} , \cdots , x^d }\right) \, . 
\eqno ({\mathop{\rm IV.9 }\nolimits})$$

\bigskip
\nn
The corresponding Lie derivative is given by
$$\La_{X_{(\a )} } f = 
\l {\part f \over \part x^\a }
\eqno ({\mathop{\rm IV.10 }\nolimits})$$
and the parameter $\l$ is chosen equal to $(h/m)^{{1 \over 2}}$.
The general differential equation $dx = X_{(\a)} (x) \cdot d z^\a$ reduces to 
$d x^\a = \l d z^\a$ (for $\a$ in $\{1, \cdots, d \}$). Hence the solution
$$x(t, z) = x_b + \l z(t) 
\eqno ({\mathop{\rm IV.11 }\nolimits})$$
describes the parametrization of the space $\Pa_b$ of paths $x$ with $x(t_b) = x_b$ by the space
${\bf Z}_b$ of paths $z$ with $z(t_b) = 0$. Accordingly, we obtain for the action
$$- {1 \over s \hbar } 
S (x(\cdot, z)) = - 
{\pi \over s} \int_{ {\bf T} }^{ } dt \, 
\left\vert{ \dot z (t) }\right\vert_{}^{2} +
{1 \over s \hbar }  \int_{ {\bf T} }^{ } dt \, V (x(t, z)) \, . 
\eqno ({\mathop{\rm IV.12 }\nolimits})$$

\medskip
\nn
Substituting $z$ for the integration variable $x$ in equation (IV.8), 
we obtain the path integral (remember that $s = i$ in quantum mechanics)
$$\Psi \! \left({ t_b, x_b }\right) =
\int_{ {\bf Z}_b }^{ } \Da z \cdot \exp \!
\left({- { 1\over s \hbar }  S (x(\cdot, z)) }\right)
\cdot \phi (x (t_a, z)) \, . 
\eqno ({\mathop{\rm IV.13 }\nolimits})$$

\medskip
\nn
>From our general results in section II, the function $\Psi$ is a solution of the differential equation
$$
 {\part \Psi \over \part t} =
 {s \over 4 \pi } \sum_{ \a }^{ }
\La_{X_{(\a)} }^2 \Psi + { 1\over s \hbar }  V \Psi \, , 
\eqno ({\mathop{\rm IV.14 }\nolimits})$$
that is 
$$s \hbar  {\part \Psi \over \part t} =
{s^2 \hbar^2 \over 2m} \D \Psi + V \Psi \, . 
\eqno ({\mathop{\rm IV.15 }\nolimits})$$

\medskip
\nn
For $s = i$, this is the sought-for Schr\"odinger equation. The integrator $\Da z$ in the space
${\bf Z}_b$ is invariant under translations and is normalized by\footnote{$^{11}$}
{\nn The scaling factor $\l = (h/m)^{1/2}$ has been chosen in such a way that no physical constant enters
in this normalization. With the notations of paragraph A.3.6, we have the dimensional equations
$$
\left\lbk{ z^\a }\right\rbk = \Ta^{1/2}
\quad , \quad [t] = \Ta \, .$$
}
$$ 
\int_{ {\bf Z}_b }^{ } \Da z \cdot \exp \!
\left({ - { \pi \over s} \int_{ {\bf T} }^{ } dt 
\left\vert{ \dot z (t)  }\right\vert_{}^{2}
 }\right) = 1
\, . 
\eqno ({\mathop{\rm IV.16 }\nolimits})$$

\bigskip
We derive now the {\it semiclassical approximation}\footnote{$^{12}$}{
\nn Since $\l = (h/m)^{1/2}$, the semiclassical expansion will proceed according to the powers of
$h^{1/2}$. }.
We use the initial wave function 
$\phi (x) = \d \left({ x - x_a }\right)$ and reparametrize the paths {\it around the classical
path} $x_{{\mathop{\rm  cl }\nolimits}}$, that is
$$x(t, \z)  = x_{{\mathop{\rm  cl }\nolimits}} (t) + \l \z (t) \, , 
\eqno ({\mathop{\rm IV.17 }\nolimits})$$
with $\z$ in ${\bf Z}_b$.
Since the integrator in ${\bf Z}_b$ is invariant under translations, and 
$x_{{\mathop{\rm  cl }\nolimits}} (t_a) = x_a $, we transform equation (IV.13) into
$$\Psi \! \left({ t_b, x_b }\right) =
\int_{ {\bf Z}_b }^{ } \Da \z \cdot \exp \!
\left({ - { 1 \over s \hbar} 
S \! \left({ x_{{\mathop{\rm  cl }\nolimits}} + \l \z }\right)
 }\right) \d (\l \z (t_a))
\, . 
\eqno ({\mathop{\rm IV.18 }\nolimits})$$

\medskip
\nn
In expanding $S \! \left({ x_{{\mathop{\rm  cl }\nolimits}} + \l \z }\right)$ in powers of $\l$, there
is no term linear in $\l$, since $x_{{\mathop{\rm  cl }\nolimits}}$ is a critical point 
of the action
functional $S$ {\it and} because $\z (t_a) = \z(t_b) = 0$.
Hence we obtain
$$
\eqalignno{
- { 1 \over s \hbar} 
S \! \left({ x_{{\mathop{\rm  cl }\nolimits}} + \l \z }\right)
&=-
{ 1 \over s \hbar} 
S \! \left({ x_{{\mathop{\rm  cl }\nolimits}}  }\right)
- {\pi  \over s } \int_{ {\bf T}}^{ } 
dt \ \vert  \dot \z (t) \vert_{}^{2}
& ({\mathop{\rm IV.19 }\nolimits}) \cr
& \ \ \  + {\pi  \over ms  } \int_{ {\bf T}}^{ }  dt \, 
\bgt_\a \bgt_\b V \! 
\left({ x_{{\mathop{\rm  cl }\nolimits}} (t)  }\right)
\z^\a (t) \z^\b (t) + O(\hbar^{1/2} ) \, . \cr}
$$

\medskip
\nn
The action 
$S \! \left({ x_{{\mathop{\rm  cl }\nolimits}} }\right)$ corresponding to the classical path
$x_{{\mathop{\rm  cl }\nolimits}}$ with endpoints 
$x_{{\mathop{\rm  cl }\nolimits}} (t_a) = x_a$, 
$x_{{\mathop{\rm  cl }\nolimits}}(t_b) = x_b$ is nothing else than the 
{\it classical action function }
$\Sa \! \left({ t_b, x_b ; t_a, x_a }\right)$.

\bigskip 
Omitting terms of order $\hbar^{1/2}$, we obtain the WKB approximation \break
$\Psi_{{\mathop{\rm WKB }\nolimits}} 
\! \left({ t_b, x_b }\right)$ to $\Psi \! \left({ t_b, x_b }\right)$. 
Using formulas (IV.18) and (IV.19), we derive
$$
\Psi_{{\mathop{\rm WKB }\nolimits}} 
\! \left({ t_b, x_b }\right) =
\exp \!
\left({ - {1 \over s \hbar }
\Sa \! \left({ t_b, x_b ; t_a , x_a }\right) 
 }\right) \cdot I \, , 
\eqno ({\mathop{\rm IV.20 }\nolimits})$$
with the integral
$$I =
\int_{ {\bf Z}_b }^{ } \Da \z \cdot \exp \!
\left({ - {\pi  \over s}  Q_1 (\z) }\right)
\d \! \left({ \l \z (t_a) }\right)
\, . 
\eqno ({\mathop{\rm IV.21 }\nolimits})$$

\medskip
\nn
Besides the quadratic form
$$Q_0 (\z) = \int_{ {\bf T} }^{ } dt \, 
\, \vert \dot \z (t) \vert^2
\eqno ({\mathop{\rm IV.22 }\nolimits})$$
corresponding to the free particle, we need the quadratic form
$$Q_V (\z) = - {1 \over m} \int_{ {\bf T} }^{ } dt \, 
h_{\a \b } (t) \, \z^\a (t) \z^\b (t) 
\eqno ({\mathop{\rm IV.23 }\nolimits})$$
with 
$h_{\a \b } (t) = \bgt_\a \bgt_\b V \!
\left({  x_{ {\mathop{\rm cl }\nolimits} } (t)  }\right)$.
In (IV.21), we use the quadratic form $Q_1 = Q_0 + Q_V$.

\bigskip
The easiest method to calculate the functional integral $I$ consists of the following steps.

\medskip
a) {\it Changing the integrator}: according to formulas (A.126), (A.132) and (A.134) in
Appendix A, we obtain $I = I_1 I_2 $ where
$$
\eqalignno{
I_1 &=
\left\vert{ \Det \!
\left({ Q_0 / Q_1 }\right)
 }\right\vert_{}^{1/2}
s^{- {\mathop{\rm Ind  }\nolimits} \! \left({  Q_1 }\right) } \, ,
&({\mathop{\rm IV.24 }\nolimits})\cr
I_2 &=
\int_{ {\bf Z}_b  }^{ }
\Da_{s, Q_1} \z \cdot \exp
\! \left({-  {\pi  \over s } Q_1 (\z)  }\right)
\d \! \left({ \l \z (t_a) }\right) \, .
&({\mathop{\rm IV.25 }\nolimits})\cr} $$
Here ${\mathop{\rm Ind  }\nolimits} \! \left({  Q_1 }\right)$ is the number of negative
directions for the quadratic form $Q_1$.

\medskip
b) {\it Restricting the domain of integration}:
to treat the $\d$ factor in $I_2$, we use the linear change of variables 
$\z \mps \z (t_a) $ from ${\bf Z}_b$ to $\RR^d$. By a method similar to the one used in
paragraph A.3.8, we obtain 
$$I_2 =
\left({ \det s W_1 \! \left({ \d_{t_a}  }\right) }\right)_{}^{- 1/2}
\l^{- d } \, , 
\eqno ({\mathop{\rm IV.26 }\nolimits})$$
where $W_1$ is the quadratic form on ${\bf Z}'_b$ inverse to $Q_1$.
Using the Green's function $G$ given by equation (B.26), we evaluate the matrix
 $W_1 \! \left({ \d_{t_a}  }\right)$ as follows
$$W_1 \! \left({ \d_{t_a}  }\right) =
G \! \left({ t_a, t_a }\right) =
N \! \left({ t_a, t_b }\right)
J \! \left({ t_b, t_a }\right) \, .
\eqno ({\mathop{\rm IV.27 }\nolimits})$$

\medskip
c) {\it Reducing a functional determinant to a finite determinant}:
using the same strategy as in paragraph III.2 we obtain
$$\Det \! \left({ Q_0 / Q_1 }\right) =
\det \!
\left({ K \! \left({ 1 ; t_b, t_a  }\right)^{- 1 }
 K \! \left({ 0 ; t_b, t_a }\right)  }\right) \, . 
\eqno ({\mathop{\rm IV.28 }\nolimits})$$
Here 
$K \! \left({ \nu ; t, t_a }\right)$ is the Jacobi field defined by
$$Q_\nu K \! \left({ \nu ; t, t_a }\right) = 0 \quad , 
\quad 
K \! \left({ \nu ; t_a, t_a }\right)
= \un 
\eqno ({\mathop{\rm IV.29 }\nolimits})$$
for $\nu$ equal to $0$ or $1$.

\bigskip
Finally, collecting the previous equations (IV.20) to (IV.29) and using equations (B.11),
(B.12), (B.7), (B.8), (B.19) and (III.42), we obtain {\it the} WKB {\it approximation to the
point-to-point transition amplitudeJ} (in the case $s = i$):
$$
\left\langle{ t_b, x_b \bigm|  t_a, x_a }\right\rangle_{
 WKB 
}^{}
= c \, 
h^{- d/2}
\left\vert{ \det \part^2 \Sa /
\part x_a^\a \cdot \part x_b^\b }\right\vert_{}^{1/2}
e^{i \Sa / \hbar} \, . 
\eqno ({\mathop{\rm IV.30 }\nolimits})$$
Here 
$\Sa = \Sa \! 
\left({ t_b , x_b ; t_a, x_a }\right)$ is the classical action function and the phase factor 
$c$ is $e^{\pi i (p - q)/4}$ where $p$ $(q)$ is the number of positive (negative) eigenvalues of
the van Vleck-Morette matrix
$\left({ \part^2 \Sa / \part x_a^\a \cdot \part x_b^\b }\right)_{
{1 \leq \a \leq d \atop 1 \leq \b \leq d }}$.

\bigskip
For the higher-order terms in the semiclassical approximation, whether 
$x_{{\mathop{\rm cl}\nolimits}} $ is or is not a degenerate critical point of the action, we
refer the reader to the literature [5] or the references at the end of Appendix B. The reader
can easily transfer these old results using the simpler and more general formalism of this
paper.

\bigskip
\nn
{\bf 2. Polar coordinates.}\footnote{$^{13}$}{Contributed by John La Chapelle.}
\smallskip
Our basic integral is formulated in terms of a transformation 
$\Si(t, z)$ on a manifold $N$, in a form valid for an arbitrary system of coordinates.
Before considering the case of a general Riemannian manifold in paragraph IV.3, we consider
polar coordinates in a plane. The case of cylindrical coordinates in $\RR^3$ is very similar.

\bigskip
``Path integrals in polar coordinates'', the  1964 paper by S.F. Edwards and Y.V.
Gulyaev [17] has been, and still is, at the origin of many investigations. To the best of our
knowledge, all the papers on the subject deal with the discretized version of the path
integral, and propose various path integral prescriptions when the discretized values of the
paths are expressed in coordinates other than cartesian.

\bigskip
The basic manifold is $N = \RR^2 \bl \{ 0 \}$ with coordinates 
$x^1 $ and $x^2$. We also consider the manifold $\ww N =
]0 , + \Inf [ \ts \RR$ with coordinates $r, \t$ (sole restriction $r > 0$) and the covering map 
$\Pi : \ww N \ra N$ taking $(r, \t) $ into 
$ \left({ x^1, x^2 }\right)$ with
$$x^1 = r \cos \t \quad , \quad  x^2 = r \sin \t \, . 
\eqno ({\mathop{\rm IV.31 }\nolimits})$$
Two points of $\ww N$ map onto the same point of $N$ if, and only if, their $\t$-coordinates differ
by an integral multiple of $2 \pi$.

\bigskip
Let $\D = 
{\part^2 \over (\part x^1)^2 } + 
{\part^2 \over (\part x^2)^2 }$ be the Laplacian in cartesian coordinates.
We know how to solve the Schr\"odinger equation
$${\part \Psi \over \part t} =
{s \over 4 \pi } \D \Psi
\eqno ({\mathop{\rm IV.32 }\nolimits})$$
by means of the path integral
$$\Psi \! \left({ t_b, x_b }\right)
= \int_{{\bf Z}_b }^{ }
\Da z \cdot e^{ - \pi Q_0 (z) / s} 
\phi \! \left({  x(t_a, z)  }\right)
\eqno ({\mathop{\rm IV.33 }\nolimits})$$
provided that $x(t, z)$ is a solution of the differential system
$$dx^1 = dz^1 \quad , \quad dx^2 = dz^2
\eqno ({\mathop{\rm IV.34 }\nolimits})$$
with boundary condition 
$x \! \left({ t_b, z }\right) = x_b$.
To solve the Schr\"odinger equation in polar coordinates, we need only to transform the system
(IV.34) in polar coordinates with the help of the transformation equations (IV.31). Hence
$$
\left\{\matrix{
dr 
= \cos \t \cdot d z^1 + \sin \t \cdot dz^2 =
X_{(1)}^1 dz^1 + X_{(2)}^1 dz^2
  \hfill\cr
\noalign{\medskip}
d \t = - 
 {\ds \sin \t  \over \ds r} \cdot dz^1 +
   {\ds \cos \t  \over \ds r} \cdot dz^2 =
X_{(1)}^2 dz^1 + X_{(2)}^2 dz^2 \, . 
  \hfill \cr}
\right. \eqno ({\mathop{\rm IV.35 }\nolimits})
$$
The vector fields $X_{(1)}$ and $X_{(2)}$ can be read off from the above equations:
$$\La_{X_{(1)}}
=
\cos \t \cdot {\part  \over \part r} -
{\sin \t \over r} \cdot {\part  \over \part \t} \,  ,
\eqno ({\mathop{\rm IV.36 }\nolimits})$$
$$\La_{X_{(2)}}
=
\sin \t \cdot {\part  \over \part r} +
{\cos \t \over r} \cdot {\part  \over \part \t} \, .
\eqno ({\mathop{\rm IV.37 }\nolimits})$$
Hence, if 
$(r(t, z), \t(t, z))$ is the solution of the differential system 
(IV.35) such that $r(t_b, z) = r_b$,
$\t(t_b, z) = \t_b$, the path integral
$$\Psi
\! \left({ t_b, r_b, \t_b }\right) : = 
\int_{ {\bf Z}_b }^{ } \Da z
\cdot e^{- \pi Q_0 (z) / s}
\phi (r(t_a, z), \t(t_a, z)) 
\eqno ({\mathop{\rm IV.38 }\nolimits})$$
solves the Schr\"odinger equation (IV.32). As expected, the operator
$$\D = 
\La_{X_{(1)}}^2 + \La_{X_{(2)}}^2
=
{ \part^2\over \part r^2}
+ {1 \over r^2} {\part^2 \over \part \t^2} +
{1 \over r} {\part  \over \part r}
$$
is the Laplacian on $\RR^2 \bl \{ 0 \}$ in polar coordinates.

\bigskip
We want to evaluate the point-to-point  transition amplitudes in polar coordinates, denoted
by $\left\langle{ t_b, r_b, \t_b \bigm| t_a , r_a, \t_a  }\right\rangle$. To obtain them, it suffices to put
$\phi (r, \t) = \d (r - r_a) \d (\t - \t_a)$ in equation (IV.38).
 Solving the system (IV.35) is easy by reverting
to cartesian coordinates, hence
$$
\left\{\matrix{
r(t, z) \cos \t (t, z) = r_b \cos \t_b + z^1 (t)
 \hfill\cr
\noalign{\medskip}
  r(t, z) \sin \t (t, z) = r_b \sin \t_b + z^2 (t) \, . 
  \hfill \cr}
\right. \eqno ({\mathop{\rm IV.39 }\nolimits}) 
$$
For given $x^1$, $x^2$, the equations (IV.31) in $r, \t$ have infinitely many solutions, and we derive easily
$$\d \!
\left({ r \cos \t - r_a \cos \t_a  }\right)
\d \! \left({ r \sin \t - r_a \sin \t_a }\right)
= 
{1 \over r_a}
\sum_{ n \in \ZZ}^{ }  
\d \! \left({ r - r_a  }\right) 
\d \! \left({ \t - \t_a - 2n \pi  }\right) \, .
\eqno ({\mathop{\rm IV.40 }\nolimits}) $$
Substituting $r(t, z)$ to $r$ and $\t (t, z)$ to $\t$ and taking into account equations (IV.39), we obtain the
path integral representation
$$
{1 \over r_a}
\sum_{ n \in \ZZ}^{ }  
\left\langle{t_b, r_b, \t_b \bigm| t_a, r_a , \t_a + 2 n \pi  }\right\rangle
= I \, . 
\eqno ({\mathop{\rm IV.41 }\nolimits}) $$

\medskip
\nn
Here we use the definitions
$$\phi_1 (u) =
\d \! \left({ r_b \cos \t_b - r_a \cos \t_a - u^1 }\right)
\d \! \left({ r_b \sin \t_b - r_a \sin \t_a - u^2 }\right)
\eqno ({\mathop{\rm IV.42 }\nolimits}) $$
for $u = \left({ u^1, u^2 }\right)$ in $\RR^2$, and
$$ I = \int_{ {\bf Z}_b }^{ }
\Da z \cdot e^{- \pi Q_0 (z) / s }
\phi_1 (z(t_a)) \, . 
\eqno ({\mathop{\rm IV.43 }\nolimits}) $$
This integral is nothing else than a point-to-point transition amplitude in cartesian coordinates.
We give now a derivation of the well-known formula for this amplitude, by relying on our methods.

\bigskip
The integrand 
in equation (IV.43) depends only on 
$z^1 (t_a) $, $z^2 (t_a)$, hence as explained in paragraph A.2.2, we introduce the linear map
$L : {\bf Z}_b \ra \RR^2 $ mapping $z$ into $z(t_a)$.
The image of the integrator 
$\Da z \cdot e^{- \pi Q_0 (z) / s }$ on ${\bf Z}_b $ is a Gaussian integrator on $\RR^2$; we
need only the value of the corresponding covariance matrix
$$W^{\a \b } =
\int_{ {\bf Z}_b }^{ }
\Da z \cdot e^{- \pi Q_0 (z) / s } \, 
z^\a (t_a) z^\b (t_a) \, .  
\eqno ({\mathop{\rm IV.44 }\nolimits})$$

\medskip
\nn
According to formulas (A.57) and (A.62), we obtain
$$W^{\a \b } =
{s \over 2 \pi }
G_b^{\a \b }\! \left({ t_a, t_a }\right)
=
{s \over 2 \pi } \d^{\a \b }
\! \left({ t_b - t_a }\right)
\, .
\eqno ({\mathop{\rm IV.45 }\nolimits})
$$

\medskip
\nn
According to formulas (A.27) and (A.32), we transform the integral (IV.43) into
$$I =
\left (s \left({ t_b - t_a }\right)\right)^{-1}
\int_{\RR }^{ } du^1 \, \int_{\RR }^{ } du^2 
\exp \! 
\left({- {\pi  \over s} 
{ \left\vert{ u }\right\vert_{}^{2} \over t_b - t_a }
}\right)
\phi_1 (u) \, . 
\eqno ({\mathop{\rm IV.46 }\nolimits})$$

\medskip
\nn
The function $\phi_1 (u) $ is a $\d$-factor and the final result is
$$
\eqalign{I 
&= \left(s\left({ t_b - t_a }\right)\right)^{-1} 
\exp \! 
\left({ -{ \pi \over s}
{r_a^2 + r_b^2 \over t_b - t_a} }\right)
\exp \! 
\left({ { 2\pi \over s}
{r_a  r_b \cos (\t_b - \t_a) \over t_b - t_a} }\right)
\cr
&
=
\left (s \left({ t_b - t_a}\right)\right) ^{-1}  
\exp \! 
\left({- {\pi  \over s} 
{ \left\vert{ x_b - x_a }\right\vert_{}^{2} \over t_b - t_a }
}\right) \, . \cr} 
\eqno ({\mathop{\rm IV.47 }\nolimits})$$

\bigskip
>From the equations (IV.41) and (IV.47) we cannot derive directly the point-to-point transition amplitudes
in polar coordinates. We defer to paragraph IV.7 a further discussion of this point.

\bigskip
\nn
{\bf 3. Frame bundles over Riemannian manifolds.}\footnote{$^{14}$}{\nn
For a discussion of functional integrals when paths take their values in a Riemannian space, see
[5], [14] and [18]. }

\smallskip
We consider a Riemannian (or pseudo-Riemannian) manifold $M$ of dimension $d$, with metric
$g$, and the orthonormal frame bundle $N = O(M)$ over $M$, with projection 
$\pi : O(M) \ra M$. We want to choose vector fields 
$X_{(1)}, \cdots, X_{(d)}$ on $O(M)$ such that the {\it equation on the bundle}
$O(M)$ 
$${ \part \over \part t_b} 
\Psi \! \left({ t_b, \r_b }\right) =
{s \over  4 \pi} h^{\a \b } X_{(\a)} X_{(\b)}
\Psi \! \left({ t_b, \r_b }\right) \, ,
\eqno ({\mathop{\rm IV.48 }\nolimits})$$
where $\r_b$ belongs to $O(M)$ and the first-order differential operators are given by
$X_{(\a)} = X_{(\a)}^\l (\r_b) 
{\part  \over \part \r_b^\l}$, gives by projection {\it an equation on the base space $M$},
namely:
$$
\eqalign{
{\part  \over \part t_b } \psi \! 
\left({ t_b, x_b  }\right) 
&= { s \over  4 \pi } g^{\l \mu} D_\l D_\mu \psi \! 
\left({ t_b, x_b  }\right)  \cr
& 
= { s \over  4 \pi } \D \psi \! 
\left({ t_b, x_b  }\right)  \, . \cr}
\eqno ({\mathop{\rm IV.49 }\nolimits})$$
Here $D_\l$ is the covariant derivative defined by the Riemannian connection, and 
$\Psi = \psi \circ \pi$; moreover $\D$ is the Laplace-Beltrami operator on $M$.

\bigskip
It has been shown in [13] that the covariant Laplacian $\D$ at a point $x_b$ of $M$ can be lifted
to a sum of products of Lie derivatives $h^{\a \b } X_{(\a)} X_{(\b)}$ at the frame 
$\r_b $ in $O(M)$; the integral curves of the set of vector fields
$\left\lbc{ X_{(\a)} }\right\rbc$ starting from $\r_b$ at time $t_b$ are the horizontal lifts
of a set of geodesics at $x_b$, tangent to the basis 
$\left\lbc{ e_{\a} }\right\rbc$ of $T_{x_b} M$ corresponding to the frame $\r_b$.
The constant matrix $\left({ h_{\a \b} }\right)$ has been chosen with the same signature as
the metric $g$ on $M$, and $g\! \left({ e_\a, e_\b }\right) = h_{\a \b }$.

\bigskip
An explicit construction of $\left\lbc{ X_{(\a)} }\right\rbc$ goes as follows. Let $\r(t) $ be
the horizontal lift of a path $x(t) $ in $M$ defined by the Riemannian connection map\footnote{$^{15}$}{
More explicitly, for a given frame $\r$ at a point $x$ of $M$, $\s (\r)$ is a linear map from 
$T_x M$ to $T_\r O(M) $, mapping $u$ into $\s (\r) \cdot u$ for $u$ in $T_x M$.}
$\s : O (M) \ra L (TM, TO (M))$; it satisfies the differential equation
$$\dot \r (t) = \s (\r (t)) \cdot \dot x (t) \quad, \quad
\r(t_b) = \r_b \, .
\eqno ({\mathop{\rm IV.50 }\nolimits})$$

\medskip
\nn
If we put $\r(t) = 
(x(t), u(t))$ and $\r_b = \left({ x_b, u_b }\right)$, a solution $\r(t)$ of the previous
equation corresponds to the frame $u(t) $ obtained by parallel transport of $u_b$,  along the
path $x$ from $x_b$ to $x(t)$. The frame $u(t)$ is also an admissible map
$$u(t) : \RR^d \ra T_{x(t)} M$$
i.e. $u(t)$ maps a $d$-tuple into a vector whose components in $u(t)$ are the chosen $d$-tuple;
equivalently
$\left({ u(t)^{- 1} \dot x (t)  }\right)_{}^{\a}$
is the $\a$-coordinate of $\dot x (t) $ in the $u(t) $ frame.
Set 
$$\dot z (t) : = u(t)^{-1}
\dot x (t) = \dot z^\a (t) \oo e_\a \, , 
\eqno ({\mathop{\rm IV.51 }\nolimits})$$
where $\left\lbc{ \oo e_\a }\right\rbc$ denotes the canonical basis of the model space $\RR^d$.
Then we can express (IV.50) in the canonical form (II.21)
$$\dot \r (t) = X_{(\a)} \!
\left({ \r(t)   }\right) \dot z^\a (t)
\quad, \quad \r (t_b) = \r_b 
\eqno ({\mathop{\rm IV.52 }\nolimits})$$
where $X_{(\a)}$ is defined by 
$$X_{(\a)} (\r) =
\s(\r) \cdot e_\a = 
\left({ \s (\r) \circ u }\right) \cdot \oo e_\a \, . 
\eqno ({\mathop{\rm IV.53 }\nolimits})$$
Here $x$ is a point of $M$ and $\r = (x, u) $ a frame, where
$u : \RR^d \ra T_x M$ is an admissible map.

\bigskip
If $\dot z(t) = \oo e_\a $, i.e. $\dot z^\b (t) = \d_\a^\b $, the coordinates of $\dot x(t) $ are
constant in the frame $u(t) $ parallel transported along $x(t)$. Therefore $x(t)$ is the
geodesic defined by 
$$x(t_b) = x_b \quad , 
\quad \dot x (t_b) = e_\a \, . 
\eqno ({\mathop{\rm IV.54 }\nolimits})$$
With 
$\dot z (t) = \oo e_\a$, equation (IV.52) reads 
$$\dot \r_{(\a)} (t) =
X_{(\a) } \! \left({ \r_{(\b)} (t) }\right)
\d_\a^\b = 
X_{(\a) } \! \left({ \r_{(\a)} (t) }\right)\, . $$

\medskip
\nn
{\it The horizontal lift $\r_{(\a)} (t)$ of the geodesic} (IV.54) {\it is, as desired, the
integral curve of 
$X_{(\a)}$ going through $\r_b$
at time $t_b$.
}
With $X_{(\a)}$ defined by (IV.53), $\r$ can be expressed in terms of the {\it Cartan
development map
}
$$(\pi \circ \r)(t) = x(t) =
(\Dev z)(t) \, . 
\eqno ({\mathop{\rm IV.55 }\nolimits})$$
The Cartan development is a bijection from a space of pointed paths (paths with a fixed end
point) on $T_{x_b} M$ (identified to $\RR^d$ via the frame $\r_b$) into a space of pointed paths on $M$ -- or
vice versa.  Here $$\Dev : \Pa_0 T_{x_b} M \ra  \Pa_{x_b} M \quad \by \quad z \mps x \, .
\eqno ({\mathop{\rm IV.56 }\nolimits})$$
The path $x$ is said to be the development of $z$, if $\dot x (t)$ parallel transported along $x$ from
$x(t)$ to $x_b$ is equal to $\dot z(t) $ trivially transported to the origin of $T_{x_b} M$,
for every $t \in {\bf T}$.

\bigskip
The path integral solution of (IV.48) is 
$$\Psi
\! \left({ t_b, \r_b  }\right) =
\int_{ {\bf Z_b}J}^{ }
\Da z \cdot \exp
\! \left({ - {\pi  \over s } Q_0 (z)
}\right)  \Phi 
\! \left({ \r (t_a, z) }\right)
\eqno ({\mathop{\rm IV.57 }\nolimits})$$
with $\r (t, z)$ solution of (IV.52). The path integral solution of (IV.49) is 
$$
\psi
\! \left({ t_b, x_b  }\right)
=
\int_{ {\bf Z_b}J}^{ }
\Da z \cdot \exp
\! \left({- {\pi  \over s} Q_0 (z)}\right)
\phi ((\Dev z)(t_a)) \, . 
\eqno ({\mathop{\rm IV.58 }\nolimits})$$
If a scalar potential is desired in (IV.48) or (IV.49), one can proceed along either of the
methods outlined in paragraph II.4.

\bigskip
The semiclassical approximation [5] of $\Psi$, with or without scalar potential, is considerably
more complicated to compute than the semiclassical approximation of $\Psi$ given in section
III. For an initial wave function $\phi$ of type (III.2) the blue-print given in Section III.2
is complete. If one wishes to compute the point-to-point propagator on a Riemannian space $M$,
one chooses the initial wave function on $M$ to be
$$\phi (x) = \d_{x_a } (x) \, . 
\eqno ({\mathop{\rm IV.59 }\nolimits})$$
The detailed  calculation can be found in reference [5, pp. 309-311].
The development map cannot parametrize spaces of paths with two fixed points\footnote{$^{16}$}{
\nn For instance, two geodesics on $S^2$ intersect at two antipodal points; they are the
developments of two halflines with one common origin.}
but (IV.59) with $(\Dev z) (t_a)$ substituted to $x$ restricts the domain of integration appropriately.

\bigskip
The following remarks simplify the calculations. 

\bigskip
\nn
{\it Remark 1.}
According to the formulas
$$
\eqalign{
\int_{{\bf T} }^{ } dt \, 
h_{\a \b } \dot z^\a (t) \dot z^\b (t) 
&= 
\int_{{\bf T} }^{ } dt \, 
g_{\mu \nu } (x(t)) \dot x^\mu (t) \dot x^\nu (t)\cr
& 
= \int_{{\bf T} }^{ } dt \, 
g_{\mu \nu} (x(t)) \dot \r^\mu (t) \dot \r^\nu (t) \, ,\cr}
\eqno ({\mathop{\rm IV.60 }\nolimits})$$
development map and horizontal lift preserve lengths and angles.

\bigskip
\nn
{\it Remark 2.}
If $z$ develops into a classical paths 
$x_{{\mathop{\rm cl }\nolimits}}$, the determinant of the derivative mapping 
$\Dev'(z)$ is unity. For the proof see, for instance, reference [5], p. 308.

\bigskip
\nn
{\bf 4. A multiply connected manifold.}

\smallskip
The domain of integration, a space of pointed paths $\Pa_x N$, is the union of disjoint sets
made of paths in different homotopy classes. We recall in paragraphs a) and b) earlier
calculations of propagators on multiply connected spaces.
 Then, in paragraph c), we explain how these methods fit into our general framework.

\bigskip
a) It was shown [3] in 1971 that, for a system with a multiply connected configuration space
$N$, the propagator $K$ is a linear combination of propagators $K_{(\a)}$
$$
\left\vert{ K^A }\right\vert
=
\left\vert{ \sum_{ g_\a \in \pi_1 (N) }^{ }  
\chi^A (g_\a) K_{(\a)} }\right\vert \, . 
\eqno ({\mathop{\rm IV.61 }\nolimits})$$

\medskip
\nn
Each $K_{(\a)}$ is obtained by summing over paths in the same homotopy class, say $\a$; the set
$\left\lbc{ \chi^A (g_\a }\right\rbc_{_\a}^{}$ forms a representation, labeled $A$, of the
fundamental group $\pi_1 (N)$. Since a homotopy class cannot be identified uniquely with an
element of $\pi_1 (N)$, the propagator is defined modulo an overall phase factor. There are as
many propagators $K^A$ as there are inequivalent representations of $\pi_1 (N)$.

\medskip
The proof of (IV.61) uses two facts:

\smallskip
i) the superposition principle of quantum propagators implies the linear combination of
partial propagators;

\smallskip
ii) the fundamental group based at a point is isomorphic to the fundamental group based at
another point, but not canonically so. Therefore the pairing 
$\left({ g_\a, K_{(\a)} }\right)$ is done by choosing an homotopy mesh (choosing a point for
the fundamental group, and pairing one group element with one homotopy class), then requiring
that the result be independent of the homotopy mesh.

\bigskip
b) Later on [14, see also 18, p. 65] the same result (worked out for a different example) was
obtained from stochastic processes on fibre bundles. 
The basic steps are as follows\footnote{$^{17}$}{We refer the reader to paragraph IV.7 for an explicit example
where we use this strategy.}.

\smallskip
i) A universal covering $\ww N$ is a principal $G$-bundle with projection \break $\Pi : \ww N\ra N$,
where $N = \ww N / G$ and $G$ is a discrete group of automorphisms of $\ww N$ isomorphic to
the fundamental group of $N$. For example, $N = S^1 = \RR / \ZZ$ and $\ww N = \RR$ is a
$\ZZ$-principal bundle over $S^1$.

\smallskip
ii) The wave function for a system with configuration space $N$ is a section of a bundle {\it
weakly} associated to $\ww N$. i.e. a bundle whose typical fibre is associated to a not
necessarily faithful representation of $G$, hence to a representation of a group homomorphic
to $G$, says $G_0$. 
For example if $N = S^1$, a vector bundle over $S^1$ with structure group $U(1)$ is weakly
associated to a $\ZZ$-principal bundle over $S^1$ by a homomorphism 
$h_\a$ of $\ZZ$ into $U(1)$ mapping $n$ into $e^{in\a}$.

\smallskip
iii) There is a unique connection on $\ww N$: the horizontal lift $\ww x$ of a pointed path
$x$ with fixed point $\ww x (t_0) = \ww x_0$ is uniquely defined by the lift $\ww x_0 $ of
$x_0$. This unique connection defines the parallel transport of the wave function $\phi (x
(t))$ back to $x_0$.

\smallskip
Consider for instance the case where $\ww N = \RR$, $N = S^1$ and 
$\Pi (x) = e^{ix}$. Take for $\phi$ a section of a $U(1)$-bundle defined by the homomorphism
$h_\a$ as above. Let $x(t)$ be a map into $S^1$, lifted to a map $\ww x_0 (t)$ into $\RR$ in such a way that 
$x(t) = e^{i \ww x_0 (t)}$. The other liftings are given by
$$\ww x_k (t) = \ww x_0 (t) + 2 \pi k
\eqno ({\mathop{\rm IV.62 }\nolimits})$$
for $k$ in $\ZZ$.
The parallel transport of $\phi 
\! \left({ \ww x_k (t) }\right)
$ to $x_0$ is given by the formulas
$$\tau_{t_0}^t \phi 
\! \left({ \ww x_k (t) }\right) =
\ww h_\a (\ww x_0) \ww h_\a
\! \left({ \ww x_k (t) }\right)_{}^{-1}
\phi (x(t)) 
\eqno ({\mathop{\rm IV.63 }\nolimits})$$
and
$$
\ww h_\a
\! \left({ \ww x_k (t) }\right)_{}^{-1} =
\ww h_\a
\! \left({ \ww x_0 (t) }\right)_{}^{-1}
h_\a (k)^{-1}
= e^{- i k \a }
\ww h_\a
\! \left({ \ww x_0 (t) }\right)_{}^{}\, . 
\eqno ({\mathop{\rm IV.64 }\nolimits})$$
The map $\ww h_\a$ is the bundle map from the $\ZZ$-bundle $\RR $ over $S^1$ to the
$U(1)$-bundle corresponding to the map $h_\a : \ZZ \ra U (1)$.

\smallskip
iv) Summing the wave function $\tau_{t_0}^t \phi \! \left({ \ww x_k (t) }\right)$
over  all $k$'s gives the solution $\Psi (t, x_0)$ of a parabolic
equation with initial value $\Psi (t_0, x_0) = \phi (x_0)$.

\bigskip
In both computations outlined above, one chooses a representation of the fundamental group,
$\left\lbc{ \chi^A (g_\a) }\right\rbc$, or the homomorphism $h : G \ra G_0$.

\bigskip
c) We specialize our basic formula (II.1)
$$
\Psi \! \left({ t_b, x_b }\right)
= \int_{ {\bf Z }_b }^{ }
\Da z 
\cdot e^{- \pi Q_0 (z) / s}
\phi \! \left({ x_b \cdot \Si (t_a, z) }\right) 
\eqno ({\mathop{\rm IV.65 }\nolimits} )$$
to the case where $N$ is multiply-connected.
To calculate the point-to-point transition amplitudes, we select $\phi (x)$ of the form $\d_{x_a} (x)$
with a $\d$-factor centered at a point $x_a$ of $N$.
Denote the evaluation map taking $z$ into $x_b \cdot \Si (t_a, z)$ by 
$\ve : {\bf Z}_b \ra N$.
Since ${\bf Z}_b $ is contractible, we can lift $\ve$ to a map
$\ww \ve : {\bf Z}_b \ra \ww N$ into the universal covering $\ww N $ of $N$ and hence
$\ve = \Pi \circ \ww \ve$. In the path integral (IV.65), the domain of integration is
restricted by the $\d$-factor 
$\d_{x_a} (x) $ to the inverse image 
$\ve^{-1}(x_a)$. It consists of paths such that
$x_b \cdot \Si (t_a, z) = x_a$ and {\it splits as the union of domains} 
${\bf Z}_\G = \ww \ve^{-1} \! \left({ \ww x_\G }\right)$ where $\ww x_\G$ runs over the various
points of $\ww N$ mapping to $x_a$ by $\Pi$. The labels $\G$ correspond to the various
homotopy classes of paths 
$x : {\bf T} \ra N$ such that 
$x(t_a) = x_a$, $x(t_b) = x_b$.
Hence the integral (IV.65) splits into a sum of integrals over the various subdomains
$$\Psi \! \left({ t_b, x_b }\right)
= \sum_{ \G }^{ } \int_{ {\bf Z }_\G }^{ }
\Da z 
\cdot e^{- \pi Q_0 (z) / s} \d_{x_a}
 \! \left({ x_b \cdot \Si (t_a, z) }\right) \, .
\eqno ({\mathop{\rm IV.66 }\nolimits} )$$

\medskip
This equation is the justification of the heuristic idea used in (IV.61) that the building
blocks of $K$ are the propagators 
$K_{(\a)}$ obtained by summing over paths in the same homotopy class.
To explain the coefficients $\chi^A (g_\a)$ we can proceed as follows: using the previous
notations 
$\ww N $, $\Pi$, $G$, consider a homomorphism $\chi^A$ of $G$ into $U(1)$.
The corresponding wave functions are functions $\ww \phi $ on $\ww N$ such that
$\ww \phi (  \ww x  g ) = \chi^A (g) \ww \phi (\ww x)$ for $g$ in $G$ and 
$\ww x $ in $\ww N$. We denote by 
$\ww \Si (t, z) $ the lifting of $\Si (t, z) $ to $\ww N$ and generalize equation (IV. 65)
by
$$\ww \Psi \! \left({ t_b, \ww x_b }\right)
=  \int_{ {\bf Z }_b }^{ }
\Da z 
\cdot e^{- \pi Q_0 (z) / s} \ww \phi
 \! \left({\ww  x_b \cdot \ww \Si (t, z) }\right) \, .
\eqno ({\mathop{\rm IV.67 }\nolimits} )$$
The function $\ww \Psi$ will satisfy the same transformation property as $\ww \phi$ and when $\ww
\phi$ is a $\d$-factor, we can split the integration domain into subdomains ${\bf Z}_\G$
as above. We recover then the results derived in paragraphs IV 4a and b.

\bigskip
\nn
{\bf 5. Gauge fields.}

\smallskip
We begin with the case of an abelian gauge group. 
In physical terms, we consider a particle of mass $m$ and electric charge $e$ moving  under the influence of a
magnetic potential $A$, with components $A_\a (x)$ at the point $x$. We consider generally a
$d$-dimensional space $\RR^d$ in cartesian coordinates $x^\a$ ($\a \in \{1, \cdots, d \}$) and 
metric $\left\vert{ x }\right\vert^2 = \d_{\a \b}
 x^\a x^\b $. The classical Lagrangian is given by
$$ L (x, \dot x) =
{m \over 2} \left\vert{ \dot x }\right\vert^2
+ e \, A_\a (x) \dot x^\a \, ,
\eqno ({\mathop{\rm IV.68 }\nolimits} )$$
hence the action functional
$$S(x) = 
\int_{{\bf T} }^{ } dt \ 
 L(x , \dot x) = {m \over 2} \int_{{\bf T} }^{ }
{ \left\vert{ dx }\right\vert^2 \over dt}  + e 
\int_{{\bf T} }^{ } A_\a d x^\a \, . 
\eqno ({\mathop{\rm IV.69 }\nolimits} )$$

\medskip
\nn
The equation of motion can be derived from this Lagrangian, and can be put into the
Hamiltonian form with the following definitions
$$p_\a = m \dot x^\a + e A_\a \quad , \quad 
H = \left\vert{ p - eA }\right\vert^2 / 2m \, . 
\eqno ({\mathop{\rm IV.70 }\nolimits} )$$
The corresponding Schr\"odinger equation is obtained in the standard way by replacing $p_\a $ by the operator 
${\hbar  \over i} 
{\part  \over \part x^\a} $
in the definition of $H$, and reads as
$$i \hbar { \part \psi \over \part t} 
= {1 \over 2m } 
\sum_{ \a }^{ }
\left({ { \hbar \over i}
{ \part \over \part x^\a} - e A_\a (x)
 }\right)^2 \psi \, . 
\eqno ({\mathop{\rm IV.71 }\nolimits} )$$
Our goal in this paragraph is to fit the well-known path integral solution of this equation
into our general framework.

\bigskip
It has long been recognized [19] that it is desirable to treat $\psi$ as a section of a complex
line bundle (here over $\RR^d$) associated to a principal $U(1)$-bundle via the canonical
representation of $U(1)$ acting on $\CC$ by multiplication. We describe the main steps of this
construction.

\vfill\eject
\nn
5.1. {\it The invariant formalism.}
\smallskip

The {\it base space} is $M = \RR^d$.
The {\it gauge group} $G$ is the set $U(1)$ of complex numbers of modulus one
$g = e^{i \T}$.
We consider a {\it principal bundle } $P$, with $G$ acting from the right via
$(p, g) \mps p \cdot g$, and projection $\Pi : P \ra M$.

\bigskip
The {\it connection } is a differential form $\om$ on $P$ with the transformation rule
$$\om (p \cdot g) = \om (p) + g^{- 1} dg \, . 
\eqno ({\mathop{\rm IV.72 }\nolimits} )$$
With the angular coordinate $\T$ such that 
$g = e^{i \T}$, one obtains $g^{- 1} dg = id \T $ for the invariant differential form on 
$U(1)$, hence the Lie algebra  ${\hbox{\tenfm g}} $ of $U(1) $ is naturally identified with
the set of pure imaginary numbers, and $\om$ is pure imaginary.

\bigskip
For any path $x : {\bf T} \ra M$ and any point $p_b$ in $P$ with 
$\Pi (p_b) = x(t_b)$, the {\it horizontal lifting}
is a curve $\xi : {\bf T} \ra P$ satisfying the following conditions:
$$\Pi \xi (t) = x(t) \quad, \quad \xi (t_b) = p_b 
\eqno ({\mathop{\rm IV.73 }\nolimits} )$$
$$\left\langle{  \om_{\xi (t) }, \dot \xi (t) }\right\rangle
= 0 \, ,
\eqno ({\mathop{\rm IV.74 }\nolimits} )$$
where $\om_{\xi (t) }$ is the value of $\om$ at the point
$\xi (t) $ of $P$, that is a linear form on the tangent space
$T_{\xi (t) } P$.

\bigskip
Let $L$ be the {\it associated line bundle}.
For any point $x$ in $M$, a point $p$ of $P$ with $\Pi (p) = x$ corresponds to an admissible
map
$\hh p : \CC \ra L_x $ where $L_x $ is the fiber of $L$ above the point $x$.
If $\psi$ is a {\it section} of $L$ over $M$, its value at point $x$ is an element $\psi (x) $
of $L_x$, hence $\hh p^{- 1} (\psi (x))$ is a complex number $\Psi (p)$. In this way  (see
e.g. [20], vol. I, p. 404), we identify the section $\psi $ of $L$ to a function $\Psi : P
\ra \CC$ with the transformation rule
$$\Psi (p \cdot g) = g^{- 1} \cdot \Psi (p) 
\eqno ({\mathop{\rm IV.75 }\nolimits} )$$
for $p$ in $P$ and $g$ in $G$.

\bigskip
\nn
5.2. {\it Fixing the gauge.}

\smallskip
A fixing of the gauge corresponds to a section $s : M \ra P$ of the principal bundle. When $s$
is chosen, we may identify $P$ with $M \ts G$ in such a way that $p =
s(x) \cdot g $ in $P$ correspond to $(x, g)$ in $M \ts G$.
A section $\psi$ of the line bundle $L$ corresponds now to  a {\it wave function},
that is to a complex-valued function on $M$ and the corresponding function $\Psi $ on $M \ts
G$ is given by $\Psi (x, g) = g^{- 1} \psi (x) $ and conversely 
$\psi (x) = \Psi (x, 1)$, or intrinsically $\psi = \Psi \circ s$. 

\bigskip
The differential form $\om$ on $P$ gives by pull-back via $s : M \ra P$ a purely imaginary
differential form on $M$, to be written as 
$- { ie\over \hbar} A$ to fit with standard physical dimensions. Hence the differential form $A$
on $M$ can be written as $A_\a (x) dx^\a$ and the functions $A_\a (x)$ are the components of
the magnetic potential. On $M \ts G$ the differential form $\om$ is given by
$$\om = g^{- 1} dg - 
{ ie\over \hbar} A =
i 
\left({ d \T - {e \over \hbar} A_\a dx^\a }\right)
\eqno ({\mathop{\rm IV.76 }\nolimits} )$$
(for $g = e^{i \T}$).

\bigskip
Let $x : {\bf T} \ra M$ be a path.
The horizontal lifting $\xi $ of $x$ is now described by 
$ \xi (t) = 
\left({ x(t) , e^{i \T (t) } }\right)$ and since $\om$ induces a zero form on the image $\xi ({\bf T})
\sbs M \ts G$, we obtain the differential equation
$$\dot \T =  {e \over \hbar}
A_\a (x(t)) \dot x^\a (t)
\eqno ({\mathop{\rm IV.77 }\nolimits} )$$
(see e.g. [18, pp. 64-65]).

\bigskip
Changing the gauge corresponds to choosing another section 
$s_1 : M \ra P$. There exists then a function $R : M \ra \RR$ such that 
$s_1 (x) = s(x) \cdot e^{- i e R (x) / \hbar}$.
In the new gauge, the section of the line bundle $L$ corresponds to a new wave function 
$$\psi_1 (x) = e^{ i e R (x) / \hbar} \psi (x) 
\eqno ({\mathop{\rm IV.78 }\nolimits} )$$
and the new components of the magnetic potential are given by 
$$A_\a^1 (x) = A_\a (x) +
{\part R (x) \over \part x^\a}
\, . 
\eqno ({\mathop{\rm IV.79 }\nolimits} )$$

\bigskip
\nn
5.3. {\it Path integrals.}
\smallskip
We revert to the notations in paragraph IV.1. The constant $\l$ is again $(h / m)^{1/2}$ and
we parametrize the paths $x$ with $x (t_b) = x_b$ by
$$x(t, z) = x_b + \l z (t) 
\eqno ({\mathop{\rm IV.80 }\nolimits} )$$
where $z$ runs over the space ${\bf Z}_b$. The horizontal lift of the previous path is given
by 
$$\xi (t, z) = 
\left({ x(t, z), e^{i \T (t, z) } }\right) \, . 
\eqno ({\mathop{\rm IV.81 }\nolimits} )$$
Taking into account the differential equation (IV.77), we obtain the following differential
system
$$\left\{\matrix{
d x^\a (t) = \l  d z^\a (t) \quad 
{\mathop{\rm for }\nolimits}
\  \a \ {\mathop{\rm in }\nolimits} \ \{ 1, \cdots, d \}
  \hfill\cr
\noalign{\medskip}
 d \T (t) = 
{e \l  \over \hbar} A_\a (x(t)) d z^\a (t) \, . 
  \hfill \cr}
\right.
\eqno ({\mathop{\rm IV.82 }\nolimits} )$$

\medskip
\nn
This system has the canonical form (II.21) where the vector fields 
$X_{(1)}, \cdots , X_{(d)}$, $Y$ are given by
$$\La_{X_{(\a)}} = \l \! 
\left({ {\part  \over \part x^\a} + {e \over \hbar} A_\a (x) {\part \over \part \T } }\right)
\quad, \quad Y = 0 \, .
\eqno ({\mathop{\rm IV.83 }\nolimits} )$$
Notice that for $\Psi (x, \T) = e^{- i \T} \psi (x) $, we have
$$\La_{X_{(\a)}}
\Psi (x, \T) = 
e^{- i \T } \cdot \l D_\a \psi (x) 
\eqno ({\mathop{\rm IV.84 }\nolimits} )$$
with the differential operator (see e.g. [20, p. 405]):
$$D_\a = 
 {\part  \over \part x^\a } 
-  {ie \over \hbar}  A_\a (x) \, . 
\eqno ({\mathop{\rm IV.85 }\nolimits} )$$

\bigskip
Our general partial differential equation
$$
 {\part \Psi \over \part t} 
=
 {s \over 4 \pi} \sum_{\a }^{ }
\La_{X_{(\a)}}^2 \Psi
\eqno ({\mathop{\rm IV.86 }\nolimits} )$$
translates now as
$$
{\part \psi \over \part t} 
=
 {s \over 4 \pi} {h \over m}
\sum_{\a }^{ } D_\a^2 \psi 
\eqno ({\mathop{\rm IV.87 }\nolimits} )$$
for the $\psi$-compoment of $\Psi$.
For $s = i$, this equation coincides with the Schr\"odinger equation (IV.71).

\bigskip
Feynman path integral solution to this equation reads as follows:
$$\psi \!  \left({ t_b, x_b  }\right) =
\int_{\Pa_b  }^{ }
\Da x \cdot e^{i S (x) / \hbar} \phi (x(t_a)) 
\eqno ({\mathop{\rm IV.88 }\nolimits} )$$
where the action $S(x)$ is given by equation (IV.69). Parametrizing the paths $x$ in $\Pa_b$ by the
paths $z$ in ${\bf Z}_b$ (see equation (IV. 80)), we can rewrite the previous path integral
as\footnote{$^{18}$}{
\nn We give the formula in the oscillatory case $s = i$. The reader is invited to work out the
formulas for the case $s = 1$. The definition of $Q_0 (z)$ is given in equation (IV.22).
}
$$
\psi  \! \left({ t_b, x_b  }\right) =
\int_{{\bf Z }_b }^{ }
\Da z \cdot
e^{ \pi i Q_0 (z)}
\exp \!
\left({ { ie\over \hbar  } 
\int_{{\bf T } }^{ }
A_\a (x(t, z)) d x^\a (t, z)
}\right) \phi (x(t_a, z)) \, . 
\eqno ({\mathop{\rm IV.89 }\nolimits} )$$

\medskip
\nn
Replacing $\psi (x) $ by $\Psi (x, \T) = e^{- i \T} \psi (x)$ and similarly $\phi (x) $ by
$\Phi (x, \T) =e^{- i \T} \phi (x)$, we can absorb the phase factor and obtain
$$\Psi \! \left({ t_b, x_b, \T_b }\right) =
\int_{{\bf Z }_b }^{ }
\Da z \cdot
e^{ \pi i Q_0 (z)}
\Phi \! \left({ \left({ x_b , \T_b }\right) \cdot \Si \!
\left({ t_a, z }\right) }\right) \, .
\eqno ({\mathop{\rm IV.90 }\nolimits} )$$
The transformation $\Si (t, z)$ of the bundle space $M \ts G$ takes 
$\left({x_b, \T_b  }\right)$ \break into 
$(x_b +\l z(t), \ \T_b - {e \l  \over \hbar }
\int_{t }^{ t_b} A_\a 
\left({ x_b + \l z (t) }\right) d z^\a (t))$ and corresponds to the integration of the
differential system (IV.82).

\bigskip
\nn
5.4. {\it Various generalizations.}

\smallskip
a) It is easy to incorporate an {\it electric  potential } $V$.
The complete action functional is now
$$S(x) = {m \over 2} \int_{{\bf T } }^{ }
{ \left\vert{ dx }\right\vert^2 \over dt} + e  \int_{{\bf T } }^{ }
 A_\a dx^\a - V dt \, . 
\eqno ({\mathop{\rm IV.91 }\nolimits} )$$
Feynman solution (IV.88) is still valid and can be made explicit as
$$\psi \! \left({ t_b, x_b }\right) =
\int_{{\bf Z }_b }^{ }
\Da z \cdot
e^{ \pi i Q_0 (z)} 
\exp \!
\left({ {ie \over \hbar } \int_{{\bf T } }^{ }  A_\a dx^\a - V dt 
 }\right) \phi\! \left({ x (t_a, z) }\right)
\eqno ({\mathop{\rm IV.92 }\nolimits} )$$
where the line integral 
$\int_{{\bf T } }^{ }  A_\a dx^\a - V dt $ is calculated along the path
$x(\cdot, z)$. The Schr\"odinger equation reads as follows:
$$i \hbar { \part \psi \over \part t}
= {1 \over 2m } \sum_{ \a }^{ } 
\left({ {\hbar \over i } {\part  \over \part x^\a}
- e A_\a (x)  }\right)_{}^{2} \psi + e V \psi \, . 
\eqno ({\mathop{\rm IV.93 }\nolimits} )$$

\bigskip
b) A non-abelian gauge group is for instance 
$G = U(N)$; a more general compact gauge group can always be realized as a closed subgroup of
some unitary group $U(N)$. We mention a few of the required changes. The Lie algebra 
${\hbox{\tenfm g}}$ is the set of matrices of the form $iE$ where $E$ is an $N$-by-$N$
hermitian matrix. The connection form $\om$ on the principal bundle $P$ takes its values in 
${\hbox{\tenfm g}}$ and the transformation rule is now
$$\om(p \cdot g) =
g^{-1} \om (p) g + g^{- 1}  dg
\eqno ({\mathop{\rm IV.94 }\nolimits} )$$
where $\om (p)$, $g$, $dg$ are 
 $N$-by-$N$ matrices and the products are given by matrix multiplication.

\bigskip
For the associated vector bundle $L$, we consider the natural action of $U(N)$ on the complex
vector space $\CC^N$. Hence a section $\psi$ of $L$ corresponds to a function 
$\Psi :P \ra\CC^N $ such that $\Psi (p \cdot g) = g^{- 1} \cdot \Psi (p)$ for $p$ in $P$ and
$g$ in $G = U(N)$.

\bigskip
In a given gauge the form  $\om $ is given by $\om (x, g) = g^{- 1} dg - 
{ ie \over \hbar} A(x) $, where $A(x)  = A_\a (x) dx^\a$ is a hermitian  $N$-by-$N$ matrix of
differential forms on $M$. A gauge transformation is given by the formulas
$$
\eqalignno{
\psi_1 (x) 
&= U(x)^{- 1} \psi (x) 
&({\mathop{\rm IV.95 }\nolimits} )\cr
A_\a^1 (x) 
&= U(x)^{- 1} A_\a (x) U(x) +  {i \hbar  \over e} U(x)^{- 1}
{\part \over \part x^\a } U(x),
&({\mathop{\rm IV.96 }\nolimits} )\cr}$$
where $U(x)$ is a unitary matrix depending on  the point $x$ of $M$. 

\bigskip

The horizontal lift of a path $x : T \ra M$ is of the form
$\xi (t) = (x(t), U(t))$ where the unitary matrix $U(t)$ satisfies the differential equation
$$\dot U (t) =  {i e \over \hbar} U(t) \cdot 
A_\a (x(t)) \dot x^\a (t) \, .
\eqno ({\mathop{\rm IV.97 }\nolimits} )$$
We can solve this equation in the standard way using time-ordered exponentials $\Ta \exp$.

\bigskip
The Schr\"odinger equation is still written in the form (IV.87), but $D_\a$ is now a
matrix of differential operators, namely ${\part  \over \part x^\a }
\cdot \un -  {ie \over \hbar}  A_\a (x)$ where $A_\a (x) $ is an $N$-by-$N$ hermitian matrix.
In the path integral (IV.89) replace the exponential factor by 
$$\Ta \exp \! 
\left({ {i e  \l  \over \hbar  }
\int_{ t_a}^{t_b }
A_\a \! \left({ x_b + \l  z(t)  }\right) dz^\a (t)
 }\right) \, . 
\eqno ({\mathop{\rm IV.98 }\nolimits} )$$

\bigskip
c) We could consider gauge groups over a curved manifold and combine the results of paragraph
IV.3 with those of the present paragraph. 

\bigskip
\nn
{\bf 6. A symplectic manifold.}\footnote{$^{19}$}{\nn Contributed by John LaChapelle.}

\smallskip
Let $N$ be a symplectic manifold $\Ma$ of dimension $d = 2n$. The manifold $\Ma$ represents
the classical phase space of a physical system -- usually the cotangent bundle $T^* Q$ of a
configuration space $Q$ but not necessarily. Paths in $\Ma$ have $n$ initial and $n$ final boundary conditions. These boundary conditions
are consistent with the requirement of quantum uncertainty, and they imply a choice of
polarization\footnote{$^{20}$}{\nn Roughly speaking, a polarization is a foliation of $\Ma$ whose
leaves are Lagrangian submanifolds of dimension $n$. }.
 Hence, it is possible to cover $\Ma$ with a family of open subsets $
\left\lbc{ U_i }\right\rbc$ such that each $U_i$ is diffeomorphic to a product of two
transverse Lagrangian submanifolds $L_i \ts L'_i$. For simplicity consider a symplectic
manifold which admits global transverse Lagrangian submanifolds $L$ and $L'$,
and identify $\Ma$ with $L \ts L'$. Set $x(t) $ in the form 
$\left({ x_L (t), x_{L'} (t)}\right)$ with $x_L (t)$ in $L$ and $x_{L'}(t)$ in $L'$, in such a way
 that
$x_L (t_b) = x_b$ and  that 
$x_{L'}(t_a) = x_a$.

\bigskip
We require the group generated by the transformations $\Si ({\bf T}, z)$ to be a subgroup of the group of
symplectomorphisms which leave the polarization invariant. This implies that the set
$ \left\lbc{ X_{(\a)} }\right\rbc$ is of the form
$\left\lbc{ X_{(a)},  X_{(a')} }\right\rbc$ such that
$ X_{(a)} \! \left({ x_{L'} (t) }\right) = 0$,
$X_{(a')} \! \left({ x_{L} (t) }\right) = 0$, and 
$\left\lbk{ X_{(a)}, X_{(a')} }\right\rbk = 0$. Here 
$a \in \left\lbc{ 1, \cdots, k }\right\rbc$,
$a' \in \left\lbc{k +  1, \cdots, d }\right\rbc$, and $k$  is a fixed integer between $1$ and
$d$,  possibly, but not necessarily, equal to $d/2$.
Consequently, a path satisfies the differential equations
$$
\left\{\matrix{
dx_L (t)
=  X_{(a)} \! \left({ x_{L} (t) }\right) dz^a + 
Y \! \left({ x_{L} (t) }\right) dt
  \hfill\cr
\noalign{\medskip}
 dx_{L'} (t)
=  X_{(a')} \! \left({ x_{L'} (t) }\right) dz^{a'} + 
Y \! \left({ x_{L'} (t) }\right) dt  \, ,
  \hfill \cr}
\right.
\eqno ({\mathop{\rm IV.99 }\nolimits} )$$
and the general formula (II.1) becomes
$$
\left.\matrix{ 
\left({ U_{t_b, t_a} \phi }\right)
\! \left({ x_b, x_a }\right) := 
\int_{{\bf Z}_{L} }^{ } \int_{{\bf Z}_{L'} }^{ }
\Da z_L \, \Da z_{L'} \cdot 
e^{
 - \pi Q_0 \! \left({ z_L, z_{L'} }\right) /s
}
\cr
\noalign{\medskip}
\ts \phi
\! \left({ x_b \cdot \Si
\! \left({ t_a, z_L }\right), x_a \cdot \Si
\! \left({ t_b, z_{L'} }\right)
 }\right) \, , 
  \cr}\right.
\eqno ({\mathop{\rm IV.100 }\nolimits} )
$$
where now $h^{\a \b}
= \pmatrix{
h^{ab} & 0 \cr
0 & h^{a'b'} \cr
}$. Here ${\bf Z}_L$ is the space of paths 
$z_L : {\bf T} \ra \RR^k$ such that $z_L (t_b) = 0$, and 
${\bf Z}_{L'}$ is the space of paths $z_{L'} : {\bf T} \ra \RR^{d-k}$
such that $z_{L'}(t_a) = 0$. Each functional integral separately satisfies a partial differential equation:
$$
\left\{\matrix{
 {\ds\part \Psi_L \over \ds \part t_a} 
=
 {\ds s \over \ds 4 \pi } 
h^{ab} \La_{X_{(a)}} \La_{X_{(b)}} \Psi_L + \La_{Y} \Psi_L
  \hfill\cr
\noalign{\medskip}
  {\ds\part \Psi_{L'} \over \ds \part t_b} 
=
 {\ds s \over \ds 4 \pi } 
h^{a'b'} \La_{X_{(a')}} \La_{X_{(b')}} \Psi_{L'} + \La_{Y} \Psi_{L'} \, , 
  \hfill \cr}
\right. \eqno ({\mathop{\rm IV.101 }\nolimits} )
$$
where $\Psi_L := 
\left.{ \Psi }\right\vert_{L }$
and 
$\Psi_{L'} := 
\left.{ \Psi }\right\vert_{L' }$.

\vfill\eject
\nn
6.1. {\it The case of a cotangent bundle.}
\smallskip
Choosing an initial function $\phi $ is choosing a transition amplitude.
We consider the case where $\Ma$ is the cotangent bundle of a flat configuration space $Q$,
hence $\Ma = Q \ts P$ where $Q$ and $P$ are finite-dimensional vector spaces in duality.
In order to define {\it position-to-position transition amplitudes}, we choose
 the initial function of the form
$\phi (q, p) = \d \! \left({q  - q_a }\right)
$.
Equation (IV.100) yields
$$\Ka \! \left({  q_b, t_b ; q_a, t_a 
 }\right) =
\int_{{\bf Z}_Q }^{ } \int_{{\bf Z}_P }^{ }
\Da z_Q \, \Da z_P \cdot 
\exp \! 
\left({ - {\pi  \over s} Q_0  \! \left({ z_Q, z_P }\right) }\right)
h \! \left({ z_Q, z_P  }\right)
\eqno ({\mathop{\rm IV.102 }\nolimits} )$$
where the integrand is given by 
$$
h \! \left({ z_Q, z_P  }\right)
=
\d \! 
\left({ q_b \cdot \Si \! \left({ t_a, z_Q }\right) - q_a }\right)
\, . \eqno ({\mathop{\rm IV.103 }\nolimits} )$$

\medskip
In order for the transition amplitudes to be consistent with the initial wave function, we require
$\build { \lim }_{ t_a \ra t_b }^{ }
\Ka \! \left({ q_b, t_b ; q_a, t_a }\right)$ to be equal to $\d \! 
\left({ q_b - q_a }\right)$. But the integrand 
$h \! \left({ z_Q, z_P }\right)$ tends to 
$\d \! \left({ q_b \cdot \Si \! \left({ t_b, z_Q }\right) - q_a}\right)
= \d \! \left({ q_b - q_a }\right)$ when $t_a$ tends to $t_b$, a limit independent of 
$z_Q$, $z_P$. Hence, provided we can interchange limit and integration, we get
$$\build { \lim }_{ t_a \ra t_b }^{ }
\Ka \! \left({ q_b, t_b ; q_a, t_a }\right) =
\int_{{\bf Z}_Q }^{ } \int_{{\bf Z}_P }^{ }
\Da z_Q \, \Da z_P \cdot
\exp \! 
\left({ - {\pi  \over s} Q_0  \! \left({ z_Q, z_P }\right) }\right)
\d \! \left({ q_b - q_a }\right) \, .
$$
Hence this is equal to $\d \! \left({ q_b - q_a }\right)$ by the normalization of our integrator.

\medskip
We handle the other cases in a similar way.

\smallskip
a) {\it Momentum-to-position amplitude}
$\Ka \! \left({ q_b, t_b ; p_a, t_a }\right)$: use the initial function $
 h^{-n/2} e^{i q_b \cdot p / \hbar }$
and the integrand
$$ h^{-n/2}
\exp \! 
\left({ { i \over \hbar} q_b \cdot 
\left({ p_a \cdot \Si \! \left({ t_b, z_P }\right) }\right)
 }\right)
\, . \eqno ({\mathop{\rm IV.104 }\nolimits}_{{\mathop{\rm a }\nolimits}} )
$$

\smallskip
b) {\it Position-to-momentum amplitude}
$\Ka \! \left({ p_a, t_a ; q_b, t_b }\right)$: use the initial function $
h^{-n/2} e^{- i q \cdot p_a / \hbar }$
and the integrand
$$ h^{-n/2}
\exp \! 
\left({ - { i \over \hbar} 
\left({ q_b \cdot \Si \! \left({ t_a, z_Q }\right) }\right) \cdot p_a
 }\right)
\, . \eqno ({\mathop{\rm IV.104 }\nolimits}_{{\mathop{\rm b }\nolimits}} )$$

\smallskip
c) {\it Momentum-to-momentum amplitude}
$\Ka \! \left({ p_b, t_b ; p_a, t_a }\right)$: use the initial function $\d \! \left({ p - p_b }\right)$
and the integrand
$$\d \! 
\left({ p_a \cdot \Si \! \left({ t_b, z_P }\right) - p_b}\right) 
\, . \eqno ({\mathop{\rm IV.104 }\nolimits}_{{\mathop{\rm c }\nolimits}} )
$$

\medskip
\nn
It follows easily from these definitions that
$$
\eqalignno{
\Ka \! \left({ p_b, t_b ; p_a, t_a }\right) 
&=
h^{- n/2 }
\int_{Q }^{ } dq_b \  e^{- i q_b \cdot p_b / \hbar }  \Ka \! \left({ q_b, t_b ; p_a, t_a }\right)
& ({\mathop{\rm IV.105 }\nolimits} )  \cr
\Ka \! \left({ q_b, t_b ; q_a, t_a }\right) 
&=
h^{- n/2 }
\int_{P }^{ } dp_b \  e^{ i q_b \cdot p_b / \hbar } \Ka \! \left({ p_b, t_b ; q_a, t_a }\right)  \, . &
({\mathop{\rm IV.106 }\nolimits} ) \cr} $$

\bigskip
\nn
6.2. {\it Coherent states.}

\smallskip
More general transition amplitudes, which cannot be interpreted as posi-tion-to-momentum
transitions, are possible by choosing more complicated initial wave functions, different
polarizations, and/or by having non-trivial phase spaces. For instance, coherent state
transitions can be calculated -- given a suitable characterization of coherent states.

\medskip

Choose a K\"ahler polarization on a non-trivial symplectic manifold, so that $\Ma$ is a
K\"ahler manifold $M$ of complex dimension $n$. It is convenient (though not necessary) to take the phase space
to be the product manifold $ M \ts \oo M$ with coordinates $(\z, \oo z)$ and complex dimension $2n$.
Following the work of Berezin [21] and Bar-Moshe and Marinov [22], we use generalized coherent
state wave functions
$\phi_{\z '} (\z) = \exp  K(\z, \oo \z') $
where 
$K(\z, \oo \z')$ is the K\"ahler potential.
Consider the function $\phi : M \ts \oo M \ra \CC  $ given by 
$\phi (\z, \oo \z') = 
e^{K \! \left({ \z, \oo \z_a }\right) }
$
and take  $z_M : {\bf T} \ra \CC^n$ such that $z_M (t_b) = 0$ and 
$z_{\oo M} : {\bf T} \ra \CC^n$ such that $z_{\oo M}(t_a) = 0$.
Then
$$ \Ka \!  
\left({ \z_b, t_b ; \oo \z_a, t_a  }\right): = \int_{ {\bf Z}_M }^{ }  \int_{ {\bf Z}_{\oo M} }^{ }
\Da z_M \, \Da z_{\oo M} \cdot  
e^{
 - \pi   Q_0  \left({ z_M, z_{\oo M } }\right) /s
}
e^{K
\! \left({  \z_b \cdot \Si  \left({ t_a, z_{ M} }\right), \oo \z_a  }\right)
}  
$$
 are generalized coherent state transition amplitudes.
Note that
$\Ka 
 \!  
\left({\oo \z_b, t_b ;  \z_a, t_a  }\right) $ \break $=
\oo { \Ka  \!  
\left({ \z_b, t_b ; \oo \z_a, t_a  }\right) } $,
and 
$\build {\lim }_{ t_a \ra t_b  }^{ }
\Ka 
 \!  
\left({ \z_b, t_b ; \oo \z_a, t_a  }\right)
=
\exp  K \! \left({ \z_b, \oo \z_a }\right) 
$ provided the limit can be taken inside the integral.

\bigskip
\nn
{\bf 7. A simple model of the Bohm-Aharonov effect.}

\smallskip
We illustrate  on an elementary example various techniques expounded in this section.
Our basic manifold $N$ is the euclidean plane $\RR^2$ with the origin removed. As in paragraph IV.2, we denote
the cartesian coordinates by $x^1, x^2$. The universal covering of $N$ is the set $\ww N$ of pairs of real
numbers $r, \t$ with $r > 0$, and the covering map $\Pi : \ww N \ra N $ is described by equation (IV.31), namely
$x^1 = r \cos \t$, $x^2 = r \sin \t$.

\medskip
We introduce a {\it magnetic potential} $A$ with components 
$$A_1 = - {Fx^2 \over 2 \pi \left\vert{ x }\right\vert^2} 
\quad , \quad  A_2 = 
{Fx^1 \over 2 \pi \left\vert{ x }\right\vert^2} \, .
\eqno ({\mathop{\rm IV.107 }\nolimits} )$$
The magnetic field has one component $B$ perpendicular to our plane, given by
$$B = {\part A_1 \over \part x^2} - {\part A_2 \over \part x^1} \, ,
\eqno ({\mathop{\rm IV.108 }\nolimits} )$$
and an obvious calculation gives $B = 0 $ {\it outside of the origin.}
But the circulation $\oint  A_1 dx^1 + A_2 dx^2 
$ on any loop around the origin is equal to $F$, hence
$$ B = F \d (x) \, .
\eqno ({\mathop{\rm IV.109 }\nolimits} )$$
Physically, we have a wire perpendicular to the plane through its
origin, carrying a magnetic flux equal to $F$.

\medskip
Denote by $m$ the mass and by $e$ the electric charge of a particle. For a path $x : {\bf T} \ra N$, with
components
$x^1 (t)$, $x^2(t)$ at time $t$, the action is given by 
$$S(x) = S_0 (x) + S_M (x) \, ,
\eqno ({\mathop{\rm IV.110 }\nolimits} )$$
with the {\it kinetic action}
$$S_0 (x) = {m \over 2}
\int_{ {\bf T}}^{ }
 { \left\vert{  dx }\right\vert_{}^{2 } \over dt}
\eqno ({\mathop{\rm IV.111 }\nolimits} )$$
and {\it the magnetic action}
$$S_M (x) = {eF \over 2 \pi } \int_{ {\bf T}}^{ }
{ x^1 dx^2 - x^2 d x^1 \over (x^1)^2 + (x^2)^2}
\, . 
\eqno ({\mathop{\rm IV.112 }\nolimits} )$$
We fix a point $x_b$ in $N$, and denote by ${\bf X}_b$ the set of paths $x : {\bf T} \ra N$ for which the
kinetic action $S_0 (x)$ is finite (that is the $L^{2,1}$ paths) with the boundary condition 
$x(t_b) = x_b$. On ${\bf X}_b$, we denote by $\Da x$ the translation-invariant integrator normalized by 
$$\int_{{\bf X}_b}
\Da x \cdot \exp \! 
\left({ { i \over \hbar}  S_0 (x) }\right) = 1 \, . 
\eqno ({\mathop{\rm IV.113 }\nolimits} )$$

\medskip
The free position-to-position transition amplitudes are defined by
$$
\left\langle{ t_b, x_b \bigm| t_a , x_a }\right\rangle_{0}^{}
= \int_{ {\bf X}_b}^{ }
\Da x \cdot \exp \! 
\left({ { i \over \hbar}  S_0 (x) }\right)
\d \!
\left({ x(t_a) - x_a }\right) \, .
\eqno ({\mathop{\rm IV.114 }\nolimits} )$$
According to the calculation made at the end of paragraph IV.2, we get\footnote{$^{21}$}{We express our formulas
directly in terms of the path $x$, and have no need for the scaling $x = \l z$ introduced earlier.} 
$$\left\langle{ t_b, x_b
\bigm| t_a , x_a }\right\rangle_{0}^{} = {- mi \over h(t_b - t_a)} 
\exp \!
\left({ { i \over \hbar}  S_0 (x_{{\mathop{\rm cl }\nolimits}}) }\right)
\, , \eqno ({\mathop{\rm IV.115 }\nolimits} )$$
where $x_{{\mathop{\rm cl }\nolimits}}$ is the classical path of a free particle, namely:
$$x_{{\mathop{\rm cl }\nolimits}} (t) =
{ x_a (t_b - t) + x_b (t - t_a) \over t_b - t_a} 
\, . 
\eqno ({\mathop{\rm IV.116 }\nolimits} )$$
Hence, we obtain explicitely
$$S_0 (x_{{\mathop{\rm cl }\nolimits}}) = 
 {m \over 2} 
{ \left\vert{  x_b - x_a}\right\vert_{}^{2} \over t_b - t_a} \, , 
\eqno ({\mathop{\rm IV.117 }\nolimits} )$$
that is, the WKB approximation (IV.30) is exact in our case.

\medskip
We consider now the transition amplitudes in the given magnetic potential:
$$\left\langle{ t_b, x_b
\bigm| t_a , x_a }\right\rangle_{F}^{}
=
\int_{ {\bf X}_b}^{ }
\Da x \cdot \exp \! 
\left({ { i \over \hbar}  S (x) }\right)
\d \!
\left({ x(t_a) - x_a }\right) \, .
\eqno ({\mathop{\rm IV.118 }\nolimits} )$$
According to (IV.87), we consider the Schr\"odinger equation
$$
{ \part \psi_F \over \part t} =
{ i \hbar \over 2m}
 \left({ D_1^2 + D_2^2 }\right)
\psi_F 
\eqno ({\mathop{\rm IV.119 }\nolimits} )$$
with the differential operators 
$D_\a = {\part  \over \part x^\a} - { ie \over \hbar} A_\a$.
Explicitly, we obtain
$${\part \psi_F \over \part t}
=
{ i \hbar \over 2m} L \psi_F \, ,
\eqno ({\mathop{\rm IV.120 }\nolimits} )$$
with
$$L =
\left({ {\part \over \part x^1} }\right)_{}^{2} +
\left({ {\part \over \part x^2} }\right)_{}^{2} -
 {2 ci \over \left\vert{ x }\right\vert^2 }
\left({ x^1{ \part \over \part x^2} - x^2{ \part \over \part x^1} }\right)
- {c^2 \over \left\vert{ x }\right\vert^2 } \, . 
\eqno ({\mathop{\rm IV.121 }\nolimits} )$$
The constant $c$ is equal to $eF / h$; it is dimensionless. The solution to the equation (IV.119) is given by
$$\psi_F \! \left({ t_b, x_b  }\right) =
\int_{ N }^{ } dx_a 
\left\langle{ t_b, x_b \bigm| t_a, x_a }\right\rangle_F \psi_F \! 
\left({  t_a, x_a}\right) \, ,
\eqno ({\mathop{\rm IV.122 }\nolimits} )$$
a formula equivalent to our more familiar one
$$
\psi_F \!\left({ t_b, x_b  }\right) = 
\int_{ {\bf X}_b}^{ }
\Da x \cdot \exp \! 
\left({ { i \over \hbar}  S (x) }\right)
\phi \!
\left({ x(t_a)  }\right)
\eqno ({\mathop{\rm IV.123 }\nolimits} )$$
if we take into account the initial wave function $\phi(x) = 
\psi_F \!\left({ t_a, x }\right)$ at time $t_a$.

\medskip
We lift now everything to the universal covering $\ww N$. By means of the formulas (IV.31), the wave function is
now a function $\ww \psi_F (t, r, \t)$ with the restriction
$$\ww \psi_F (t, r, \t + 2 \pi) = \ww \psi_F (t, r, \t)
\eqno ({\mathop{\rm IV.124 }\nolimits} )$$
that is, a function on ${\bf T} \ts \ww N$ invariant under the group $\ZZ$ acting on $\ww N$ by 
$(r , \t)  \cdot n = (r , \t + 2 \pi n)$. The Schr\"odinger equation (IV.120) keeps its form with
 $L$ changed into the new
operator
$$
\ww L = 
 {\part^2 \over \part r^2} +  {1 \over r^2}   {\part^2\over \part \t^2} +
{1 \over r}   {\part \over \part r}
-  {2ci \over r^2}  {\part  \over \part \t} -  {c^2 \over r^2} \, .
\eqno ({\mathop{\rm IV.125 }\nolimits} )$$
This operator can be written as $\ww L = e^{c i \t} \, \D e^{- ci \t}$ where
$\D $ is the Laplacian in polar coordinates, namely:
$$
\D = {\part^2 \over \part r^2} +  {1 \over r^2}   {\part^2\over \part \t^2} +
{1 \over r}   {\part \over \part r} \, .
\eqno ({\mathop{\rm IV.126 }\nolimits} )$$
Hence the new wave function
$$\ww \psi_0 (t, r, \t) = e^{- ci \t }
\ww \psi_F (t, r, \t)
\eqno ({\mathop{\rm IV.127 }\nolimits} )$$
satisfies the equation
$$
 {\part \ww \psi_0 \over \part t} =
 {i \hbar \over 2m}  \D \ww \psi_0 \, ,
\eqno ({\mathop{\rm IV.128 }\nolimits} )$$
that is the Schr\"odinger equation for a free particle, written in polar coordinates. As we shall see in a
moment, equation (IV.127) expresses a {\it gauge transformation} which apparently removes the magnetic
potential. But the periodicity condition (IV.124) transforms into
$$\ww \psi_0 (t, r , \t + 2 \pi) = e^{- 2 \pi i c}
\ww \psi_0 (t, r, \t) \, ,
\eqno ({\mathop{\rm IV.129 }\nolimits} )$$
hence the wave function $\ww \psi_0$ on $\ww N$ is not the lifting of a function $\psi_0 $ on $N$.
Rather, the linear representation $n \mps e^{ 2 \pi i n c }$ of the group $\ZZ$ into $U(1)$ defines a line
bundle on $N = \RR^2 \bl \{ 0 \} $ and $\ww \psi_0$ {\it corresponds to a section $\psi$ of this line
bundle.}
Otherwise stated, despite the fact that the magnetic field $B = \part A_1 / \part x^2 - \part A_2 / \part x^1$
is identically zero on $N$, there is no single-valued function $R$ on $N$ such that $A_\a = - \part R / \part
x^\a$: {\it the space $N$ is multiply connected.}

\bigskip
To explain gauge transformations, let us introduce the principal bundles 
$P = N \ts U(1)$ and $\ww P = \ww N \ts U(1)$. According to the general formulas the connection form $\om_F$ on
$P$ is given by
$$\om_F =
i 
\left({d \T -   {e \over \hbar } A}\right)
= i \left({d \T - c { x^1 dx^2 - x^2 d x^1\over 
\left({ x^1 }\right)^2 + \left({ x^2 }\right)^2
} }\right)
\, ,
\eqno ({\mathop{\rm IV.130 }\nolimits} )$$
where $\T$ is the angular coordinate on $U(1)$ (taken modulo $2\pi$). When there is no magnetic field, it reduces
to 
$\om_0 = id\T$. We cannot transform $\om_F$ into $\om_0$
 by a gauge transformation, but if we lift these differential forms to $\ww N$, we
obtain
$$\ww \om_F = i(d \T - cd\t) \quad , \quad \ww \om_0 = id \T \, .
\eqno ({\mathop{\rm IV.131 }\nolimits} )$$
The transformation $\ww U$ of $\ww P$ into $\ww P$ taking $(r, \t, \T)$ into $(r, \t, \T + c\t)$ is an
automorphism of $U(1)$-bundle, and 
$\ww U^* \ww \om_F = \ww \om_0$, but $\ww U$ is not the lifting of an automorphism $U$ of $P$, 
{\it except when $c$ is an integer.}

\bigskip
To the wave function $\psi_F (x) $ on $N$ we associate the function
$$\Psi_F (x, \T) = e^{- i \T} \psi_F (x)
\eqno ({\mathop{\rm IV.132 }\nolimits} )$$
on $P$. The Schr\"odinger equation (IV.119) translates into our standard form
$$
 {\part \Psi_F \over \part t} 
=  {i \over 4 \pi} h^{\a \b}
\La_{X_{(\a)}^F} \La_{X_{(\b)}^F} \Psi_F
\eqno ({\mathop{\rm IV.133 }\nolimits} )$$
provided we take $h_{\a \b}$ equal to $m \d_{\a \b}/h$, with the vector fields
$$
\left\{\matrix{
\La_{X_{(1)}^F} = {\ds \part  \over \ds\part x^1} - {\ds cx^2 \over \ds\left\vert{ x }\right\vert^2 }
{ \ds\part \over \ds\part \T} 
  \hfill\cr
\noalign{\medskip}
\La_{X_{(2)}^F} = {\ds\part  \over \ds\part x^2} + {\ds cx^1 \over \ds\left\vert{ x }\right\vert^2 }
{\ds \part \over \ds\part \T}   
\, .
  \hfill \cr}
\right.
\eqno ({\mathop{\rm IV.134 }\nolimits} )$$
We lift now everything to $\ww P$, so $\psi_F$ lifts to $\ww \psi_F, \ \Psi_F $ to $\ww \Psi_F$ and the lifted
vector fields are given by\footnote{$^{22}$}
{For $F = 0$, they reduce to the vector fields given in equations (IV.36) and (IV.37).}
$$
\left\{\matrix{
\La_{\ww X_{(1)}^F} = \cos \t \cdot {\ds\part  \over \ds\part r} - {\ds\sin \t \over\ds r } \cdot 
\left({ {\ds \part \over \ds\part \t} + c {\ds \part \over \ds\part \T} }\right) 
  \hfill\cr
\noalign{\medskip}
\La_{\ww X_{(2)}^F} = \sin \t \cdot {\ds\part  \over \ds\part r} + {\ds\cos \t \over \ds r } \cdot 
\left({ {\ds \part \over \ds\part \t} + c {\ds \part \over \ds\part \T} }\right)   
\, .
  \hfill \cr}
\right.
\eqno ({\mathop{\rm IV.135 }\nolimits} )$$
The Schr\"odinger equation (IV.133) remains the same with $\Psi_F$ replaced by $\ww \Psi_F$ and
$X_{(\a)}^F $ by $\ww X_{(\a)}^F$.

\medskip
We describe the effect of the gauge transformation $\ww U$ given by 
$\ww U (r, \t, \T) = (r, \t, \T + c \t)$.
We noticed that 
$\ww U^* \ww \om_F $ is equal to $\ww \om_0$, and similarly, we get 
$\ww \Psi_F \circ \ww U = \ww \Psi_0$ where $\ww \Psi_0$ is defined by analogy to $\ww \Psi_F$, namely
$$\ww \Psi_0 (r, \t, \T) =
e^{- i \T }
\ww \psi_0 (r, \t) \,.
\eqno ({\mathop{\rm IV.136 }\nolimits} )$$
Dually, $\ww U$ transforms the vector field $\ww X_{(\a)}^0$ into 
$\ww X_{(\a)}^F$ for $\a = 1$ or $2$.
Hence $\ww \Psi_F$ is a solution of the Schr\"odinger equation (IV.133) if, and only if,
$\ww \Psi_0$ is a solution of the corresponding equation for $F = 0$.

\medskip
We conclude this paragraph by a discussion of path integrals.
We parame-trize paths in ${\bf X}_b$ by paths in ${\bf Z}_b$, using the correspondence
$$x(t) = x_b + z(t) \, . 
\eqno ({\mathop{\rm IV.137 }\nolimits} )$$
Fix $r_b$ and $\t_b$ in such a way that
$$x_b = 
\left({ r_b \cos \t_b, r_b \sin \t_b }\right)
\, . 
\eqno ({\mathop{\rm IV.138 }\nolimits} )
$$
Then $x$ can be lifted in a unique way to a path $\ww x : {\bf T} \ra \ww N$ of the form 
$\ww x (t) = (r(t), \t (t))$ such that 
$\ww x (t_b) = 
\left({ r_b, \t_b }\right)$,
that is 
$r(t_b) = r_b $ and $\t (t_b) = \t_b$.
The initial point $x_a = x(t_a)$ is lifted to $\ww x (t_a) =
\left({ r_a, \t_a }\right)$ where 
$r(t_a) = r_a$, $\t(t_a) = \t_a$.
Moreover, the magnetic action of the original path is given by 
$$S_M (x) / \hbar = c \cdot 
\left({ \t_b - \t_a }\right)
\, . 
\eqno ({\mathop{\rm IV.139 }\nolimits} )$$

\smallskip
\nn
Consider the space ${\bf X}_{a, b}$ of paths $x : {\bf T} \ra N$ such that 
$x(t_a) = x_a$, $x(t_b) = x_b$. {\it Then two such paths are in the same homotopy class if, and only if,
the lifted paths in $\ww N$ correspond to the same determination of the angular coordinate
$\t_a$ of $x_a$.
An equivalent condition is that they have the same magnetic action.}

\medskip
Using the connection of the principal bundle $\ww P$ over $\ww N$, we can define the horizontal lifting of
$\ww x$; it is the unique path of the form 
$y(t) = (r(t), \t(t), \T(t))$ on which $\ww \om_F$ induces a zero form, that is 
$\T(t) - c \t(t)$ is a constant in time. For the phase factors, we get
$$e^{i \T (t_b)} =
e^{i \T (t_a)} e^{i S_M (x) / \hbar}
\, . 
\eqno ({\mathop{\rm IV.140 }\nolimits} )$$
Moreover, the horizontal lifting is a solution of our standard differential equation $dy =
X_{(\a)} (y) \cdot dz^\a$.

\medskip
>From all this, it follows that the transition amplitudes as defined in the Feynman way (IV.118) agree with our
standard expression (see for instance (IV.38)). The result obtained in paragraph IV.2 (see formulas (IV.41) and
(IV.47)) can be restated as
$$
\left\langle{ t_b, x_b \bigm| t_a, x_a  }\right\rangle_0
= { 1 \over r_a}
\sum_{ n\in \ZZ }^{ } 
\left\langle{ t_b, r_b,\t_b \bigm| t_a, r_a, \t_a + 2 n \pi}\right\rangle
\, . 
\eqno ({\mathop{\rm IV.141 }\nolimits} )$$
A similar reasoning yields a more general formula:
$$
\left\langle{ t_b, x_b\bigm| t_a, x_a  }\right\rangle_F
= { 1 \over r_a}
\sum_{ n\in \ZZ }^{ }  e^{ic (\t_b - \t_a - 2n \pi) }
\left\langle{ t_b, r_b,\t_b \bigm| t_a, r_a, \t_a + 2 n \pi}\right\rangle
\, . 
\eqno ({\mathop{\rm IV.142 }\nolimits} )$$
We can invert this formula\footnote{$^{23}$}{Notice that the transition probabilities are functions of $F$ which
admit a period $h/e$, a well-known physical effect.}
and obtain
$$\left\langle{ t_b, r_b, \t_b \bigm| t_a, r_a, \t_a  }\right\rangle
= r_a \int_{ 0}^{ 1} dc\,  e^{- ic \left({ \t_b- \t_a }\right) }
\left\langle{ t_b, x_b\bigm| t_a, x_a }\right\rangle_{hc/e} \, .
\eqno ({\mathop{\rm IV.143 }\nolimits} )$$
Hence, {\it from the knowledge of the transition amplitude 
$\left\langle{ t_b, x_b\bigm| t_a, x_a }\right\rangle_{F}$
as a function of the magnetic flux $F$, we can infer the value of the transition amplitudes in polar
coordinates.  }

\vfill \eject
\centerline {{\bf Acknowledgements.} }

\bigskip
Our collaboration has been made possible by the financial support of the NATO collaborative
Research Grant \# 910101, the University (of Texas) Research Institute grant R-270, the Jane
and Roland Blumberg Centennial Professorship in Physics, and the Institut des Hautes Etudes
Scientifiques (I.H.E.S.) at Bures-sur-Yvette.

\medskip
The Center for Relativity of the University of Texas, the Ecole Normale Sup\'erieure, Paris
and the Institut des Hautes Etudes Scientifiques, Bures-sur-Yvette have given us ideal
working conditions.

\medskip
The manuscript which went through  several revisions has been diligently prepared by
Marie-Claude Vergne at I.H.E.S. and Debbie Hajji at the Center for Relativity.

\vfill\eject
\vglue 1cm
\centerline{J\bf Appendix A } 
\medskip
\centerline {\titre   Functional integration}
\bigskip
\bigskip

In this appendix, we develop the basic properties of our integrators. In the oscillating case $(s
= i)$, our theory is, up to some inessential changes in notation, the same as the one expounded by
Albeverio and H\o egh-Krohn in [9]. The introduction of a parameter $s$ equal to $1$ or $i$
enables us to treat in a unified way the oscillating integrators of Fresnel type
$e^{\pi i Q(x) } \Da x$ and the Gaussian integrators of type 
$e^{- \pi Q(x)} \Da x$ (with $Q(x) > 0$ for $x \ne 0$ in the latter case).
We content ourselves by giving here the {\it basic formulas} and the {\it computational tools}.
In the Gaussian case, we would have to justify these formal manipulations, since, as it is
well-known, the $L^{2, 1}$ {\it functions are a set of measure $0$ for the Wiener measure}, and 
{\it we
are considering functions on this set of measure $0$.} Our claims can be fully vindicated, but we
defer the complete justification to another publication.

\bigskip
\nn
{\bf 1. Gaussian integrators.}
\smallskip
\nn
1.1. {\it The setup.}
\smallskip

We denote by  ${\bf X}$ a real, separable, Banach space, by ${\bf X'}$ its dual and by 
$\left\langle{ x', x }\right\rangle$ (or sometimes
$\left\langle{ x', x }\right\rangle_{{\bf X} }$) the duality between 
${\bf X}$ and ${\bf X'}$. We suppose given a continuous linear map
$D : {\bf X} \ra {\bf X'}$ with the following properties:

\smallskip
-- (symmetry) $\left\langle{ Dx, y }\right\rangle =
\left\langle{ Dy, x }\right\rangle$ for $x, y $ in ${\bf X}$;

\smallskip
-- (invertibility) there exists a continuous linear map
$G : {\bf X'} \ra {\bf X}$ inverse of $D$, that is $DG = \un$ and $GD = \un$. 

\smallskip
\nn
Out of these data one constructs two quadratic forms, $Q$ on ${\bf X}$ and $W$ on ${\bf X'}$, by
the rules
$$Q(x)  = 
\left\langle{ Dx, x }\right\rangle \quad, \quad
W(x') = \left\langle{ x', Gx' }\right\rangle. 
\eqno ({\mathop{\rm A.1 }\nolimits})$$

\smallskip
\nn
They are related to each other as follows:
$$Q(x)  = 
W (Dx) \quad, \quad
W(x') = Q(Gx')
\eqno ({\mathop{\rm A.2 }\nolimits})$$
for $x$ in ${\bf X}$
and $x'$ in ${\bf X'}$.

\bigskip
We denote by $s$ a parameter equal to $1$ or $i$. The function 
$e^{- \pi s W}
$ on ${\bf X'}$ is continuous. It is bounded in the following cases:

\smallskip
-- $s = i$;

\smallskip
-- $s = 1$ and $W(x') > 0$ for 
$x' \ne 0$ in ${\bf X'}$. Equivalently, by (A.2), $s = 1$ and $Q(x) > 0$
 for $x \ne 0$ in ${\bf X}$.

\bigskip
\nn
1.2. {\it The oscillatory case $(s = i)$.}
\smallskip
The integrator $\Da x$ is characterized by the following integration formula:
$$
\int_{ {\bf X} }^{ }
\Da x \cdot \exp \!
\left({ \pi i Q(x) - 2 \pi i 
\left\langle{ x', x }\right\rangle
 }\right) =
\exp (- \pi i W (x'))
\eqno ({\mathop{\rm A.3 }\nolimits})$$
for every $x'$ in ${\bf X'}$. This relation should be interpreted as follows.

\bigskip
Since the metric
space ${\bf X'}$ is complete and separable, we know the notion of a complex bounded Borel
measure\footnote{$^{25}$}{That is: a $\s$-additive functional from the $\s$-algebra of Borel
subsets of ${\bf X'}$ into the complex numbers (see [23] for instance). We use a simplified
terminology by refering to $\mu$ as a ``measure'' on ${\bf X'}$.} $\mu$ on ${\bf X'}$.
The {\it Fourier-Stieltjes transform} $\Fa \mu$ of $\mu$ is given the customary definition:
$$\left({ \Fa \mu }\right)
(x) = \int_{{\bf X'}}^{ } d \mu (x') \,  e^{- 2 \pi i 
\left\langle{ x', x }\right\rangle } .
\eqno ({\mathop{\rm A.4 }\nolimits})$$ 
It is a continuous bounded function on ${\bf X}$.

\bigskip
We denote by\footnote{$^{26}$}{The initial $\Fa$ stands for Fresnel or Feynman according to the
worshipping habits of the reader.} $\Fa({\bf X}) $ the set of functions 
on ${\bf X}$ of the form $\Phi_\mu =
e^{\pi i Q} \cdot \Fa \mu$, where $\mu$ runs over the  measures on
${\bf X'}$. Since the map $\mu \mps \Fa \mu$ is injective, $\Fa({\bf X})$ is a Banach space, the
norm of the function $\Phi_\mu$ being taken equal to the total variation\footnote{$^{27}$}{According
to the standard definition, this is the l.u.b. of the set of numbers
$\ds \sum_{i =1 }^{ p} \left\vert{ \mu (A_i)  }\right\vert$ where 
$\left({ A_1, \cdots, A_p }\right)$ runs over the set of all partitions of ${\bf X'}$ into a
finite number of Borel subsets.
} 
 $\Var (\mu)$ of $\mu$.
On this Banach space, one defines a continuous linear form, denoted as an integral, by
$$\int_{ {\bf X}  }^{ } \Da x \cdot  \Phi_\mu (x) =
\int_{ {\bf X'}  }^{ }
d \mu  (x') \, e^{- \pi i W (x') }.
\eqno ({\mathop{\rm A.5 }\nolimits})$$

\smallskip
\nn
Hence by definition, we have
$$
\int_{ {\bf X}  }^{ } \Da x 
\int_{ {\bf X'}  }^{ }
d \mu (x') \exp \!
\left({  \pi i Q(x) - 2 \pi i 
\left\langle{ x', x }\right\rangle
}\right)
=
\int_{ {\bf X'}  }^{ }
d \mu  (x') \exp (- \pi i W(x')) \, .
\eqno ({\mathop{\rm A.6 }\nolimits})$$

\smallskip
\nn
Formally, this is the equation obtained by integrating equation (A.3) w.r.t. $d \mu (x') $ and
then interchanging the integrations:
$$\int_{ {\bf X'}  }^{ }
d \mu  (x') \int_{ {\bf X}  }^{ } \Da x  =
\int_{ {\bf X}  }^{ } \Da x  
\int_{ {\bf X'}  }^{ }
d \mu  (x') \, .$$

\bigskip
The space $\Fa ({\bf X})$ of Feynman-integrable functionals on ${\bf X}$
 is invariant under translations by elements of ${\bf X}$, and so is the integral, namely:
$$\int_{ {\bf X}  }^{ } \Da x  \cdot F (x) =
\int_{ {\bf X}  }^{ } \Da x \cdot 
F \! \left({ x + x_0 }\right)
\eqno ({\mathop{\rm A.7 }\nolimits})$$ 
or, in shorthand notation, 
$\Da x = \Da \! \left({ x + x_0 }\right)$, for any fixed element $x_0$ of ${\bf X}$.

\bigskip
According to equations (A.1) and (A.2), formula (A.3) can be rewritten as
$$
\int_{ {\bf X}  }^{ } \Da x  \cdot e^{\pi i Q(x - Gx')}
= 1. 
\eqno ({\mathop{\rm A.8 }\nolimits})$$ 

\smallskip
\nn
Hence {\it assuming  the invariance under translations
 of the integral}, the normalization of $\Da x$ is achieved by
$$
\int_{ {\bf X}  }^{ } \Da x \cdot 
 e^{\pi i Q(x )}
= 1. 
\eqno ({\mathop{\rm A.9 }\nolimits})$$ 

\bigskip
\nn
1.3. {\it The positive case $(s = 1)$.}
\smallskip
Henceforth, we assume that $s = 1$ and that $Q$ is positive-definite (that is $Q(x) > 0$ for $x
\ne 0$ in ${\bf X}$).
We take as basic integration formula
$$\int_{ {\bf X}  }^{ } \Da x 
\cdot  \exp \!
\left({ - \pi  Q(x) - 2 \pi i 
\left\langle{ x', x }\right\rangle
}\right) = \exp (- \pi W (x')) \, .
\eqno ({\mathop{\rm A.10 }\nolimits})$$ 
We can interpret this relation as above.
The conclusion is the integration formula
$$
\int_{ {\bf X}  }^{ } \Da x 
\int_{ {\bf X'}  }^{ }
d \mu (x') \exp \!
\left({  - \pi  Q(x) - 2 \pi i 
\left\langle{ x', x }\right\rangle
}\right)
=
\int_{ {\bf X'}  }^{ }
d \mu  (x') \exp (- \pi  W(x')) \, .
\eqno ({\mathop{\rm A.11 }\nolimits})$$ 

\smallskip
\nn
But the space $\Fa ({\bf X})$ is no longer invariant under translations. 
To restore this invariance, we have to replace Fourier-Stieltjes transforms by
{\it Laplace-Stieltjes transforms}. 
For that purpose, we have to consider the complex dual space
${\bf X}'_\CC$
consisting of the continuous real-linear maps
$x' : {\bf X}\ra \CC$ and measures $\mu$ on 
${\bf X}'_\CC$.
We leave the details to a forthcoming publication.

Formulas (A.3) and (A.10) are obtained as the specializations of formula (I.7) for $s = i$ and $s
= 1$ respectively.

\bigskip
\nn
{\bf 2. Linear changes of variables.}
\smallskip
\nn
2.1. {\it The case of Lebesgue integrals.}
\smallskip
Let ${\bf Y}$ be another real, separable, Banach space, and let 
$L : {\bf X} \ra {\bf Y}$ be a Borel-measurable map.
If $\om_{\bf X}$ is a measure on ${\bf X}$, its image under $L$ is the measure
$\om_{\bf Y}$ on ${\bf Y} $ defined by 
$\om_{\bf Y} (B) = \om_{\bf X} \!
\left({ L^{- 1} (B) }\right)$ for any Borel subset 
$B$ of ${\bf Y}$. In functional terms, this definition is tantamount to
$$
\int_{ {\bf X}  }^{ }  d \om_{\bf X} (x) \, 
g(L(x)) =
\int_{ {\bf Y}  }^{ }  d \om_{\bf Y} (y) \, g(y)
\eqno ({\mathop{\rm A.12 }\nolimits})$$ 
for any bounded (or non-negative) Borel-measurable function $g$ on ${\bf Y}$.

\bigskip
Assume now that $L$ is linear and continuous.
We consider the transpose $\ww L$ of $L$, that is the map
$\ww L : {\bf Y'} \ra {\bf X'}$ such that
$$
\langle \ww L y', x \rangle_{\bf X}
=
\langle y', Lx \rangle_{\bf Y}
\eqno ({\mathop{\rm A.13 }\nolimits})$$
for $x$ in ${\bf X}$ and $y'$ in ${\bf Y'}$.
The measures $\om_{\bf X}$ on ${\bf X}$ and 
$\om_{\bf Y}$ on ${\bf Y}$ being as above, introduce their Fourier-Stieltjes transforms 
$\Fa \om_{\bf X}$ and $\Fa \om_{\bf Y}$ respectively.
Hence $\Fa \om_{\bf X}$ is a bounded continuous function on 
${\bf X'}$, and similarly for 
$\Fa \om_{\bf Y}$ on ${\bf Y'}$.
By definition, we have
$$
\eqalignno{
\Fa \om_{\bf X} (x')
&
= \int_{ {\bf X}  }^{ } d \om_{\bf X} (x) 
\exp \!
\left({  - 2 \pi i 
\left\langle{ x', x }\right\rangle_{{\bf X}}
}\right)&({\mathop{\rm A.14 }\nolimits})\cr
\Fa \om_{\bf Y} (y')
&
= \int_{ {\bf Y}  }^{ } d \om_{\bf Y} (y) 
\exp \!
\left({  - 2 \pi i 
\left\langle{ y', y }\right\rangle_{{\bf Y}}
}\right) \, . &({\mathop{\rm A.15 }\nolimits})\cr}$$

\bigskip
By specializing 
$g(y) = \exp \!
\left({  - 2 \pi i 
\left\langle{ y', y }\right\rangle_{{\bf Y} } }\right)$
into formula (A.12) and taking into account formulas (A.13) to (A.15), we get
$\Fa \om_{\bf Y} (y') =
\Fa \om_{\bf X} (\ww L y')$ for every $y'$ in ${\bf Y'}$. Hence the composition formula:
$$\Fa \om_{\bf Y}  =
\Fa \om_{\bf X} \circ \ww L \, .
\eqno ({\mathop{\rm A.16 }\nolimits})$$ 
[The figure can be found in the Journal of Mathematical Physics {\bf 36} p.~2284]
\vglue 5cm 

\bigskip
The case of a translation is similar.
Assume now that $x_0$ is a given vector in ${\bf X} $
and denote by $T$ the translation taking $x$ into 
$x + x_0$ in ${\bf X}$.
If $\om$ is any measure on ${\bf X}$
and $\om_{x_0}$ its image under the translation $T$, a suitable specialization of formula (A.12)
gives the Fourier transform of $\om_{x_0}$, namely
$$\left({ \Fa \om_{x_0} }\right) (x') =
\exp \! \left({  - 2 \pi i 
\left\langle{ x', x_0 }\right\rangle
}\right)
\cdot \Fa \om (x') \, . 
\eqno ({\mathop{\rm A.17 }\nolimits})$$ 

\smallskip
\nn
2.2. {\it Infinite-dimensional integrators.}
\smallskip
We go back to the setup of paragraph A.1.1. We denote by $\Da x$ the integrator characterized by
$$\int_{ {\bf X}  }^{ } \Da x \cdot 
\exp \! \left({- { \pi \over s } Q(x)  }\right) \cdot 
\exp \! \left({  - 2 \pi i 
\left\langle{ x', x }\right\rangle
}\right)
=
\exp (- \pi s W (x')) \, .
\eqno ({\mathop{\rm A.18 }\nolimits})$$

\medskip
\nn
Formally, this means that the integrator $\Da \om $ defined by 
$$\Da \om (x) = \exp \!
\left({- { \pi \over s } Q(x)  }\right) 
\cdot \Da x
\eqno ({\mathop{\rm A.19 }\nolimits})$$ 
has a Fourier transform equal to 
$\exp (- \pi sW)$, namely:
$$\int_{ {\bf X}  }^{ } \Da \om (x) 
\exp \! \left({  - 2 \pi i 
\left\langle{ x', x }\right\rangle
}\right) = \exp 
\left({  -  \pi s 
W(x')
}\right). 
\eqno ({\mathop{\rm A.20 }\nolimits})$$ 

\smallskip
\nn
This can be interpreted as follows: the integrator $\Da \om$ is a continuous linear form on the
Banach space of Fourier-Stieltjes transforms $\Fa \mu$, given by
$$
\int_{ {\bf X}  }^{ } \Da \om (x) \Fa \mu (x) =
\int_{ {\bf X'}  }^{ } d \mu (x') \exp
\left({  -  \pi s 
W(x')
}\right). 
\eqno ({\mathop{\rm A.21 }\nolimits})$$  
This is another form of the Parseval relation.

\bigskip
We proceed to study the image of the integrator $\Da \om$ under a linear map.
We use the following notations:

\smallskip
-- $Q_{\bf X}$ is a quadratic form on ${\bf X}$, with inverse $W_{{\bf X}}$ 
on ${\bf X'}$;

\smallskip
-- 
$Q_{\bf Y}$ is a quadratic form on ${\bf Y}$, with inverse $W_{{\bf Y}}$ 
on ${\bf Y'}$;

\smallskip
-- $L$ is a continuous linear map from ${\bf X}$ into ${\bf Y}$.

\medskip
\nn
We assume that $L$ is surjective, hence $\ww L$ is injective, and that the quadratic forms
$W_{{\bf X}}$ and $W_{{\bf Y}}$ are related by
$$W_{{\bf Y}} = W_{{\bf X}} \circ \ww L \, .
\eqno ({\mathop{\rm A.22 }\nolimits})$$  
 
\bigskip
Consider now the integrators\footnote{$^{28}$}{
In the space ${\bf X} $, the integrator $\Da x$ depends on $Q_{{\bf X} } $
and $s$, and should be written more explicitly as 
$\Da_{s, Q_{{\bf X} } } x$. Similarly for $\Da y$ in space ${\bf Y}$.}
$$\eqalignno{
\Da \om_{{\bf X}} (x)
&= \exp \! 
\left({ - {\pi  \over s}  Q_{{\bf X} } (x) }\right)\cdot \Da x
&({\mathop{\rm A.23 }\nolimits})\cr
\Da \om_{{\bf Y}} (y)
&= \exp \! 
\left({ - {\pi  \over s}  Q_{{\bf Y} } (y) }\right) \cdot \Da y \, . 
&({\mathop{\rm A.24 }\nolimits})\cr
}
$$

\smallskip
\nn
The Fourier transforms are given respectively by 
$$\Fa \om_{{\bf X}} 
= \exp \! \left({ - \pi s  W_{{\bf X} }  }\right) \quad ,
\quad
\Fa \om_{{\bf Y}} 
= \exp \! \left({ - \pi s  W_{{\bf Y} }  }\right)
\eqno ({\mathop{\rm A.25 }\nolimits})$$  
and according to formula (A.22), we obtain
$$
\Fa \om_{{\bf Y}} 
= \Fa \om_{{\bf X}} 
\circ \ww L \, . 
\eqno ({\mathop{\rm A.26 }\nolimits})$$  

\smallskip
\nn
This is the same as formula (A.16), hence 
$\om_{{\bf Y}} $ is the image of $\om_{{\bf X}} $ under the linear mapping $L$.
Explicitly stated, we obtain the integration formula:
$$\int_{ {\bf X} }^{ } \Da x \cdot \exp \! 
\left({ - {\pi  \over s}  Q_{{\bf X} } (x) }\right)
\cdot g (Lx) =
\int_{ {\bf Y} }^{ } \Da y \cdot \exp \! 
\left({ - {\pi  \over s}  Q_{{\bf Y} } (y) }\right)
\cdot g (y) \, . 
\eqno ({\mathop{\rm A.27 }\nolimits})$$  

\smallskip
\nn
This holds if $g$ is a Fourier-Stieltjes transform
$\Fa \nu$ for some measure $\nu$
on
${\bf Y'}$.

\bigskip
The relationship between the quadratic forms 
$Q_{{\bf X}}$ on ${\bf X}$ and $Q_{{\bf Y}}$ on ${\bf Y}$
is expressed by the relation (A.22) between their inverses. Let
$D_{{\bf X}} : {\bf X} \ra {\bf X'}
$ be the continuous linear map such that 
$Q_{{\bf X}} (x) = 
\left\langle{ D_{{\bf X}} x, x }\right\rangle_{ {\bf X} }$
and 
$\left\langle{ D_{{\bf X}} x_1, x_2 }\right\rangle_{ {\bf X} }=
\left\langle{ D_{{\bf X}} x_2, x_1 }\right\rangle_{ {\bf X} }$,
 and define similarly 
 $ D_{{\bf Y}}$. Between the inverses $G_{\bf X}$ of $D_{{\bf X}}$ and 
$G_{{\bf Y}}$ of $D_{{\bf Y}}$ there holds the relation:
$$G_{{\bf Y}} = L \circ G_{{\bf X}} \circ \ww L \, . 
\eqno ({\mathop{\rm A.28 }\nolimits})$$  

\medskip
\nn
We mention a few particular cases:

\smallskip
-- if $L$ is invertible, with inverse $L^{- 1}$, then
$Q_{{\bf Y}} = Q_{{\bf X}} \circ L^{- 1}$;

\smallskip
-- if $Q_{{\bf X}}$ is positive-definite, then
$Q_{{\bf Y}} (y)$ is the infimum of $Q_{{\bf X}}$ over the set of elements $x$ in 
${\bf X}$ such that  $Lx = y$;

\smallskip
--
for any given $y$ in ${\bf Y}$ the equation
$$D_{{\bf X}} x =
\left({ \ww L \circ D_{{\bf Y}} }\right) (y)
\eqno ({\mathop{\rm A.29 }\nolimits})$$  
has a unique solution $x = x(y)$ in ${\bf X}$ and we obtain
$$Q_{{\bf Y}} (y) = Q_{{\bf X}} (x(y)) \, .
\eqno ({\mathop{\rm A.30 }\nolimits})$$

\bigskip
\nn
{\bf 3. Examples and applications.}

\smallskip
\nn
3.1. {\it The finite-dimensional case.} 
\smallskip
We record here the basic formulas. 
Assume that ${\bf X}$ is of finite dimension $d$; after choosing a linear frame, we represent a
vector $x$ by a column matrix $(x^\a)$ where
$\a \in 
\left\lbc{ 1, \cdots, d }\right\rbc$.
The elements of ${\bf X'}$ are represented by row matrices and the duality is given by
$\left\langle{ x', x }\right\rangle
= x'_\a x^\a$ (Einstein's summation convention).
The volume element is given by 
$d x = dx^1 \cdots dx^d$. 

\bigskip
The quadratic forms $Q$ and $W$ correspond to  symmetric matrices, namely:
$$Q(x) = h_{\a \b} x^\a x^\b \quad , \quad 
W(x') = h^{\a \b} x'_\a x'_\b
\eqno ({\mathop{\rm A.31 }\nolimits})$$  
with 
$h_{\a \b} h^{\b \g } = \d_\a^\g$.
The integrators are given by 
$$
\eqalignno{
\Da x
& = \left\vert{  \det h_{\a \b } }\right\vert^{1/2}
dx^1 \cdots dx^d 
& ({\mathop{\rm A.32 }\nolimits}) \cr
\Da \om (x)
& = 
\exp \! 
\left({ - \pi h_{\a \b} x^\a x^\b }\right) \cdot \Da x
& ({\mathop{\rm A.33 }\nolimits}) \cr}$$
when $s = 1$ and the matrix 
$\left({  h_{\a \b}  }\right)$ is positive-definite.
In the oscillating case, we have to multiply 
$\Da x$ by $e^{\pi i (q - p) /4}$ where the symmetric matrix
$\left({  h_{\a \b}  }\right)$ has $p$ positive and $q$ negative eigenvalues.

\bigskip
The quadratic form $W$ gives the covariance. More precisely, one obtains
$$\int_{ {\bf X} }^{ } \Da \om (x) 
\left\langle{ x', x }\right\rangle^2 =
{ s \over 2 \pi  } W(x')  
\eqno ({\mathop{\rm A.34 }\nolimits})$$
hence
$$\int_{ {\bf X} }^{ } \Da \om (x) 
x^\a x^\b =
{ s \over 2 \pi  } h^{\a \b }  \, .
\eqno ({\mathop{\rm A.35 }\nolimits})$$

\bigskip
\nn
3.2. {\it Image under a linear form.}
\smallskip
Let $x'$ in the dual ${\bf X'}$ of  ${\bf X}$. By specializing the results of paragraph A.2.2 to
the linear map $L : x \mps \left\langle{ x', x }\right\rangle$ from 
${\bf X}$ into $\RR$, we obtain the following  result. We identify $\RR$ with its dual, hence
$\ww L $ takes a number $\l$ to $\l x'$ in ${\bf X'}$. With $Q_{{\bf X}}$ equal to $Q$, hence
$W_{{\bf X}}$ to $W$ , we obtain
$$W_\RR (\l) = \l^2 W (x'), 
\eqno ({\mathop{\rm A.36 }\nolimits})$$
hence
$$Q_\RR (u) = u^2  / W (x') . 
\eqno ({\mathop{\rm A.37 }\nolimits})$$

\medskip
\nn
>From formula (A.27), we deduce the following integration formula:
$$\int_{{\bf X} }^{ }
\Da x \cdot e^{- \pi Q (x) / s } 
g 
\left({ \left\langle{ x', x }\right\rangle }\right)
 =
C \int_{ \RR }^{ } du \, e^{- \pi u^2  / s W (x')} 
g(u) 
\eqno ({\mathop{\rm A.38 }\nolimits})$$
where the normalization constant $C$ is 
$1 / (s W (x'))^{1/2}$ (principal branch of the square root). More explicitly:

\smallskip
-- for $s = 1$, hence $W(x') > 0$, then 
$ C = {\ds  1\over \ds\sqrt{ W (x')  }  } $;

\smallskip
-- for $s = i$, and $W(x') > 0$, then 
$ C = { \ds e^{- \pi i / 4} \over \ds  \sqrt{ W (x')  }  } $;

\smallskip
--  for $s = i$, and $W(x') < 0$, then 
$ C = {\ds  e^{\pi i / 4}  \over \ds \sqrt{ |W (x')|  }  } $.

\medskip
\nn
If we take in particular $g(u) = u^2 $, we get
$$\int_{{\bf X} }^{ }
\Da x \cdot e^{- \pi Q (x) / s } 
 \left\langle{ x', x }\right\rangle^2 =
 { s \over 2 \pi  } W(x') \, . 
\eqno ({\mathop{\rm A.39 }\nolimits})$$

\medskip
\nn
By polarization, we obtain the more general formula
$$
\int_{{\bf X} }^{ }
\Da x \cdot e^{- \pi Q (x) / s } 
 \left\langle{ x', x }\right\rangle
\left\langle{ y', x }\right\rangle
= 
{s  \over 2 \pi }
\left\langle{ x', Gy' }\right\rangle
\eqno ({\mathop{\rm A.40 }\nolimits})$$
for $x'$ and $y'$ in ${\bf X'}$, thus giving the covariance of our integrator. 

\bigskip
\nn
{\it Remark.}
In the application IV.1 to point-to-point transitions, we need to restrict the domain of integration (see also
paragraph A.3.8). This can be achieved by inserting a delta function
$\d \left({ \left\langle{ x'_0, x }\right\rangle }\right)$ and then using equation 
(A.38): we restrict the domain of integration from the space ${\bf X}$ to the hyperplane
${\bf X}_0$ with equation $\left\langle{ x'_0, x }\right\rangle  = 0$.

\bigskip
In general, any integral of the form
$$
\int_{{\bf X} }^{ }
\Da x \cdot e^{- \pi Q (x) / s } 
f \!
\left({  \left\langle{ x'_1, x }\right\rangle , \cdots , 
\left\langle{ x'_n, x }\right\rangle }\right)
$$
with fixed elements 
$x'_1, \cdots , x'_n $ in ${\bf X}'$ can be reduced to an $n$-dimensional integral.
For instance, combining formulas (A.38) and (A.40), the Green's function 
$G_{a,b} (t, u)$ on ${\bf Z}_{a, b}$ given by formulas (A.58) and (B.28), can be evaluated using an integral
$$
\int_{{\bf Z}_b }^{ }
\Da_b z \cdot \exp 
\left({ - { \pi \over s }
\int_{{\bf T} }^{ }dt \, z(t)^2
 }\right)
\d \! \left({ z (t_a) }\right) z(t) z(u)$$
and reducing it to a 3-dimensional integral.

\vfill\eject
\nn
3.3. {\it Space of paths of finite action.}
\smallskip
The time interval ${\bf T} = 
\left\lbk{ t_a,t_b }\right\rbk$ being given, we denote by $L^{2,1}$
(or more accurately $L^{2,1} ({\bf T})$) the space of real-valued functions 
$z(\cdot)$ on ${\bf T}$ with square-integrable derivative
$\dot z (\cdot)$ and we define the quadratic form $Q_0$ on $L^{2,1}$ by
$$Q_0 (z) =
\int_{ {\bf T} }^{ } dt \, 
\dot z (t)^2 \, .
\eqno ({\mathop{\rm A.41 }\nolimits})$$

\medskip
\nn
The definition of $L^{2,1}$ can be rephrased as follows:
the function $z$ belongs to $L^{2,1}$ if and only if there exists a function $\dot z$ in 
$L^2 ({\bf T})$ such that
$$z (t') - z(t) =
\int_{ t }^{ t' } du\, \dot z (u)  
\eqno ({\mathop{\rm A.42 }\nolimits})$$
whenever $t, t' $ are epochs in ${\bf T}$ such that $t < t'$.
The function $\dot z$ is then defined up to a null set; indeed by Lebesgue's derivation theorem,
one gets
$$\dot z(t) = 
\build { \lim }_{ \tau = 0 }^{ }
\left({ z(t + \tau) - z(t) }\right)/ \tau
\eqno ({\mathop{\rm A.43 }\nolimits})$$
for almost all $t$ in ${\bf T}$. By Cauchy-Schwarz inequality, one deduces from (A.42) the
inequality
$$
\left\vert{ z(t') - z(t) }\right\vert_{}^{2}
\leq Q_0 (z) \cdot 
\left\vert{ t' - t }\right\vert
\eqno ({\mathop{\rm A.44 }\nolimits})$$
for $t$, $t'$ in ${\bf T}$. Hence any function $z$ in $L^{2,1}$ satisfies a Lipschitz condition
of order $1/2$, and {\it a fortiori} it is a continuous function. 

\bigskip
The quantity $Q_0 (z)$ in equation (A.41) can be calculated using a {\it discretization of time}.
Consider a subdivision $\Ta$ of the time interval ${\bf T}$ by epochs
$$t_a \leq  t_0 < t_1 < \dots < t_{N-1} < t_N \leq  t_b\, .$$

\medskip
\nn
Set $\D t_i = t_i - t_{i-1}$ for 
$0 \leq i \leq N + 1$ with the convention $t_{-1} = t_a$, $t_{N+1}=t_b$  and denote by 
$\d (\Ta) $ the largest among the increments $\D t_i$, that is the {\it mesh} of the subdivision
$\Ta$. For any function $z : {\bf T} \ra \RR$ set
$$z_i = z(t_i) \quad , \quad \D z_i = z_i - z_{i-1}
\eqno ({\mathop{\rm A.45 }\nolimits})$$
for $0 \leq i \leq N + 1 $.
The {\it quadratic variation of} $z$ w.r.t. the subdivision $\Ta$ is defined as
$$Q_\Ta (z) = 
\sum_{ i = 1}^{ N } 
\left({ \D z_i }\right)_{}^{2}
/
\D t_i \, .
\eqno ({\mathop{\rm A.46 }\nolimits})$$

\medskip
\nn
Then the function $z$ belongs to $ L^{2,1}$ if and only if the set of quadratic variations 
$Q_{\Ta} (z)$ is bounded when $\Ta$ runs over all subdivisions of $ {\bf T} $. Then
$$Q_0 (z) = {\mathop{\rm l.}\nolimits} 
\build { {\mathop{\rm u}\nolimits}  }_{ \Ta}^{ }
 {\mathop{\rm .b. }\nolimits}\ Q_\Ta (z) \, .
\eqno ({\mathop{\rm A.47 }\nolimits})$$

\medskip
\nn
More precisely, for any sequence of subdivisions $\Ta (n)$ (for $n \in \NN$) whose mesh $\d (\Ta
(n))$ tends to $0$, one gets
$$Q_0(z) = 
\build { \lim }_{ n = \Inf  }^{ }
Q_{\Ta (n)}  (z) \, .
\eqno ({\mathop{\rm A.48 }\nolimits})$$

\medskip
\nn
It is therefore justified to write  $Q_0 (z) $ in the form 
$\ds \int_{{\bf T} }^{ }
{(dz)^2 \over  dt} .$

\bigskip
\nn
3.4. {\it Green's functions.}
\smallskip
Fix an element $t_0$ in ${\bf T}$ and denote by 
${\bf Z}_0$ the space of functions $z$ in $L^{2,1}$ such that 
$z (t_0) = 0$. This is a (real) Hilbert space, with scalar product
$$
\left\langle{ z_1 \bigm|  z_2 }\right\rangle
= \int_{{\bf T} }^{ } dt\, 
\dot z_1 (t) \dot z_2 (t)  ,
\eqno ({\mathop{\rm A.49 }\nolimits})$$
hence $Q_0 (z)$ is equal to 
$\left\langle{ z \bigm|  z }\right\rangle$.
>From (A.44), one deduces the inequality
$$\left\vert{ z(t) }\right\vert^2 \leq 
\left\vert{ t - t_0 }\right\vert
\cdot
\left\langle{ z \bigm|  z }\right\rangle,
\eqno ({\mathop{\rm A.50 }\nolimits})$$
and therefore there is an element $\d_t$ in the dual 
${\bf Z}'_0 $ of ${\bf Z}_0$ such that $z(t) =
\left\langle{ \d_t, z }\right\rangle$
(that is: $\d_t $ is a Dirac ``function'' centered at $t$).
According to the general theory, one introduces
a linear continuous and invertible map
$G_0 : {\bf Z}'_0 \ra {\bf Z}_0$ corresponding to the quadratic form $Q_0$ on ${\bf Z}_0$. 
The 
{\it Green's function} is defined as follows:
$$G_0(t, u) = 
\left\langle{ \d_t , G_0 \d_u }\right\rangle.
\eqno ({\mathop{\rm A.51 }\nolimits})$$
That is, the function $t \mps G_0(t, u)$ belongs to ${\bf Z}_0$ and is equal to 
$G_0 \d_u$. By definition, one gets
$$z(u) = 
\left\langle{ G_0 \d_u, z }\right\rangle
\eqno ({\mathop{\rm A.52 }\nolimits})$$
for every $z$ in ${\bf Z}$. More explicitly, these conditions on $G_0(t, u)$ can be expressed
as follows:
$$ \quad \qquad \qquad \qquad\left\{\matrix{
  G_0 \!\left({ t_0 , u }\right) = 0
   \hfill \cr
\noalign{\medskip}
\ds \int_{ {\bf T} }^{ } dt \, {\ds \part  \over \ds \part t} G_0(t, u) {\ds \part  \over \ds \part
t} z(t)  = z(u). 
  \hfill \cr}
\right. \qquad    \quad    \quad \qquad \qquad
\left.\matrix{
   ({\mathop{\rm A.53 }\nolimits})\cr
\noalign{\medskip}
 ({\mathop{\rm A.54}\nolimits})
  \cr}
\right. 
$$

\medskip
\nn
The unique solution to these equations is given by
$$
G_0(t, u) = 
\left\{\matrix{
\inf \left({ t - t_0, u - t_0 }\right)\hfill &
{\mathop{\rm for }\nolimits}
\ t \geq t_0 , \ u \geq t_0 , 
  \hfill\cr
\noalign{\medskip}
\inf \left({  t_0 - t ,  t_0 - u  }\right)\hfill &
{\mathop{\rm for }\nolimits}
\ t \leq t_0 , \ u \leq t_0 , 
   \hfill \cr
\noalign{\medskip}
0\hfill &
{\mathop{\rm otherwise }\nolimits}
\, .
  \hfill \cr}
\right. \eqno ({\mathop{\rm A.55 }\nolimits})$$

\medskip
We consider some special cases:

\smallskip
-- space ${\bf Z}_a$ of functions $z$ in $L^{2,1}$ with 
$z \!\left({ t_a }\right) = 0$;

\smallskip
-- space ${\bf Z}_b$ of functions $z$ in  $L^{2,1}$ with 
$z \!\left({ t_b }\right) = 0$;

\smallskip
-- space ${\bf Z}_{a, b} = {\bf Z}_a \cap {\bf Z}_b $ with boundary conditions 
$z \!\left({ t_a }\right) = z \!\left({ t_b }\right) = 0.$

\medskip
\nn
The corresponding Green's functions are as follows:
$$\eqalignno{
G_a (t, u) 
&=
\t(t - u) 
\left({ u - t_a }\right) + \t (u - t) \left({ t - t_a }\right)
& ({\mathop{\rm A.56 }\nolimits}) \cr
G_b (t, u) 
&=
\t(t - u) 
\left({ t_b - t }\right) + \t (u - t) \left({ t_b - u }\right)
& ({\mathop{\rm A.57 }\nolimits}) \cr
G_{a, b} (t, u) 
&=
 \t (t - u) \left({ t - t_b }\right) \left({ t_b - t_a }\right)^{-1 } 
\left({ t_a - u }\right)  
 \cr
& \ \  \  - \t (u - t) \left({ t - t_a }\right) \left({ t_a - t_b }\right)^{-1 } 
\left({ t_b - u }\right)  
 \, . 
& ({\mathop{\rm A.58 }\nolimits})\cr
}$$

\bigskip
\nn
{\it Remark. }
Equation (A.58) can be obtained by integrating
$$
\int_{{\bf Z}_b }^{ }
\Da \om (z) \, \d \! \left({ z(t_a) }\right)
\, z(t) \, z(u) \, .
$$

\bigskip
\nn
3.5. {\it Vector-valued functions.}
\smallskip
As in paragraph II.1.1, we consider vector functions
$z = \left({ z^1, \cdots, z^d }\right)
$ with components in $L^{2,1}$.
We introduce a real symmetric matrix
$\left({ h_{\a \b} }\right)$ with an inverse denoted by 
$\left({ h^{\a \b} }\right)$. The basic quadratic form is given by
$$Q_0 (z) = 
\int_{ {\bf T} }^{ } dt\,  h_{\a\b }\, 
\dot z^\a (t) \, \dot z^\b (t) \, . 
\eqno ({\mathop{\rm A.59 }\nolimits})$$

\medskip
\nn
The auxiliary condition 
$z (t_0) = 0$ defines the space ${\bf Z}_0$. Similarly, the spaces
${\bf Z}_a$, ${\bf Z}_b$ and ${\bf Z}_{a,b}$ are described by the respective boundary conditions:

\smallskip
-- $z(t_a) = 0 $ for ${\bf Z}_a$;

\smallskip
-- $z(t_b) = 0 $ for ${\bf Z}_b$;

\smallskip
-- $z(t_a) = z(t_b) = 0 $ for ${\bf Z}_{a,b}$.

\medskip
\nn
In each case, the linear form
$\d_t^\a$ taking $z$ into $z^\a (t)$ belongs to the dual of the corresponding of paths, and the
Green's function is characterized by 
$$G_0^{\a  \b} (t, u) =
\left\langle{ \d_t^\a, G_0 \d_u^\b  }\right\rangle \, . 
\eqno ({\mathop{\rm A.60 }\nolimits})$$

\medskip
\nn
Explicitly, we obtain
$$ G_0^{\a  \b} (t, u)  =
h^{\a \b } G_0(t, u)
\eqno ({\mathop{\rm A.61 }\nolimits})$$
where the Green's function $G_0$ refers to the boundary condition of the relevant space of paths
 (see formulas (A.55) to (A.58)).
The Green's function can also be obtained from formula (A.40) by specialization:
$$\int_{ {\bf Z} }^{ }
\Da z \cdot e^{- \pi Q_0 (z) / s}
\, z^\a (t) \, z^\b (u) =
{s \over 2 \pi } G_0^{\a \b }\,  (t, u) \, .
\eqno ({\mathop{\rm A.62 }\nolimits})$$
Here $\Da z $ stands for $\Da_{s, Q_0} z$.

\bigskip
\nn
{\it Remark}.
In general, the Green's functions $G$ satisfying a second-order differential equation $DG = \un$ are uniquely
determined by $d$ conditions at $t = t_a$ and $d$ conditions at $t = t_b$. These conditions are obvious from equation
(A.62) for the space ${\bf Z}_{a,b}$, namely $G^{\a \b } (t, u) = 0 $ if $t$ (or $u$) is one of the end points $t_a,
t_b$.

\medskip

In the case of the space ${\bf Z}_a$, the boundary conditions are
$$
\left.{ G }\right\vert_{t = t_a} = 0 \quad , \quad 
\left.{ \part G / \part t }\right\vert_{t = t_b} = 0 \, . $$
According to equation (A.62), these conditions can be expressed as 
$$
\eqalign{
\int_{ {\bf Z}_a }^{ } \Da \om (z) \cdot z^\a (t_a) z^\b (u) &= 0, \cr
\int_{ {\bf Z}_a }^{ } \Da \om (z) \cdot \dot z^\a (t_b) z^\b (u) &= 0\, .\cr}$$
The interpretation is as follows: any path $z$ in ${\bf Z}_a$ satisfies
$z^\a (t_a) = 0$; the time derivative 
$\dot z^\a (t_b) $ at $t_b$ is totally unspecified, nevertheless it vanishes in a statistical sense being
uncorrelated to the position and velocity at any other time.

\bigskip
On the space ${\bf Z}'_0$ dual to ${\bf Z}_0$, we have defined the quadratic form
$W_0$ inverse to $Q_0$. The Greens's function can be expressed as follows:
$$G_0^{\a \b} (t, u) =
{ {}_{1} \over {}^{2} }
\left\lbk{ W_0 \!
\left({ \d_t^\a + \d_u^\b }\right) -
W_0\!
\left({ \d_t^\a  }\right) -
W_0\! 
\left({  \d_u^\b }\right)
 }\right\rbk \, .
\eqno ({\mathop{\rm A.63 }\nolimits}) $$

\medskip
\nn
Conversely, given an element $z'$ of ${\bf Z'}$ represented by 
$$
\left\langle{ z', z }\right\rangle =
\int_{ {\bf T}J}^{ }dt \,  z'_\a (t) z^\a (t) , 
\eqno ({\mathop{\rm A.64 }\nolimits})$$
the quadratic form is given by
$$W_0(z') = \int_{ {\bf T}J}^{ } dt \int_{ {\bf T}J}^{ }
dt \, du \,
G_0^{\a \b} (t, u) 
z'_\a (t) \, z'_\b (u) 
\, .
\eqno ({\mathop{\rm A.65 }\nolimits})$$

\bigskip
\nn
3.6. {\it Scaling the paths.}

\smallskip
Here is a concise {\it dimensional analysis} of our quantities.

\nn
{\it Conventions} : $\La \ (\Ta)$ stands for the dimension of length (time) and $[X]$ for the
dimensional content of a quantity $X$.

\medskip
Since $Q_0 (z)$ appears in 
$e^{- \pi Q_0 (z) /s}$, it has to be a pure number, hence we can use formula (A.59) to deduce
$\left\lbk{ h_{\a \b } }\right\rbk
= \La^{- 2} \Ta$. From formulas (A.55) to (A.58) we infer the dimension of 
$G_0 (t, u)$ to be $\Ta$. 
>From (A.62) we infer that
$\left\lbk{G_0^{\a \b } }\right\rbk = \La^2$, hence 
$\left\lbk{h^{\a \b } }\right\rbk = \La^2 \Ta^{-1}$ from
(A.61). This is in accordance with the matrix relation 
$h_{\a \b } h^{\b \g }
= \d_\a^{\ \g}$. Notice also that for a particle of mass $m$ in a flat space, we have
$$
h_{\a \b } = m \d_{\a \b } /h 
\quad , \quad
h^{\a \b } = h \d_{\a \b } / m 
\eqno ({\mathop{\rm A.66 }\nolimits})$$
(where $h$ is the Planck constant and $\d_{\a \b}$ the Kronecker delta), and that the dimension of
$h/m $ is $\La^2 \Ta^{-1}$. We summarize our findings :

\bigskip
\centerline{TABLE 1}
$$\vbox{
\offinterlineskip
\halign
{\tv#&#  &\tv#&#  &#  &#  &#  &#  &#  &\tv#  \cr
\noalign{\hrule}
\tvi
\cc{Quantity}
&&\cc {$z^\a$} &\cc {$t$}
 &\cc {$h_{\a \b }$} &\cc  {$G_0 (t, u)$} &{$ G_0^{\a \b} (t, u)$} &
\cc {$h^{\a \b }$}  &&\cr
\tvi 
\cc{Dimension}
&&\cc {$\La$} &\cc {$\Ta$}
 &\cc {$\La^{-2} \Ta$} &\cc  {$\Ta$} &\cc {$ \La^2$} &
\cc {$\La^2 \Ta^{-1}$}  &&\cr
\noalign{\hrule} }}$$

\medskip
We can confirm these results by {\it scaling } our paths.
Let $\L > 0$ be a numerical constant, and denote by $\L  z$ the path with components
$\L z^\a (t) $ at time $t$. The basic integrator being written as 
$\Da \om (z) = e^{- \pi Q_0 (z) / s } \Da z$, we scale it into an integrator 
$\Da \om^\L (z) $ according to the formula
$$
\int_{ {\bf Z}_0 }^{ } \Da \om^\L (z) \, F(z) =
\int_{ {\bf Z}_0 }^{ } \Da \om (z) \, F (\L z) \, . 
\eqno ({\mathop{\rm A.67 }\nolimits})$$

\medskip
\nn
By suitably specializing $F(z)$, we obtain the new Green's function:
$$
\int_{ {\bf Z} }^{ } \Da \om^\L (z) \, 
z^\a (t) z^\b (u) =
\L^2 G_0^{\a \b} (t, u)\, ,
\eqno ({\mathop{\rm A.68 }\nolimits})$$
in accordance with 
$\left\lbk{ G_0^{\a \b } }\right\rbk = \La^2$. Similarly any $n$-point function scales as
$\La^n$. 
The map taking $z$ into $\L z$ is a linear map from ${\bf Z}_0 $ to ${\bf Z}_0 $.
By the general theory in paragraph A.2.2, $\Da \om^\L $ is another Gaussian integrator. From the
formulas (A.65) and (A.68), it corresponds to the quadratic form
$\L^2 W_0 (z')$ on ${\bf Z}'_0$, with an inverse quadratic form on ${\bf Z}_0$ given by
$$\L^{-2} Q_0 (z) = 
\int_{ {\bf T} }^{ } dt \, \L^{-2} \cdot
h_{\a \b} \, \dot z^\a (t) \dot z^\b (t) \, .
\eqno ({\mathop{\rm A.69 }\nolimits})$$

\medskip
\nn
This relation confirms $\left\lbk{  h_{\a \b} }\right\rbk = \La^{- 2} $.
A similar analysis applies to the scaling of time.

\bigskip
We can rewrite formula (A.67) in the following form:
$$
\int_{ {\bf Z}_0 }^{ }
\Da^\L z \cdot e^{- \pi \L^{-2} Q_0 (z) / s  }
H \!\left({ \L^{- 1} z }\right)
= \int_{ {\bf Z}_0 }^{ } \Da z \cdot e^{- \pi Q_0 (z) / s } 
H (z) \, ,
\eqno ({\mathop{\rm A.70 }\nolimits})$$
with a new integrator 
$\Da^\L z$  which is invariant under translations.
It is justified to summarize the previous formula by 
$\Da^\L z = \Da \left({ \L^{- 1} z}\right) $, hence
$\Da \om^\L (z) = \Da \om \!\left({ \L^{- 1} z}\right) $ according to formula (A.67).
In the standard heuristic derivations, one writes $\Da z$ in the form
$C \cdot \prod_{ t, \a  }^{ }
d z^\a (t) $. According to our normalization
$$
\int_{ {\bf Z}_0 }^{ } \Da z \cdot e^{- \pi Q_0 (z) / s } 
= 1, 
\eqno ({\mathop{\rm A.71 }\nolimits})$$
we should write
$$
\Da z =
{ \ds \prod_{t, \a  }^{ } d z^\a (t)  \over
\ds \int_{ }^{ } \ds \prod_{t, \a  }^{ } d z^\a (t) \cdot
\exp \left({ - \pi Q_0 (z) / s }\right)
 } \, , 
\eqno ({\mathop{\rm A.72 }\nolimits})$$
and similarly
$$
\Da^\L z =
{ \ds \prod_{t, \a  }^{ } d z^\a (t)  \over 
\ds \int_{ }^{ } \ds \prod_{t, \a  }^{ } d z^\a (t) \cdot
\exp \left({ - \pi Q_0 \left({ \L^{- 1} z }\right)  / s }\right)
} \, . 
\eqno ({\mathop{\rm A.73 }\nolimits})$$
Replacing $z$ by $\L^{- 1} z$ amounts to replacing 
$d z^\a (t) $ by $\L^{-1} dz^\a (t)$, hence the volume element 
$\prod_{t, \a  }^{ } d z^\a (t) $ is multiplied by $\L^{- N}$ where 
$N$ is the (infinite) number of degrees of freedom $t, \a$.
If we calculate $\Da \left({ \L^{- 1} z }\right)$ in accordance with (A.72), both numerator and 
denominator acquire a factor $\L^{- N}$ which drops out, and the correct formula
$\Da^\L z = \Da \left({ \L^{- 1} z }\right)$ is obtained from the heuristic formulas
(A.72) and (A.73).

\bigskip
The heuristic constant $\L^{- N}$ is equal to $\Inf$, 1 or 0 according to the three cases
$0 < \L < 1$, $\L = 1$, $\L > 1$. This is reflected in the rigorous theory by the following fact :
{\it for $\L \ne 1$  no functional $F (z)$, except the constant $0$, is such that both $F(z)$ and 
$F (\L z)$ can be simultaneously integrated } w.r.t. $\Da z$ In a finite-dimensional space $\RR^N$, the
volume element $dx^1 \cdots dx^N$ is the only one, up to a multiplicative constant, which is
invariant under translations. In our infinite-dimensional setup there exist many
translation-invariant integrators, but {\it they act in different functional sectors}.

\bigskip
\nn
3.7. {\it White noise representations.}

\smallskip
We consider now the Hilbert space $L^2 ({\bf T}) $ (or $L^2$) consisting of the square-integrable functions
$\xi : {\bf T} \ra \RR$, with the standard quadratic form
$$H(\xi) = \int_{ {\bf T} }^{ } dt \, \xi (t)^2 \, . 
\eqno ({\mathop{\rm A.74 }\nolimits})$$
As usual we identify 
$L^2 ({\bf T})$ with its own dual by means of the following scalar product
$$
\left\langle{ \xi', \xi  }\right\rangle
=
\int_{ {\bf T} }^{ } dt \, \xi' (t)  \xi (t) \, . 
\eqno ({\mathop{\rm A.75 }\nolimits})$$
Let us denote for a while the space
$L^2 ({\bf T})$ by ${\bf H}$. Under our general conventions, both
$D : {\bf H} \ra {\bf H' } $ and 
$G : {\bf H'} \ra {\bf H } $ are identity operators, hence the quadratic form 
${\bf H}$ is its own inverse.
The corresponding covariance is given by the kernel of the identity operator, namely
$\d (t - u)$. Hence denoting by $\Da \xi$ the basic integrator, we get from equation (A.40):
$$
\int_{ L^2 ({\bf T})  }^{ }
\Da \xi \cdot \exp 
\left({ 
- { \pi \over s}
\int_{ {\bf T} }^{ } dt \, \xi (t)^2 
 }\right)
\xi (t) \xi (u) =
{s \over 2 \pi}
\d (t - u) \, . 
\eqno ({\mathop{\rm A.76 }\nolimits})$$

\bigskip
Let 
$L_0^{2,1}({\bf T})$ (or $L_0^{2,1}$) be the subspace of functions $z$ in 
$L^{2,1} ({\bf T})$ such that $z(t_0) = 0$.
The spaces $L^2 ({\bf T})$ and $L_0^{2,1}({\bf T})$ are isomorphic under a correspondance 
$\xi \lra z$ where
$\xi = \dot z $ and conversely 
$z(t) = \int_{ t_0}^{ t} dt' \, \xi (t')$. This can be expressed by the formulas
$$z(t) =
\int_{ {\bf T} }^{ } dt'\, \T (t', t) \xi (t') \, , 
\eqno ({\mathop{\rm A.77 }\nolimits})$$
$$
\T (t', t) =
\left\{\matrix{
 1\hfill &  {\mathop{\rm for }\nolimits} \hfill & t_0 \leq t' \leq t,  \hfill\cr
\noalign{\medskip}
   - 1\hfill &  {\mathop{\rm for }\nolimits} \hfill & t \leq t' \leq t_0,  \hfill\cr
\noalign{\medskip}
 0\hfill &  {\mathop{\rm otherwise .}\nolimits} \hfill &   \hfill\cr}
\right. \eqno ({\mathop{\rm A.78 }\nolimits})$$

\medskip
\nn
The following transformation rule holds for functional integrals:
$$\int_{ L_0^{2,1} }^{ }
\Da z \cdot e^{- \pi Q_0 (z) / s} F(z) =
\int_{ L^{2} }^{ }
\Da \xi \cdot e^{- \pi H (\xi) / s}
\Phi (\xi) \, , 
\eqno ({\mathop{\rm A.79 }\nolimits})$$
where $\Phi (\xi)$ is equal to $F(z)$ if $\xi$ corresponds to $z$ by (A.77).
The covariance in 
$L_0^{2,1}$ is obtained by specializing $F(z)$ to $z(t) z(u) $ in (A.79).
Using (A.76), we obtain
$$G_0 (t, u) =
\int_{ {\bf T} }^{ } dt' \, \T (t', t) \T (t', u)
\eqno ({\mathop{\rm A.80 }\nolimits})$$
in agreement with formula (A.55).

\bigskip
To make contact with the heuristic definitions, we introduce the 
``coordinates'' 
$X_t = \xi (t) \sqrt{ dt } $ for the function 
$\xi$, hence
$$H(\xi) = 
\sum_{ t }^{ } X_t^2 \quad , 
\quad \Da \xi = \prod_{ t }^{ }
\left({ dX_t / \sqrt{ s}  }\right) \, .
\eqno ({\mathop{\rm A.81 }\nolimits})$$
We urge the reader to substantiate these claims by resorting to a subdivision 
$\Ta$ of the time interval ${\bf T}$, like in paragraph A.3.3.

\bigskip
Let ${\bf X}$ be any closed vector subspace of
$L^2 ({\bf T})$, and let $\Pi$ be the orthogonal projector from
$L^2 ({\bf T} )$ onto ${\bf X}$, represented as an integral operator with kernel $\Pi (t, u)$.
The quadratic form $H_{\bf X}$ on ${\bf X}$ is obtained by restriction of $H$, namely:
$$H_{\bf X} (\xi) = 
\int_{ {\bf T} }^{ } dt \, \xi (t)^2
\eqno ({\mathop{\rm A.82 }\nolimits})$$
for $\xi$ in ${\bf X}$.
Denote by 
$\Da_{{\bf X}} \xi$ the corresponding integrator on ${\bf X}$.
Then the covariance is expressed as follows:
$$\int_{ {\bf X} }^{ }
\Da_{{\bf X} } \xi \cdot \exp 
\left({ -{ \pi \over s}
\int_{ {\bf T} }^{ }
dt \, \xi (t)^2
 }\right)
\cdot \xi (t) \xi (u) =
{ s \over 2 \pi } \Pi (t, u) \, .
\eqno ({\mathop{\rm A.83 }\nolimits})$$

\medskip
\nn
For instance, if we denote by $L_0^2 ({\bf T})$ the subspace of $L^2 ({\bf T})$ defined by the
condition $\int_{ {\bf T}}^{ } dt \, \xi (t) = 0$, the orthogonal projector is given by
$$(\Pi \xi) (u) = \xi (u) - 
{1 \over t_b - t_a} \int_{ {\bf T}}^{ } dt \, \xi (t),
\eqno ({\mathop{\rm A.84 }\nolimits})$$
hence
$$\Pi (t, u) = \d(t - u) - 
{1 \over t_b - t_a} \, . 
\eqno ({\mathop{\rm A.85 }\nolimits})$$
The map $\xi \mps z$ where $z(t) =
\int_{t_a }^{t } dt' \, \xi (t')$ takes isomorphically the space
$L_0^2 ({\bf T})$ onto the space denoted by 
${\bf Z}_{a,b}$ in paragraph A.3.4. It follows that the Green's function 
$G_{a,b}$ corresponding to the space ${\bf Z}_{a,b}$ is given by
$$G_{a,b}(t, u) = \int_{t_a }^{t } dv \int_{t_a }^{u} dv' \, \Pi (v, v') \, . 
\eqno ({\mathop{\rm A.86 }\nolimits})$$
Using (A.85), we obtain easily (A.58).

\bigskip
Explicit formulas for transformations, mapping a given quadratic form on
$L^2 ({\bf T}) $ into any given quadratic form on 
$L^{2,1} ({\bf T})$, can be found in 
[5, p. 274]. The reader is urged to extend the results of this paragraph to vector-valued
functions.

\bigskip
\nn
3.8. {\it Fixing the endpoints.}
\smallskip
We denote by 
$\Da_a z$, $\Da_b z$ and $\Da_{a, b} z$ respectively the integrators on the spaces
${\bf Z}_a$, ${\bf Z}_b$ and ${\bf Z}_{a, b}$ (see paragraph A.3.5.). These integrators are
related by the following formula:
$$
\int_{ {\bf Z}_a }^{ } \Da_a z \cdot e^{- \pi Q_0 (z) / s }
F(z) \d \!\left({ z \!\left({ t_b }\right) }\right)
= C^{- 1 }
\int_{ {\bf Z}_{a,b} }^{ } \Da_{a,b} \, z \cdot e^{- \pi Q_0 (z) / s }
F(z) \, ,
\eqno ({\mathop{\rm A.87 }\nolimits})$$
where the constant $C$ is given by
$$\eqalign{
C &= 
\left({ \det s G_a^{\a, \b } \left({ t_b, t_b  }\right) }\right)_{}^{1/2}\cr
&= \left({ \det sh^{\a \b } }\right)_{}^{1/2}
\left({ t_b - t_a }\right)_{}^{d/2}
\, . \cr} 
\eqno ({\mathop{\rm A.88 }\nolimits})$$
A similar formula (where the roles of $t_a$ and $t_b$ are exchanged) can be found in 
[5, pp. 284, 357].

\bigskip
We give a new derivation of this formula to illustrate our methods.
We simplify the notations by putting $\Da \om_a (z)$ equal to 
$\Da_a z \cdot e^{- \pi Q_0 (z) / s }$, and defining 
$\Da \om_{a, b} (z)$ similarly. According to our general strategy, we need only to prove formula
(A.87) in the special case
$$F(z) = e^{- 2 \pi i \left\langle{  z', z}\right\rangle }
\eqno ({\mathop{\rm A.89 }\nolimits})$$
where $z'$ is an element of ${\bf Z}'_a$. That is, we have to establish the following equality:
$$\int_{ {\bf Z}_a }^{ } \Da \om_a (z) \cdot e^{- 2\pi i
\left\langle{  z', z}\right\rangle  }
\d \! \left({ z \! \left({ t_b }\right) }\right)
= C^{- 1} 
e^{- \pi s W_{a, b } (z') } \, ,  
\eqno ({\mathop{\rm A.90 }\nolimits})$$
where $ W_{a, b }$ is the quadratic form on 
${\bf Z}'_{a,b}$ inverse of the restriction of $Q_0$ to the subspace 
${\bf Z}_{a,b}$ of ${\bf Z}_{a}$.

\bigskip
To prove (A.90), we use the well-known formula 
$$\d(x) =
\int_{\RR  }^{ }
e^{- 2 \pi i u x} \, du \, .
\eqno ({\mathop{\rm A.91 }\nolimits})$$
We shall consider in detail the scalar case $d = 1$ and leave the general case to the reader.
The left-hand side $L$ of (A.90) can be rewritten as
$$L = \int_{ \RR }^{ } du 
\int_{ {\bf Z}_a }^{ } \Da \om_a (z) \, 
\exp 
\left({ - 2 \pi i 
\left\langle{ z' + u \d_{t_b}, z }\right\rangle
 }\right), 
\eqno ({\mathop{\rm A.92 }\nolimits})$$
hence
$$L = \int_{ \RR }^{ } du  \, \exp 
\left({ - \pi s W_a 
\left({ z' + u \d_{t_b} }\right)
 }\right)
\eqno ({\mathop{\rm A.93 }\nolimits})$$
by (A.18). We develop now the exponent
$$W_a  \left({ z' + u \d_{t_b }  }\right)
= u^2 G_a \! \left({ t_b, t_b }\right)
+ 2 u 
\int_{ {\bf T} }^{ } dt \, z'(t) G_a 
\! \left({ t_b, t  }\right)
+ W_a (z') 
\eqno ({\mathop{\rm A.94 }\nolimits})$$
and use the integration formula
$$
\int_{ \RR }^{ } du \, 
e^{- \pi s  
\left({ au^2 + 2 bu + c  }\right) }=
{1 \over \sqrt{ as }  }
\, \exp 
\left({ - \pi s 
\left({ c - b^2/a }\right)
}\right)
\eqno ({\mathop{\rm A.95 }\nolimits})$$
to obtain 
$$L = 
{1  \over \sqrt{s G_a \! \left({ t_b, t_b  }\right) } } \cdot 
\eqno ({\mathop{\rm A.96 }\nolimits})
$$
$$\exp 
\left({ - \pi s \! \left({  W_a (z')  
- 
G_a \! \left({ t_b, t_b  }\right)^{- 1 } 
\int_{ {\bf T} }^{ } dt 
\int_{ {\bf T} }^{ } du \,
z'(t) z'(u) 
G_a \!\left({ t_b, t  }\right)
G_a \! \left({ t_b, u  }\right)  }\right) }\right) \, . 
$$

\medskip
\nn
This can be transformed in the desired right-hand side of (A.90) with 
$C = \sqrt{s G_a \!\left({ t_b, t_b  }\right) }$ provided we establish the following identity:
$$G_{a, b} (t, u) = 
G_a(t, u) - G_a \!\left({ t, t_b }\right) 
G_a \!\left({ t_b, t_b }\right)^{- 1 }
G_a \!\left({ t_b, u }\right) \, .
\eqno ({\mathop{\rm A.97 }\nolimits})$$
We leave it as an exercice (see formulas (A.56) and (A.58)).

\bigskip
We can prove a more general formula by a similar reasoning. We fix points 
$z_a$ and $z_b$ and consider $L^{2,1}$ paths $z(\cdot)$ with 
$z(t_a) = z_a$. We obtain the affine space 
$z_a + {\bf Z}_a$, and transport to it the integrator 
$\Da_a z$. Then, with a suitable constant $C_{a, b}$ one gets
$$
\int_{ z_a + {\bf Z}_a }^{ }
\Da_a z \cdot 
e^{- \pi Q_0 (z) / s } F(z) 
\d \left({   z(t_b) - z_b}\right)$$
$$= C_{a, b}^{-1}
\int_{ {\bf Z}_{a, b} }^{ } \Da_{a,b}  \z \cdot 
e^{- \pi Q_0 (\z) / s } F
\left({ z_{{\mathop{\rm cl }\nolimits}}  + \z }\right)
\, .
\eqno ({\mathop{\rm A.98 }\nolimits})$$
Here $ z_{{\mathop{\rm cl }\nolimits}}$ is the affine-linear map
$ z_{{\mathop{\rm cl }\nolimits}} (t) = \l t + \mu $ with endpoints 
$z(t_a) = z_a $, $z(t_b) = z_b$, and the path $\z$ in 
${\bf Z}_{a,b}$, that is 
$\z (t_a) = \z (t_b) = 0$, is a {\it quantum fluctuation.}

\bigskip
\nn
{\bf 4. Infinite-dimensional determinants.}
\smallskip
For this section, we refer the reader to [24], [25, chapter 5] and [26].

\bigskip
\nn
4.1. {\it The case of operators.}
\smallskip
We return to our Banach space ${\bf X}$.
We assume that there exists on ${\bf X}$ an invertible\footnote{$^{29}$}{
A positive-definite continuous quadratic form is not necessarily invertible.
For instance, let ${\bf X}$ be the space $\ell^2$ of sequences
$\left({ x_1, x_2, \cdots }\right)$ of real numbers with
$\sum_{n = 1 }^{ \Inf } \left({ x_n }\right)^2 < \Inf
$ and define the norm by
$\left\Vert{ x }\right\Vert_{}^{2}=
\sum_{ n=1}^{ \Inf }
\left({ x_n }\right)^2
$.
We can identify ${\bf X}$ with its dual ${\bf X'}$, the scalar product being given by 
$\sum_{ n=1}^{ \Inf } x'_n x_n$.
The quadratic form 
$Q(x) = \sum_{ n=1}^{ \Inf } \left({ x_n / n }\right)^2$ corresponds to the map
$D : {\bf X } \ra {\bf X' }$ taking 
$\left({  x_1, x_2, \cdots}\right)$ into
$\left({ x_1/1, x_2/2, \cdots }\right)$. The inverse of $D$ does not exist as  a map from
$\ell^2$ into $\ell^2$ since the sequence $1, 2, 3, ...$ is unbounded.
 }
positive-definite quadratic form $H$ (for instance
$H(z) = \sum_{ \a }^{ }
\int_{ {\bf T } }^{ }
dt\, \dot z^{\a } (t)^2$ on the space ${\bf Z}_0$). 
We use the norm defined on ${\bf X}$ by
$\left\Vert{ x }\right\Vert
= H(x)^{1/2}$, and derive from it the dual norm
$\left\Vert{ x' }\right\Vert$ on ${\bf X'}$ as usual.

\bigskip
According to Grothendieck, an operator $T$ in ${\bf X}$ is called {\it nuclear} if it admits a
representation in the form
$$Tx = \sum_{ n \geq 0 }^{ }
\left\langle{ x'_n, x }\right\rangle
\cdot x_n
\eqno ({\mathop{\rm A.99 }\nolimits})$$
with elements 
$x_n$ in ${\bf X}$ and $x'_n$ in ${\bf X'}$ such that
$\sum_{ n \geq 0 }^{ }
\left\Vert{ x_n }\right\Vert
\cdot
\left\Vert{ x'_n }\right\Vert$ be finite.
The g.l.b. of all such sums 
$\sum_{ n  }^{ }
\left\Vert{ x_n }\right\Vert
\cdot
\left\Vert{ x'_n }\right\Vert$ is called the {\it nuclear norm} of $T$,
to be denoted by
$\left\Vert{ T }\right\Vert_1$. The nuclear operators in ${\bf X}$ form a Banach space, denoted by
$\La^1 ({\bf X})$, with norm 
$\left\Vert{ \cdot }\right\Vert_1$.
On $\La^1 ({\bf X})$, there exists a continuous linear form, the {\it trace}, such that
$$\Tr (T) = \sum_{ n \geq 0 }^{ }
\left\langle{ x'_n, x_n }\right\rangle 
\eqno ({\mathop{\rm A.100 }\nolimits})$$
for an operator $T$ given by (A.99).

\bigskip
We introduce a power series in $\l$, namely:
$$
\sum_{p \geq 0  }^{ } 
\s_p (T) \l^p : = 
\exp 
\left({ \l \Tr (T) -
{\l^2 \over 2}  \Tr (T^2) + {\l^3 \over 3} \Tr (T^3) - \cdots
 }\right) \, . 
\eqno ({\mathop{\rm A.101 }\nolimits})$$
By using Hadamard's inequality about determinants, we obtain the {\it basic estimate}:
$$\left\vert{ \s_p (T) }\right\vert
\leq p^{p/2}
\left\Vert{ T }\right\Vert_{1}^{p}
/ p!\, .
\eqno ({\mathop{\rm A.102 }\nolimits})$$

\medskip
\nn
It follows that the power series $\sum_{ p \geq 0 }^{ }  \s_p (T) \l^p $ has an infinite radius of
convergence. We can therefore define the determinant as follows:
$$\Det(1 + T):= \sum_{ p \geq 0 }^{ } \s_p (T)
\eqno ({\mathop{\rm A.103 }\nolimits})$$
for any nuclear operator $T$. By homogeneity, we obtain more generally
$$\Det (1 + \l T) = 
\sum_{ p \geq 0 }^{ } \s_p (T) \l^p \, .
\eqno ({\mathop{\rm A.104 }\nolimits})$$
The fundamental property of determinants is, as expected, the
{\it multiplicative rule}:
$$\Det \left({ U_1 \circ U_2 }\right)
= \Det \left({ U_1 }\right)
\Det \left({ U_2 }\right) \, , 
\eqno ({\mathop{\rm A.105 }\nolimits})$$
where $U_i$ is of the form $1 + T_i$, with $T_i$ nuclear
(for $i = 1, 2)$. From this, and the relation $\s_1 (T) = \Tr (T)$, we get a
{\it variation formula} (for $U - 1$ and $\d U$ nuclear):
$$
{ \Det (U + \d U) \over \Det (U) } 
= 1 + \Tr 
\left({ U^{- 1} \cdot \d U}\right)
+ O
\left({ \left\Vert{ \d U }\right\Vert_{1}^{2} }\right) \, .
\eqno ({\mathop{\rm A.106 }\nolimits})$$
Otherwise stated, if $U(\nu)$ is an operator of the form 
$1 + T(\nu)$, where $T(\nu)$ is nuclear, depending smoothly on the parameter $\nu$,
we get the {\it derivation formula}:
$${ d \over d \nu } 
\ell n \Det (U (\nu)) =
\Tr 
\left({ U (\nu)^{-1} {d \over d \nu} U (\nu) }\right) \, .
\eqno ({\mathop{\rm A.107 }\nolimits})$$

\bigskip
\nn
{\it Remark.}
For any other norm 
$\left\Vert{ \cdot }\right\Vert^1$ defining the topology of ${\bf X}$, we have an estimate
$$C^{- 1}
\left\Vert{ x }\right\Vert
\leq \left\Vert{ x }\right\Vert^1 \leq C 
\left\Vert{ x }\right\Vert \, , 
\eqno ({\mathop{\rm A.108 }\nolimits})$$
with a  finite numerical constant $C > 0$.
It follows easily that the previous definitions are independent of the choice of the particular
norm
$\left\Vert{ x }\right\Vert = H (x)^{1/2}$ in ${\bf X}$.

\bigskip
\nn
4.2. {\it Explicit formulas.}
\smallskip
Introduce a basis
$\left({  e_n}\right)_{n \geq 1}^{}$ of ${\bf X}$ orthonormal for the quadratic form 
$H$, hence $H ( \sum_{ n }^{ } t_n e_n) = \sum_{ n }^{ } t_n^2$.
An operator $T$ in ${\bf X}$ has a matrix $(t_{mn})$ such that
$$T e_n = \sum_{ m }^{ } e_m \cdot t_{mn} \, . 
\eqno ({\mathop{\rm A.109 }\nolimits})$$

\bigskip
Assume that $T$ is nuclear. 
Then the series $\sum_{ n }^{ }  t_{nn}$ of diagonal terms in the matrix converges absolutely and 
the trace
$\Tr (T)$ is equal to $\sum_{ n }^{ }  t_{nn}$, as it should be.
Furthermore $\s_p (T)$ is the sum of the series made of the principal minors of order $p$:
$$\s_p(T) = 
\sum_{ i_1 < \cdots < i_p}^{ }
\det \left({ t_{i_\a, i_\b} }\right)_{
{1 \leq \a \leq p \atop 1 \leq \b \leq p} }
\, . 
\eqno ({\mathop{\rm A.110 }\nolimits})$$

\medskip
\nn
For the operator $U = 1 + T$, with the matrix with elements 
$u_{mn} = \d_{mn} + t_{mn}$, we obtain the determinant as a 
{\it limit of finite-size determinants}:
$$\Det(U) = 
\build { \lim }_{ N = \Inf }^{ }
\det 
\left({ u_{mn} }\right)_{ {1 \leq m \leq N \atop 1 \leq n \leq N}  }^{}
\, . 
\eqno ({\mathop{\rm A.111 }\nolimits})$$

\bigskip
As a special case, suppose that the basic vectors $e_n$ are 
{\it eigenvectors} for $T$, namely
$$T e_n = \l_n e_n \, . 
\eqno ({\mathop{\rm A.112 }\nolimits})$$
Then we get
$$\Tr (T) = \sum_{ n }^{ } \l_n \quad , 
\quad \Det (1 + T) = 
\prod_{ n}^{ }
\left({ 1 + \l_n }\right)$$
where both the series and the infinite product converge absolutely.

\bigskip
The nuclear norm 
$\left\Vert{ T }\right\Vert_1$ can also be computed as follows: there exists an orthonormal basis
$(e_n)$ such that the vectors $T e_n$ are mutually orthogonal (for the quadratic form $H$) and
then
$\left\Vert{ T }\right\Vert_1 =
\sum_{ n }^{ } 
\left\Vert{ T e_n }\right\Vert$.

\bigskip
\nn
{\it Remark.}
Let $T$ be a continuous linear operator in ${\bf X}$. Assume that the series of diagonal terms 
$\sum_{ n }^{ } t_{nn}$ converges absolutely for {\it every} orthonormal basis.
Then $T$ is nuclear. When $T$ is symmetric and positive, it is enough to assume that this
statement holds for {\it one} given orthonormal basis, and then it holds for all. There are
counterexamples when $T$ is not symmetric and positive. 

\bigskip
\nn
4.3. {\it The case of quadratic forms.}
\smallskip
Contrary to a widespread misbelief, {\it there is no such thing like the determinant of a
quadratic form.}
Consider for instance a quadratic form $Q$ on some finite-dimensional space with coordinates
$x^1, \cdots , x^d$, namely
$$Q(x) = h_{\a \b } x^\a x^\b \, . 
\eqno ({\mathop{\rm A.113 }\nolimits})$$
If we introduce a new system of coordinates 
$\oo x^1, \cdots, \oo x^d$ such that 
$x^\a = u_\l^\a \oo x^\l $, then we obtain
$$Q(x) = \oo h_{\l \mu} \oo x^\l \oo x^\mu
\eqno ({\mathop{\rm A.114 }\nolimits})$$
with a new matrix
$$\oo h_{\l \mu } = u_\l^\a u_\mu^\b h_{\a \b }
\, . 
\eqno ({\mathop{\rm A.115 }\nolimits})$$
The determinants 
$D = \det \left({ h_{\a \b } }\right)$ and 
$\oo D = \det \left({ \oo h_{\l \mu } }\right)$ are connected by the 
{\it scaling relation}:
$$\oo D = D \cdot 
\left({ \det (u_\l^\a) }\right)_{}^{2} \, . \eqno ({\mathop{\rm A.116 }\nolimits})$$

\medskip
\nn
Hence, what makes sense is the {\it ratio of determinants}
$$\det 
\left({ h_{\a \b }^{(1)} }\right)
/
\det 
\left({ h_{\a \b }^{(0)} }\right)
\eqno ({\mathop{\rm A.117 }\nolimits})$$
associated to two quadratic forms
$$Q_0 (x) =
h_{\a \b }^{(0)} x^\a x^\b \quad , \quad 
Q_1 (x) = h_{\a \b }^{(1)} x^\a x^\b
\eqno ({\mathop{\rm A.118 }\nolimits})$$
on the same space. We shall denote it by 
$\det (Q_1 / Q_0)$.

\bigskip
To obtain an intrinsic definition, let us consider two continuous quadratic forms $Q_0$ and $Q_1$
on  a Banach space ${\bf X}$ and assume that $Q_0$ is invertible.
We associate to $Q_0$ and $Q_1$ two continuous linear maps
$D_i : {\bf X} \ra {\bf X'}$ such that
$$
\left\langle{ D_i x_1, x_2 }\right\rangle =
\left\langle{ D_i x_2, x_1 }\right\rangle
\quad , \quad
Q_i (x) =
\left\langle{ D_i x, x  }\right\rangle \, ,
\eqno ({\mathop{\rm A.119 }\nolimits})$$
for $x$, $x_1$ and $x_2$ in ${\bf X}$.
By assumption, $D_0$ is invertible, hence there exists a unique continuous 
operator $U$ in ${\bf X}$ such that 
$D_1 = D_0 \circ U$.
We denote $U$ by\footnote{$^{30}$}{A better notation should perhaps be $Q_0 \bl Q_1$. }
$Q_1/Q_0$. In case the determinant of $U$ is defined,
namely when $U - 1$ is nuclear, we denote it  by 
$\Det \left({ Q_1 / Q_0 }\right)$. Of course, when ${\bf X}$ is finite-dimensional, this
definition agrees with the previous one and
$$\Det \left({ Q_1 / Q_0 }\right) =
\det \left({ Q_1 / Q_0 }\right)
\eqno ({\mathop{\rm A.120 }\nolimits})$$
in this case.

\bigskip
A procedure to calculate this determinant is as follows. Let $V$ be a finite-dimensional subspace
of ${\bf X}$. By restricting $Q_0$ and $Q_1$ to $V$ we obtain two quadratic forms 
$Q_{0, V}$ and $Q_{1, V}$ on $V$. 
Assume now that $V$ runs through an increasing sequence of subspaces, whose union is dense in ${\bf X}$,
and  that $Q_{0, V}$ be invertible for every $V$.
Then 
$$\Det \left({ Q_1 / Q_0 }\right) =
\build {\lim }_{ V}^{ }
\det 
\left({ Q_{1,V} / Q_{0,V} }\right) \, . 
\eqno ({\mathop{\rm A.121 }\nolimits})$$

\bigskip
\nn
4.4. {\it A pencil of quadratic forms.}

\smallskip
Here are our assumptions:

\smallskip
-- {\it there exists an invertible positive-definite quadratic form on ${\bf X}$};

\smallskip
-- {\it $Q_0$ and $Q_1$ are continuous quadratic forms on ${\bf X}$};
\smallskip

-- {\it $Q_0$ is invertible};

\smallskip
-- {\it if we set $Q = Q_1 - Q_0$, the operator $T = Q / Q_0$ on ${\bf X}$ is nuclear.}

\smallskip
The last condition is an intrinsic property of the quadratic form $Q$, namely the existence of a
representation like
$$Q(x) =
\sum_{ n \geq 1 }^{ } 
\a_n 
\left\langle{ x'_n, x }\right\rangle^2, 
\eqno ({\mathop{\rm A.122 }\nolimits})$$
where the series
$\sum_{ n }^{ } \a_n$ converges absolutely, and the sequence of numbers
$\left({ 
\left\langle{ x'_n, x }\right\rangle
 }\right)_{n \geq 1}^{} $ is bounded for every vector $x$ in ${\bf X}$.
We express this property by saying that the {\it quadratic form $Q$ is nuclear.}

\bigskip
We say  that $\l$ is an {\it eigenvalue } of a quadratic form $Q'$ 
w.r.t. $Q_0$ if it is an eigenvalue of the operator $Q' / Q_0$.
This is tantamount to saying that the quadratic form $Q' - \l Q_0 $ 
is not invertible.

\bigskip
We interpolate between the quadratic forms $Q_0$ and 
$Q_1 = Q_0 + Q$ by putting $Q_\nu = Q_0 + \nu Q$ where $\nu$ is a real or complex parameter.
The determinant 
$$\D (\nu) = \Det (Q_\nu / Q_0)
\eqno ({\mathop{\rm A.123 }\nolimits})$$
is defined, being equal to $\Det (1 + \nu T)$.
This is an entire function of the complex variable $\nu$, hence it vanishes for a discrete
set of values of $\nu$ (possibly empty). Furthermore, $\l$ is an eigenvalue of $Q$ w.r.t. $Q_0$ if and only if the
quadratic form $Q_{- 1/\l}$ is non invertible, that is if and only if $\D (- 1 /\l) = 0$.

\medskip

As a consequence, the function of real variable $\nu \mps \D (\nu)$ has, at most, a finite number
of zeroes $\nu_{(1)}, \cdots , \nu_{(p)}$ in the interval $[0, 1]$ with 
$0 < \nu_{(1)} < \cdots < \nu_{(p)} \leq 1$. According to the relation (A.107), between two such
zeroes, the following differential equation holds:
$${ d \over d \nu } \ell n \Det (Q_\nu / Q_0) = \Tr (Q / Q_\nu) 
\eqno ({\mathop{\rm A.124 }\nolimits})$$
where $Q_\nu = Q_0 + \nu Q$.

\bigskip
\nn
4.5. {\it Some integration formulas.}
\smallskip
With the previous notations, consider any $\nu$ in $[0, 1]$ distinct from
$\nu_{(1)}, \cdots ,$ $ \nu_{(p)}$. Hence the quadratic form $Q_\nu$ on ${\bf X}$ is invertible,
with an inverse quadratic form $W_\nu$ on ${\bf X}'$. There exists an integrator 
$\Da_{\nu} \, x$ on ${\bf X}$ characterized by the formula
$$
\int_{ {\bf X} }^{ } \Da_\nu x\cdot \exp 
\left({ - { \pi \over s} Q_\nu (x) - 2 \pi i 
\left\langle{ x', x }\right\rangle
}\right) = \exp 
\left({ - \pi s \, W_\nu (x') }\right) \, . 
\eqno ({\mathop{\rm A.125 }\nolimits})$$
Our claim is that $\Da_\nu \, x$ {\it is proportional to } $\Da_0 x$, namely that there exists a
constant $I(\nu)$ such that $\Da_\nu \, x = I (\nu) \Da_0 x$. 
More explicitly, we assert the formula
$$
\int_{ {\bf X} }^{ } \Da_\nu x \cdot F (x) = I (\nu) 
\int_{ {\bf X} }^{ } \Da_0 x \cdot F(x) \, , 
\eqno ({\mathop{\rm A.126 }\nolimits})$$
and that a functional $F(\cdot)$ on ${\bf X}$ is Feynman-integrable for 
$\Da_\nu x$ if and only if it is Feynman-integrable for $\Da_0 x$.
The constant $I(\nu)$ can be obtained by putting 
$F(x) =
\exp \left({- { \pi \over s}  Q_\nu (x) }\right)$ into equation (A.126), hence
$$I(\nu)^{- 1 } = \int_{ {\bf X} }^{ }  
\Da_0 x \cdot \exp 
\left({ - { \pi \over s} 
\left({ Q_0 (x) + \nu Q (x)  }\right)
}\right) \, . 
\eqno ({\mathop{\rm A.127 }\nolimits})$$

\bigskip
\nn
{\bf First case $s = 1$:}
\smallskip
According to our conventions, the quadratic forms $Q_0$ and $Q_1$ are positive-definite and
invertible. Since the quadratic form $Q = Q_1 - Q_0$ is nuclear, it follows from the spectral
theory that $Q_0$ and $Q$ can be simultaneously diagonalized. Hence there exists a basis
$\left({ e_n }\right)_{n\geq 1}^{}$ for ${\bf X}$ such that
$$
\eqalignno{
Q_0 \left({ \sum_{ n=1}^{\Inf } t_n e_n  }\right)
&= \sum_{ n=1}^{\Inf } \left({ t_n }\right)^2
& ({\mathop{\rm A.128 }\nolimits}) \cr
Q \left({ \sum_{ n=1}^{\Inf } t_n e_n  }\right)
&= \sum_{ n=1}^{\Inf } \l_n \left({ t_n }\right)^2
& ({\mathop{\rm A.129 }\nolimits})\cr}$$
with real constants $\l_n$ such that 
$\sum_{ n=1}^{\Inf } \left\vert{ \l_n }\right\vert < \Inf $.
Since $Q_1 = Q_0 + Q$ is positive-definite, we have $1 + \l_n > 0$,
hence
$$Q_\nu  \left({ \sum_{ n=1}^{\Inf } t_n e_n  }\right)
= \sum_{ n=1}^{\Inf }
\left({ (1 - \nu) + \nu (1 + \l_n) }\right)
\cdot \left({ t_n }\right)^2
\eqno ({\mathop{\rm A.130 }\nolimits})$$
is again positive-definite and invertible for 
$\nu $ in $[0, 1]$.
Hence the determinant of 
$Q_\nu / Q_0$ is defined and
$$\Det (Q_\nu / Q_0) =
\prod_{n=1 }^{ \Inf}
\left({ 1 + \nu \l_n }\right) > 0 \, .
\eqno ({\mathop{\rm A.131 }\nolimits})$$

\medskip

The main result is given by the following formula:
$$I(\nu) = \Det (Q_\nu / Q_0)^{{1 \over 2}} 
\, . 
\eqno ({\mathop{\rm A.132 }\nolimits})$$
The proof is obtained without difficulty using equations (A.127) to (A.131) and the 
{\it approximation formula}:
$$\int_{ {\bf X} }^{ } \Da_0 x \cdot F (x) =
\build { \lim }_{ N = \Inf }^{ }
\int_{ \RR^N }^{ } d^N t \, 
F 
\left({ t_1 e_1 + \cdots + t_N e_n   }\right)  \, . 
\eqno ({\mathop{\rm A.133 }\nolimits})$$

\bigskip
\nn
{\bf Second case $s = i$:}
\smallskip
Here the basic formula is given by 
$$I (\nu) =
\left\vert{ \Det
\left({ Q_\nu / Q_0  }\right)
 }\right\vert^{1/2} 
i^{
{\mathop{\rm Ind  }\nolimits}
\left({ Q_\nu / Q_0  }\right)
}
\eqno ({\mathop{\rm A.134 }\nolimits})$$
where the {\it index} ${\mathop{\rm Ind  }\nolimits}
\left({ Q_\nu / Q_0  }\right)$ is the number of negative eigenvalues of 
$Q_\nu$ w.r.t. $Q_0$. To simplify the statements, we shall assume that $Q_0$ is
positive-definite. When $\nu$ runs over the interval $[0, 1]$, this index remains constant except
when $\nu$ is crossing a singular value $\nu_{(k)}$, where it experiences a jump.
Hence the index ${\mathop{\rm Ind  }\nolimits}
\left({ Q_\nu / Q_0  }\right)$ is a sum of local contributions from the exceptional values
$\nu_{(1)} , \cdots, \nu_{(p)}$, a phenomenon reminiscent of caustics.

\medskip
The easiest proof of formula (A.134) is obtained by following a strategy initiated by Nelson and
Sheeks in [6]. It works in both cases $s = 1$ and $s = i$.
Starting from equation (A.127), we obtain by derivation
$$
\eqalign{
{d  \over d \nu } I(\nu)^{-1}
&= \int_{{\bf X} }^{ } \Da_0 x \cdot \exp 
\left({ - {\pi  \over s}
\left({ Q_0 (x) + \nu Q (x)  }\right)
 }\right)
\left({ - {\pi  \over s} Q (x)  }\right) \cr
&= I (\nu)^{- 1 } 
\int_{{\bf X} }^{ } \Da_\nu x \cdot \exp 
\left({ - {\pi  \over s}
 Q_\nu (x)  
 }\right)
\left({ - {\pi  \over s} Q (x)  }\right) , \cr}
\eqno ({\mathop{\rm A.135 }\nolimits})$$
that is
$$
{d  \over d \nu } \, \ell n \, I(\nu)
= {\pi  \over s} \int_{ {\bf X} }^{ } 
\Da \om_\nu (x) \, Q(x) \, .
\eqno ({\mathop{\rm A.136 }\nolimits})$$

\medskip
\nn
Expanding $Q(x) $, according to formula (A.122) and using formula (A.39) for the covariance, we
obtain
$${d  \over d \nu } \, \ell n \, I (\nu) = 
{\pi  \over s} \cdot 
\sum_{n=1 }^{\Inf  }
\a_n \cdot {s  \over 2 \pi }
W_\nu (x'_n) 
\eqno ({\mathop{\rm A.137 }\nolimits})$$
which can be transformed easily into
$$
{d \over d \nu } \, \ell n \, I (\nu) = 
{ {}_{1} \over {}^{2} } \Tr \left({ Q / Q_\nu }\right)
\, . \eqno ({\mathop{\rm A.138 }\nolimits})$$

\medskip
\nn
According to equation (A.124), we conclude
$${d  \over d \nu } \, \ell n \, I (\nu)
= { {}_{1} \over {}^{2} } {d  \over d \nu } \, \ell n \, 
\Det \left({ Q_\nu / Q_0 }\right) \, . 
\eqno ({\mathop{\rm A.139 }\nolimits})$$

\medskip
\nn
It remains to study the shift in phase when $\nu$ goes through an exceptional value
$\nu_{(k)}$.

\bigskip
\nn
{\it Remark.}
To simplify matters, assume again  that $Q_0$ is positive-definite. Let $L$ be an invertible operator in
${\bf X}$, of the form $1 + T$ where $T$ is nuclear.
Putting $Q_1 = Q_0 \circ L$, it can be shown that the quadratic form 
$Q = Q_1 - Q_0$ is nuclear and that 
$\Det \left({ Q_1 / Q_0 }\right) $ is equal to 
$\Det (L)^2$. 
The quadratic from $Q_1$ is positive-definite, hence
${\mathop{\rm Ind }\nolimits}
\left({ Q_1 / Q_0 }\right)  = 0$.
>From our main formulas (A.132) and (A.134), we derive
$$\Da_{s, Q_0} (Lx) =
\left\vert{ \Det (L)  }\right\vert
\cdot \Da_{s, Q_0}  x \, .
\eqno ({\mathop{\rm A.140 }\nolimits})$$

\bigskip

Coming back to the general situation, where $Q_\nu = Q_0 + \nu Q$, the reader will find in [5, p. 277-279] examples
of operators $L(\nu)$ such that $Q_\nu = Q_0 \circ L (\nu)$. In such a case, we obtain 
$\Da_\nu x = \left\vert{ \Det L(\nu) }\right\vert
\cdot \Da_0 x $.

\bigskip
Let us summarize the main results obtained in this paragraph:

\medskip
{\it Let $Q_0$ and $Q_1$ be continuous quadratic forms on the Banach space ${\bf X}$.
 We assume that $Q_0$ and $Q_1$
are invertible and that $Q_1-Q_0$ is nuclear. 

\smallskip
(A) Assuming that both $Q_0$ and $Q_1$  be positive-definite, we obtain
$$\int_{{\bf X}}{} \Da_{Q_0} x \cdot e^{- \pi Q_1 (x) } =
\Det \! \left({Q_1/Q_0 }\right)_{}^{-1/2} \, .
\eqno ({\mathop{\rm A.141 }\nolimits})$$

\smallskip
(B) For the the oscillatory integral, we obtain
$$\int_{{\bf X}}{} \Da_{i, Q_0} x \cdot e^{ \pi i Q_1 (x) } =
\left\vert{ \Det \! \left({Q_1/Q_0 }\right) }\right\vert_{}^{-1/2}
i^{-{\mathop{\rm Ind  }\nolimits} \! \left({Q_1/Q_0 }\right)  } 
\eqno ({\mathop{\rm A.142 }\nolimits})$$
where the index ${\mathop{\rm Ind  }\nolimits} \! \left({Q_1/Q_0 }\right)$ counts 
the number of negative eigenvalues
of $Q_1 $ w.r.t. $Q_0$.

\bigskip
Equations } (A.141) {\it and } (A.142) {\it  justify the basic formula } (A.126)
{\it
$$ 
\int_{{\bf X}}{} \Da_{\nu} x \cdot F(x) =
I(\nu) 
\int_{{\bf X}}{} \Da_{0} x \cdot F(x)
$$
with $I(\nu) $ given by } (A.132) {\it or}  (A.134).

\vfill\eject
\vglue 1cm
\centerline{J\bf Appendix B } 
\medskip
\centerline {\titre   Functional determinants of Jacobi operators}
\bigskip
\bigskip

Many functional determinants of Jacobi operators have been computed in works on
semiclassical expansions. They are conveniently expressed in terms of finite-dimensional
determinants of Jacobi matrices. We give here the backbone of the method, and
applications to Jacobi operators along a classical path $x_{{\mathop{\rm cl }\nolimits}}$ 
characterized by $d$ initial
conditions (position or momentum at $t_a$) and $d$ final conditions (position or momentum
at $t_b$). Details, proofs, applications and generalizations of the equations presented
here are scattered in the literature and we give a few selected references at the end of
this appendix. A general presentation on the properties of Jacobi operators will be found
in the Ph.D dissertation of John La Chapelle [University of Texas (Austin) Ph. D. expected in 1995].

\bigskip
\nn
{\bf 1. Jacobi fields and Jacobi matrices.}
\smallskip
A {\it Jacobi field} is a vector field along a classical path obtained by variation through
classical paths. A $d \ts d$ Jacobi matrix is built up from Jacobi fields; each column
consists of the components of a Jacobi field. We give explicit constructions when the action
functional is
$S(x) = 
\int_{ {\bf T} }^{ }
dt \, L(x(t), \dot x (t)) $.
 
\medskip
Let $\{ x_{{\mathop{\rm cl }\nolimits}} (\mu) \} $ be a $2d$-parameter
family of classical paths (critical paths of the functional $S$)  with values in a
manifold  $M^d$ 
$$x_{{\mathop{\rm cl }\nolimits}}(\mu): {\bf T} \ra M^d \quad , \quad {\bf T} = [0, T]\, .$$

\smallskip
\nn
We introduce the notations 
$$
\eqalignno{
 x_{{\mathop{\rm cl }\nolimits}} (t ; \mu )&: = (x_{{\mathop{\rm cl }\nolimits}} (\mu)) (t) 
& ({\mathop{\rm B.1 }\nolimits})\cr
\dot x_{{\mathop{\rm cl }\nolimits}} (t ; \mu )&: = \part x_{{\mathop{\rm cl }\nolimits}} (t ; \mu
) / \part t  & ({\mathop{\rm B.2 }\nolimits}) \cr
x'_{{\mathop{\rm cl }\nolimits}, \a }(t ; \mu )&: =
 \part x_{{\mathop{\rm cl }\nolimits}} (t ; \mu ) / \part  \mu^\a.
& ({\mathop{\rm B.3 }\nolimits}) \cr} $$

\medskip
\nn
We choose $\mu = (\mu^1, \cdots , \mu^{2d})$ to stand for $2d$ initial
conditions that characterize the classical path $x_{{\mathop{\rm cl }\nolimits}} 
(\mu)$, assumed to be unique for the
time being; namely,
$$S' (x_{{\mathop{\rm cl }\nolimits}} (\mu)) = 0 \quad
\eqno ({\mathop{\rm B.4 }\nolimits})$$
has a unique solution for a given set $\mu$ of parameters.

\medskip
By varying successively the $2d$ initial conditions, one obtains $2d$ Jacobi fields
$\{ \part x_{{\mathop{\rm cl }\nolimits}} / \part \mu^\a \} $ for 
$\a \in \{ 1, \cdots, 2 d \} $. It is easy to show that
$$S'' (x_{{\mathop{\rm cl }\nolimits}} (\mu)) \cdot 
{\part x_{{\mathop{\rm cl }\nolimits}} \over \part \mu^\a } = 0 \, .
\eqno ({\mathop{\rm B.5 }\nolimits})$$

\medskip
\nn
We assume for the time being that the quadratic form 
$S'' (x_{{\mathop{\rm cl }\nolimits}} (\mu))\cdot \xi \xi  $ is not degenerate for variations $\xi
$ which respect the boundary conditions specified by $\mu$.  Therefore the $2d$ Jacobi fields are
linearly independent. 

\medskip
There are several convenient basis for the $2d$-dimensional space of Jacobi
fields. Choose as initial condition
$$ \mu = (x_a, p_a)
\eqno ({\mathop{\rm B.6a }\nolimits})$$
where
$$ x_a: = x (t_a)
\eqno ({\mathop{\rm B.6b }\nolimits})$$
$$ p_a: = \part L / \!\!
\left.{ \part \dot x (t) }\right\vert_{t = t_a} \, . 
\eqno ({\mathop{\rm B.6c }\nolimits})$$
The corresponding Jacobi fields are
$$j^{\bu \b } (t) = \part x_{{\mathop{\rm cl }\nolimits}}^\bu 
\left({ t; x_a, p_a }\right)  / \part p_{a ,\b}
\eqno ({\mathop{\rm B.7 }\nolimits})$$
$$k^\bu{}_\b (t) = \part x_{{\mathop{\rm cl }\nolimits}}^\bu 
\left({ t; x_a, p_a }\right) / \part x_a^\b \, .
\eqno ({\mathop{\rm B.8 }\nolimits})$$
Having introduced $p_a$ as one of the initial conditions, the reader suspects that we shall
also use\footnote{$^{31}$}{The classical momentum 
$p_{ {\mathop{\rm  cl }\nolimits}, \bu }  \left({ t; x_a, p_a }\right)  $
is obtained by evaluating $\part L/ \part \dot x (t) $ at the point
$x (t) = x_{ {\mathop{\rm  cl }\nolimits} } 
\left({ t; x_a, p_a }\right) $.
}
$$\ww k_\bu{}^\b (t) = 
\part p_{{\mathop{\rm cl }\nolimits}, \bu } 
\left({ t; x_a, p_a }\right)  / \part p_{a ,\b}
\eqno ({\mathop{\rm B.9 }\nolimits})$$
$$\ell_\bu{}_\b (t) = 
\part p_{{\mathop{\rm cl }\nolimits}, \bu  } 
\left({ t; x_a, p_a }\right) / \part x_a^\b \, .
\eqno ({\mathop{\rm B.10 }\nolimits})$$

\medskip
\nn
Let the {\it Jacobi matrices} $J, K, \ww K, L$ be made of the Jacobi fields $j, k, \ww k, \ell$
respectively
$$J^{\a \b} (t, t_a) = j^{\a \b } (t)
\eqno ({\mathop{\rm B.11 }\nolimits})$$
$$
K^\a {}_\b (t, t_a) = k^\a {}_\b (t)
\eqno ({\mathop{\rm B.12 }\nolimits})$$
$$\ww K_\a {}^\b (t, t_a) = \ww k_\a {}^\b (t)
\eqno ({\mathop{\rm B.13 }\nolimits})$$
$$L_{\a \b } (t, t_a) = \ell_{\a \b } (t) \, .
\eqno ({\mathop{\rm B.14 }\nolimits})$$

\medskip
\nn
The $2d \ts 2d $ matrix constructed from the four $d \ts d$ blocks
$J, K, \ww K, L$ is a solution of the equation defined by the Jacobi operator in phase
space. The properties of the Jacobi matrices are numerous. We note only their properties at
$t = t_a$, which can be read off from their definitions :
$$J (t_a, t_a) = 0 \quad , \quad L (t_a, t_a) = 0 
\eqno ({\mathop{\rm B.15 }\nolimits})$$
$$K(t_a, t_a) = \un \quad , \quad \ww K (t_a, t_a) = \un \, .
\eqno ({\mathop{\rm B.16 }\nolimits})$$

\medskip
\nn
The matrices 
$K$ and $\ww K$ are indeed the transposed of each other
$$\ww K_\a {}^\b (t, t_a) = K^\b {}_\a (t_a, t) .
\eqno ({\mathop{\rm B.17 }\nolimits})$$

\medskip
\nn
Let $M, N, \ww N, P$ be the matrix inverses of $J, K, \ww K, L$ respectively.
Hence we have
$$J^{\a \b } (t, t_a) M_{\b \g } (t_a, t) = \d_\g^\a
\eqno ({\mathop{\rm B.18 }\nolimits})$$
and three similar equations expressed in condensed form as follows
$$KN = \un, \quad \ww K \ww N = \un, \quad LP = \un \, .
\eqno ({\mathop{\rm B.19 }\nolimits})$$
The matrices $M, N, \ww N, P$ are the hessians of the corresponding action functions (Van
Vleck matrices):
$$M_{\b \a } (t_b, t_a) =
{\part^2 \Sa  \over \part x_{{\mathop{\rm cl }\nolimits}}^\b (t_b) 
\part x_{{\mathop{\rm cl }\nolimits}}^\a (t_a)} ,
\eqno ({\mathop{\rm B.20 }\nolimits})$$
and similarly $N$ and $P$ are the hessians of 
$\Sa (x_{{\mathop{\rm cl }\nolimits}} (t_b), p_{{\mathop{\rm cl }\nolimits}} (t_a)) $ and 
$\Sa (p_{{\mathop{\rm cl }\nolimits}} (t_b), p_{{\mathop{\rm cl }\nolimits}} (t_a))$ respectively.

\bigskip
\nn
{\bf 2. Jacobi operators and their Greens's functions.}
\smallskip
The Green's functions of the Jacobi operators in phase space can be expressed in terms of
the Jacobi matrices $J, K, \ww K, L$; they include the Green's functions of the Jacobi
operators in configuration space. For the sake of brevity, but at the cost of elegance, we
consider here only the Green's functions of the Jacobi operators in configuration space.

\medskip
In the previous paragraph, $\ $ Jacobi fields were obtained by variation \break through classical
paths. In this paragraph, we consider a one-parameter family of paths $ x (\l)  $
satisfying $d$ fixed boundary conditions at $t_a$ and $d$ fixed boundary conditions at 
$t_b$, 
referred in brief as ``$a$'' and ``$b$''. The corresponding path space is
 denoted by $\Pa_{a, b} M$, hence
$$x(\l) \in \Pa_{a, b} M \sbs \Pa M \, . 
\eqno ({\mathop{\rm B.21 }\nolimits})$$

\medskip
\nn
We assume that for $\l = 0$, the path 
$x_{{\mathop{\rm cl }\nolimits}} = x(0)$ is the unique critical point of the action functional 
$S$ restricted to $\Pa_{a, b} M$. We set
$$\dot x (\l, t): = \part x (\l , t) / \part t \ , \quad  x' (\l, t): = \part x (\l, t) / \part 
\l \ ,  \quad x' (0, t) = \xi (t) .$$

\medskip
\nn
The second variation of the action functional gives the Jacobi operator 
at $x_{{\mathop{\rm cl }\nolimits}}$:
$$\eqalign{
S'' (x_{{\mathop{\rm cl }\nolimits}}) \cdot \xi \xi
& : = {d^2 \over d \l^2 } 
\left.{  S (x (\l)) }\right\vert_{\l = 0} 
\cr
& = 
\int_{{\bf T}}^{ } dt 
\left({ L_{1 \a, 1 \b} \xi^\a (t) + L_{2 \a, 1 \b } \dot \xi^\a (t) }\right)
\xi^\b (t) \cr
&\ \ \ + \int_{{\bf T} }^{ } dt 
\left({ L_{1 \a, 2 \b} \xi^\a (t) + 
L_{2 \a, 2 \b} \dot \xi^\a (t) }\right)
\dot \xi^\b (t) \cr} \eqno ({\mathop{\rm B.22 }\nolimits}) $$
where
$L_{2 \a, 2 \b}$ is $  \part^2 L / \part \dot x^\a (t) \part \dot x^\b (t)$ evaluated at
 $x_{{\mathop{\rm cl }\nolimits}}$,
etc...
\smallskip

\medskip
If we integrate (B.22) by parts, we obtain the Jacobi operator in its differential garb
(second order differential operator on the space of vector fields $\xi$ 
along $x_{{\mathop{\rm cl }\nolimits}}$
together with boundary terms at $t_a$ and $t_b$). It is often simpler to work with the
quadratic form 
$S'' (x_{{\mathop{\rm cl }\nolimits}}) \cdot \xi \xi$ as written in (B.22). We call functional
Jacobi operator
 the kernel corresponding to the quadratic form 
$S'' (x_{{\mathop{\rm cl }\nolimits}}) \cdot \xi \xi$, namely:
$$ {1 \over 2}
 {\d^2 S'' (x_{{\mathop{\rm cl }\nolimits}}) \cdot \xi \xi \over \d \xi^\a (s) \d \xi^\b (t)}
=: \Ja_{\a \b } (x_{{\mathop{\rm cl }\nolimits}}, s, t ) \, .
\eqno ({\mathop{\rm B.23 }\nolimits})$$

\medskip
\nn
Its inverse will be called its {\it Green's function }
$G^{\bu \bu } \left({ x_{{\mathop{\rm cl }\nolimits}}, t, u }\right) $, namely
$$\int_{\bf T }  dt \, \Ja_{\a \b}
(x_{{\mathop{\rm cl }\nolimits}}, s, t)  G^{\b \g}
(x_{{\mathop{\rm cl }\nolimits}}, t, u) =
\d_\a^\g \d (s - u ) \, .
\eqno ({\mathop{\rm B.24 }\nolimits})$$

\medskip
\nn
Provided the quadratic form 
$S'' (x_{{\mathop{\rm cl }\nolimits}})\cdot  \xi \xi$ is non degenerate, the functional Jacobi
operator has a unique inverse. We say $S'' (x_{{\mathop{\rm cl }\nolimits}})$ is  {\it non
degenerate} if for any $\xi \ne 0$ in $T_{x_{{\mathop{\rm cl }\nolimits}} } \Pa_{a, b} M$
there exists $\eta$ in this vector space with 
 $$S'' (x_{{\mathop{\rm cl }\nolimits}}) \cdot \xi \eta \ne 0 \, .
\eqno ({\mathop{\rm B.25 }\nolimits})$$

\medskip
\nn
This equation says that there are no (nonzero) Jacobi field in the tangent space
 $T_{x_{{\mathop{\rm cl }\nolimits}}}
\Pa_{a, b} M$ to the space of paths with boundary conditions $a$ and $b$.

\bigskip
We list below the Green's functions of Jacobi operators at classical paths with different
boundary conditions. A more abstract formula could encode all cases. 
Explicit formulas may be more
useful for applications.

\medskip
i) 
{\it The classical path is characterized by }
$ p_{{\mathop{\rm cl }\nolimits}} (t_a) = p_a, \ x_{{\mathop{\rm cl }\nolimits}}
 (t_b) = x_b$ 
$$
\eqalign{ G(t, s) &= \t (s - t) K (t, t_a) N (t_a, t_b) J (t_b , s)\cr
&\ \ \ - 
\t (t- s) J (t, t_b) \ww N (t_b, t_a) \ww K (t_a, s) \, .\cr}
\eqno ({\mathop{\rm B.26 }\nolimits})$$

\medskip
ii)
{\it The classical path is characterized by } 
 $x_{{\mathop{\rm cl }\nolimits}} (t_a) = x_a$, $p_{{\mathop{\rm cl }\nolimits}} (t_b) =p_b$
$$\eqalign{ G(t, s) &=
\t (s - t) J (t, t_a) \ww N (t_a, t_b) \ww K (t_b, s) \cr
&\ \ \ - \t (t - s) K (t, t_b) N (t_b, t_a) J (t_a, s) \, .\cr}
\eqno ({\mathop{\rm B.27 }\nolimits})$$

\medskip
iii)
{\it The classical path is characterized by }
$x_{{\mathop{\rm cl }\nolimits}} (t_a) = x_a$, $x_{{\mathop{\rm cl }\nolimits}} (t_b) = x_b$
$$\eqalign{ G(t, s) &= \t (s - t) J (t, t_a) M (t_a, t_b) J (t_b, s) \cr
&\ \ \  - \t (t - s) J
(t, t_b) M (t_b, t_a) J (t_a, s) \, .\cr}
\eqno ({\mathop{\rm B.28 }\nolimits})$$

\medskip
iv)
{\it The classical path is characterized by 
} $p_{{\mathop{\rm cl }\nolimits}} (t_a) = p_a$, 
$p_{{\mathop{\rm cl }\nolimits}} (t_b) = p_b$
$$\eqalign{ G(t, s) &=
\t (s - t) K (t, t_a) P (t_a, t_b) \ww K (t_b, s) \cr
&\ \ \ -
\t (t - s) K (t, t_b) \ww P (t_b, t_a) \ww K (t_a, s) \, .\cr}
\eqno ({\mathop{\rm B.29 }\nolimits})$$

\bigskip
\nn
{\bf 3. Semiclassical expansions.}
\smallskip
We have given examples of WKB approximations in two cases:

i) The classical path of reference is characterized by initial momentum, final position (paragraph
III.2).

\medskip
ii) The classical path is characterized by initial position, final position (paragraph IV.1).

\medskip
Two other cases are often needed:

iii) The classical path is characterized by initial position and final momentum. The transposition
from case i) is straightforward.

\medskip
iv) The classical path is characterized by initial momentum and final momentum.
The critical points of the action functional are degenerate when the classical system is
constrained by a conservation law. We refer the reader to the reference [15] for this case.

\medskip
The computation of semiclassical expansions along the lines of paragraphs III.2 and IV.1, but
with general action functionals $S$, brings in a second variation
$$S'' (x_{{\mathop{\rm cl }\nolimits}}) \cdot \xi \xi = Q_0 (\xi) + Q (\xi)
 \eqno ({\mathop{\rm B.30 }\nolimits})$$
where
$$Q_0 (\xi) =  
\int_{ {\bf T} }^{ } dt \,  L_{2 \a, 2 \b }
\dot \xi^\a (t) \dot \xi^\b (t) \, .
\eqno ({\mathop{\rm B.31 }\nolimits})$$

\medskip
\nn
Provided the Legendre matrix
$$L_{2 \a, 2 \b } 
(x_{{\mathop{\rm cl }\nolimits}} (t), \dot x_{{\mathop{\rm cl }\nolimits}} (t)) =
\part^2 L / \part \dot x_{{\mathop{\rm cl }\nolimits}}^\a \part 
\dot x_{{\mathop{\rm cl }\nolimits}}^\b
\eqno ({\mathop{\rm B.32 }\nolimits})$$
is invertible, the Gaussian integrator defined by the quadratic form(B.30) can be handled
 by the same
techniques as the Gaussian integrator defined by the quadratic form 
$$\int_{{\bf T }} dt \,  h_{\a \b} \dot \xi^\a (t) \dot \xi^\b (t) \, .
\eqno ({\mathop{\rm B.33 }\nolimits})$$

\medskip
The contribution of the second variation 
$S'' (x_{{\mathop{\rm cl }\nolimits}}) \cdot \xi \xi $ to the semiclassical expansion of
$\Psi(t_b, x_b)$ is
$$\int \Da_{ s, Q_0}^{ } \xi\, \,   \exp 
\left({ - { \pi \over s } (Q_0 (\xi) + Q(\xi) }\right)=
 \Det\left({ Q_0 / (Q_0 + Q) }\right)_{}^{1/2}.
\eqno ({\mathop{\rm B.34 }\nolimits})$$
Then one uses the Green's functions (B.26-29) to identify the ratio of Jacobi matrices whose
determinant is equal to $\Det (Q_0 / (Q_0 + Q))$. The results are identical to the results obtained
by discretizing the functional determinants, and reported in [16].

\bigskip
\nn
{\bf 4. References for Appendix B.}
\smallskip
The most complete {\it summary}, to date, of Jacobi fields in phase space, including
degenerate critical points of the action, can be found in Appendix A of 

\smallskip
- C. DeWitt-Morette ``Feynman path integrals'' Acta Physica Austriaca, Suppl. XXVI,
101-170 (1984).

\medskip
For the {\it proofs}  of the results presented in the present  paper, and for properties of Jacobi
matrices used in the present  paper, see [2, 5, 15, 16] and 

- C. DeWitt-Morette and T.-R. Zhang ``Feynman-Kac formula in phase space with application to
coherent state transitions'' Phys. Rev. {\bf D 28}, 2517-2525 (1983).

\medskip
\nn
More on degenerate critical points can be found in 

- C. DeWitt-Morette, B. Nelson and T.-R. Zhang ``Caustic problems in quantum mechanics
with applications to scattering theory'' Phys. Rev. {\bf D 28}, 2526-2546 (1983).

- C. DeWitt-Morette and B.L. Nelson ``Glories and other degenerate points of the action''
 Phys. Rev.
{\bf D 29}, 1663-1668 (1984).

\medskip
A short introduction to the Hamiltonian techniques underlying the above results can be found in 

- P. Cartier ``Some fundamental techniques in the theory of integrable systems''
in {\it Lectures on Integrable Systems} (O. Babelon, P. Cartier and Y. Kosmann-Schwarzbach edit.),
World Scientific, Singapore (1994).

\vfill\eject
\vglue 1cm
\centerline{J\bf Appendix C } 
\medskip
\centerline {\titre   A new class of ordinary differential equations}
\bigskip
\bigskip

The purpose of this Appendix is to extend to the $L^{2,1}$ case the familiar theorems about the existence and
uniqueness of solutions of differential equations, and to describe the parametrization of paths in a curved
space by means of paths in a flat space.

\bigskip
\nn
{\bf 1. Solutions of differential equations: t$\!$he classical case.}

\smallskip
We shall follow the usual strategy, as described in any standard textbook, for instance Bourbaki  [27]. More
precisely, consider a domain $\O$ in the Euclidean space $\RR^n$ and a vector field in $\O$ associating to
every epoch $t$ in the time interval ${\bf T}  = 
\left\lbk{ t_a, t_b }\right\rbk
$ and to every point $x$ in $\O$ a velocity vector $v(t, x) $ in $\RR^n$. 
We assume that $v(t, x)$ is a continuous function of $t$ and $x$.
The two basic remarks are as follows:

\smallskip
a) {\it any trajectory of the vector field $v$ can be prolongated as long as it does not reach the boundary of
$\O$;

}

\smallskip
b)J$\,${\it by the mean value theorem, if the absolute velocity $|v|$ is bounded by a constant $V$ along a given
trajectory leading from $x_a$ at time $t_a$ to $x_b$ at time $t_b$, then the mean velocity
$ { \left\vert{x_b - x_a  }\right\vert \over t_b - t_a } $ is bounded by $V$. 

}

\smallskip
\nn
>From these remarks follows the following existence theorem (Peano):

\smallskip
{\it Assume that $\O$ is a closed ball centered at $x_a$ of radius $L$, and that 
$$\left\vert{ v (t, x) }\right\vert
< L / T
\eqno ({\mathop{\rm C.1 }\nolimits})$$
holds uniformly for $t$ in  ${\bf T}$ and $x$ in $\O$, where $T = t_b - t_a$ is the length of ${\bf T}$.
Then there exists a solution $x : {\bf T} \ra \O$ of the differential equation $\dot x = v (t, x)$ with the
initial condition $x(t_a) = x_a$. }

\medskip
For the proof, we construct first an approximate solution by the Euler method. We select epochs 
$t_1, \cdots , t_{N-1}$ such that $t_a < t_1 < \cdots < t_{N-1} < t_b $ and set
$t_0 = t_a$ , $t_N = t_b$. Then we define inductively points $x_0, x_1, \cdots, x_N$ in $\O$ by $x_0 = x_a$ and
$$x_i = x_{i-1} + \left({ t_i - t_{i-1} }\right) v \! \left({ t_{i-1}, x_{i-1} }\right)
\eqno ({\mathop{\rm C.2 }\nolimits})$$
for $1 \leq i \leq N$; the estimate (C.1) guarantees that the points $x_0, x_1,\cdots, x_N$ are in $\O$.
We then interpolate linearly in each subinterval 
$\left\lbk{ t_{i-1}, t_i }\right\rbk$, and generate a function 
$x_\Ta : {\bf T} \ra \O$ depending on the subdivision 
$\Ta = \left({ t_0 < t_1 < \cdots < t_N }\right)$ of ${\bf T}$.
By the mean value theorem and the estimate (C.1) we obtain
$$\left\vert{ x_\Ta (t) - x_\Ta (t') }\right\vert
\leq V \left\vert{ t - t' }\right\vert
\eqno ({\mathop{\rm C.3 }\nolimits})$$
for $t, t'$ in ${\bf T}$, where $V = L/ T$.  We then use Ascoli's theorem: it asserts 
the existence of a sequence of subdivisions $\Ta (n) $ of ${\bf T}$, whose mesh $\D(n) $ tends to $0$, such that
$x_{\Ta (n)} (t)$ tends to a limit $x(t)$ uniformly for $t$ in ${\bf T}$. Then, one shows that 
$x_{\Ta (n)} $ satisfies an approximate integral equation
$$
\left\vert{ x_{\Ta (n)}(t) - x_a - \int_{ t_a }^{  t}
ds \ v (s, x(s))
 }\right\vert \leq \ve_n \, ,
\eqno ({\mathop{\rm C.4 }\nolimits})$$
where $\build { \lim }_{ n = \Inf  }^{ } \ve_n = 0$.
By uniform convergence, the limit function $x(\cdot) $ satisfies the integral equation
$$x(t) = x_a + \int_{ t_a}^{ t}
ds \  v(s, x(s)) \, ,
\eqno ({\mathop{\rm C.5 }\nolimits})$$
fully equivalent to the differential equation
$$\left\{\matrix{
\dot x (t) = v(t, x (t))
 \hfill\cr
\noalign{\medskip}
x(t_a) = x_a \, .   
 \hfill \cr}
\right.
\eqno ({\mathop{\rm C.6 }\nolimits})$$

\medskip
Both the {\it uniqueness } of a solution to the previous equation and the {\it continuous dependence of the
unique solution on the initial position} $x_a$ are established by an analysis of the {\it stability} of
trajectories. Suppose given another trajectory $y : {\bf T} \ra \O$ with initial position $y_a = y (t_a)$; such
a trajectory exists if $y_a$ is sufficiently close to $x_a$. Hence we have the integral relation
$$y(t) = y_a + \int_{t_a }^{ t} ds \, v(s, y (s)) \, .
\eqno ({\mathop{\rm C.7 }\nolimits})$$
Furthermore, assume that the velocity field satisfies a {\it Lipschitz condition}:
$$\left\vert{ v(t, x) - v(t, x') }\right\vert
\leq k \left\vert{ x - x' }\right\vert
\eqno ({\mathop{\rm C.8 }\nolimits})$$
for $x, x'$ in $\O$ and $t$ in ${\bf T}$, with a fixed constant $k > 0$. Denote by $\d(t) $ the distance
between $x(t) $ and $y(t)$. From (C.5), (C.7) and (C.8) one derives
$$\left\vert{ \d(t) - \d(t_a)  }\right\vert
\leq k \int_{t_a }^{ t} ds \, \d (s) \, .
\eqno ({\mathop{\rm C.9 }\nolimits})$$

\nn
By using an iteration procedure reminiscent of Picard's method, one derives  the following estimate:
$$\d (t) \leq \d (t_a) e^{k (t - t_a)}
\eqno ({\mathop{\rm C.10 }\nolimits})$$
(``{\it  Gronwall's lemma}''). 
If $x_a = y_a$, then $\d (t_a) = 0$, hence $\d(t) = 0$ for all $t$ and $x(t) = y(t)$: uniqueness. Furthermore
if $\d(t_a) = \left\vert{ y_a - x_a }\right\vert$ tends to $0$, then $y(t)$ tends
uniformly to $x(t)$ for  $t$ in the finite interval ${\bf T}$.

\medskip
The results obtained so far are local.
To get a global existence theorem, assume that $N$ is a compact manifold, of dimension $n$, and that $X$ is a
time-dependent vector field on $N$, namely $X(t, x) $ belongs to $T_x N$ for $t $ in ${\bf T}$, and $x$ in $N$.
Assume furthermore that $X(t, x) $ is continuous in $t, x$, with continuous first derivatives in $x$.
Choose a point $x_a$ in $N$ and a coordinate system around $x_a$. Then $X$ is expressed in this coordinate
system by a function $v : {\bf T} \ts \O \ra \RR^n$, where $\O$ is a  closed ball centered at $x_a$, with radius
$L$. Since a continuous function is bounded on any compact set, there exists a constant $\tau > 0$ such that
$| v(t, x) | \leq L / 2 \tau$ uniformly for $t$ in ${\bf T}$ and $x$ in $\O$.
>From the local existence theorem, it follows that for any point $x_0$ in the open ball centered at $x_a$ with
radius $L/2$, and any $t_0$ in ${\bf T}$, there exists a unique trajectory 
$x : \left\lbk{ t_0 , t_0 + \tau }\right\rbk \ra N$ such that 
$x(t_0) = x_0$. By compacity of $N$, we can cover $N$ by a finite number of such balls.
Defining $\s $ as the minimum of the life-times $\tau$, we conclude that the existence and uniqueness of a
trajectory $x : \left\lbk{ t_0 , t_0 + \s }\right\rbk \ra N$ with $x(t_0) = x_0$ hold for every point $x_0$ in
$N$.

\medskip
Consider now the time interval ${\bf T} = \left\lbk{ t_a , t_b }\right\rbk $ and subdivide it into subintervals
$\left\lbk{ t_i , t_{i+1} }\right\rbk $ of lengths $\leq \s$ (for $i$ in 
$\{ 0, 1, \cdots, N \}$)
where $t_0 = t_a$, $t_{N+1} = t_b$. So given any point $x_a$ in $N$, there exists a trajectory $x_0(t)$
 (for $t_0
\leq t \leq t_1$) such that $x_0(t_0) = x_a$.
Set $x_1 = x_0(t_1)$ and consider the trajectory $x_1(t) $ (for $t_1 \leq t \leq t_2$) with initial point $x_1$.
By repeating this procedure, we construct step by step the required trajectory $x : {\bf T}\ra N$ such that 
$x(t_a) = x_a$.

\medskip
Proceeding backwards rather than forwards, we conclude that given any point $x_b$ in $N$, there exists a unique
solution $x(t) = x(t ; x_b)$ of the differential equation 
$\dot x (t) = X (t, x(t)) $ with final position $x(t_b) = x_b$.
For fixed $t$, we define the transformation $\Si (t) $ in $N$ taking $x_b$ to $x(t ; x_b)$, hence
$$x (t; x_b) = x_b \cdot \Si (t)\, .
\eqno ({\mathop{\rm C.11 }\nolimits})$$

\bigskip
\nn
{\bf 2. Solutions of differential equations: the $L^{2,1}$ case.}

\smallskip
We sketch the necessary modifications. A more detailed account will appear elsewhere.

\smallskip
Consider an $L^{2,1}$ path $x : {\bf T} \ra \O$; its action is defined by
$$A = \int_{ {\bf T} }^{ }
{|dx|^2\over dt} \, .
\eqno ({\mathop{\rm C.12 }\nolimits})$$
The mean value theorem is replaced by the following estimate
$$
{ \left\vert{ x_b - x_a }\right\vert^2 \over t_b - t_a }
\leq A \, , 
\eqno ({\mathop{\rm C.13 }\nolimits})$$
which follows easily from Cauchy-Schwarz inequality.

\medskip
We make the following assumptions about the velocity field $v(t, x)$.

\smallskip
(A) {\it The function $v(t, x)$ is jointly-measurable in $(t, x)$, and continuous in $x$ for a given $t$, and
furthermore there exists a numerical $L^2$ function \break $V : {\bf T} \ra [0, + \Inf[$ such that
$$
\left\vert{ v(t, x) }\right\vert \leq V(t)
\eqno ({\mathop{\rm C.14 }\nolimits})$$
holds for $x$ in $\O$ and $t$ in ${\bf T}$.
}

\smallskip
(B) {\it
Denoting by $A$ the action $\int_{ {\bf T}}^{ } dt \, V(t)^2$, the radius $L$ of the ball $\O$ centered at
$x_a$, and the length $T$ of the time interval ${\bf T}$ obey the estimate
$$A < L^2 / T \, . \eqno ({\mathop{\rm C.15 }\nolimits})$$
}
\smallskip
Under these conditions, one proves an existence theorem for the differential equation (C.6) with an $L^{2,1}$
solution $x : {\bf T} \ra \O$. As a first step, one replaces Euler's approximation (C.2) by the following
approximate solution
$$x_\Ta (t) = x_{\Ta} (t_{i-1}) +
\int_{t_i }^{t }ds \, v \! 
\left({ s, x_\Ta (t_{i-1}) }\right)
\eqno ({\mathop{\rm C.16 }\nolimits})$$
for $t$ in the subinterval $\left\lbk{ t_{i-1}, t_i }\right\rbk$.
The inequalities (C.13) to (C.15) guarantee that this trajectory $x_\Ta (t)$ remains in the closed ball $\O$.
Furthermore, the velocity $\dot x_\Ta (t)$ of this trajectory satisfies the estimate
$\left\vert{ \dot x_\Ta (t) }\right\vert
\leq V(t)$ for $t$ in ${\bf T}$, and by (C.13), we obtain the estimate
$$\left\vert{ x_\Ta (t) - x_\Ta (t') }\right\vert^2
\leq A \left\vert{ t - t' }\right\vert \, .
\eqno ({\mathop{\rm C.17 }\nolimits})$$
We can then invoke Ascoli's theorem, and find a sequence of approximate solutions $x_{\Ta (n)}$ converging
uniformly on ${\bf T}$ towards an $L^{2,1}$ function $x$. Each approximation $x_{\Ta (n)}$ satisfies an
approximate integral equation, and in the limit the integral equation (C.5) is obtained using Lebesgue's
dominated convergence.

\medskip
For the uniqueness, we need a Lipschitz condition of the type
$$
 \left\vert{ v(t, x) - v(t, x') }\right\vert
\leq k(t) \left\vert{ x - x' }\right\vert
\eqno ({\mathop{\rm C.18 }\nolimits})$$
where the integral $\int_{ {\bf T} }^{ } dt \, k(t)^2$ is finite. 
We use again a variant of Gronwall's lemma.

\medskip
The global results are obtained for a compact manifold $N$ by the reasonings used at the end of paragraph C.1.
Only minor modifications are needed.

\bigskip
\nn
{\bf 3. Parametrization of paths.}
\smallskip
We consider a manifold $N$ and $d$ vector fields 
$X_{(1)}, \cdots , X_{(d)}$ on $N$.
We assume that they are of class $C^1$ and linearly independent at each point of $N$.
For $x$ in $N$ denote by $H_x$ the vector subspace of $T_x N$ generated by 
$X_{(1)}(x), \cdots , X_{(d)}(x)$. The collection of these vector spaces is a subbundle\footnote{$^{32}$}{The
letter $H$ stands for ``horizontal''.} $H$ of the tangent bundle $TN$ to $N$. We consider also a
real  symmetric invertible
matrix
$\left({ h_{\a \b} }\right)$ of size $d \ts d$, and define a field of quadratic forms $h_x$ on $H_x$ by
$$h_x 
\left({ X_{(\a)} (x), X_{(\b)} (x) }\right)
= h_{\a \b} \, .
\eqno ({\mathop{\rm C.19 }\nolimits})$$

\medskip
Fix a point $x_b$ in $N$.
We denote by $\Pa_{x_b}^H N$ the set of $L^{2,1}$ paths 
$x : {\bf T} \ra N$ which satisfy the following conditions:

\smallskip
(A) {\it The endpoint $x(t_b) $ is equal to $x_b$.}

\smallskip
(B) {\it For each epoch $t$ in ${\bf T}$, the velocity vector $\dot x (t)$ lies in the subspace 
$H_{x(t)}$ of $T_{x(t)} N$.

}

\medskip
\nn
We define the {\it action}
 of such a path by 
$$A(x) = 
\int_{{\bf T} }^{ } dt \, h_{x(t)}
\! \left({ \dot x(t), \dot x(t) }\right)\, .
\eqno ({\mathop{\rm C.20 }\nolimits})$$
Since the vectors $X_{(\a)} (x(t))$ (for $1 \leq \a \leq d$) form a basis of 
$H_{x(t)}$, and from the hypothesis that the path $x$ is of class $L^{2,1}$, we infer that there exist functions
$\dot z_\a $ in $L^2 ({\bf T})$ such that
$$\dot x(t) = X_{(\a)} (x(t)) \dot z^\a (t) \, . 
\eqno ({\mathop{\rm C.21 }\nolimits})$$
Furthermore, the function $\dot z^\a$ is the derivative of a function $z^\a$ in $L^{2,1}({\bf T})$ normalized
by 
$z^\a (t_b) = 0$. The vector function $z = \left({ z^1, \cdots z^d }\right)$ is an element of the space denoted
by ${\bf Z}_b$ in paragraph A.3.8. {\it This construction associates to a path $x$ in $\Pa_{x_b}^H N$ a path $z$
in ${\bf Z}_b$ with conservation of action:} 
$$A(x) = Q_0 (z) \, , 
\eqno ({\mathop{\rm C.22 }\nolimits})$$
where as usual $Q_0 (z)$ is equal to 
$\int_{{\bf T} }^{ } dt \, h_{\a \b}
\dot z^\a (t) \dot z^\b (t) $. 

\medskip
Assume now that $N$ is compact. By using the theory of $L^{2, 1}$ differential equations sketched in paragraph
C.2, it can be shown that we can invert the transformation $x \mps z$. Given any $z$ in ${\bf Z}_b$, the
differential equation (C.21) has a unique solution $x$ in $\Pa_{x_b}^H N $, and we obtain {\it a
parametrization
 $$P : {\bf Z}_b \ra \Pa_{x_b}^H N$$
of a space of paths in a curved space $N$ by a space of paths in a flat space $\RR^d$.}
For $z$ in ${\bf Z}_b$, we denote by $x(t, z)$ the solution of the differential equation (C.21) with endpoint 
$x(t_b, z) = x_b$. We can also introduce a global transformation $\Si (t, z) : N \ra N$ taking $x_b$ into $x(t,
z)$.
If necessary, we include 
$t_b$ in the notation and denote this transformation by
$\Si ({\bf T} ; z)$ or
 $\Si \! \left({  t_b, t_a ; z}\right)$ in the case $t = t_a$. The chain rule
$$\Si \! \left({  t_b, t_a ; z_{ba}J}\right)
=
\Si \! \left({  t_b, t_c ; z_{bc}}\right)
\cdot
\Si \! \left({  t_c, t_a ; z_{ca}}\right)
\eqno ({\mathop{\rm C.23 }\nolimits})$$
is a consequence of the uniqueness of the solution of the differential equation (C.21).
Here $z_{ba}$ is $z$, $t_c$ is an intermediate epoch, 
and the paths $z_{bc} : \left\lbk{  t_c, t_b }\right\rbk \ra \RR^d$
and $z_{ca} : \left\lbk{  t_a, t_c }\right\rbk \ra \RR^d$
are given by 
$$z_{bc}(t) = z(t) \quad , \quad
z_{ca} (t) = z(t) - z(t_c) \, . 
\eqno ({\mathop{\rm C.24 }\nolimits})$$

\medskip
\nn
{\it Remark.}
The more general differential equation
$$\dot x (t) = X_{(\a)} 
(x(t)) \dot z^\a (t)  + Y (x(t)) 
\eqno ({\mathop{\rm C.25 }\nolimits})$$
can be handled in a similar way.
In this case, we replace $\Pa_{x_b}^H N$ by the space
$
\Pa_{x_b}^{H, Y} N$ of $L^{2,1}$ paths 
$x : {\bf T} \ra N$ such that
$x (t_b) = x_b $ and $\dot x(t) - Y(t) $ belong to 
$H_{x (t)} $ for every $t$ in ${\bf T}$.

\medskip
These constructions are related to the {\it Cartan development map}. Take for $M$ a compact Riemannian manifold,
and let $N$ be the corresponding bundle of orthonormal frames; it is a compact manifold.
Fix $x_b$ in $M$, and a frame $\r_b = 
\left({ x_b, u_b }\right)$ at $x_b$.
Then, by the Riemannian connection, there is defined a  ``horizontal'' subspace 
$T_{\r_b}^H N$ of the tangent space $T_{\r_b} N$,
and the projection $\Pi : N \ra M$ induces an identification of $T_{\r_b}^H N$
with $T_{x_b} M$.
Since $u_b$ defines an orthonormal basis of $T_{x_b} M$, we obtain a basis
$X_{(1)} (\r_b) , \cdots , X_{(d)} (\r_b)$
of $T_{\r_b}^H N$. This construction is valid for every point in $N$, hence we define vector fields
$X_{(1)} , \cdots , X_{(d)}$ on
 $N$.
By using the previous construction, we get a parametrization
$$P : {\bf Z}_b \ra \Pa_{\r_b}^H N \, . 
$$
But the paths in $\Pa_{\r_b}^H N$ are nothing else than the horizontal liftings of the 
$L^{2,1}$ paths in $M$.
More precisely denote by $\Pa_{x_b} M$ the set of $L^{2,1}$
paths $x : {\bf T} \ra M$, such that $x(t_b) = x_b$. Then, the projection 
$\Pi : N \ra M$ induces a mapping 
$\ww x \mps \Pi \circ \ww x$ of $\Pa_{\r_b}^H N$ into $\Pa_{x_b} M$, and this map is a bijection.
To conclude, we get a diagram
$${\bf Z}_b \ra \Pa_{\r_b}^H N \ra \Pa_{x_b} M$$
and by composition a parametrization of 
$\Pa_{x_b} M$ by ${\bf Z}_b$.
This is the Cartan development map for $L^{2,1}$ paths.
The standard theory works for $C^1$ paths.

\vfill \eject
\centerline {{\bf REFERENCES} }

\bigskip
\bigskip

\item{[1]}
C. Morette 
``On the definition and approximation of Feynman's path integral''
Phys. Rev. {\bf 81}, 848-852 (1951).
 
\medskip

\item{[2]}
C. DeWitt-Morette, B. Nelson, and T.R. Zhang
``Caustic problems in quantum mechanics with applications to scattering theory''
Phys. Rev. D {\bf 28}, 2526-2546 (1983).
 
\medskip
\item{[3]}
M.G.G. Laidlaw and C. Morette-DeWitt 
``Feynman functional integrals for systems of indistinguishable particles'' 
Phys. Rev. D {\bf  3}, 1375-1378 (1971).

\medskip
\item{[4]}
R.H. Cameron and W.T. Martin
``Transformation of Wiener integrals under a general class of linear transformations''
Trans. Amer. Math. Soc. {\bf 58}, 184-219 (1945).

\medskip
\item{[5]}
C. DeWitt-Morette, A. Maheshwari, and B. Nelson 
``Path Integration in Non-Relativistic Quantum Mechanics''
Phys. Rep. {\bf 50}, 266-372 (1979).
This article includes a summary of earlier articles.

\medskip
\item{[6]}
B. Nelson and B. Sheeks 
``Fredholm determinants associated with Wiener integrals''
J. Math. Phys. {\bf 22},
2132-2136 (1981).

\medskip
\item{[7]}
K.D. Elworthy
``{\it Stochastic Differential Equations on Manifolds}''
Cambridge University Press, Cambridge U.K. (1982).

\medskip
\item{[8]}
C. Morette-DeWitt 
``Feynman's Path Integral, definition without limiting procedure''
Commun. Math. Phys. {\bf 28}, 47-67 (1972), and ``Feynman Path Integrals, I. Linear and
affine techniques, II. The Feynman-Green function'' 
Commun. Math. Phys. {\bf 37}, 63-81 (1974).

\medskip
\item{[9]}
S.A. Albeverio and R.J. H\o egh-Krohn
{\it ``Mathematical Theory of Feynman Path Integrals''J}
Springer Verlag Lecture Notes in Mathematics 523 (1976).

\medskip
\item{[10]}
P. Kr\'ee 
``Introduction aux th\'eories des distributions en dimension infinie''
Bull. Soc. Math. France {\bf 46}, 143-162 (1976), and references therein, in particular, {\it
Seminar P. Lelong} Springer Verlag Lecture Notes in Mathematics {\bf 410} and {\bf 474}
(1972-1974).

\medskip
\item{[11]}
P. Cartier and C. DeWitt-Morette
``Int\'egration fonctionnelle; \'el\'ements d'axiomatique''
 C.R. Acad. Sci. Paris, t.316, S\'erie II,  733-738 (1993).

\medskip
\item{[12]}
A. Young and C. DeWitt-Morette ``Time substitutions in stochastic Processes as a Tool in Path
Integration'' Ann. of Phys. {\bf 69}, 140-166 (1986).
 
\medskip
\item{[13]}
Y. Choquet-Bruhat and C. DeWitt-Morette
``Supplement to {\it Analysis, Manifolds and Physics}''
Armadillo preprint, Center for Relativity, University of Texas, 
Austin TX 78712.

\medskip
\item{[14]}
C. DeWitt-Morette, K.D. Elworthy, B.L. Nelson, and G.S. Sammelman
``A stochastic scheme for constructing solutions of the Schr\"odinger equation''
Ann. Inst. H. Poincar\'e A {\bf 32}, 327-341 (1980).

\medskip
\item{[15]}
C. DeWitt-Morette and T.-R. Zhang 
``Path integrals and conservation laws'' 
Phys. Rev. D {\bf  28}, 2503-2516 (1983).

\medskip
\item{[16]}
C. DeWitt-Morette
``The semiclassical expansion''
Ann. of Phys. {\bf 97}, 367-399 (1976).
(Correct a misprint p. 385, l. 5, the reference is [5].)

\medskip
\item{[17]}
S.F. Edwards and Y.V. Gulyaev 
``Path integrals in polar coordinates''
Proc. Roy. Soc. London A{\bf 279}, 224 (1964).

\medskip
\item{[18]}
C. DeWitt-Morette 
``Quantum mechanics in curved spacetimes; stochastic processes on frame bundles''
pp. 49-87 in {\it Quantum Mechanics in Curved Space-Time} 
Eds J. Audretsch and V. de Sabbata, Plenum Press, New York (1990).

\medskip
\item{[19]}
W. Greub and H.-R. Petry 
``Minimal coupling and complex line bundles''
J. Math. Phys. {\bf 16}, 1347-1351 (1975).

\medskip
\item{[20]}
Y. Choquet-Bruhat and C. DeWitt-Morette {\it ``Analysis, Manifolds and Physics Part I:
Basics, Part II: 92 applications''}, North Holland, Amsterdam (1989).

\medskip
\item{[21]}
F.A. Berezin ``Quantization'' Math. USSR, Izvestija {\bf 8}, 1109-1165 (1974).

\medskip
\item{[22]}
D. Bar-Moshe and M.S. Marinov
``Berezin quantization and unitary representations of Lie groups'',
preprint, to appear in Berezin Memorial volume (1994).

\medskip
\item{[23]}
K.R. Parthasarathy
{\it `` Probability measures on metric spaces''}
Academic Press, New York (1967).

\medskip
\item{[24]}
N. Bourbaki
{\it ``Int\'egration, chapitre 9''}
Masson, Paris (1982).

\medskip
\item{[25]}
N. Bourbaki
{\it ``Espaces vectoriels topologiques''}
Masson, Paris (1981).

\medskip
\item{[26]}
P. Cartier ``A course on determinants''
in {\it Conformal invariance and string theory} Eds P. Dita and V. Georgescu, Academic Press, New York (1989).

\medskip
\item{[27]}
N. Bourbaki {\it ``Fonctions d'une variable r\'eelle''}
Masson, Paris (1982).

\vfill\eject
\vglue 	1cm
\centerline {Figure captions}
\vglue 	2cm
{\it Figure} A.1: ``Linear change of variables''

\bye